\documentclass[11pt]{amsart}
\usepackage[utf8]{inputenc}
\usepackage{amsmath}
\usepackage{float}

\usepackage[shortlabels]{enumitem}
\usepackage[makeroom]{cancel}
\usepackage{tikz-cd} 
\usepackage{quiver}
\usepackage{todonotes}
\usepackage{amsfonts, amsthm}

\usepackage{lmodern}
\usepackage{dsfont}

\usepackage{subfig}
\usepackage{mathtools}
\usepackage{amssymb, hyperref, color, graphicx, caption}
\usepackage[top=2cm,bottom=3cm,textwidth=460pt]{geometry}
\usepackage{setspace}
\usepackage{hyperref}
\usepackage{tikzsymbols}
\usepackage{braket}
\usepackage{tcolorbox}
\usepackage{xcolor}   
\usepackage{hyperref}
\hypersetup{
    colorlinks=true, 
    linktoc=all,     
    linkcolor=blue,  
    citecolor=blue
}
\usepackage{calligra,mathrsfs}

\usepackage{subfiles}
\usepackage{macros} 
\usepackage[final,nopatch=footnote]{microtype}
\usepackage{internal_macros}
\usepackage{subfiles}
\usepackage{tikz}
\usepackage{enumitem}

\usetikzlibrary{decorations.markings}
\usepackage{stmaryrd}
\SetSymbolFont{stmry}{bold}{U}{stmry}{m}{n}

\title[Complex Rank Parabolic Category \strO]{Parabolic Category \strO in Complex Rank via Fock Space Tensor Product Categorifications}
\author{Hamilton Wan}
\address{Department of Mathematics \\
Massachusetts Institute of Technology \\
77 Massachusetts Avenue \\
Cambridge, MA 
02139}

\email{wanh24@mit.edu}


\calclayout
\begin{document}
\pagenumbering{gobble}

\begin{abstract}
    We initiate the study of complex rank analogues of parabolic categories $\cO$ for general linear Lie algebras defined via Deligne's interpolating categories. We regard these categories as a family varying over an affine parameter space and conjecture that their structure is controlled by a countable locally finite hyperplane arrangement, that is, they are constant along facets. We prove this conjecture on \textit{admissible} facets using the theory of $\fsl_{\ZZ}$-categorification. The main technical ingredient is a uniqueness theorem for highest weight categories equipped with a categorical type A action categorifying an ordered tensor product of highest and lowest Fock space representations of $\fsl_{\ZZ}$. Under some combinatorial conditions on the parameters, this rigidity result allows us to compare complex rank category $\cO$ with stable limits of classical parabolic categories $\cO$. These equivalences yield character formulas for simple objects in terms of stable limits of parabolic Kazhdan--Lusztig polynomials, answering a problem posed by Etingof. For the case of two Levi blocks of non-integral size, the admissibility assumption is unnecessary, giving a complete description in terms of stable representation theory. As an application, we obtain multiplicities for parabolic analogs of hyperalgebra Verma modules in the large rank and large characteristic limit. 
\end{abstract}

\maketitle


\newpage

{\small
\tableofcontents
}

\clearpage

\newpage

\pagenumbering{arabic}

\section*{Introduction}\label{sec:intro}

The Bernstein--Gelfand--Gelfand category $\cO$ and its parabolic analogues have played pivotal roles in the representation theory of reductive Lie algebras over algebraically closed fields of characteristic zero. For the general linear Lie algebra $\fgl_m$, standard parabolic subalgebras are classified by compositions $(m_1,\ldots,m_n)$ of $m$. Hence, we can view the parabolic categories $\cO$ for general linear Lie algebras as a family of categories indexed by arbitrary compositions. One motivation for this paper is to understand how the structure of these categories varies with the ``block size" parameters $(m_1,\ldots,m_n)$. To understand these categories as an algebraic family, we work in the setting of \textit{complex rank representation theory}, where we can study analogues of parabolic category $\cO$ parameterized by arbitrary \textit{complex} numbers.

The study of representation theory in complex rank aims to understand representation theoretic phenomena for families of groups indexed by some natural number $n$ that, in some sense, depend polynomially on the parameter $n$. Deligne and Milne \cite{deligne_milne} observed that certain structural properties of the category $\Rep{\GL_n}$ of rational representations of the general linear group $\GL_n$ over an algebraically closed field $\fk$ of characteristic zero depend polynomially on $n$ for sufficiently large $n$. By interpolating these polynomials to non-integer $t \in \fk$, they introduced the \textit{Deligne categories} $\rep{\GL_t}$, a family of symmetric tensor categories (STCs) not equivalent to the STCs of rational representations of any algebraic group. 

Etingof \cite{etingof_complex_rank_1,Etingof_complex_rank_2} observed that many representation theoretic objects have natural complex rank analogs and initiated the study of representation theory \textit{internal} to the Deligne interpolating categories. Complex rank representation theory thus illuminates classical representation theoretic phenomena for $\GL_n$ that depend only polynomially on the rank. In fact, Harman \cite{harman16} showed that, in a precise model-theoretic sense, Deligne categories are large rank and large characteristic limits of the classical categories $\Rep{\GL_n}$ over algebraically closed fields of positive characteristic. Hence, representation theory in Deligne categories also sheds light on stability phenomena in usual representation theory. For previous results in complex rank representation theory, see Entova-Aizenbud's work on rational Cherednik algebras \cite{ea_rational_cherednik} and Schur-Weyl duality \cite{ea_schur_weyl}; Kalinov's work on Yangians in complex rank \cite{kalinov_yangians}; Motorin's \cite{motorin} and Pakharev's \cite{pakharev} works on affine Kac--Moody algebras; Utiralova's work on Harish-Chandra bimodules \cite{utiralova_hc,utiralova_ktype}; and Riesen's recent work \cite{riesen} on Feigin--Frenkel duality.

This paper contributes to this growing body of literature by studying natural complex rank analogues of parabolic category $\cO$. The Deligne category $\rep{\GL_t}$ contains a tautological representation $V_t$ of categorical dimension $t$ which gives us the Lie algebra $\fgl_{t} := V_t \ox V_t^*$. The choice of complex numbers $t_1,\ldots,t_n$ satisfying $t = t_1 + \cdots + t_n$ yields a ``block triangular'' decomposition of $\fgl_{t}$ and a parabolic subalgebra of ``block upper triangular matrices.'' Writing $\bt = (t_1,\ldots,t_n)$, we define an associated \textit{parabolic category} $\cO^{\bt}$. This category decomposes as a product of subcategories $\ots$ depending additionally on a sequence $\bs := (s_1,\ldots,s_n) \in \fk^n$ describing a character of the center of the Levi subalgebra of block diagonal matrices (see Definition \ref{defn:ots}).

Many constructions in the classical category $\cO$ generalize to $\ots$. In particular, we define Verma modules $\mts(\blambda)$ with ``highest weight'' given by a multipartition $\blambda = (\lambda_1,\ldots,\lambda_{2n})$, where each $\lambda_i$ is a partition. These Verma modules have unique simple quotients $\lts(\blambda)$ classifying all simple objects in $\ots$. In \cite{Etingof_complex_rank_2}, Etingof posed the following basic problem regarding the structure of $\ots$.

\begin{introproblem}[Etingof]\label{introproblem:etingof}
    Compute $[\mts(\blambda):\lts(\bmu)]$ for any $2n$-multipartitions $\blambda,\bmu$.
\end{introproblem}

We study these categories as a family over the affine space \[\scr{X}_n = \Spec \fk[t_1,\ldots,t_n,s_1,\ldots,s_n].\] For sufficiently generic points in $\scr{X}_n$ (see Definition \ref{defn:interpolatable}), we use the ultraproduct construction of Deligne categories to show that multiplicities in $\ots$ coincide with certain stable limits of multiplicities in usual parabolic categories $\cO$. However, this approach cannot access multiplicities away from these ``sufficiently generic'' points, e.g., at points whose coordinates $(t_1,\ldots,t_n)$ are algebraic numbers. In principle, one could adapt the ultraproduct technique to positive characteristic and compute multiplicities in $\ots$ in terms of stable multiplicities from modular representation theory. Multiplicities in modular representation theory are already very complicated, so this approach is somewhat roundabout. 

We formulate a conjecture to compute multiplicities more generally. We believe that multiplicities are controlled by a countable hyperplane arrangement in the parameter space $\scr{X}_n$ (see Subsection \ref{subsec:stratification}). That is, multiplicities at points with non-integer $t_i$-coordinates are constant on each facet in this hyperplane arrangement. In fact, we conjecture that this numerical behavior is an artifact of a deeper structural fact: the categories $\ots$ themselves ought to be constant along each facet. The precise statement of this conjecture requires a certain abelian subcategory $\extots \subset \ots$, to be introduced in Definition \ref{defn:extots}.

\begin{introconj}\label{introconj:structural}
    For any facet $W$ with $\dim W \geq 1$, the category $\extots$ does not depend on $(\bt,\bs) \in W$ when $t_i \not \in \ZZ$ for $i = 1,\ldots,n$. That is, for any $(\bt,\bs),(\bt',\bs') \in W$ with non-integer $t$-coordinates, there is an equivalence $\extots \simeq (\cO^{t'}_{s'})^{\mathrm{res}}$ intertwining labels of simple objects.
\end{introconj}


\subsection*{Summary of Results}

Our main result is a proof of Conjecture \ref{introconj:structural} under some conceptually illuminating combinatorial conditions. The validity of Conjecture \ref{introconj:structural} in these cases lends credence to the idea that the structure of $\ots$ is tied to the geometry of the parameter space $\scr{X}_{n}$. When $n = 2$, i.e., for a ``maximal parabolic," any facet containing a point $(\bt,\bs)$ with $t_1,t_2 \not \in \ZZ$ satisfies the aforementioned conditions, so our proof establishes Conjecture \ref{introconj:structural} unconditionally then. 

The main tool is the theory of type A categorical actions. We axiomatize the notion of a categorification of an external tensor product of $\fsl_{\ZZ}$-representations given by a tensor product of highest weight and lowest weight Fock spaces (Definition \ref{defn:multi_fock}). This definition extends Losev and Webster's definition of a tensor product categorification from \cite{losev_webster} and incorporates highest weight categories that are unions of upper finite highest weight categories. The presence of both highest and lowest weight factors differentiates our complex rank setting from previous work.

Under suitable assumptions on the ordering of the highest and lowest weight factors, we prove that our tensor product categorifications are essentially unique (Theorem \ref{thm:mixed_admissible_uniqueness}). Our proof extends the approach of Losev and Webster in \textit{loc. cit.} The posets underlying our categories are products of lower and upper finite posets, so their proof, built on a Soergel-esque paradigm of understanding projectives under a quotient functor, does not readily generalize here. Instead, we study tilting objects with respect to a compatible standardly stratified structure.

In turn, we establish complex rank interpolations of classical representation theoretic results, e.g., the Jantzen determinant formula (Theorem \ref{thm:jantzen_determinant}) and sum formula (Theorem \ref{thm:jantzen_sum}), to prove that the categories $\extots$ have the natural structure of Fock space tensor product categorifications (Theorem \ref{thm:ots_mftpc}). Using our uniqueness theorem along with an explicit construction in terms of stable limits of classical parabolic categories $\cO$ (Theorem \ref{thm:construct_mftp}), we obtain a proof of Conjecture \ref{introconj:structural} in the admissible case. 

In particular, our approach shows that the categories $\extots$ are not only constant along admissible facets, but they are also equivalent to explicitly constructed stable limits of parabolic categories $\cO$. As a corollary, we obtain explicit multiplicity formulas for simple objects in Verma modules in terms of stable limits of parabolic Kazhdan--Lusztig polynomials. As an application, we obtain some stability results in modular representation theory.

\subsection*{Detailed Discussion and Main Theorems} 

The guiding philosophy behind our work is to leverage the rich internal symmetry of $\extots$ afforded by the compatible structures of a \textit{categorical $\fsl_{\ZZ}$-action} (in the sense of Rouquier \cite{rouquier_2kacmoody}) and of a \textit{semi-infinite highest weight category}. In particular, we also build on the literature on categorical representations of Kac--Moody algebras.

In general, one should not expect an equivalence of highest weight categories with the same poset. Remarkably, the existence of a categorical Kac--Moody action suitably compatible (e.g., in the sense of Definition \ref{defn:multi_fock}) with an underlying highest weight structure is an incredibly rigid requirement. Most notably, Losev and Webster \cite{losev_webster} axiomatically defined categorified tensor products of irreducible representations of Kac--Moody algebras and proved a strong uniqueness result for them, generalizing work of Rouquier \cite{rouquier_2kacmoody} on the uniqueness of minimal categorifications of irreducible highest weight representations. 



Given that many representation-theoretic categories carry a natural Kac--Moody action compatible with their highest weight structures, it is not surprising that this framework has already borne ample fruit. For example, Brundan, Losev, and Webster \cite{brundan_losev_webster} used the uniqueness results of Losev and Webster to prove multiplicity formulas for integral blocks of parabolic category $\cO$ for the general linear Lie superalgebra. Similarly, Losev \cite{losev_vv_conjecture} used categorical actions on truncations of affine parabolic categories $\cO$ to study categories $\cO$ for rational Cherednik algebras. Another related work is a preprint by Elias and Losev \cite{elias_losev2017} relating modular Kazhdan--Lusztig combinatorics to diagrammatic singular Soergel bimodules using categorical actions of $\widehat{\fsl}_p$. 

Our complex rank setting differs from the aforementioned work in a few ways, and we now highlight what we see as the deviations from the existing theory. The first issue is the decategorified $\fsl_{\ZZ}$-module underlying $\extots$. We will observe that $\extots$ categorifies the tensor product of $2n$ copies of integrable highest \textit{and} lowest weight Fock space representations of $\fsl_{\ZZ}$. Let's call these \textit{Fock space tensor product categorifications} (FTPCs). We remark that Webster \cite{webster_canonical_bases} produced explicit diagrammatic constructions for these categorifications when the factors come in a suitable order, but he did not prove any uniqueness result for these categorifications. 


Interestingly, we can only understand FTPCs when the factors come in an \textit{admissible} order, i.e., all lowest weight factors come first or all highest weight factors come first. We call the former case \textit{upper admissible} and the latter case \textit{lower admissible}. We define a natural highest weight structure on FTPCs. Roughly speaking, our admissible orderings correspond to upper or lower finite highest weight structures (in the language of \cite{brundan_stroppel_semiinfinite}), respectively. That is, the subcategories of projective or tilting objects, respectively, are well-behaved for admissible orderings.

The presence of mutually commuting actions of multiple copies of $\fsl_{\ZZ}$ on $\extots$ poses a more significant obstacle. These actions jointly categorify external tensor products of FTPCs. We call these categorifications \textit{multi-Fock tensor product categorifications (MFTPCs)}. Writing $\cF^0_{\sigma}$ for the highest weight Fock space representation of $\fsl_{\ZZ}$ (where the fundamental weight $\varpi_{\sigma}$ is its highest weight) and $\cF_{\sigma}^1$ for the lowest weight representation with lowest weight $-\varpi_{\sigma}$, an MFTPC categorifies a $\fsl_{\ZZ}^{\oplus r}$-module of the form
\[
    \bigboxtimes_{\ell=1}^r \bigotimes_{i=1}^{a_\ell} \cF_{\sigma_{\ell,i}}^{c_{\ell,i}},
\]
where each $c_{\ell} := (c_{\ell,1},\ldots,c_{\ell,a_\ell})$ is a $01$-sequence with either all zeros appearing first or all ones appearing first (the admissibility condition), and each $\sigma_\ell = (\sigma_{\ell,1},\ldots,\sigma_{\ell,a_\ell})$ is a sequence of integers. The combinatorial data $\Xi := ((\b{\sigma}_1,\b{c}_1),\ldots,(\b{\sigma}_r,\b{c}_r))$ is called the \textit{type} of the MFTPC.

The main complication for proving uniqueness of MFTPCs arises from the potential for ``mixed'' admissibility, i.e., where the external factors can be either upper or lower admissible. In our axiomatization of MFTPCs, these cases correspond to neither upper nor lower finite highest weight categories, but rather, ``interval-finite" highest weight categories with a filtration by upper finite highest weight categories. 

Projective and tilting objects are no longer well-behaved with respect to this structure, so the ``double centralizer" approach used by Losev and Webster does not go through. Instead, we produce a compatible standardly stratified structure on mixed admissible MFTPCs and use the tilting objects with respect to this compatible structure as a substitute for projective objects.

We actually prove our uniqueness theorem for \textit{restricted} MFTPCs, which are \textit{finite} highest weight categories that categorify a finite direct sum of weight spaces; they generalize the ``restricted tensor product categorifications" introduced by Losev \cite{losev_vv_conjecture}. We work with restricted categories because they allow us to ``glue" subquotients of classical parabolic category $\cO$ to obtain genuine MFTPCs. The following statement is a simplified version of Theorem \ref{thm:mixed_admissible_uniqueness}.

\begin{introthm}\label{intro_thm:uniqueness_restricted}
    If $\cC$ and $\cD$ are mixed admissible restricted MFTPCs of the same type, then there is a strongly equivariant equivalence $\cC \simeq \cD$ intertwining the labels of simple objects.
\end{introthm}

To our knowledge, the existing literature does not address categorifications of ``mixed'' external tensor products of tensor products of highest and lowest weight representations. In an unpublished note \cite{losev_note}, Losev proves the uniqueness of FTPCs of $\hat{\fsl}_p$ with one highest weight factor and one lowest weight factor; however, his result does not generalize directly to our setting, for instance, because of the aforementioned obstacles presented by ``mixed admissible" tensor factors. A related paper by Entova-Aizenbud \cite{ea_deligne_categorical_actions} studies categorical actions in Deligne categories. In our language, she proves uniqueness of FTPCs in the $n = 1$ case when $r = 1$ and $a_1 = 2$ or $r = 2$ and $a_1 = a_2 = 1$, so our work in this article is a generalization of her results. Her methods use properties specific to the Deligne category and thus do not generalize to our setting.

Using a categorical gluing construction described in Subsection \ref{subsec:restricted_vs_full}, we prove the uniqueness theorem for genuine MFTPCs (see Corollary \ref{cor:full_uniqueness} for a detailed statement).

\begin{introthm}
    If $\cC$ and $\cD$ are mixed admissible MFTPCs of the same type, then there is a strongly equivariant equivalence $\cC \simeq \cD$ intertwining the labels of simple objects.
\end{introthm}

We construct examples of restricted MFTPCs arising from classical representation theory using subquotients of parabolic categories $\cO$ for the usual general linear Lie algebra. Precisely, for any type $\Xi$ and integer $m$, we construct restricted categorifications $\cO^{\Xi}_m$ of level $m$ in terms of usual parabolic categories $\cO$. This step uses the combinatorics of linkage for parabolic categories $\cO$.

\begin{introthm}\label{intro_thm:construction}
    Any mixed admissible restricted MFTPC is equivalent to the Deligne tensor product of stable limits of subquotients of usual parabolic category $\cO$ for the general linear Lie algebra.
\end{introthm}



Finally, we apply this technology to study the complex rank categories $\ots$. By interpolating the classical Jantzen sum formula (Subsection \ref{subsec:jantzen_sum}), we produce a highest weight structure on the \textit{extended} category $\extots$ compatible with the natural categorical type A action given by translation functors. These data give $\extots$ the structure of an MFTPC under the assumption that $t_1,\ldots,t_n$ are non-integers. This condition is required to guarantee the existence of projectives in certain truncations of $\extots$.

As always, we require that this MFTPC is \textit{mixed admissible}. This amounts to condition that $(\bt,\bs)$ lies away from the intersections of some affine hyperplanes in $\fk^{2n}$, which precisely define the \textit{admissible} facets $W$ described earlier. Using Theorem \ref{intro_thm:construction}, we can produce an equivalence of $\extots$ with a limit of subquotients of classical parabolic categories $\cO$, making explicit the connection between complex rank representation theory and stable classical representation theory. 

\begin{introthm}\label{introthm:fruit}
If $\bt$ and $\bs$ are admissible parameters with $t_i \not \in \ZZ$ for any $i = 1,\ldots,n$, then the multiplicity $[\mts(\blambda):\lts(\bmu)]$, for any $2n$-multipartitions $\blambda,\bmu$ can be explicitly computed as a product of parabolic Kazhdan--Lusztig polynomials evaluated at $1$. In fact, the category $\extots$ is equivalent to a stable limit of the categories from Theorem \ref{intro_thm:construction}. 

{\textit{In particular, Conjecture \ref{introconj:structural} holds for admissible strata.}}
\end{introthm}

In the case of two blocks (i.e., when $n = 2$), the admissibility assumption is extraneous given the underlying assumption that $t_1,t_2$ are not integers. Thus, when $n = 2$ and $t_1,t_2 \not \in \ZZ$, we can completely describe the structure of $\extots$ in terms of stable representation theory. Our strategy also provides a heuristic for why we expect Conjecture \ref{introconj:structural} to hold in general. Namely, the type $\Xi$ of the MFTPC $\extots$ depends only on the facet that $(\bt,\bs)$ belongs to. Thus, an extension of Theorem \ref{intro_thm:uniqueness_restricted} to inadmissible MFTPCs would yield Conjecture \ref{introconj:structural} in general.

Finally, using ultraproduct constructions in positive characteristic, we derive results on multiplicities for parabolic analogs of the hyperalgebra Verma modules introduced by Haboush \cite{haboush} in the large rank and large characteristic limit. 

\vspace{8pt}

\noindent\textbf{Outlook.} At minimum, we view the validity of Conjecture \ref{introconj:structural} on admissible facets as evidence that the structure of the categories $\ots$ ought to be understood through the geometry of $\scr{X}_n$. To this end, we hope to further study the hyperplane arrangement introduced in this paper in the future. For example, our categorification approach indirectly yields equivalences between categories on the same admissible facets, but it is unclear whether these equivalences can be realized through ``translation functors'' defined intrinsically (i.e., without reference to categorical actions). We hope these hypothetical translation functors give an avenue toward proving Conjecture \ref{introconj:structural} in general.

\subsection*{Outline}

In Section \ref{sec:main_conjecture}, we introduce the complex rank parabolic categories $\ots$ as a family of categories on an affine parameter space $\scr{X}_n$ and illustrate, using ultraproducts, how to compute multiplicities at ``sufficiently generic'' points. We then introduce the stratification of $\scr{X}_n$ and our main Conjectures \ref{conj:main_conjecture} and \ref{conj:categorical_conjecture}. 

In Section \ref{sec:fock_tensor_prod_categorification}, we lay out technical preliminaries needed to prove Conjecture \ref{introconj:structural} in the admissible case. We recall the necessary background on highest weight categories, Fock space representations of $\fsl_{\ZZ}$, and categorical type A actions. We also axiomatize multi-Fock tensor product categorifications (MFTPCs) (Definition \ref{defn:multi_fock}) along with their restricted variants (Definition \ref{defn:restricted_mftpc}).

Section \ref{sec:uniqueness} is the technical centerpiece of this article. We prove the uniqueness theorem for restricted mixed admissible MFTPCs (Theorem \ref{thm:mixed_admissible_uniqueness}). We also study the relationship between restricted and full MFTPCs, resulting in a uniqueness statement of full MFTPCs (Corollary \ref{cor:full_uniqueness}).

In Section \ref{sec:construction}, we provide explicit constructions of restricted MFTPCs in terms of parabolic categories $\cO$ for the general linear Lie algebra (Theorem \ref{thm:construct_mftp}).

In Section \ref{sec:complex_rank_category_o}, we situate $\ots$ in the setting of MFTPCs. We interpolate Jantzen's formula for the determinant of the Shapovalov form on a parabolic Verma module to complex rank (Theorem \ref{thm:jantzen_determinant}) and derive a Jantzen sum formula (Theorem~\ref{thm:jantzen_sum}). We derive a linkage principle (Subsection~\ref{subsec:interpolated_linkage_order}) and define the extended category $\extots$ (Definition \ref{defn:extots}).

In Section \ref{sec:categorical_actions}, we use the MFTPC technology to produce concrete results about $\extots$. We study a natural categorical type A action on $\extots$ (Subsection \ref{subsec:categorical_action}) that gives it the structure of an MFTPC. When $(\bt,\bs)$ lies on an admissible facet, we produce an equivalence of $\extots$ with a stable limit of parabolic categories $\cO$.

We conclude our work in Section \ref{sec:applications} with applications to modular representation theory.

\subsection*{Notation and Conventions}

\begin{itemize}
    \item Throughout this paper, we let $\fk$ denote an algebraically closed field of characteristic zero. 
    \item We will let $\cP$ denote the set of all partitions, and for any integer $r$, we write $\cP^{r}$ for the set of $r$-multipartitions, that is, $r$-tuples of partitions.
    \item We represent partitions $\lambda \in \cP$ as Young diagrams in the ``French style.'' That is, lengths of rows decrease going up. Given a box with coordinates $(x,y)$ in a partition $\lambda$, we say that the content of this box is given by $x - y$. \textit{Warning: this convention is the opposite of the usual convention for content}. We will frequently conflate partitions with their Young diagrams.
    \item We write $|\lambda|$ for the number of boxes in any partition $\lambda \in \cP$. Moreover, for any $c \in \fk$, we write $|\lambda|_c$ for the number of boxes of content $c$ in $\lambda$. Finally, $\ell(\lambda)$ denotes the length of the partition, i.e., the number of rows in its Young diagram.
    \item We use \textbf{bold font} to represent sequences of numbers or partitions. For instance, we write $\blambda \in \cP^n$ for a multipartition $\blambda = (\lambda_1,\ldots,\lambda_n)$, where each $\lambda_i \in \cP$.
    \item We use \ul{\textbf{underline}} to denote sequences of multipartitions. For instance, we write $\ublambda \in \cP^{a_1} \times \cdots \times \cP^{a_r}$ to denote a sequence $(\blambda_1,\ldots,\blambda_r)$, where each $\blambda_\ell$ is an $a_\ell$-multipartition.
    \item An \textit{upper finite} poset is a poset $(\Lambda,\leq)$ where the set $\{\mu \in \Lambda \ |\ \mu \geq \lambda\}$ is finite for any $\lambda \in \Lambda$. Similarly, we say that $(\Lambda,\leq)$ is \textit{lower finite} if the set $\{\mu \in \Lambda \ |\ \mu \leq \lambda\}$ is finite for any $\lambda \in \Lambda$. Finally, an \textit{interval finite} poset is one where intervals $\{\mu \in \Lambda\ |\ \lambda_1 \leq \mu \leq \lambda_2\}$ are finite for any $\lambda_1,\lambda_2 \in \Lambda$. 
    \item Let $(\Lambda,\leq)$ be a poset. A \textit{poset ideal} (resp. \textit{poset coideal}) in $\Lambda$ is a subset $S \subset \Lambda$ such that $\mu \leq \lambda$ (resp. $\mu \geq \lambda$) and $\lambda \in S$ implies $\mu \in S$. 
    \item Let $\cC$ be an abelian category. Suppose $M \in \cC$ is an object and $L \in \cC$ is a simple object. For any finite filtration $\cF$ of $M$, we can consider the multiplicity $[M:L]_{\cF}$. In turn, we define the \textit{composition multiplicity} $[M:L]$ to be the supremum of $[M:L]_{\cF}$ across all finite filtrations $\cF$ of $M$.
    \item We write $A\modcat_{\cC}$ (resp. $\rmodcat{A}_{\cC}$) for the category of left (resp. right) $A$-modules in a symmetric tensor category $\cC$, where $A$ is an algebra object in $\cC$. We omit the subscript $\cC$ if the category $\cC$ is obvious from context.
    \item A finite abelian category is a $\fk$-linear abelian category equivalent to the category of left $A$-modules for some finite-dimensional $\fk$-algebra $A$. Given finite abelian categories $\cC := A_1\modcat$ and $\cD := A_2\modcat$, their Deligne tensor product is the finite abelian category $\cC \boxtimes \cD := (A_1 \ox_{\fk} A_2)\modcat$. 
    \item The \textit{additive Karoubi envelope} of a $\fk$-linear category $\cC$ is the $\fk$-linear category obtained by formally adjoining images of all idempotents to the additive completion of $\cC$.
    \item The Serre subcategory \textit{spanned} by a subset $S$ of simple objects in an abelian category $\cC$ is the Serre subcategory of $\cC$ consisting of object $X \in \cC$ such that $[X:L] \neq 0$ for a simple object $L \in \cC$ implies $L \in S$.
\end{itemize}

\subsection*{Acknowledgements}

I would like to thank Pavel Etingof for suggesting this problem to me, for many valuable conversations, and for helpful suggestions that greatly improved the exposition in this paper. I am also grateful to Ivan Losev for his invaluable guidance (in particular, introducing me to his work on type A categorifications) and for his generous sharing of unpublished notes. Finally, I want to thank Jon Brundan for a stimulating discussion of my work. This work was partially supported by NSF grant DMS-2001318.

\section{Parabolic Category \strO in Complex Rank and the Main Conjecture}\label{sec:main_conjecture}

In this section, we recall some background on Deligne categories and formally define the parabolic categories $\cO$ in complex rank. Viewing these categories as a family over an affine parameter space, we then compute multiplicities for ``sufficiently generic'' values of these parameters using stable representation theory. The precise details of this computation, using ultraproducts, are postponed to Appendix \ref{sec:ultraproducts}.  In turn, we introduce a hyperplane arrangement on the parameter space and formulate conjectures on the constancy of multiplicities along each facet.

\subsection{Background on Deligne Categories}

A \textit{symmetric tensor category} is a rigid symmetric monoidal category with a compatible structure of a $\fk$-linear abelian category, subject to the requirement that the monoidal bifunctor $\ox$ is additive and $\fk$-linear in both arguments. Moreover, we ask that the monoidal unit be simple, that all Hom spaces be finite-dimensional, and that all objects have finite length. We refer the reader to \cite{EGNO} for more details.

By abuse of notation, we will often omit the associativity and unit constraints as well as the braiding for a symmetric tensor category. Moreover, we will write $\ev_X: X^* \ox X \to \unit$ and $\coev_X: \unit \to X \ox X^*$ for the evaluation and coevaluation maps for the object $X$, dropping the subscript $X$ if it is clear from context.

Deligne's interpolating categories $\cD_t$ are certain symmetric tensor categories that cannot be realized as the symmetric tensor category of representations of any algebraic supergroup scheme. They are obtained by interpolating the representation categories of the general linear groups $\GL_m$. For a detailed exposition, we refer the reader to \cite{hu2024introductiondelignecategories}.  

\begin{defn} 
For any $t \in \fk$, the \textit{oriented Brauer category} $\cD_t^0$ is the rigid $\fk$-linear symmetric monoidal category generated by an object $\uparrow$ of dimension $t$. We write $\downarrow$ for the dual of $\uparrow$.
\end{defn}

For $r,s \geq 0$, the endomorphism algebra of $\uparrow^{\ox r} \ox \downarrow^{\ox s}$ in $\cD_t^0$ is the \textit{walled Brauer algebra} of dimensions $r,s$ and parameter $t$. Consider the Karoubi envelope $\cD_t$ obtained from $\cD_t^0$ by taking its additive completion and formally adjoining images of all idempotents. If $t \not \in \ZZ$, then Deligne and Milne \cite{deligne_milne} showed that the category $\cD_t$ is a semisimple symmetric tensor category. 

\begin{defn}
    We define the Deligne category as the category $\cD_t$ obtained above. We emphasize that $\cD_t$ is a symmetric tensor category only when $t \not \in \ZZ$. 
\end{defn}

\begin{rem}
When $t \in \ZZ$, the usual representation category $\ul{\mathrm{Rep}}(\GL_t)$ (for $t < 0$, the group $\GL_t$ should be taken to mean the supergroup $\GL(0|t)$) is the \textit{semisimplicification} of $\cD_t$, i.e., the quotient of $\cD_t$ by the tensor ideal of negligible morphisms. Moreover, Entova-Aizenbud, Hinich, and Serganova \cite{EAHS} showed that $\cD_t$ admits an abelian envelope $\cD_t^\mathrm{ab}$ that can be constructed as a limit of representation categories of the \textit{supergroups} $\GL(t+n|n)$ as $n \to \infty$. 
\end{rem}

The Deligne category satisfies the following universal property. For any $\fk$-linear symmetric tensor category $\cC$ and any object $V \in \cC$ of categorical dimension $t$, there exists a unique additive symmetric monoidal functor $\cD_t \to \cC$ sending $\uparrow$ to $V$.

We say that an indecomposable object in $\cD_t$ is \textit{primitive} if it is an indecomposable direct summand of $\uparrow^{\ox r} \ox \downarrow^{\ox s}$ for some $r,s \geq 0$. Building on work of Koike on primitive idempotents of the walled Brauer algebra for integer parameters, Comes and Wilson showed that primitive indecomposable objects in $\cD_t$ are classified by bipartitions of arbitrary size \cite{comes_wilson}. In particular, when $t \not \in \ZZ$, the simple objects in $\cD_t$ are indexed by bipartitions. 

\subsection{Lie Algebras in Symmetric Tensor Categories} Given any object $X \in \cC$, we may consider the object $\fgl(X) := X \ox X^*$, which is an associative algebra in $\cC$ with unit given by $\coev: \unit \to X \ox X^*$ and multiplication given by the composition \[m: (X \ox X^*) \ox (X \ox X^*) \simeq X \ox (X^* \ox X) \ox X^* \xto{\id_X \ox \ev \ox \id_{X^*}} X \ox X^*.\] Hence, we may define a Lie bracket $[-,-]: \fgl(X)\ox \fgl(X) \to \fgl(X)$ by the commutator 
\[
[-,-] := m -m \circ b_{\fgl(X),\fgl(X)} : \fgl(X) \ox \fgl(X) \to \fgl(X).
\]
This construction realizes $\fgl(X)$ as a \textit{Lie algebra object} in the category $\cC$. One may then define representations of $\fgl(X)$ as objects $M \in \cC$ equipped with a morphism $\rho: \fgl(X) \ox M \to M$ compatible with Lie bracket. 
More generally, these representations $M$ are allowed to be ind-objects in $\cC$, i.e., filtered colimits of objects of $\cC$. Morphisms between representations of $\fgl(X)$ are morphisms in $\Ind(\cC)$ that intertwine the action of $\fgl(X)$. Moreover, the tensor product of representations can be equipped with a $\fgl(X)$-action in the usual way. 

\subsection{Complex Rank Parabolic Category \texorpdfstring{$\cO$}{O}}\label{subsec:categoryO}

We now introduce parabolic category $\cO$ in complex rank first defined by Etingof \cite[Section 4]{Etingof_complex_rank_2}. This category depends on a \textit{generalized composition} of a complex number $t$, i.e., a sequence of complex numbers $\bt := (t_1,\ldots,t_n)$ satisfying $t = t_1 + \cdots + t_n$. We assume throughout that $t_1,\ldots,t_n \not \in \ZZ$, though $t$ is allowed to be an integer.

Objects in this category are certain ind-objects of the symmetric tensor category \[\cL_{\bt} := \cD_{t_1} \boxtimes \cD_{t_2} \boxtimes \cdots \boxtimes {\cD_{t_n}},\] which we may think of as the representation category for the ``Levi subgroup" of $GL_{t}$ associated with the composition $\bt$. For any object $Y \in \cD_{t_i}$, we abuse notation and write $Y \in \cL_{\bt}$ to refer to the object
$\unit^{\boxtimes (i-1)} \boxtimes Y \boxtimes \unit^{\boxtimes (n-i)}$. Then, if we let $V_{t_i}$ denote the basic object of $\cD_{t_i}$, we see that the object $V := V_{t_1} \oplus \cdots \oplus V_{t_n} \in \cL_{\bt}$ has categorical dimension $t$. Then, we can consider the general linear Lie algebra $\gl_t := V \ox V^*$ in $\cL_{\bt}$. We define the \textit{universal enveloping algebra} of $\gl_t$ as the ind-object
\[
U_{\bt} := T(\gl_t)/(x \ox y - y \ox x - [x,y]),
\]
where $T(X)$ denotes the tensor algebra of $X$. We can think of $U_{\bt}$-modules in $\Ind \cL_{\bt}$ as $\gl_t$-modules with a prescribed action of the \textit{Levi subalgebra} \[\fl_{\bt} := \gl_{t_1} \oplus \cdots \oplus \gl_{t_n} = \bigoplus_{i=1}^{n} V_{t_i} \boxtimes V_{t_i}^*.\]
The ordered decomposition of $V$ as a direct sum of subobjects $V_{t_i}$ yields a block triangular decomposition of $\fgl_{\bt}$. Namely, we first set
\[
\mathfrak{p}_{\bt} := \bigoplus_{i \leq j} V_{t_i} \boxtimes V_{t_j}^*.
\]
We should think of $\fp_{\bt}$ as a \textit{parabolic subalgebra} of $\fgl_t$. We also write 
\[
\mathfrak{u}^+_{\bt} := \bigoplus_{i < j} V_{t_i} \boxtimes V_{t_j}^*, \quad \mathfrak{u}^-_{\bt} := \bigoplus_{i > j} V_{t_i} \boxtimes V_{t_j}^*.
\]
Then, observe that we have the decomposition
\[
\fgl_t = \mathfrak{u}^+_{\bt} \oplus \fl_{\bt} \oplus \mathfrak{u}^-_{\bt}
\]
as objects of $\cL_{\bt}$.  Moreover, there is a splitting $\fl_{\bt} = \mathfrak{z}_{\bt} \oplus \fl_{\bt}^{ss}$ where $\fz_{\bt} := \Hom(\unit, \fl_{\bt})$ is a complex $n$-dimensional abelian Lie algebra and $\fl_{\bt}^{ss} := [\fl_{\bt},\fl_{\bt}]$ is semisimple. We can identify $\fz_{\bt}^*$ with $\fk^n$ and thus equip it with the usual type $A_{n-1}$ dominance order. 

We are now ready to introduce the parabolic category $\cO$ that will be central to our story. We start with an extended variant. The main category of interest will be a Serre subcategory of this extended category $\cO$.

\begin{defn}\label{defn:completion_ot}
    We define the \textit{completed parabolic category $\widehat{\cO}^{\bt}$} associated to $\bt$ as the full subcategory of $U_{\bt}$-modules $M$ in $\Ind \cL_{\bt}$ such that the following conditions hold.
    \begin{enumerate}[(O+1)]
        \item The action of $\fl_{\bt}^{ss}$ on $M$ coincides with its defining action as an object of $\Ind \cL_{\bt}$.
        \item The subalgebra $\fz_{\bt}$ acts semisimply.
        \item The parabolic subalgebra $\fp_{\bt}$ acts locally finitely on $M$.
    \end{enumerate}
\end{defn}

This definition allows weight subobjects in $\widehat{\cO}^{\bt}$ to be possibly ind-objects in $\cL_{\bt}$. In particular, the category $\widehat{\cO}^{\bt}$ is closed under all small coproducts. 
 
It is easy to see that $\widehat{\cO}^{\bt}$ is closed under taking subquotients and is thus an abelian subcategory of $U_{\bt}\modcat$. Following \cite[Section 4]{Etingof_complex_rank_2}, we now construct some basic objects in this category. Recall that the simple objects in $\cD_{t_i}$ are labelled by bipartitions. Then, for any $2n$-tuple of partitions $\blambda = (\lambda_1,\ldots,\lambda_{2n}) \in \cP^{2n}$, we adopt the notation \[X^{\bt}(\blambda) := X^{t_1}_1(\lambda_1, \lambda_{2}) \boxtimes \cdots \boxtimes X^{t_n}_n(\lambda_{2n-1}, \lambda_{2n}) \in \cL_{\bt},\] where $X^{t_i}_i(\lambda_{2i-1},\lambda_{2i})$ is the simple object in $\cD_{t_i}$ corresponding to the bipartition $(\lambda_{2i-1},\lambda_{2i})$. In particular, if $t_i$ is not an integer, then $X^{t_i}_i(\lambda_{2i-1},\lambda_{2i})$ is irreducible. For fixed $\bs = (s_1,\ldots,s_n) \in \fk^n$, we consider the $\fl_{\bt}$-module
\[
X^{\bt}_{\bs}(\blambda) := X^{\bt}(\blambda) \ox \unit_{\bs},
\]
where $\unit_{\bs}$ is the unit object of $\cL_{\bt}$ equipped with the action of $\fl_{\bt}$ where $\fz \simeq \fk^n$ acts by the character $\bs$ and $[\fl_{\bt},\fl_{\bt}]$ acts trivially.

By inflating $X^{\bt}_{\bs}(\blambda)$ to a module over $\fp$, we define the \textit{parabolic Verma module}
\[
M^t_{\bs}(\blambda) := U_{\bt} \ox_{U(\fp)} X^{\bt}_{\bs}(\blambda),
\]
which is a $U_{\bt}$-module in $\operatorname{Ind} \cL_{\bt}$. Observe that $M^t_{\bs}(\blambda)$ is indeed an object of $\widehat{\cO}^{\bt}$.

\begin{lem}
    The Verma module $M^{\bt}_{\bs}(\blambda)$ has a unique simple quotient $L^{\bt}_{\bs}(\blambda)$. Moreover, every simple object in $\widehat{\cO}^{\bt}$ (and hence in $\cO^{\bt}$) has the form $\lts(\blambda)$ for some $\blambda \in \cP^{2n}$ and $\bs \in \fk^n$. 
\end{lem}

\begin{proof}
    The sum of all proper $U_{\bt}$-submodules is the unique maximal proper submodule of $\mts(\blambda)$, which proves that $\mts(\blambda)$ has a unique simple quotient. Moreover, every simple object in $\widehat{\cO}^{\bt}$ is generated as a $U_{\bt}$-module by a simple $\fl_{\bt}$-submodule annihilated by $\f{u}^+_{\bt}$, proving the claim.
\end{proof}

Etingof posed the following problem regarding these simple modules. Strictly speaking, he asked for the \textit{characters} of the simple modules $\lts(\blambda)$, but as usual (it is easy to see that) this problem reduces to the following problem.

\begin{prob}{\cite[Problem 4.7]{Etingof_complex_rank_2}}\label{question:etingof}
    Compute the multiplicities $[\mts(\blambda):L_{\b{r}}^{\bt}(\b{\mu})]$ for any multipartitions $\blambda,\b{\mu}$ and any sequences $\bs,\b{r}$ of complex numbers. 
\end{prob}

To make the problem more amenable to our eventual categorification approach, it helps to restrict our attention to a fixed sequence $\bs$. 

\begin{defn}\label{defn:ots}
    For each $\bs \in \fk^n$, we define $\widehat{\cO}^{\bt}_{\bs}$ as the Serre subcategory consisting of objects $M$ in $\widehat{\cO}^{\bt}$ such that $[M:L_{\bs'}^{\bt}(\blambda)] \neq 0$ implies $\bs' = \bs$.
\end{defn}

\begin{lem}
    Fix $\blambda \in \cP^{2n}$ and $\bs \in \fk^n$. Suppose $X \in \widehat{\cO}^{\bt}$ is a simple object with $[\mts(\blambda):X] \neq 0$. Then, there exists $\bmu \in \cP^{2n}$ such that $X \simeq \lts(\bmu)$.  
\end{lem}

\begin{proof}

As an $\fl_{\bt}$-module, the PBW theorem for $U_{\bt}$ gives us an isomorphism 
\[
    \mts(\blambda) \simeq \Sym(\f{u}^-_{\bt}) \ox X_{\bs}^{\bt}(\blambda),
\]
and thus, every $\fl_{\bt}$-isotypic component of $\mts(\blambda)$ is a finite direct sum of objects of the form $X_{\bs}^{\bt}(\bmu)$ for some $\bmu \in \cP^{2n}$. It is not hard to see that this implies the lemma, since the objects $X_{\bs}^{\bt}(\bmu)$ and $X_{\bs'}^{\bt}(\bmu')$ are non-isomorphic unless $\bs = \bs'$ and $\bmu = \bmu'$.
\end{proof}

\begin{rem}
    Note that every object in $\widehat{\cO}^{\bt}$ has two actions of the central subalgebra $\fz_{\bt}$: one obtained by differentiating the tautological action of the fundamental group $L_{\bt} = \GL_{t_1} \times \cdots \times \GL_{t_n}$ of $\cL_{\bt}$ (this $\fz_{\bt}$-action has integer eigenvalues), and the other action coming from the defining structure of a $U_{\bt}$-module (this action has arbitrary eigenvalues). The difference between these actions on a simple $\fl_{\bt}$-module is precisely encoded by the character $\b{s}$. 
    
    Classically, the second action completely determines the isomorphism class of the underlying representation of the Levi subalgebra, since the Levi subgroup has a determinant character -- in particular, if the difference between the two actions is integral, then the representation was integrable to begin with (i.e., we may use the determinant character to relabel the actions so that they coincide). However, the complex rank group scheme $\GL_t$ does not have a determinant character, which means we must remember both actions. Concretely, the modules $X_{\bs}^{\bt}(\blambda)$ and $X_{\bs'}^{\bt}(\blambda')$ are non-isomorphic unless $\bs = \bs'$ and $\blambda = \blambda'$. 
\end{rem}

Thus, to address Problem \ref{question:etingof}, we may restrict our attention to the subcategories $\widehat{\cO}^{\bt}_{\bs}$, whose simple objects are labeled by $2n$-multipartitions. In fact, we shall restrict our attention to the following subcategory of $\widehat{\cO}^{\bt}_{\bs}$, which will play a crucial role throughout this paper. As we shall see in Section \ref{sec:complex_rank_category_o}, this subcategory is a union of \textit{Schurian} upper finite highest weight categories. 

\begin{defn}\label{defn:complex_rank_o}
    We define the \textit{parabolic category} $\ots$ associated to $\bt$ and $\bs$ as the Serre subcategory of the completed category $\widehat{\cO}^{\bt}_{\bs}$ consisting of objects with finite composition multiplicities.
\end{defn}

It is easy to see that the Verma modules $\mts(\blambda)$ must belong to $\ots$, as the $\fz_{\bt}^*$ weight subobjects in $\mts(\blambda)$ are all finite length (so the length of the $X_{\bs}^{\bt}(\bmu)$ isotypic component and thus the multiplicity of $\lts(\bmu)$ in $\mts(\blambda)$ are bounded).

\begin{rem}
    Another plausible (and perhaps more classical) definition for the parabolic category $\cO$ would be to consider the full subcategory of $U_{\bt}$-modules $M$ in $\Ind \cL_{\bt}$ such that (1) the action of $\fl^{ss}_{\bt}$ on $M$ coincides with its defining action as an object of $\Ind \cL_{\bt}$, (2) the subalgebra $\fz_{\bt}$ acts semisimply with weight objects belonging to $\cL_{\bt}$, and (3) all $\fz_{\bt}$-weights in $M$ are bounded above by a finite set (depending on $M$) of weights in $\fk^n$ (with respect to the dominance order). However, this definition fails to be ``large'' enough for our purposes; in particular, it will not be a union of Schurian categories (a property that will prove crucial for the highest weight structures we consider in this paper).
\end{rem}

\begin{rem}
    One could also require, as in \cite[Section 4]{Etingof_complex_rank_2}, that the objects in the category $\cO$ be finitely generated. However, it is not immediate obvious that this condition results in an abelian category. Moreover, dual Verma modules (defined in Section \ref{subsec:duality}) are not necessarily finitely generated in this setting -- any generating subobject for a dual Verma module must include the set of all ``singular subobjects'' in some Verma module $\mts(\blambda)$ (subobjects annihilated by $\fp_{\bt}$), and ultraproduct arguments (see Appendix \ref{sec:ultraproducts}) involving Verma's embedding theorem that there can be infinitely many of these singular subobjects in general. 

     A related question is whether finitely generated objects form a Serre subcategory of $\cO^{\bt}_{\bs}$, i.e., whether submodules of finitely generated modules are finitely generated themselves. For \textit{admissible parameters} $(\bt,\bs)$, the existence of projectives in certain truncations of $\ots$ established in Subsection \ref{subsubsec:projectives_truncations} show that this statement is true for finitely generated finite length objects of $\ots$ if and only if any submodule of any Verma module $\mts(\blambda)$ is also finitely generated.
\end{rem}


\subsection{Multiplicities at Interpolatable Points}\label{subsec:interpolatable}

Throughout the remainder of this section, we let $\fk = \CC$. We regard our categories $\ots$ as a family of $\CC$-linear abelian categories over the space $\scr{X}_n = \AA^n_{\CC} \times \AA^n_{\CC}$ parameterizing points \[\{(\bt,\bs) \in \CC^n \times \CC^n\}.\] We equip $\scr{X}_n$ with the diagonal action of the automorphism group $G := \mathrm{Aut}(\CC)$. For any point $(\bt,\bs) \in \scr{X}_n$, let $\ov{G(\bt,\bs)}$ denote the Zariski closure of its $G$-orbit $G(\bt,\bs)$.

We will use an ultraproduct construction to compute multiplicities in $\ots$ at sufficiently generic points in $\scr{X}_n$. First, we must specify what ``sufficiently generic'' means.

\begin{defn}
Fix some $(\bt,\bs) \in \scr{X}_n$. A sequence of points $(\mathbf{m}_j,\b{\varsigma}_j) \in \ov{G(\bt,\bs)} \cap (\ZZ^n \times {\ov{\QQ}}^n)$ is said to be \textit{$(\bt,\bs)$-interpolating} if the following condition holds.
\begin{itemize}
    \item For any non-principal ultrafilter $\scr{F}$ on $\NN$, a classical result of Steinitz guarantees the existence of an isomorphism $\phi_{\mathscr{F}}: \prod_{\mathscr{F}} \ov{\QQ} \to \CC$. We require that this isomorphism can be chosen so that $\phi_{\scr{F}}((\mathbf{m}_j, \b{\varsigma}_j)_{j \in \NN}) = (\bt,\bs)$. 
\end{itemize}
\end{defn}

\begin{defn}\label{defn:interpolatable} 
     A closed point $(\bt,\bs) \in \scr{X}_n$ is \textit{interpolatable} whenever there exists a $(\bt,\bs)$-interpolating sequence \[(\mathbf{m}_1, \b{\varsigma}_1),(\mathbf{m}_2,\b{\varsigma}_2),(\mathbf{m}_3,\b{\varsigma}_3),\ldots \in \ov{G(\bt,\bs)} \cap (\ZZ^n \times \ov{\QQ}^n).\] 
\end{defn}

\begin{rem}
    Equivalently, $(\bt,\bs) \in \scr{X}_n$ is interpolatable if the set \[\{(\mathbf{m}, \b{\varsigma}) \in \ov{G (\bt,\bs)}\ |\ \mathbf{m} \in \ZZ^n\}\] is Zariski dense in $\ov{G (\bt,\bs)}$.
\end{rem}

Let $(\bt,\bs) \in \scr{X}_n$ be an interpolatable point with $t_i \neq 0$ for each $i = 1,\ldots,n$ (this will be the only relevant case anyway), and pick a $(\bt,\bs)$-interpolating sequence \[(\mathbf{m}_j,\b{\varsigma}_j) \in \ov{G(\bt,\bs)} \cap (\ZZ^n \times \ov{\QQ}^n).\] 
For $i = 1,\ldots,n$, we have $t_i \neq 0$, so $\mathbf{m}_{j,i} \neq 0$ for infinitely many $j \geq 0$. Thus, by passing to a suitable subsequence, we may as well assume $\mathbf{m}_{j,i} \neq 0$ for each $i = 1,\ldots,n$ and $j \geq 0$. 

Now, for each point $\mathbf{m}_j = (\mathbf{m}_{j,1},\ldots,\mathbf{m}_{j,n})$, we consider the vector superspace
\begin{align}\label{eqn:super_decomposition}V(j) := \ov{\QQ}^{\mathbf{m}_{j,1}} \oplus \cdots \oplus \ov{\QQ}^{\mathbf{m}_{j,n}},\end{align} where we use the convention that $\ov{\QQ}^{m}$ is the purely odd superspace of dimension $|m|$ whenever $m < 0$. 
Let $\fgl(j)$ denote the general linear Lie superalgebra of endomorphisms of $V(j)$. Let $\fl(j)$ denote the Levi sub-superalgebra stabilizing the decomposition (\ref{eqn:super_decomposition}). In particular, we have
\[
    \fl(j) = \fgl_{m_{j,1}} \oplus \cdots \oplus \fgl_{m_{j,n}},
\]
where we use the convention that $\fgl_{m_{j,i}} := \fgl_{0|-m_{j,i}}$ whenever $m_{j,i} < 0$. Finally, the \textit{ordered} decomposition (\ref{eqn:super_decomposition}) gives an ascending filtration on $V(j)$, and we let $\fp(j) \subset \fgl(j)$ denote the parabolic sub-superalgebra corresponding to the decomposition \[\dim V(j) = \mathbf{m}_{j,1} + \ldots + \mathbf{m}_{j,n}.\] 

By choosing $j \gg 0$ (i.e., so the coordinates of $\mathbf{m}_j$ are arbitrarily large), we can regard a multipartition $\blambda \in \cP^{2n}$ as an $\fl(j)$-dominant weight $\theta_j(\blambda) = (\theta_j(\blambda)_1,\ldots,\theta_j(\blambda)_n)$, where \[\theta_j(\blambda)_k := (\lambda'_{2k-1,1},\ldots,\ \lambda'_{2k-1,\ell(\lambda'_{2k-1})},\ 0, \ldots, 0,\ -\lambda'_{2k,\ell(\lambda_{2k})},\ldots,\ -\lambda'_{2k,1}) + (\varsigma_{j,k},\ldots,\varsigma_{j,k}) \in \ov{\QQ}^{\mathbf{m}_{j,k}}. \] Above, for each $r = 1,\ldots,2n$, we write $\lambda'_{r}$ for the partition $\lambda_r$ if $\b{m}_{j,\lceil r/2 \rceil} > 0$ and the \textit{conjugate} partition $\lambda_r^\dagger$ if $\b{m}_{j,\lceil r/2 \rceil} < 0$. 
In turn, we may define the parabolic Verma supermodule $\Delta^{j}_{\bt,\bs}(\blambda)$ with highest weight $\theta_j(\blambda)$ and its simple quotient $L^j_{\bt,\bs}(\blambda)$. 

The following proposition is conceptually motivating but not particularly germane to the focus of our work in this paper. Thus, we postpone the proof, which requires some technical setup using ultraproducts, to Appendix \ref{sec:ultraproducts}. Roughly speaking, we can compute multiplicities at ``sufficiently generic" points in $\scr{X}_n$ using large rank stability phenomena. 

\begin{prop}\label{prop:ultraproduct_multiplicity}
    Continue to use the notation from the previous paragraph. For any $\blambda,\bmu \in \cP^{2n}$, the multiplicities $[\Delta^j(\blambda):L^j(\bmu)]$ stabilize as $j \to \infty$. In fact, we have 
    \[
        [M^{\bt}_{\bs}(\blambda): L^{\bt}_{\bs}(\bmu)] = \lim_{j \to \infty} [\Delta^j_{\bt,\bs}(\blambda):L^j_{\bt,\bs}(\bmu)].
    \]
    In particular, the right-hand side is independent of the $(\bt,\bs)$-interpolating sequence $(\mathbf{m}_j,\b{\varsigma}_j)$.
\end{prop}

Evidently, any interpolatable point has transcendental $t_i$ coordinates, so Proposition \ref{prop:ultraproduct_multiplicity} can yield multiplicities only at points whose $t_i$-coordinates are transcendental.

\begin{rem}
Section \ref{sec:applications} demonstrates how a positive characteristic adaptation of the ultraproduct procedure described in Appendix \ref{sec:ultraproducts} can, in principle, yield multiplicities at points with algebraic coordinates. However, we see this kind of approach as an end and not as a means (i.e., computing multiplicities in $\ots$ should tell us something about modular representation theory).
\end{rem}

\subsection{Stratification of Parameter Space and the Main Conjecture}\label{subsec:stratification}

For $1 \leq i < j \leq 2n$ and $k \in \ZZ$, consider the affine hyperplane $H_k(i,j) \subset \scr{X}_n$ whose closed points $(\bt,\bs)$ satisfy \[(t_1 + \cdots + t_{\lceil j/2 \rceil}) - (t_1 + \cdots + t_{\lceil i/2\rceil}) + s_{\lceil i/2\rceil} - s_{\lceil j/2\rceil} = k.\]
Let $\cH$ denote the set of all triples $(i,j,k)$ as above. For each subset $I \subset \cH$, we define \[V_I := \bigcap_{(i,j,k) \in I} H_k(i,j).\] If $V_I \neq \emptyset$, note that $I$ is necessarily finite. For any $I \subset \cH$ with $V_I \neq \emptyset$, we define 
\[
W_I := V_I \setminus \left(\bigcup_{(i,j,k) \in \cH \setminus I} H_k(i,j)\right) \subset V_I.
\]
We refer to each locally closed subset $W_I \subset \scr{X}_n$ as a \textit{facet}. Note that facets $W_I$ and $W_{I'}$ may coincide even when $I \neq I'$, but this is not an issue for our purposes.

\begin{defn}
    Fix a facet $W_I$. A closed point $(\bt,\bs) \in W_I$ is \textit{$W_I$-generic} if its $G$-orbit is Zariski dense in $W_I$. We say that $W_I$ is \textit{interpolatable} if all of its $W_I$-generic points are interpolatable.
\end{defn}

\begin{ex}
    Any facet $W_I$ where $|I| \leq 1$ is interpolatable. On the other hand, in the case $n = 2$, the facets $W_I$ with $I = \{(1,4,k), (2,3,k')\}$ are not interpolatable for any $k,k' \in \ZZ$.
\end{ex}

In particular, Proposition \ref{prop:ultraproduct_multiplicity} asserts that multiplicities for $W_I$-generic points on interpolatable facets can be computed via stable limits of parabolic Kazhdan-Lusztig polynomials. Generally, we conjecture that multiplicities on any facet $W_I$ can be computed through the multiplicities at any $W_I$-generic point. Let $W_I^{\natural} \subset W_I$ denote the subset of all points $(\bt,\bs) \in W_I$ such that $t_i \not \in \ZZ$ for any $i = 1,\ldots,n$, where $\bt = (t_1,\ldots,t_n)$.

\begin{conj}\label{conj:main_conjecture}
    Let $W_I$ be a facet with $\dim W_I \geq 1$. The multiplicity $[\mts(\blambda):\lts(\bmu)]$ does not depend on the choice of $(\bt,\bs) \in W_I^{\natural}$ and can be computed at any $W_I$-generic point $(\bt,\bs)$. In particular, if $W_I$ is interpolatable, then these multiplicities are computed via Proposition \ref{prop:ultraproduct_multiplicity}.
\end{conj}

Assuming Conjecture \ref{conj:main_conjecture}, taking any point $(\bt,\bs) \in W_I^{\natural}$ whose coordinates are algebraic numbers should yield, via model theory, stability results for multiplicities in modular representation theory. We make this connection precise in Section \ref{sec:applications}. 

In fact, we believe that Conjecture \ref{conj:main_conjecture} is a numerical ``shadow" of a deeper categorical statement. The precise statement of this conjecture requires an abelian subcategory $\extots$ of $\ots$, to be introduced in Definition \ref{defn:extots}. As we will see in Subsection \ref{subsec:extended_category_o} and Section \ref{sec:categorical_actions}, it is more natural to study $\extots$ from the perspective of categorification. In any case, the following conjecture implies Conjecture \ref{conj:main_conjecture}.

\begin{conj}\label{conj:categorical_conjecture}
    Let $W_I$ be a facet with $\dim W_I \geq 1$. The structure of the category $\extots$ does not depend on the choice of $(\bt,\bs) \in W_I^{\natural}$. More precisely, for any $(\bt,\bs),(\bt',\bs') \in W_I^{\natural}$, there exists an equivalence of categories $\extots \simeq (\cO^{\bt'}_{\bs'})^{\mathrm{res}}$ intertwining the labels of Verma modules.
\end{conj}

We will prove Conjecture \ref{conj:categorical_conjecture} in some special, yet conceptually illuminating, cases. 

\begin{defn}\label{defn:admissible_facets}
A facet $W_I$ is \textit{inadmissible} if there exist $(i_1,i_2,k_1), (i_2,i_3,k_2) \in I$ such that
\begin{itemize}
\item[(a)] $1 \leq i_1 < i_2 < i_3 \leq 2n$,
\item[(b)] $i_1 \equiv i_3 \pmod{2}$ and yet $i_1 \not \equiv i_2 \pmod{2}$.
\end{itemize}
Otherwise, we say that $W_I$ is \textit{admissible}.
\end{defn}

\begin{rem}
    The seemingly ad hoc definitions of facets and admissibility will be justified in Section \ref{sec:categorical_actions}. We shall see that these definitions correspond to some natural conditions on the structure of certain $\fsl_{\ZZ}$-categorifications arising from $\extots$.
\end{rem}

\begin{ex}
    Any facet $W_I$ with $|I| \leq 1$ is admissible. For a more complicated example, let $n = 3$. Then, the subset $I = \{(1,3,k_1),(3,6,k_2),(2,4,k_1'), (4,5,k_2')\}$ (for any $k_1,k_1',k_2,k_2' \in \ZZ$) gives rise to an admissible facet $W_I \subset \scr{X}_3$. On the other hand, for any $k,k' \in \ZZ$, the subset $I' = \{(1,4,k),(4,5,k')\}$ gives rise to an \textit{inadmissible} facet $W_{I'}$.
\end{ex}

In particular, we prove Conjecture \ref{conj:categorical_conjecture} when the facet $W_I$ is admissible. In the case $n = 2$, every facet $W_I$ with $W_I^{\natural} \neq \emptyset$ is admissible, so our work completely resolves Conjecture \ref{conj:categorical_conjecture}. 

Using the $\fsl_{\ZZ}$-categorification strategy outlined in the introduction, we actually show that each category $\extots$ on an admissible facet is equivalent, as highest weight categories, to an explicitly constructed stable limit of classical parabolic categories $\cO$. This equivalence gives us a completely concrete and purely algebraic description of the category $\extots$ that does not rely on model theoretic arguments. A key feature of our result is that it does not depend on whether the admissible facet $W_I$ is interpolatable. Thus, our results not only lead to a more explicit interpretation of the multiplicities $[\mts(\blambda):\lts(\bmu)]$ than the multiplicities computed in Proposition \ref{prop:ultraproduct_multiplicity}, but they also yield multiplicities on non-interpolatable facets, where this ultraproduct approach fails to say anything about multiplicities even at generic points.

Our categorification approach also provides a heuristic for why we expect that Conjecture \ref{conj:categorical_conjecture} should hold in general. Roughly speaking, our strategy hinges on the uniqueness of certain type A categorifications, which we call multi-Fock tensor product categorifications (MFTPCs), defined in a way that depends only on the set $I \subset \cH$. We can prove uniqueness of the MFTPCs defined from admissible facets. The main obstruction to proving Conjecture \ref{conj:categorical_conjecture} is that our uniqueness theorem fails to handle the case of inadmissible MFTPCs. Nonetheless, we expect that a uniqueness theorem should hold in the inadmissible case as well, which would yield Conjecture \ref{conj:categorical_conjecture} in general.

\section{Fock Space Tensor Product Categorifications}\label{sec:fock_tensor_prod_categorification}

In this section, we establish some preliminaries on highest weight categories and type A categorical actions. We also introduce the central definition in this paper, that of a multi-Fock tensor product categorification. 

\subsection{Semi-Infinite Standardly Stratified Categories}\label{subsec:standardly_stratified}

We give some background on highest weight categories with possibly infinite poset, following work of Brundan and Stroppel \cite{brundan_stroppel_semiinfinite}.

We work with sufficiently ``nice'' abelian categories where semi-infinite analogues of highest weight categories are well-behaved. An associative (not necessarily unital) $\fk$-algebra $A$ is \textit{locally unital} if it has a system of mutually orthogonal idempotents $\{e_i\}_{i \in I}$ such that $A = \bigoplus_{i,j \in I} e_iAe_j$. Moreover, $A$ is \textit{locally finite-dimensional} if each $e_iAe_j$ is finite-dimensional. A (left) module over a locally unital algebra $A$ is a (left) $A$-module $M$ satisfying $M = \bigoplus_{i \in I} e_iM$. We say that $M$ is \textit{locally finite-dimensional} if $\dim_{\fk} e_iM<\infty$ for all $i \in I$. We write $A\modlfd$ for the full subcategory of locally finite-dimensional $A$-modules. 
 
\begin{defn}
A $\fk$-linear abelian category $\cC$ is \textit{Schurian} if it is equivalent to the $A\modcat_{\lfd}$ for some locally finite-dimensional locally unital algebra $A$. A Schurian category with finitely many isomorphism classes of simple objects is called a \textit{finite} abelian category.
\end{defn}

Note that every simple object in a Schurian category has a projective cover. Moreover, a Schurian category is a Grothendieck abelian category and thus has enough injectives.

\begin{defn}
    A \textit{stratification} on a $\fk$-linear abelian category $\cC$ is a tuple $(\Theta,\Lambda,\leq,\rho)$, where $\Theta$ is an indexing set for the isomorphism classes of simple objects of $\cC$, $(\Lambda,\leq)$ is a poset, and $\rho: \Theta \to \Lambda$ is a map with finite fibers.
\end{defn}

For each $\lambda \in \Lambda$, let $\cC_{\leq \lambda}$ (resp. $\cC_{<\lambda}$) denote the Serre subcategory spanned by simples $L(\mu)$ for $\mu \leq \lambda$ (resp. $\mu < \lambda$). In turn, the \textit{stratum} corresponding to $\lambda$ is the Serre quotient $\cC_{\lambda} := \cC_{\leq\lambda}/\cC_{<\lambda}$, which is an abelian category with finitely many simple objects $L_\lambda(\theta)$ labelled by $\theta \in \rho^{-1}(\lambda)$. We will require the category $\cC$ to satisfy the following ``recollement'' condition.
\begin{equation} \tag{Q} \label{cond:Q}
\parbox{0.85\textwidth}{
The quotient functor $\pi^\lambda: \cC_{\leq \lambda} \to \cC_{\lambda}$ has a left adjoint $\pi^\lambda_!: \cC_{\lambda} \to \cC_{\leq\lambda}$ and a right adjoint $\pi^\lambda_*: \cC_\lambda \to \cC_{\leq \lambda}$. Moreover, $\cC_\lambda$ is a finite abelian category.
}
\end{equation}

\begin{rem}
    For instance, (\ref{cond:Q}) holds if $\cC_{\leq \lambda}$ is Schurian \cite[Lemma 2.24]{brundan_stroppel_semiinfinite}.
\end{rem}

If $\cC$ satisfies (\ref{cond:Q}), then we write $\Delta_\lambda:\cC_{\lambda} \to \cC$ and $\nabla_\lambda: \cC_{\lambda} \to \cC$ for the compositions of the left and right adjoints, respectively, of $\pi^\lambda$ with the inclusion $\cC_{\leq\lambda} \to \cC$. 

\begin{defn}
Define the \textit{associated graded category} $\operatorname{gr} \cC:= \bigoplus_{\lambda \in \Lambda} \cC_{\lambda}$ and the \textit{standardization and costandardization functors} \[\Delta := \bigoplus_{\lambda \in \Lambda} \Delta_\lambda: \operatorname{gr}\cC \to \cC, \quad \nabla := \bigoplus_{\lambda \in \Lambda} \nabla_\lambda: \operatorname{gr}\cC \to \cC,\] respectively. Letting $L_\lambda(\theta)$ denote the simple in $\cC_\lambda$ with label $\theta \in \rho^{-1}(\lambda)$ and $P_\lambda(\theta)$ denote its projective cover, define the \textit{standard} and \textit{proper standard} objects in $\cC$ to be 
\[
\Delta(\theta) := \Delta(P_\lambda(\theta)),\quad \ov{\Delta}(\theta) := \Delta(L_\lambda(\theta)).
\]
Similarly, we can define the costandard and proper costandard objects, respectively, by
\[
\nabla(\theta) := \nabla(I_\lambda(\theta)), \quad \ov{\nabla}(\theta) := \nabla(L_\lambda(\theta)),
\]
where $I_\lambda(\theta)$ is the injective hull of $L_\lambda(\theta)$.
\end{defn}

Note that $\Delta(\theta)$ (resp. $\nabla(\theta)$) has simple cosocle (resp. socle) $L(\theta)$. Moreover, there is a canonical surjection $\Delta(\theta) \twoheadrightarrow \ov{\Delta}(\theta)$ and a canonical inclusion $\ov{\nabla}(\theta) \hookrightarrow \nabla(\theta)$.

We will use the language of upper finite standardly stratified categories developed by Brundan and Stroppel \cite{brundan_stroppel_semiinfinite}. In general, the semi-infinite categories we encounter in this work have posets that are neither upper nor lower finite, but instead, are (roughly) products of lower finite posets with upper finite posets.

\begin{defn}\label{defn:semi_infinite}
Let $\cC$ be a $\fk$-linear abelian category with stratification $(\Theta,\Lambda,\leq,\rho)$. We say that $\cC$ is an \textit{upper finite fully stratified category} if the following axioms hold.
\begin{enumerate}[(UF1)]
    \item The poset $\Lambda$ is upper finite.
    \item The category $\cC$ is Schurian. In particular, it satisfies (\ref{cond:Q}) and has enough projectives. 
    \item {Let $P(\theta)$ denote the projective cover of $L(\theta)$. Then, the kernel of the epimorphism $P(\theta) \twoheadrightarrow \Delta(\theta)$ has a filtration by standard objects $\Delta(\theta')$ for $\rho(\theta') \geq \rho(\theta)$}. Similarly, the epimorphism $P(\theta) \twoheadrightarrow \ov{\Delta}(\theta)$ must also have a filtration by proper standard objects $\ov{\Delta}(\theta')$, where $\rho(\theta') \geq \rho(\theta)$
\end{enumerate}
If, additionally, $\Lambda$ is a \textit{finite} poset, then we say that $\cC$ is a \textit{finite fully stratified category}.
\end{defn}

Pick some upper finite ideal $J \subset \Lambda$ and $\theta \in J$. Write $\Delta_{J}(\theta)$ for the corresponding standard object in $\cC(J)$. Moreover, write $\iota_J: \cC(J) \to \cC$ for the inclusion functor. Observe that the object $\Delta(\theta):=\iota_J \Delta_{J}(\theta)$ does not depend on the choice of upper finite ideal $J$ containing $\rho(\theta)$.
The same conclusion also holds for costandard objects $\nabla_{J}(\theta)$ in $\cC$. In particular, we can define (proper) standard and costandard objects $\Delta(\theta)$ and $\nabla(\theta)$, respectively, in $\cC$ (without requiring the choice of an upper finite ideal $J$).

\begin{defn}
    Let $(\cC,\Theta,\Lambda,\leq,\rho)$ be an upper finite fully stratified category. If each stratum $\cC_\lambda$ is equivalent to the category of vector spaces, then we say $\cC$ is an \textit{upper finite highest weight category}. 
\end{defn}

For highest weight categories, the map $\rho$ is a bijection, so we may drop $\Theta$ and $\rho$ from our notation. For instance, if we write $\Delta(\lambda)$ for $\lambda \in \Lambda$, we are referring to $\Delta(\theta)$ for the unique $\theta \in \rho^{-1}(\lambda)$. Moreover, note that $\ov{\Delta}(\lambda) = \Delta(\lambda)$ and $\ov{\nabla}(\lambda) = \nabla(\lambda)$ for any $\lambda \in \Lambda$.

\subsubsection{General Properties of Standardly Stratified Categories}

We highlight some salient properties of standardly stratified categories, referring the reader to \cite[Section 3]{brundan_stroppel_semiinfinite} for details. Throughout this subsubsection, we let $(\cC,\Theta,\Lambda,\leq,\rho)$ denote an upper finite fully stratified category. 

By \cite[Theorem 3.33]{brundan_stroppel_semiinfinite}, the opposite category $\cC^{\mathrm{op}}$ is an upper finite fully stratified category with stratification $(\Theta,\Lambda,\leq,\rho)$ and \textit{standard objects} $\nabla(\theta)$.

Any poset ideal $J \subset \Lambda$ yields a finite fully stratified Serre subcategory $\cC(J)$ with stratification $(\rho^{-1}(J),J,\leq,\rho)$. We refer the reader to \cite[Theorem 3.41]{brundan_stroppel_semiinfinite} for an omnibus result on the structure of $\cC(J)$. For any poset \textit{coideal} $J \subset \Lambda$, we may consider the Serre span $\cC(\not \in J)$ of objects $L(\theta)$ with $\rho(\theta) \not \in J$. Then, the Serre quotient $\cC^J := \cC/\cC(\not \in J)$ is an upper finite fully stratified category with stratification $(\rho^{-1}(J),J,\leq,\rho)$. The Serre quotient functor $\pi: \cC \to \cC^J$ has left and right adjoints which are both fully faithful \cite[Lemma 2.24]{brundan_stroppel_semiinfinite}. In particular, they identify $\cC^J$ with some full subcategory of $\cC$. We refer the reader to \cite[Theorem 3.42]{brundan_stroppel_semiinfinite} for an omnibus theorem detailing the standard, costandard, projective, and injective objects in the quotient $\cC^J$. 

Objects admitting a filtration by standard or costandard objects play an important role in dualities between semi-infinite standardly stratified categories.

\begin{defn}
    For any upper finite fully stratified category $\cC$, we write $\Delta(\cC)$ for the full subcategory of objects with a finite $\Delta$-filtration. We similarly introduce the notation $\nabla(\cC)$.
\end{defn}

\begin{prop}\cite[Theorem 3.11]{brundan_stroppel_semiinfinite}\label{prop:homological_criterion}
    Suppose $\cC$ is a \textbf{finite} fully stratified category
    \begin{itemize}
    \item An object $V \in \cC$ belongs to $\Delta(\cC)$ if and only if $\Ext^1_{\cC}(V,\ov{\nabla}(\theta)) = 0$ for all $\theta \in \Theta$. In this case, the multiplicity of any $\Delta(\theta)$ in a standard filtration of $V$ is independent of the choice of filtration, as it can be computed by $\dim \Hom_{\cC}(V,\ov{\nabla}(\theta))$.
    
    \vspace{0.1cm}

    \item The same result holds with $\Delta(\cC)$ replaced by the full subcategory of objects with a proper standard filtration and with $\ov{\nabla}(\theta)$ replaced by $\nabla(\theta)$. 
    
    \vspace{0.1cm}

    \item An object $V \in \cC$ belongs to $\nabla(\cC)$ if and only if $\Ext^1_{\cC}(\Delta(\theta),V) = 0$ for all $\theta \in \Theta$. In this case, the multiplicity of any ${\nabla}(\theta)$ in a costandard filtration of $V$ is $\dim \Hom_{\cC}(\Delta(\theta),V)$. 
    \end{itemize}
\end{prop}

\subsubsection{Refining Standardly Stratified Structures}\label{subsubsec:compatible_stratified} We now discuss what happens when we change the stratification on a category $\cC$. For simplicity, we assume that $\Lambda$ is a finite poset.

\begin{defn}
    Suppose $\cC$ is a $\fk$-linear category with stratification $(\Theta,\Lambda,\leq,\rho)$. A \textit{refinement} of this stratification is a stratification $(\Theta,\Gamma,\preceq,\sigma)$ and an \textit{order-preserving surjective} map $q: \Gamma \twoheadrightarrow \Lambda$ with finite fibers such that $\rho = q \circ \sigma$ and so that $q(\beta) < q(\gamma)$ implies $\beta \prec \gamma$.
\end{defn}


Here is the most germane example for our purposes. Let $\Lambda_1$ and $\Lambda_2$ be (finite) posets with orders $\leq_1$ and $\leq_2$. We write $\preceq$ for the lexicographic order on $\Lambda_1 \times \Lambda_2$ and $\leq$ for the product order. Now suppose $\cC$ has stratification $(\Lambda_1 \times \Lambda_2, \Lambda_1, \leq_1, \rho)$, where $\rho: \Lambda_1 \times \Lambda_2 \to \Lambda_1$ is the projection. A stratification on $\cC$ of the form \[(\Lambda_1\times \Lambda_2, \Lambda_1 \times \Lambda_2, \preceq, \operatorname{id})\] along with the order-preserving surjection $\rho$ is thus a refinement of $(\Lambda_1 \times \Lambda_2, \Lambda_1, \leq_1, \rho)$. We record the following lemma for future use.

\begin{lem}{\cite[Lemma 3.66]{brundan_stroppel_semiinfinite}}\label{lem:refine_product_order}
If $\cC$ is a highest weight category with poset $(\Lambda_1 \times \Lambda_2, \leq)$ and each (co)standardization functor $\pi_{\lambda_1}^!, \pi_{\lambda_1}^*: \cC_{\lambda_1} \to \cC_{\leq_1 \lambda_1}$ (with respect to the coarser stratification) is exact for each $\lambda_1 \in \Lambda_1$, then $\cC$ is fully stratified with respect to the coarser stratification $(\Lambda_1 \times \Lambda_2, \Lambda_1, \leq_1, \rho)$.
\end{lem}

\subsubsection{Tilting Objects and Ringel Duality}\label{subsubsec:tilting_ringel}

We now discuss \textit{tilting objects} and \textit{Ringel duality} for \textit{finite} fully stratified categories, following \cite[Section 4]{brundan_stroppel_semiinfinite}.

\begin{defn}
    A \textit{tilting object} in a finite fully stratified category $\cC$ is an object \[T \in \Delta(\cC) \cap \ov{\nabla}(\cC).\] We write $\mathrm{Tilt}(\cC)$ for the full subcategory of tilting objects.
\end{defn}

It is well-known, see, e.g. \cite[Section 4]{brundan_stroppel_semiinfinite}, that
\begin{itemize}
    \item Any direct summand of a tilting object is a tilting object. Moreover, any tilting object decomposes uniquely as a direct sum of indecomposable tilting objects.
    \item There exists a bijection between $\theta \in \Theta$ and indecomposable tilting objects $T(\theta)$.
    \item There exists an injection $\Delta(\theta) \to T(\theta) $ whose cokernel is filtered by finitely many standard objects $\Delta(\theta')$ for $\rho(\theta') < \rho(\theta)$. Similarly, there is a surjection $T(\theta) \to {\nabla}(\theta)$ with kernel filtered by finitely many costandard objects ${\nabla}(\theta')$ for $\rho(\theta') < \rho(\theta)$.
\end{itemize}

A tilting generator for $\cC$ is an object $T \in \mathrm{Tilt}(\cC)$ for which every tilting object is isomorphic to a summand of $T$. The \textit{Ringel dual} category to $\cC$ relative to $T$ is \[\cC' := A^{\op}\modcat_{\mathrm{fd}}, \quad A := \End(T).\]
There is a \textit{Ringel duality functor} $\FF: \cC \to \cC'$ given by 
\[
M \mapsto \bigoplus_{i \in I} \Hom_{\cC}(T_i,M).
\]
In general, the Ringel dual category $\cC'$ is not fully stratified. Brundan and Stroppel \cite[Theorem 4.25]{brundan_stroppel_semiinfinite} showed that it is \textit{partially stratified} in the sense that the statement about proper standard filtrations in Definition \ref{defn:semi_infinite} still holds (but the statement about standard filtrations may no longer be true). In \cite{brundan_stroppel_semiinfinite}, partially stratified categories are called $(-)$-stratified. 

\begin{thm}\cite[Theorem 4.25]{brundan_stroppel_semiinfinite}
    The functor $\FF$ restricts to an equivalence $\nabla(\cC) \xto{\sim} \ov{\Delta}(\cC')$ and sends the indecomposable tilting $T(\theta) \in \cC$ to a projective cover of $L(\theta)$ in $\cC'$.
\end{thm}

\subsection{Partial Order on Multi-Partitions}\label{subsec:multipartition_poset}
 
We will now introduce and study a partial order on the set $\cP^{a}$ of $a$-multipartitions for some nonnegative integer $a$. This order depends on fixed sequences $\b{c} := (c_1,\ldots,c_a) \in \{0,1\}^a$ and $\b{\sigma} := (\sigma_1,\ldots,\sigma_a) \in \ZZ^a$. 

\begin{defn}\label{defn:inv_dominance}
    The \textit{inverse dominance order} $\leq$ of type $(\b{\sigma},\b{c})$ is defined as follows. Given multipartitions $\bds{\lambda} = (\lambda_1,\ldots,\lambda_a)$ and $\bds{\mu} = (\mu_1,\ldots,\mu_a)$, we declare $\bds{\lambda} \leq \bds{\mu}$ if, for each $N = 1,\ldots,a$ and each $b \in \ZZ$, we have
    \[
        \sum_{i=1}^{N} (-1)^{c_i}|\lambda_i|_{(-1)^{c_i}(b-\sigma_i)} \leq \sum_{i=1}^N (-1)^{c_i}|\mu_i|_{(-1)^{c_i}(b-\sigma_i),}
    \]
    with equality holding for $N = a$.
\end{defn}

For the following lemmas, we use the following notation: for $\blambda \in \cP^a$ and $k = 1,\ldots,a$, write \[\blambda(k) := \sum_{j=1}^k (-1)^{c_j} |\lambda_j|.\]

\begin{lem}
    The poset $(\cP^a, \leq)$ is interval-finite.
\end{lem}

\begin{proof}
    Pick $\b{\mu} \leq \b{\lambda}$ in $\cP^a$, and consider some $\b{\nu} =(\nu_1,\ldots,\nu_a)$ satisfying $\b{\mu} \leq \b{\nu} \leq \b{\lambda}$. We bound the number of boxes in $\nu_i$ for each $i = 1,\ldots,a$. To do so, we proceed by induction on $i$. If $c_i = 0$, then
    \[
    \b{\mu}(i) - \b{\nu}(i-1) \leq |\nu_i| \leq \b{\lambda}(i) - \b{\nu}(i-1).
    \]
    By induction, we have upper and lower bounds on $\b{\nu}(i-1)$, so we have a bound on $|\nu_i|$ as well. The case where $c_i = 1$ is identical.
\end{proof}

Under some combinatorial assumptions on the sequence $(c_1,\ldots,c_a)$, we can say even more about the finiteness properties of $(\cP^a,\leq)$. 

\begin{defn}\label{defn:admissible_sequence}
    A $01$-sequence $(c_1,\ldots,c_a)$ is \textit{lower admissible} (resp. \textit{upper admissible}) if there exists $p \in \{1,\ldots,a\}$ such that $c_i = 0$ (resp. $c_i = 1$) for all $i = 1,\ldots,p$ and $c_i = 1$ (resp. $c_i = 0$) for $i = p+1,\ldots,a$.
\end{defn}

\begin{lem}\label{lem:admissible_finite}
    Suppose $(c_1,\ldots,c_a)$ is lower (resp. upper) admissible. Then, the poset $(\cP^{a},\leq)$ is lower (resp. upper) finite. 
\end{lem}

\begin{proof}
    Let's prove the lower finiteness claim -- the proof of the upper finiteness claim is identical. Fix an $a$-multipartition $\bds{\lambda}$ and assume $\bds{\mu} \leq \bds{\lambda}$. Let $p$ be the swapping index. Then, the requirement $\bds{\lambda}(N) \geq \bds{\mu}(N)$ for each $N = 1,\ldots,p$ puts an upper bound on the total number of boxes appearing in the partitions $\mu_i$ for all $i = 1,\ldots,p$. On the other hand, the condition $\bds{\lambda}(a) = \bds{\mu}(a)$ puts an upper bound on the total number of boxes appearing in $\mu_i$ for all $i = p+1,\ldots,a$. 
\end{proof}

\begin{cor}\label{cor:restricted_partitions}
    Assume that $(c_1,\ldots,c_a)$ is lower (resp. upper) admissible, with $p$ being the swapping index. The subset $\cP^{a}(m_1,m_2)$ consisting of $a$-multipartitions $\bds{\lambda}$ such that 
    \[
    \sum_ {i=1}^p |\lambda_i| \leq m_1, \quad \sum_{i=p+1}^a |\lambda_i| \leq m_2
    \]
    is a poset ideal (resp. coideal) in $(\cP^a,\leq)$. When $m = m_1 = m_2$, we will write $\cP^a(m) := \cP^a(m,m)$. 
\end{cor}

\begin{proof}
 In the lower admissible case, the proof of Lemma \ref{lem:admissible_finite} shows that the condition $\b{\mu} \leq \b{\lambda}$ bounds the total number of boxes appearing in $\mu_1,\ldots,\mu_p$ by the total number of boxes appearing in $\lambda_1,\ldots,\lambda_p$, and similarly, the total number of boxes appearing in $\mu_{p+1},\ldots,\mu_a$ is bounded by the total number of boxes appearing in $\lambda_{p+1},\ldots,\lambda_a$. The upper admissible case is similar.
\end{proof}

\subsection{Fock Space Representations of \texorpdfstring{$\fsl_{\ZZ}$}{slZ}}\label{subsec:fock}

We let $\fsl_{\ZZ}$ denote the Lie algebra of finitary $\ZZ \times \ZZ$-matrices with entries in $\fk$. We write $M_{ij}$ for the elementary matrix with $1$ in row $i$ and column $j$ and zeros elsewhere. The Lie algebra $\fsl_{\ZZ}$ is a symmetriable Kac--Moody algebra with Chevalley generators $f_i := M_{i,i+1}$ and $e_i := M_{i+1,i}$ for each $i \in \ZZ$. We write $\fh_{\ZZ}$ for the Cartan subalgebra with basis given by $h_i := [e_i,f_i]$ for $i \in \ZZ$. Next, let $\varpi_i$ for $i \in \ZZ$ denote the fundamental weights of $\fsl_{\ZZ}$, i.e., the weights dual to the basis $h_i$. Finally, we define the simple roots $\alpha_i := \varpi_{i+1} + \varpi_{i-1} - 2\varpi_i$ for each $i \in \ZZ$. We can impose the usual dominance ordering on the weight lattice $\bigoplus_{s \in \ZZ} \ZZ \varpi_s$ by declaring $w_1 \geq w_2$ if $w_1 - w_2$ is a nonnegative integer linear combination of the simple roots $\alpha_i$.

\begin{defn}
    For any $s \in \ZZ$, we define the shifted Fock space $\cF_{s}$ as the irreducible highest weight representation of $\fsl_{\ZZ}$ with highest weight $\varpi_s$. We similarly define its dual $\cF_s^\vee$ as the irreducible lowest weight representation of $\fsl_{\ZZ}$ with lowest weight $-\varpi_s$. 
\end{defn}

Observe that $\cF_s$ has a weight basis $\ket{\lambda}$ indexed by partitions $\lambda \in \cP$ and 
\[
f_i\ket{\lambda} = \ket{\lambda + \square_{i-s}}, \quad e_i\ket{\lambda} = \ket{\lambda - \square_{i-s}}.
\]
In this notation, the highest weight vector is $\ket{\emptyset}$, and the weight of $\ket{\lambda}$ is given by
\[
\varpi_s - \sum_{b \in \ZZ} |\lambda|_{b-s}\alpha_b.
\]
Likewise, the dual Fock space $\cF_s^{\vee}$ has a weight basis $\bra{\lambda}$ indexed by partitions $\lambda \in \cP$ with
\[
f_i\bra{\lambda} = \bra{\lambda - \square_{-(i-s)}}, \quad e_i\bra{\lambda} = \bra{\lambda + \square_{-(i-s)}}.
\]
In this notation, the lowest weight vector is $\bra{\emptyset}$, and the weight of $\bra{\lambda}$ is given by
\[
-\varpi_s + \sum_{b \in \ZZ} |\lambda|_{-(b-s)}\alpha_b.
\]
Given sequences $(\sigma_1,\ldots,\sigma_a) \in \ZZ^a$ and $(c_1,\ldots,c_a) \in \{0,1\}^a$, consider the tensor product 
\[
\cF_{\sigma_1}^{c_1} \ox \cF_{\sigma_2}^{c_2} \ox \cdots \ox \cF_{\sigma_a}^{c_a}
\]
as a $\fsl_{\ZZ}$-representation, where \[\cF_{\sigma}^0 := \cF_\sigma, \quad \cF_{\sigma}^1 := \cF_\sigma^\vee.\] This tensor product has weight basis $\ket{\b{\lambda}} = \ket{\lambda_1} \ox \cdots \ox \ket{\lambda_a}$ indexed by multipartitions $\b{\lambda} \in \cP^a$. For each $i = 1,\ldots,a$, write $w(\b{\lambda})_i$ for the weight of $\ket{\lambda_1} \ox \cdots \ox \ket{\lambda_i}$ in \[\cF_{\sigma_1}^{c_1} \ox \cF_{\sigma_2}^{c_2} \ox \cdots \ox \cF_{\sigma_i}^{c_i}.\] Then, we can partially order the vectors $\ket{\b{\lambda}}$ by declaring $\ket{\b{\mu}} \leq_S \ket{\b{\lambda}}$ if $w(\b{\mu})_i \leq w(\b{\lambda})_i$ for all $i = 1,\ldots,a$ with equality holding for $i = a$. The induced order on $\cP^a$ is precisely the inverse dominance order from Definition \ref{defn:inv_dominance}.

\subsection{Categorical Actions of Type A}

We now recall the notion of a categorical type A action, following the seminal work of Chuang and Rouquier \cite{chuang_rouquier}. To state the definition, we first recall the definition of the degenerate affine Hecke algebra $\mathrm{dAH}_k$. This is the $\fk$-algebra with generators $x_1,\ldots,x_k$ and $t_1,\ldots,t_{k-1}$ satisfying the following relations:
\begin{itemize}
    \item the generators $x_1,\ldots,x_k$ freely generate a commutative polynomial subalgebra of $\mathrm{dAH}_k$,
    \item the generators $t_1,\ldots,t_{k-1}$ generate a subalgebra of $\mathrm{dAH}_k$ isomorphic to the group algebra of the symmetric group $\mathfrak{S}_k$ on $k$ letters,
    \item as a vector space over $\fk$, we have an isomorphism $\mathrm{dAH}_k \simeq \fk[x_1,\ldots,x_k] \ox_{\fk} \fk \mathfrak{S}_k$, where $\fk[x_1,\ldots,x_k]$ and $\fk\mathfrak{S}_k$ are the subalgebras defined in the first two relations,
    \item we have $t_ix_i = x_{i+1}t_i - 1$ and $t_{i}x_{i+1} = x_it_i + 1$; otherwise $t_ix_j = x_jt_i$. 
\end{itemize}

\begin{defn}\label{defn:categorical_type_a_action}
    Let $\cC$ be a $\fk$-linear abelian category. A \textit{categorical type A action} on $\cC$ is a tuple $(E,F,x,\tau)$, where $(E,F)$ is an adjoint pair of endofunctors of $\cC$, $x \in \End(F)$ and $\tau \in \End(F^2)$. 
    This data is required to satisfy the condition that $F$ is isomorphic to a left adjoint of $E$ and that there exists an algebra homomorphism $\dAH_k \to \End(F^k)$ given by \[x_i \mapsto F^{k-i}xF^{i-1}, \quad t_i \mapsto F^{k-i-1}\tau F^{i-1}.\] 
\end{defn}

The adjunction between $E$ and $F$ means that $x$ induces a natural transformation of $E$ as well. Decompose $F = \bigoplus_{\alpha \in I} F_\alpha$ and $E = \bigoplus_{\alpha \in I} E_\alpha$, where $I \subset \fk$ and $E_\alpha$ (resp. $F_\alpha$) is the generalized eigenspace of $x$ on $F$ (resp. $E$) with eigenvalue $\alpha \in I$. Our treatment of categorical type A actions differs slightly from that of Chuang and Rouquier in that our eigenvalues are allowed to be general complex numbers instead of integers. 

The functors $E_\alpha$ and $F_\alpha$ are necessarily biadjoint hence exact. Moreover, they induce an action of $\fsl_2$ on the complexified Grothendieck group $\CC \ox_{\ZZ} K_0(\cC)$, such that the endomorphisms $[E_{\alpha}]$ and $[F_{\alpha}]$ satisfy the Chevalley--Serre relations for the Lie algebra $\fsl_{I}$. This is the semisimple Lie algebra associated to the Dynkin diagram $D_I$ defined as follows. Let us write $I = \coprod_{\ell=1}^r I_\ell$, a partition into cosets modulo $\ZZ$. Then, $D_I$ is the union of the type A Dynkin diagrams with vertices indexed by $I_\ell$ for each $\ell =1,\ldots,r$.

\begin{defn}\label{defn:slI_categorification}
    Let $\cC$ be a $\fk$-linear abelian category equipped with a categorical type A action $(E,F,x,\tau)$. Suppose the eigenvalues of $x$ belong to $I \subset \CC$, so that the eigenfunctors $E_\alpha,F_\alpha$ give the action of $\fsl_I$ on $\CC \ox_{\ZZ} K_0(\cC)$. We say that $(\cC,E,F,x,\tau)$ is a \textit{$\fsl_I$-categorification} if
    \begin{itemize}
    \item the $\fsl_I$-module $\CC \ox_{\ZZ} K_0(\cC)$ is \textit{integrable}
    \item we can decompose $\cC = \prod_{\lambda \in \fh^*} \cC_\lambda$ so that $\CC \ox_{\ZZ} K_0(\cC_\lambda)$ is the $\lambda$-weight space of $\CC \ox_{\ZZ} K_0(\cC)$. 
    \end{itemize}
\end{defn}

We will also need a notion of equivalence to state the desired uniqueness results for our categorifications. Suppose $(\cC,E,F,x,\tau)$ and $(\cC',E',F',x',\tau')$ are categories equipped with a categorical type A action. A functor $G:\cC \to \cC'$ is said to be a \textit{strongly equivariant equivalence} if it is an equivalence of categories and there exists a natural isomorphism $\xi: F' \circ G \to G \circ F$ satisfying $\xi \circ x'G = Gx\circ \xi$ and $\xi F \circ F' \xi \circ \tau'G = G\tau \circ \xi F \circ F'\xi$.  

\begin{ex}\label{ex:categorical_type_a_glX}
For our purposes, the most important example of a categorical type A action will come from the category of representations of $\fgl(X)$ in a symmetric tensor category. This construction was already well-known in the classical setting, and Entova-Aizenbud extended it to the categorical setting in \cite[Section 6.2]{ea_deligne_categorical_actions}. Let $\cC$ be a symmetric tensor category and $X \in \cC$ any object. Suppose $\cD$ is a full subcategory of the category of $\fgl(X)$-representations in $\Ind \cC$ containing the tautological representation $X$ of $\fgl(X)$. Moreover, we assume that $\cD$ is closed under the functors $F: X \ox -$ and $E: X^* \ox -$. Then, we define $x \in \End(F)$ component-wise by setting $x_M: X \ox M \to X \ox M$ to be the endomorphism
\[
x_M := \frac{1}{2}(C_{X \ox M} - \id_X \ox C_M - C_X \ox M),
\]
where $C$ is the \textit{Casimir operator} defined as the composition
\[
C_W: W \xto{\coev_{\fgl(X)}\ox W} \fgl(X) \ox \fgl(X) \ox W \xto{\rho \circ \rho} W
\]
for any $\fgl(X)$-representation $W \in \Ind(\cC)$. On the other hand, we define $\tau \in \End(F^2)$ component-wise by setting $\tau_M: X \ox X \ox M \to X \ox X \ox M$ to be $b_{X,X} \ox \id_M$. By a straightforward computation \cite[Lemma 6.2.2]{ea_deligne_categorical_actions}, the data $(E,F,x,\tau)$ defines a categorical type A action on $\cD$. 

For a concrete example of this construction, let $\cC$ be the category of vector spaces, so that $\fgl(X) \simeq \fgl_n$ for some $n \geq 0$. Then, if we take $\cD$ to be the category $\cO$ for $\fgl_n$, the endomorphism $x_M$ described above is simply the tensor Casimir on $V \ox M$, i.e., the action of the element $\sum_{i,j} E_{ij} \ox E_{ji} \in U(\fgl_n)$. This example then recovers \cite[Section 7.4]{chuang_rouquier}.
\end{ex}

\subsection{Fock Space Tensor Product Categorifications}\label{subsec:ftp}

In this section, we consider categorifications of tensor products of highest and lowest weight Fock spaces, provided that they come in an ``admissible" order. In fact, we will encounter commuting actions of multiple copies of $\fsl_{\ZZ}$ such that the decategorified module is an external tensor product of tensor products of Fock space representations. Thus, we introduce categorifications of external tensor products of Fock space tensor products, which we refer to as \textit{multi-Fock tensor product categorifications}. Our definition is an explicit interpretation of \cite[Definition 3.2]{losev_webster}.

We begin by fixing a sequence \[\Xi := ((\b{\sigma}_1,\b{c}_1), (\b{\sigma}_2,\b{c}_2),\ldots, (\b{\sigma}_r,\b{c}_r)),\] where each $\sigma_\ell = (\sigma_{\ell,1},\ldots,\sigma_{\ell,a_\ell})$ is a sequence of integers and $\b{c}_\ell = (c_{\ell,1},\ldots,c_{\ell,a_\ell})$ is a 01-sequence. Without loss of generality, assume that $\b{c}_1,\ldots,\b{c}_n$ are \textit{lower admissible} and $\b{c}_{n+1},\ldots,\b{c}_r$ are \textit{upper admissible} for some $n \in \{1,\ldots,r\}$. 

\begin{defn}\label{defn:mixed_admissible}
We say that $\Xi$ is \textit{admissible} if each sequence $\b{c}_\ell$ is admissible. If all sequences $\b{c}_\ell$ are lower (resp. upper) admissible, then we say that $\Xi$ is \textit{lower (resp. upper) admissible}.
\end{defn}

\begin{defn}\label{defn:multi_fock}
Suppose $\cC$ is an abelian category with highest weight stratification $(\Lambda, \preceq)$. Moreover, assume $\cC$ carries a categorical type A action $(E,F,x,\tau)$. We say that $\cC$ is a \textit{multi-Fock tensor product categorification (MFTPC)}  of type $\Xi$ if the following axioms hold.
\begin{enumerate}[(MF1)]

\item The poset $(\Lambda, \preceq)$ can be identified with the product poset $\cP^{a_1} \times \cdots \times \cP^{a_r}$, where the order on each factor is the inverse dominance order of type $(\b{c}_\ell,\b{\sigma}_\ell)$ from Definition \ref{defn:inv_dominance}. Throughout the rest of this definition, we identify $\Lambda$ with $\cP^{a_1} \times \cdots \times \cP^{a_r}$.

\item For each $j\geq 0$, write $\cC^j$ for the Serre subcategory of $\cC$ corresponding to the poset ideal \[\Lambda^j := \cP^{a_1}(j) \times \cP^{a_2}(j) \times \cdots \times \cP^{a_n}(j) \times \cP^{a_{n+1}} \times \cdots \times \cP^{a_r} \subset \Lambda.\] We require that each $\cC^j$ is an upper-finite highest weight category with poset $\Lambda^j$. Moreover, we ask that $\cC$ is (equivalent to) the \textit{union} of its subcategories $\cC^j$, that is, every object in $\cC$ is isomorphic to an object belonging to some subcategory of the form $\cC^j$. 

\item We have a direct sum decomposition $F = \bigoplus_{i=1}^r F(\ell)$ (resp. $E = \bigoplus_{\ell=1}^r E(\ell)$) where the summands $F(\ell)$ (resp. $E(\ell)$) pairwise commute. Moreover, the endomorphisms $x$ and $\tau$ restrict to endomorphisms $x_\ell$ and $\tau_\ell$ of $F(\ell)$ and $F(\ell)^2$, respectively. 

\item The data $(E(\ell), F(\ell), x_\ell,\tau_\ell)$ defines a categorical type A action on $\cC$ such that the generalized eigenvalues of $x_\ell$ on $F(\ell)$ belong to $\ZZ$. We write $\eil{i}{\ell}$ and $\fil{i}{\ell}$ for the summands of $E(\ell)$ and $F(\ell)$, respectively, with eigenvalue $i \in \ZZ$. 

\item \label{cond:mftp4} For each $\ell = 1,\ldots,r$ and for any sequence $\ul{\b{\lambda}} := (\b{\lambda}_1,\ldots,\b{\lambda}_r)$, where $\b{\lambda}_j \in \cP^{a_j}$,  
\[\text{$\fil{i}{\ell}\Delta(\ul{\blambda})$ admits a filtration with sections $\Delta(\ublambda^F(\ell,j)[i])$ for $j = 1,\ldots,a_\ell$,}\] where $\ublambda^F(\ell,j)[i] = (\blambda_1,\ldots,\blambda_\ell^F(j)[i],\ldots,\blambda_r)$ and $\b{\lambda}_\ell^F(j)[i] = (\lambda_{\ell,1},\ldots,\lambda_{\ell,j}^F[i],\ldots,\lambda_{\ell,a_\ell})$, with $\lambda_{\ell,j}^F[i]$ obtained from $\lambda_{\ell,j}$ by
\vspace{1pt}
    \begin{itemize}
        \item[(0)] adding a box of content $i - \sigma_{\ell,j}$ to $\lambda_{\ell,j}$ if $c_{\ell,j} = 0$
        \item[(1)] removing a box of content $-i + \sigma_{\ell,j}$ from $\lambda_{\ell,j}$ if $c_{\ell,j} = 1$
    \end{itemize}
\vspace{3pt}
(if such a partition exists; otherwise, we agree that the corresponding section is zero). Similarly, \[\text{$\eil{i}{\ell}\Delta(\ul{\blambda})$ has a filtration with sections 
$\Delta(\ublambda^E(\ell,j)[i])$ for $j = 1,\ldots,a_\ell$,}\]where $\ublambda^E(\ell,j)[i] = (\blambda_1,\ldots,\blambda_\ell^E(j)[i],\ldots,\blambda_r)$ and $\b{\lambda}_\ell^E(j)[i] = (\lambda_{\ell,1},\ldots,\lambda_{\ell,j}^E[i],\ldots,\lambda_{\ell,a_\ell})$, with $\lambda_{\ell,j}^E[i]$ obtained from $\lambda_{\ell,j}$ by
\vspace{1pt}
    \begin{itemize}
        \item[(0)] \textit{removing} a box of content $i - \sigma_{\ell,j}$ from $\lambda_{\ell,j}$ if $c_{\ell,j} = 0$
        \item[(1)] \textit{adding} a box of content $-i + \sigma_{\ell,j}$ to $\lambda_{\ell,j}$ if $c_{\ell,j} = 1$
    \end{itemize}
\vspace{3pt}
\end{enumerate} 
\end{defn}

\begin{rem}
     Losev and Webster defined tensor product categorifications in the general setting of fully stratified categories. The only categorifications that we encounter are genuinely highest weight, so we do not lose too much by restricting to highest weight categories from the outset. 
\end{rem}

\begin{defn}
In the notation of Definition \ref{defn:multi_fock}, we say $\cC$ is \textit{mixed admissible} if $\Xi$ is admissible. If, additionally, $\Xi$ is lower (resp. upper) admissible, then we say $\cC$ is \textit{lower (resp. upper) admissible}.
\end{defn}

Since the inclusions $\cC^j \to \cC^{j+1}$ intertwine the labels of simple objects, they also intertwine the labels of standard objects. In particular, we can make sense of \textit{standard objects} $\Delta(\ublambda)$ in a MFTPC $\cC$ as well as the full subcategory $\Delta(\cC) \subset \cC$ of objects with a finite standard filtration.

We have an isomorphism of $\bigoplus_{\ell=1}^r \fsl_{\ZZ}$-modules 
\[
\CC \ox_{\ZZ} K_0(\Delta(\cC)) \xto{\sim} \bigboxtimes_{\ell=1}^r (\cF_{\sigma_{\ell,1}}^{c_{\ell,1}} \ox \cdots \ox \cF_{\sigma_{\ell,a_\ell}}^{c_\ell,a_\ell})
\]
(where the $\ell$th copy of $\fsl_{\ZZ}$ acts on the $\ell$th external tensor factor on the right-hand side) given by
\[
[\Delta(\b{\lambda}_1,\ldots,\b{\lambda}_r)] \mapsto \bigboxtimes_{\ell=1}^r (\ket{\lambda_{\ell,1}} \ox \cdots \ox \ket{\lambda_{\ell,a_\ell}}),
\]
for each $\b{\lambda}_\ell = (\lambda_{\ell,1},\ldots,\lambda_{\ell,a_\ell}) \in \cP^{a_\ell}$. In this sense, we see that $\cC$ categorifies the notion of an external tensor product of tensor products of Fock spaces. This also shows that $\Delta(\cC)$ categorifies an integrable $\fsl_{\ZZ}^{\oplus r}$-module. Moreover, note that $[\Delta(\ublambda)]$ and $[\Delta(\ubmu)]$ belong to the same weight space in $\CC \ox_{\ZZ} K_0(\Delta(\cC))$ if and only if 
\[
    \sum_{k=1}^{a_\ell} (-1)^{c_{\ell,k}}|\lambda_{\ell,k}|_{(-1)^{c_{\ell,k}}(b-\sigma_{\ell,k})} = \sum_{i=1}^{a_\ell} (-1)^{c_{\ell,k}}|\mu_{\ell,k}|_{(-1)^{c_{\ell,k}}(b-\sigma_{\ell,k})}.
\]
In particular, if $[\Delta(\ublambda)]$ and $[\Delta(\ubmu)]$ do not belong to the same weight space, then $\Delta(\ublambda)$ and $\Delta(\ubmu)$ are not linked in $\cC$.

\begin{lem}\label{lem:multifock_opposite}
    If $\cC$ can be equipped with the structure of a multi-Fock tensor product categorification of type $\Xi$, then the opposite category $\cC^{\op}$ can as well. In particular, the functors $\eil{i}{\ell}$ and $\fil{i}{\ell}$ preserve the subcategories of standardly and costandardly filtered objects.
\end{lem}

\begin{proof}
    We take the categorification data $(E,F,x,\tau)$ on $\cC^{\op}$ to be the same as the categorification data for $\cC$, so the only nontrivial axiom that we need to check is (MF5). This follows from Proposition \ref{prop:homological_criterion} and the fact that $E_i(\ell)$ and $F_i(\ell)$ are biadjoint.
\end{proof}

\subsection{Restricted Admissible Categorifications}\label{subsec:restricted_multi_fock}

Restricted categorical actions were originally introduced in \cite{losev_vv_conjecture} and subsequently used in \cite{elias_losev2017} to study categorifications of certain weight spaces in tensor products of highest weight Fock space representations of $\widehat{\fsl}_p$.

Before proceeding to the definition, we give some motivation for why we work with restricted categorifications. By definition, restricted categorifications are highest weight subquotients of full categorifications, so computing multiplicities of simple objects in standard objects can be understood locally by working in restricted categorifications. Second, restricted categorifications are genuinely finite highest weight categories, so they are easier to work with. Finally, to construct examples of full categorifications arising from ``classical" integer rank representation theory, we need to ``glue" together restricted categorifications, see Subsection \ref{subsec:restricted_vs_full}. In fact, we prove uniqueness for full categorifications by first proving them for restricted categorifications.

Let us start with notational conventions. Fix a \textit{mixed admissible} type
\[
\Xi := ((\b{\sigma}_1,\b{c}_1),\ldots,(\b{\sigma}_r,\b{c}_r)),
\]
where each $\b{c}_\ell$ is an \textit{admissible} 01-sequence and each $\b{\sigma}_\ell$ is an integer sequence. As before, we write $a_\ell$ for the length of $\b{c}_\ell$. We may assume without loss of generality that the sequences $\b{c}_1,\ldots,\b{c}_{n}$ are lower admissible and that the sequences $\b{c}_{n+1},\ldots,\b{c}_r$ are upper admissible for some $n \in \{0,\ldots,r\}$. 

Pick a positive integer $m$ and define \[(\Lambda,\leq) := \cP^{a_1}(m) \times \cP^{a_2}(m) \times \cdots \times \cP^{a_r}(m),\] where each factor $\cP^{a_\ell}(m) \subset \cP^{a_\ell}$ is equipped with the inverse dominance order of type $(\b{c}_\ell,\b{\sigma}_\ell)$ (see Definition \ref{defn:inv_dominance}), and $\Lambda$ is equipped with the product order. Finally, suppose we have a finite abelian category $\cC$ equipped with the structure of a highest weight category with poset $\Lambda$.

We will use the following notation for certain distinguished posets related to $\Lambda$:
\[
\Lambda^\bullet = \cP^{a_1}(m) \times \cdots \times \cP^{a_n}(m), \quad \Lambda^\circ = \cP^{a_{n+1}}(m) \times \cdots \times \cP^{a_r}(m).
\]
Moreover, for $j,k \in \{0,1,\ldots,m\}$, we will consider the following subsets of $\Lambda(m)$:
\begin{align*}
\Lambda^\bullet(j,k) &:= \left\{\ublambda \in \Lambda \ \bigg|\ \sum_{\substack{i=1 \\ c_{\ell,i} = 0}}^{a_\ell} |\lambda_{\ell,i}| \leq j, \quad \sum_{\substack{i=1 \\ c_{\ell,i} = 1}}^{a_\ell} |\lambda_{\ell,i}| \leq k, \quad \ell = 1,\ldots,n\right\}, 
\end{align*}
\begin{align*}
\Lambda^\circ(j,k) &:= \left\{\ublambda \in \Lambda \ \bigg|\ \sum_{\substack{i=1 \\ c_{\ell,i} = 0}}^{a_\ell} |\lambda_{\ell,i}| \leq j, \quad \sum_{\substack{i=1 \\ c_{\ell,i} = 1}}^{a_\ell} |\lambda_{\ell,i}| \leq k, \quad \ell = n+1,\ldots, r\right\},
\end{align*}
\begin{align*}
\Lambda(j,k) &:= \left\{\ublambda \in \Lambda \ \bigg|\ \sum_{\substack{i=1 \\ c_{\ell,i} = 0}}^{a_\ell} |\lambda_{\ell,i}| \leq j, \quad \sum_{\substack{i=1 \\ c_{\ell,i} = 1}}^{a_\ell} |\lambda_{\ell,i}| \leq k, \quad \ell = 1,\ldots,r\right\}.
\end{align*}
In turn, we define the following subquotients of $\cC$.
\begin{itemize}
\item[($\bullet$)] Write $\cC^\bullet(j,k)$ for the Serre subcategory of $\cC$ corresponding to the poset ideal $\Lambda^\bullet(j,k)$.
\item[($\circ$)] Write $\cC^\circ(j,k)$ for the Serre \textit{quotient} of $\cC$ corresponding to the poset coideal $\Lambda^\circ(j,k)$.
\item[($-$)] Write $\cC(j,k)$ for the quotient of $\cC^\bullet(j,k)$ given by the poset coideal $\Lambda(j,k) \subset \Lambda^\bullet(j,k)$.
\end{itemize}

As before, we use the notation $\pi_! = \pi(j,k)_!: \cC^\circ(j,k) \to \cC$ for the fully faithful left adjoint to the quotient functor $\pi = \pi(j,k): \cC \to \cC^\circ(j,k)$. 

\begin{notation}
    We adopt the following convention for labels in $\Lambda$. Namely, we will frequently conflate $\Lambda$ and the product poset $\Lambda^\bullet \times \Lambda^\circ$. For example, given $\ublambda = (\blambda_1,\ldots,\blambda_r) \in \Lambda$, we will use the notation $\ublambda^\bullet := (\blambda_1,\ldots,\blambda_n) \in \Lambda^\bullet$ and $\ublambda^\circ := (\blambda_{n+1},\ldots,\blambda_r) \in \Lambda^\circ$. Then, we interchangeably write $\ublambda \in \Lambda$ and $(\ublambda^\bullet,\ublambda^\circ) \in \Lambda^\bullet \times \Lambda^\circ = \Lambda$ (i.e., these both denote the same label). 
\end{notation}

\begin{defn}\label{defn:restricted_mftpc}
    The structure of a \textit{level $m$ restricted mixed admissible MFTPC of type $\Xi$} on $\cC$ consists of the following data.
    \begin{enumerate}[(a)]
        \item \label{data:a} A collection of functors, indexed by $i \in \ZZ$ and $\ell = 1,\ldots,n$, \[\eil{i}{\ell}: \cC^\bullet(m,m-1) \to \cC^\bullet(m-1,m), \quad \fil{i}{\ell}: \cC^\bullet(m-1,m) \to \cC^\bullet(m,m-1).\] We require that $\eil{i}{\ell}$ and $\fil{i}{\ell}$ restrict to functors 
        \[\eil{i}{\ell}: \cC^\bullet(j,k) \to \cC^\bullet(j-1,k+1), \quad \fil{i}{\ell}: \cC^\bullet(j-1,k+1) \to \cC^\bullet(j,k)\]
        for each $j = 1,\ldots,m$ and $k = 0,\ldots,m-1$. We define 
        \[F^\bullet := \bigoplus_{\ell=1}^n \bigoplus_{i \in \ZZ} \fil{i}{\ell}, \quad E^\bullet := \bigoplus_{\ell=1}^n \bigoplus_{i \in \ZZ} \eil{i}{\ell}.\] 
        
        \item \label{data:b} A collection of functors, indexed by $i \in \ZZ$ and $\ell = n+1,\ldots,r$, \begin{align*}\eil{i}{\ell}: \pi_!\cC^\circ(m-1,m) \to \pi_!\cC^\circ(m,m-1), \\[5pt] \fil{i}{\ell}: \pi_!\cC^\circ(m,m-1) \to \pi_!\cC^\circ(m-1,m).\end{align*} For $j = 0,\ldots,m-1$ and $k = 1,\ldots,m$, we require that $\eil{i}{\ell}, \fil{i}{\ell}$ restrict to functors
        \begin{align*}
        \eil{i}{\ell}: \pi_!\cC^\circ(j,k) \to \pi_!\cC^\circ(j+1,k-1), \\[5pt]
        \fil{i}{\ell} : \pi_!\cC^\circ(j+1,k-1) \to \pi_!\cC^\circ(j,k).
        \end{align*}
        Similarly define $F^\circ$ and $E^\circ$ as the sums of all $F_i(\ell)$ and $E_i(\ell)$, respectively.
        
        \item Natural transformations $x^\bullet \in \End(F^\bullet), \tau^\bullet \in \End((F^\bullet)^2)$.
        \item Natural transformations $x^\circ \in \End(F^\circ), \tau^\circ \in \End((F^\circ)^2)$. 
    \end{enumerate}

    This data is required to satisfy the following axioms.

\begin{enumerate}
    \item[(RMF1$^\bullet$)] For each $i \in \ZZ$, $\ell = 1,\ldots,n$ and all $j = 1,\ldots,m$ and $k = 0,\ldots,m-1$, the functors 
    \begin{align*}\eil{i}{\ell}: \cC^\bullet(j,k) \to \cC^\bullet(j-1,k+1) \\[5pt] \fil{i}{\ell}: \cC^\bullet(j-1,k+1) \to \cC^\bullet(j,k)\end{align*}
    are biadjoint. 
    \item[(RMF1$^\circ$)] Similarly, for each $i \in \ZZ$ and $\ell = n+1,\ldots,r$, we require biadjointness for
    \begin{align*}
        \eil{i}{\ell}: \pi_!\cC^{\circ}(j,k) \to \pi_!\cC^{\circ}(j+1,k-1), \\[5pt]
        \fil{i}{\ell} : \pi_!\cC^{\circ}(j+1,k-1) \to \pi_!\cC^{\circ}(j,k)
    \end{align*}

    \item[(RMF2$^\bullet$)]\label{axiom:RMF2_bullet} The endomorphism $x^\bullet \in \End(F^\bullet)$ acts on $F^\bullet_i(\ell)$ for $i \in \ZZ$ and $\ell = 1,\ldots,n$ with generalized eigenvalue $i \in \ZZ$. Moreover, for each $k \geq 0$, the assignments $x_i \mapsto (F^\bullet)^{k-i}x^\bullet(F^\bullet)^{i-1}$ and $t_i \mapsto (F^\bullet)^{k-i-1}\tau^\bullet (F^\bullet)^{i-1}$ give rise to a homomorphism $\dAH_k \to \End((F^\bullet)^k)$.

    \item[(RMF2$^\circ$)]\label{axiom:RMF2_circ} The same condition on $x^\circ$ and $\tau^\circ$, replacing the $\bullet$ appearing in (RMF2$^\bullet$) with $\circ$.
    
    \item[(RMF3)]\label{axiom:RMF3} For each $\ell = 1,\ldots,n$ and for any $\ublambda \in \Lambda^\bullet(m-1,m)$, or for each $\ell = n+1,\ldots,r$ and $\ublambda \in \Lambda^{\circ}(m,m-1)$,
    \[\text{$\fil{i}{\ell}\Delta(\ul{\blambda})$ admits a filtration with sections $\Delta(\ublambda^F(\ell,j)[i])$ for $j = 1,\ldots,a_\ell$,}\]
     where the notation $\ublambda^F(\ell,j)[i]$ was defined in \ref{cond:mftp4}. For $\ell \leq n$ and  $\ublambda \in \Lambda^\bullet(m,m-1)$ or for each $\ell = n+1,\ldots,r$ and $\ublambda \in \Lambda^{\circ}(m-1,m)$, \[\text{$\eil{i}{\ell}\Delta(\ul{\blambda})$ has a filtration with sections $\Delta(\ublambda^E(\ell,j)[i])$ for $j = 1,\ldots,a_\ell.$}\]

    \item[(RMF4)] The functors $\{\eil{i}{\ell}, \fil{i}{\ell}\}_{i \in \ZZ}$ pairwise commute with $\{\eil{i}{\ell'}, \fil{i}{\ell'}\}_{i \in \ZZ}$ when $\ell \neq \ell'$.
\end{enumerate}
\end{defn}

\begin{defn}
    If $n = 0$ (resp. $n = r$), then we say that the aforementioned structure is \textit{upper admissible} (resp. \textit{lower admissible}). In this case, the data labelled by $\bullet$ (resp. $\circ$) is extraneous and hence we will drop all such superscripts from our notation.
\end{defn}

We can make sense of a \textit{strongly equivariant equivalence} of restricted multi-Fock categorifications (and any of their subquotients $\cC(m_1,m_2)$) in a completely analogous way to strongly equivariant equivalences of type A categorical actions. 

\begin{defn}
To avoid confusion, we will use the terminology \textit{full multi-Fock categorification} to describe the ``un-restricted" categorifications that we introduced in Definition \ref{defn:multi_fock}.
\end{defn}

If $\cC$ is a restricted mixed admissible MFTPC of type $\Xi$, then there exists a natural identification of $\CC \ox_{\ZZ} K_0(\cC)$ with the subspace of the $\fsl_I \simeq \fsl_\ZZ^{\oplus r}$-module
\[
\bigboxtimes_{\ell=1}^r (\cF_{\sigma_{\ell,1}}^{c_{\ell,1}} \ox \cdots \ox \cF_{\sigma_{\ell,a_\ell}}^{c_{\ell,a_\ell}})
\]
spanned by the weight vectors of the form
\[
\ket{\ul{\blambda}} := \bigboxtimes_{\ell=1}^r \ket{\lambda_{\ell,1}}\ox \cdots \ox \ket{\lambda_{\ell,a_\ell}}
\]
for $\ul{\blambda} = (\blambda_1,\ldots,\blambda_r) \in \Lambda$. Namely, this identification associates any standard object $\Delta(\ublambda)$ with the corresponding weight vector. The actions of the endomorphisms $[\eil{i}{\ell}]$ and $[\fil{i}{\ell}]$ coincide with the actions of the Chevalley generators $e_i$ and $f_i$ for the $\ell$th factor of $\fsl_\ZZ$ in $\fsl_I \simeq \fsl_\ZZ^{\oplus r}$. 

\begin{lem}\label{lem:restricted_multifock_opposite}
    Suppose $\cC$ is a level $m$ restricted lower admissible MFTPC of type $\Xi$. The same categorification data gives $\cC^{\op}$ the structure of a level $m$ restricted lower admissible MFTPC of type $\Xi$, except with standard objects given by the costandard objects in $\cC$.
\end{lem}

\begin{proof}
Similar to Lemma \ref{lem:multifock_opposite}.
\end{proof}

\begin{lem}\label{lem:restricted_multifock_ringel}
    Suppose $\cC$ is a level $m$ restricted lower (resp. upper) admissible MFTPC of type $\Xi$, and let $\cC^\vee$ be its Ringel dual. Then, $\cC^\vee$ can be equipped with the structure of a level $m$ restricted upper (resp. lower) admissible MFTPC of type $\Xi^{\rev}$. 
\end{lem}

\begin{proof}
    We first consider the lower finite case. Ringel duality gives an equivalence $\Delta(\cC^\vee) \simeq \nabla(\cC)$, which allows us to define biadjoint endofunctors \[^\vee \eil{i}{\ell}: \Delta(\cC^\vee) \to \Delta(\cC^\vee), \quad {}^\vee \fil{i}{\ell}: \Delta(\cC^\vee) \to \Delta(\cC^\vee)\] corresponding to the endofunctors $\eil{i}{\ell}, \fil{i}{\ell}: \nabla(\cC) \to \nabla(\cC)$, respectively. We define \[^\vee E := \bigoplus_{i,\ell} {}^\vee \eil{i}{\ell}, \quad ^\vee F := \bigoplus_{i,\ell} {}^\vee \fil{i}{\ell}.\] 
    Note that $\cC^\vee$ is an \textit{upper finite} highest weight category with poset $\Lambda^\op$, so $\cC^\vee$ has enough projectives. Since any projective object in $\cC^\vee$ is standardly filtered and the functors $^\vee E$ and $^\vee F$ are biadjoint, we may extend them to endofunctors of $\cC^\vee$. This construction also gives us endomorphisms $x \in \End({}^\vee F)$ and $\tau \in \End({}^\vee F^2)$ satisfying the degenerate affine Hecke relations. Verifying the axioms of a MFTPC on $\cC^\vee$ is then straightforward.

    The upper finite case is handled similarly, using the injective objects in the truncated subcategories $\cC^\vee_{\leq \lambda}$ to extend the induced functors $^\vee E$ and $^\vee F$. 
\end{proof}

\section{Uniqueness of Admissible Fock Space Tensor Product Categorifications}\label{sec:uniqueness}

In this section, we establish the uniqueness of restricted mixed admissible multi-Fock tensor product categorifications (MFTPCs). Consequently, we study the relationship between full and restricted MFTPCs and prove a uniqueness result for full MFTPCs. We first establish some notational conventions. Throughout this section, we fix an admissible sequence \[\Xi = ((\b{\sigma}_1,\b{c}_1),\ldots,(\b{\sigma}_r,\b{c}_r)),\] and assume without loss of generality that the sequences $\b{c}_1,\ldots,\b{c}_n$ are lower admissible and that the sequences $\b{c}_{n+1},\ldots,\b{c}_r$ are upper admissible, where $n$ is some index between $0$ and $r$. As always, we write $a_\ell$ for the length of the sequence $\b{c}_\ell = (c_{\ell,1},\ldots,c_{\ell,a_\ell})$. 
 
For some positive integer $m$, we fix a highest weight category $\cC$ with poset \[\Lambda := \cP^{a_1}(m) \times \cdots \times \cP^{a_r}(m),\] where the poset structure on $\Lambda$ is the product order. Each factor $\cP^{a_\ell}(m)$ here is equipped with the inverse dominance order of type $(\b{c}_\ell,\b{\sigma}_\ell)$. 

We assume that $\cC$ is equipped with the structure of a level $m$ restricted mixed admissible multi-Fock tensor product categorification (MFTPC) of type $\Xi$. The reader should revisit Subsection \ref{subsec:restricted_multi_fock} to recall the related notation. Our inductive arguments in this section will require the notion of \textit{shape} for a mixed admissible MFTPC. 

\begin{notation}\label{notation:shape}
For each $\ell = 1,\ldots,r$, let us write $p_\ell$ for the number of indices $j$ for which $c_{\ell,j} = 0$. In turn, we define $q_\ell := a_\ell - p_\ell$ and \begin{align*}
&p^\bullet := p_1 + \cdots + p_n, \quad &q^\bullet := q_1 + \cdots + q_n \\[5pt] &p^\circ := p_{n+1} + \cdots + p_r, \quad &q^\circ := q_{n+1} + \cdots + q_r \\[5pt]
&p := p^\bullet + p^\circ, \quad &q := q^\bullet + q^\circ.
\end{align*}
We will say that the MFTPC $\cC$ has \textit{shape} $(p^\bullet,q^\bullet,p^\circ,q^\circ)$. 
\end{notation}

We first provide a high-level overview of our strategy. Given level $m$ restricted mixed admissible MFTPCs $\cC$ and $\cD$ of the same type $\Xi$, our goal is to produce a strongly equivariant equivalence $\cC(m,m) \simeq \cD(m,m)$ that intertwines the labels of simple objects. The core statement is a ``double centralizer'' theorem (Theorem \ref{thm:mixed_fully_faithful}) akin to \cite[Theorem 5.1]{losev_webster}. An essential prerequisite is the independence of the algebra $\End_{\cC}(\TT^\bullet)$ from the categorification $\cC$, where $\TT^\bullet := \bigoplus_{\alpha,\beta \leq m} E^\alpha F^\beta \Delta(\ubempty,\ubempty)$. We accomplish this step in Subsection \ref{subsec:heisenberg}, where we describe these endomorphism algebras in terms of cyclotomic quotients of the degenerate Heisenberg category. 

A key categorical splitting construction developed by Losev (explained in Subsection \ref{subsec:categorical_splitting}) will allow us to inductively ``peel away" Fock space factors from our categorifications and reduce our claim to Losev and Webster's double centralizer theorem. 

A major obstacle to this inductive approach is that $\TT^\bullet$ is not projective in general (unlike in \cite[Theorem 5.1]{losev_webster}, where it is in fact both projective and injective) and thus, the problem of reconstructing the category $\cC(m,m)$ from $\TT^\bullet$ is subtle. To circumvent this problem, we introduce a class of \textit{weak tilting objects} in $\cC$, which coincide with projective (resp. tilting) objects in the upper (resp. lower) admissible case. We discuss general properties of these weak tilting objects in Subsection \ref{subsec:structure}. In particular, we will show that these objects are compatible with the categorical splitting construction, allowing us to prove Theorem \ref{thm:mixed_fully_faithful}. The proof of the uniqueness theorem then amounts to matching labels of weak tilting objects.

\subsection{Categorical Splitting}\label{subsec:categorical_splitting}

A recurring theme in our uniqueness results is an inductive structure coming from subcategories or quotient categories of $\cC$ that inherit a ``truncated'' categorical action from $\cC$. The key tool here is the categorical splitting construction discovered by Losev in \cite{losev_categorical_splitting}. Strictly speaking, Losev works in the setting of ``full" categorifications, but his results easily generalize to our restricted setting. 

We use four variations on the categorical splitting construction, each relying on one of the following two subsets of $\Lambda$ (defined for some $h = 1,\ldots,r$):
\begin{align}\label{eqn:lambda_0}
  \ul{\Lambda}_{h} &:=\{\ublambda \in \Lambda \ |\ \lambda_{h,a_h} = \emptyset\} \subset \Lambda, \\[5pt]
  \label{eqn:lambda_1}
  \ov{\Lambda}_{h} &:= \{\ublambda \in \Lambda \ |\ \lambda_{h,1} = \emptyset\} \subset \Lambda.
\end{align}

The constructions coming from $\ul{\Lambda}_{h}$ (resp. $\ov{\Lambda}_{h}$) will be called \textit{tail} (resp. \textit{head}) splittings. 

On the other hand, if $h \leq n$ (resp. $h > n$), we refer to the splitting as \textit{lower admissible} (resp. \textit{upper admissible}). Lower (resp. upper) admissible splittings correspond to constructing truncated categorification data on Serre subcategories (resp. quotient categories) and will be particularly compatible with tilting (resp. projective) objects.

\begin{notation}\label{notation:add_partition}
  We will also make heavy use of the following notation throughout. First, assume $q_h \neq 0$ for some $h = 1,\ldots,r$. For any $\ublambda \in \ul{\Lambda}_h$ and any $\mu \in \cP$, we write \[\ublambda_\mu(h) := (\blambda_1,\ldots,\blambda_{h-1},(\blambda_h)_\mu,\blambda_{h+1},\ldots,\blambda_r),\] where $(\blambda_h)_\mu = (\lambda_{h,1}, \ldots, \lambda_{h,a_h-1},\mu)$. On the other hand, given $\ublambda\in \Lambda$, we write \[\ublambda^-(h) := (\blambda_1,\ldots,\blambda_{h-1},(\blambda_h)^-,\blambda_{h+1},\ldots,\blambda_r),\] where $(\blambda_h)^- = (\lambda_{h,1}, \ldots, \lambda_{h,a_h-1},\emptyset)$.

  Now, assume $p_h \neq 0$ for some $h = 1,\ldots,r$. Fixing some $\ublambda \in \ov{\Lambda}_h$ and $\mu \in \cP$, we write \[\prescript{}{\mu}{\ublambda}{}(h) := (\blambda_1,\ldots,\blambda_{h-1},\prescript{}{\mu}{(\blambda_h)},\blambda_{h+1},\ldots,\blambda_r),\] where $\prescript{}{\mu}{(\blambda_h)} = (\mu,\lambda_{h,2},\ldots,\lambda_{h,a_h})$. On the other hand, given $\ublambda \in \Lambda$, write \[\prescript{-}{}{\ublambda}{}(h) := (\blambda_1,\ldots,\blambda_{h-1},\prescript{-}{}{(\blambda_h)},\blambda_{h+1},\ldots,\blambda_r),\] where $\prescript{-}{}{(\blambda_h)} = (\emptyset,\lambda_{h,2},\ldots,\lambda_{h,a_h})$. 
\end{notation}

\subsubsection{Lower Admissible Tail Splitting}\label{subsubsec:lower_admissible_tail_splitting}

Suppose $q_h \neq 0$ for some $h \leq n$. In this case, $\ul{\Lambda}_h \subset \Lambda$ forms a poset ideal, so we can consider the corresponding Serre subcategory $\ul{\cC}_h \subset \cC$.
Write \[\ul{\Xi}_{h} := ((\b{\sigma}_1,\b{c}_1),(\b{\sigma}_2,\b{c}_2),\ldots, (\b{\sigma}_h',\b{c}_h'),\ldots,(\b{\sigma}_r, \b{c}_r)),\] where $\b{\sigma}_h' = ({\sigma}_{h,1},\ldots,{\sigma}_{h,a_h-1})$ and $\b{c}_h' = ({c}_{h,1},\ldots,{c}_{h,a_h-1})$.
Finally, for concision, write $s := -\sigma_{h,a_h}$.

We would like to equip $\ul{\cC}_h$ with a compatible restricted MFTPC of type $\ul{\Xi}_h$. Unfortunately, $E_{s}(h)$ does not restrict to a functor $\ul{\cC}^\bullet_h(m-1,m) \to \ul{\cC}^\bullet_h(m,m-1)$. On the other hand, by definition, all other categorification functors $\eil{i}{\ell}$ and $\fil{i}{\ell}$ do restrict to functors \[\uleil{i}{\ell}: \ul{\cC}^\bullet_h(m,m-1) \to \ul{\cC}^\bullet_h(m-1,m), \quad \ulfil{i}{\ell}: \ul{\cC}^\bullet_h(m-1,m) \to \ul{\cC}^\bullet_h(m,m-1)\] for $\ell \leq n$. 
Similarly, for $\ell > n$, all functors $\eil{i}{\ell}$ and $\fil{i}{\ell}$ restrict to functors \[\uleil{i}{\ell}: \ul{\cC}^\circ_h(m,m-1) \to \ul{\cC}^\circ_h(m-1,m), \quad \ulfil{i}{\ell}: \ul{\cC}^\circ_h(m-1,m) \to \ul{\cC}^\circ_h(m,m-1).\]  We apply the categorical splitting construction first discovered by Losev \cite{losev_categorical_splitting} to ``truncate'' the functor $\eil{s}{h}$ to a functor ${\ul{E}}_{s}(h): \ul{\cC}^\bullet_h(m,m-1) \to \ul{\cC}^\bullet_h(m-1,m)$. The rough idea is to
\begin{enumerate}[(1)]
  \item\label{list_item:equal_subcategory} consider the Serre subcategory $\cC_= \subset \cC$ corresponding to the poset ideal \[\Lambda_= := \{\ublambda \ |\ |\lambda_{h,a_h}| \leq 1\}\] and observe that $\eil{s}{h}$ and $\fil{s}{h}$ both preserve $\cC_=$ and $\ul{\cC}_h$ is a subcategory of $\cC_=$,  
  \item construct explicit equivalences $\ul{\cC}_h \simeq \cC_=/\ul{\cC}_h$,
  \item explicitly construct left and right adjoints to $\ulfil{s}{h}$ using adjoints to the inclusion $\ul{\cC}_h \to \cC_=$,
  \item prove that these adjoints are isomorphic to a third functor defined in terms of the aforementioned equivalences $\ul{\cC}_h \simeq \cC_=/\ul{\cC}_h$.
\end{enumerate}
We refer the reader to \cite[Section 5]{losev_categorical_splitting} for explicit details -- in Losev's terminology, the desired construction corresponds to ``freezing'' the box $(1,1)$ in the partition $\lambda_{h,a_h}$.

In particular, for any $\ublambda \in \ul{\Lambda}_h$, there is a canonical short exact sequence
\begin{align}\label{eqn:ses}
0 \to \Delta(\ublambda') \to \eil{s}{h}\Delta(\ublambda) \to \uleil{s}{h}\Delta(\ublambda) \to 0,
\end{align}
where $\ublambda' = (\blambda_1,\ldots,\blambda_h',\ldots,\blambda_r)$ and $\blambda_h' = (\lambda_{h,1},\ldots,\lambda_{h,a_h-1},\square)$. 
In turn, set \[\ul{E} := \bigoplus_{\ell=1}^r \bigoplus_{i \in \ZZ} \uleil{i}{\ell},\quad \ul{F} := \bigoplus_{\ell=1}^r \bigoplus_{i \in \ZZ} \ulfil{i}{\ell}.\] 

The natural transformations $x$ and $\tau$ descend to natural transformations $\ul{x}\in \End(\ul{F})$ and $\ul{\tau} \in \End(\ul{F}^2)$ that give actions of $\dAH_k$ on $\ul{F}^k$ for each $k \geq 0$. The following proposition is straightforward, thanks to (\ref{eqn:ses}).     

\begin{prop}\label{prop:multifock_truncation}
    The data $(\ul{E},\ul{F},\ul{x},\ul{\tau})$ realizes the structure of a level $m$ restricted MFTPC of type $\ul{\Xi}_h$. This categorification has shape $(p^\bullet,q^\bullet-1,p^\circ,q^\circ)$. 
\end{prop}

\subsubsection{Lower Admissible Head Splitting}\label{subsubsec:lower_admissible_head_splitting} Now suppose $p_h \neq 0$ for some $h \leq n$. In this case, $\ov{\Lambda}_h$ is a poset ideal in $\Lambda$ as well. We construct the corresponding Serre subcategory $\ov{\cC}_h \subset \cC$. We define 
\[\ov{\Xi}_{h} := ((\b{\sigma}_1,\b{c}_1),\ldots,(\b{\sigma}'_h,\b{c}_h'),\ldots, (\b{\sigma}_r, \b{c}_r)),\] where $\b{\sigma}_h' = ({\sigma}_{h,2},\ldots,{\sigma}_{h,a_h})$ and $\b{c}_h' = ({c}_{h,2},\ldots,{c}_{h,a_h})$.

Losev's categorical splitting construction, given by freezing the box $(1,1)$ in $\lambda_{h,1}$ this time, gives us truncated functors $\oveil{i}{\ell}$ and $\ovfil{i}{\ell}$ on $\ov{\cC}_h$. In this case, these truncated functors give rise to the structure of a level $m$ restricted MFTPC on $\ov{\cC}_h$ of type $\ov{\Xi}_h$ and shape $(p^\bullet-1,q^\bullet,p^\circ,q^\circ)$. Moreover, we have a natural isomorphism $\iota \circ \ov{F} \simeq F \circ \iota$. 

\subsubsection{Upper Admissible Tail Splitting}\label{subsubsec:upper_admissible_tail_splitting}

Suppose $p_h \neq 0$ for some $h > n$. In this case, $\ul{\Lambda}_h \subset \Lambda$ is a poset \textit{coideal}, so we consider the associated Serre quotient category $\ul{\cC}_h$. All categorification functors $\eil{i}{\ell}$ and $\fil{i}{\ell}$ naturally descend to biadjoint restricted categorification functors $\uleil{i}{\ell}$ and $\ulfil{i}{\ell}$ on the quotient $\ul{\cC}_h$ with the \textit{exception} of the functor $\eil{s}{h}$, where $s = -\sigma_{h,a_h}$. In this case, Losev truncates the functor $\eil{s}{h}$ to a (restricted) functor $\uleil{s}{h}$ on $\ul{\cC}_h$ as follows. 
\begin{enumerate}[(1)]
  \item We consider the Serre quotient $\cC_=$ corresponding to the poset \textit{coideal} $\Lambda_{=}$ from \ref{list_item:equal_subcategory}. It is preserved by $\eil{s}{h}$ and $\fil{s}{h}$ 
  \item There is a Serre subcategory $\cC_- \subset \cC_=$ corresponding to $\ublambda$ with $|\lambda_{h,a_h}| = 1$. Construct an explicit equivalence $\cC_- \simeq \ul{\cC}_h$. 
  \item Use the procedure from earlier to truncate $\eil{s}{h}$ to a (restricted) functor on $\cC_-$
  \item Use the equivalence $\cC_- \simeq \ul{\cC}_h$ to transfer the truncated functor to a functor $\uleil{s}{h}$  on $\ul{\cC}_h$.
\end{enumerate}
In this case, these truncated functors give rise to the structure of a level $m$ restricted MFTPC on $\ul{\cC}_h$ of type $\ul{\Xi}_h$ and shape $(p^\bullet,q^\bullet,p^\circ-1,q^\circ)$. 

\begin{lem}\label{lem:e_commute_quotient}
  Writing $\varphi_!: \ul{\cC}_h \to \cC$ for the left adjoint to the quotient functor $\cC \to \ul{\cC}_h$, there is a natural isomorphism $E \circ \varphi_! \simeq \varphi_! \circ \ul{E}$ when restricted to standardly filtered objects in $\ul{\cC}_h$.
\end{lem}

\begin{proof}
Completely analogous to \cite[Lemma 4.9]{losev_webster}. 
\end{proof}

\subsubsection{Upper Admissible Head Splitting}\label{subsubsec:upper_admissible_head_splitting}

Finally, suppose $q_h \neq 0$ for some $h > n$. In this case, $\ov{\Lambda}_h \subset \Lambda$ is a poset coideal, and we can consider the corresponding quotient category $\ov{\cC}_h$. The splitting construction in this case yields truncated functors $\oveil{i}{\ell}$ and $\ovfil{i}{\ell}$ giving $\ov{\cC}_h$ the structure of a level $m$ restricted MFTPC of type $\ov{\Xi}_h$ and shape $(p^\bullet,q^\bullet,p^\circ,q^\circ-1)$. 

\begin{lem}\label{lem:f_commute_quotient}
  Writing $\varphi_!: \ov{\cC}_h \to \cC$ for the left adjoint to the quotient functor $\cC \to \ov{\cC}_h$, there is a natural isomorphism $F \circ \varphi_! \simeq \varphi_! \circ \ov{F}$ when restricted to standardly filtered objects in $\ov{\cC}_h$.
\end{lem}

\subsection{Structure of Mixed Admissible Categorifications}\label{subsec:structure}

In a lower admissible categorification, i.e., when $n = r$, the object $\UU = \Delta(\ubempty,\ubempty)$ is a tilting object, while in an upper admissible categorification, $\UU$ is a projective object. In the general mixed admissible case, $\UU$ is neither tilting nor projective, but it behaves like a ``mixed'' tilting-projective object -- a \textit{weak tilting object}, which we define in this subsection. Weak tilting objects are well-behaved enough, e.g., they have no self-extensions and can be inductively constructed from the categorification functors, for us to generalize \cite[Theorem 5.1]{losev_webster} to our setting. These weak tilting objects are constructed as the tilting objects with respect to a coarser fully stratified structure on $\cC$, which we call the \textit{weak stratification}. The categorical splitting construction will feature heavily throughout this subsection -- we will exploit the fact that, when $p^\bullet + q^\bullet = 0$, the weak stratification is the same as the default highest weight structure on $\cC$. Before proceeding, we invite the reader to revisit Subsection \ref{subsec:restricted_multi_fock} for a reminder on notation.

\subsubsection{Weak Stratification}

We first produce a fully stratified structure on $\cC$ that is \textit{refined} by the default highest weight structure on $\cC$ in the sense of Subsubsection \ref{subsubsec:compatible_stratified}. This structure will be needed to introduce the aforementioned weak tilting objects in $\cC$. We remind the reader that we will often implicity identify elements $\ublambda \in \Lambda$ with pairs $(\ublambda^\bullet, \ublambda^\circ)$, where $\ublambda^\bullet \in \Lambda^\bullet$ and $\ublambda^\circ \in \Lambda^\circ$.

\begin{defn}
  The \textit{weak stratification} on $\cC$ is the stratification on $\cC$ given by the data \[(\Lambda, \Lambda^\bullet, \leq, q),\] where $q: \Lambda = \Lambda^\bullet \times \Lambda^\circ \to \Lambda^\bullet$ is the projection. For $\ublambda \in \Lambda$, we write 
  \[\Delta^\bullet(\ublambda) = \Delta^\bullet(\ublambda^\bullet, \ublambda^\circ)\]
  for the standard objects in the weak stratification. We write $\Delta^\bullet(\cC)$ for the full subcategory of $\Delta^\bullet$-filtered objects. We will refer to these objects as the \textit{weak standard} objects of $\cC$. We similarly define the \textit{weak proper costandard} objects $\ov{\nabla}^\bullet(\ublambda)$ and define $\ov{\nabla}^\bullet(\cC)$.
\end{defn}

In particular, the stratum $\cC_{\ublambda^\bullet}$ corresponding to $\ublambda^\bullet \in \Lambda^\bullet$ has the additional structure of a finite highest weight category with poset $\Lambda^\circ$. We will define the quotient categories $\cC_{\ublambda^\bullet}(j,k)$ of $\cC_{\ublambda^\bullet}$ analogously to the quotients $\cC^\circ(j,k)$ of $\cC$.

For each $\ell > n$ and $i \in \ZZ$, observe that the categorification functors $\eil{i}{\ell}$ and $\fil{i}{\ell}$ preserve the subcategories $\cC_{\leq \ublambda^\bullet}$ and $\cC_{< \ublambda^\bullet}$. Hence, they induce functors 
\begin{align*}
E_{i,\ublambda^\bullet}(\ell): \cC_{\ublambda^\bullet}(m-1,m) \to \cC_{\ublambda^\bullet}(m,m-1), \quad
F_{i,\ublambda^\bullet}(\ell): \cC_{\ublambda^\bullet}(m,m-1) \to \cC_{\ublambda^\bullet}(m-1,m).
\end{align*}
By definition, for $\ell > n$, there are functorial isomorphisms \[F_{i,\ublambda^\bullet}(\ell) \circ \pi^{\ublambda^\bullet} \simeq \pi^{\ublambda^\bullet} \circ \fil{i}{\ell}, \quad E_{i,\ublambda^\bullet}(\ell) \circ \pi^{\ublambda^\bullet} \simeq \pi^{\ublambda^\bullet} \circ \eil{i}{\ell}.\]
Then, the functors $E_{i,\ublambda^\bullet}(\ell)$ and $F_{i,\ublambda^\bullet}(\ell)$, along with the natural transformations induced by $x$ and $\tau$, equip each stratum $\cC_{\ublambda^\bullet}$ with the structure of a level $m$ restricted \textit{upper admissible} MFTPC of type $\Xi^\circ$, where \[
\Xi^{\circ} := ((\b{\sigma}_{n+1},\b{c}_{n+1}), \ldots, (\b{\sigma}_{r}, \b{c}_{r})).
\]
We now describe how to construct projective objects using the categorification functors.

\begin{defn}\label{defn:content_string}
Given $\mu \in \cP$, we say that its \textit{content string} is the sequence of integers $(\chi_1,\ldots,\chi_{|\mu|})$ determined by the contents of boxes of $\mu$, counted from left to right and bottom to top. In particular, $\chi_1 = 0, \chi_2 = -1,\ldots, \chi_{\mu_1} = 1-\mu_1$, $\chi_{\mu_1 + 1} = 1$, $\chi_{\mu_1 + 2} = 0, \ldots, \chi_{\mu_1 + \mu_2} = 1 - \mu_2$, etc. 

When used in the context of \textit{head splitting}, we define the composite functors
\begin{align*}
F^{(\ell)}(\mu) &:= \fil{\chi_{|\mu|} + \sigma_{\ell,1}}{\ell}\cdots \fil{\chi_2 + \sigma_{\ell,1}}{\ell}\fil{\chi_1 + \sigma_{\ell,1}}{\ell}, \\[5pt] E^{(\ell)}(\mu) &:= \eil{\chi_{|\mu|} + \sigma_{\ell,1}}{\ell}\cdots \eil{\chi_2 + \sigma_{\ell,1}}{\ell}\eil{\chi_1 + \sigma_{\ell,1}}{\ell}.
\end{align*}
When used in the context of \textit{tail splitting}, the indices $\sigma_{\ell,1}$ above should be replaced with $\sigma_{\ell,a_\ell}$.
\end{defn}

The reader should revisit Subsection \ref{subsec:categorical_splitting} to recall the notation in the following lemma, whose proof is straightforward and thus omitted for brevity.

\begin{lem}\label{lem:produce_proj}
  Suppose $n = 0$ and $q_h \neq 0$ for some $h = 1,\ldots,r$. Pick $\ublambda \in \ov{\Lambda}_h$ such that \[\sum_{\ell=1}^r \sum_{\substack{j=1 \\ c_{\ell,j} = 1}}^{a_\ell} |\lambda_{\ell,j}| \leq m - |\mu|.\] For any $\mu \in \cP$, the projective object $P(\ublambda_{\mu}(h))$ is the unique direct summand in $E^{(h)}(\mu)P(\ublambda)$ with minimum label. In fact, all other labels $\ubnu$ occurring in the decomposition satisfy $|\nu_{h,a_h}| < |\mu|$. 

  Now suppose $p_h \neq 0$ and fix $\mu \in \cP$. Pick $\ublambda \in \ul{\Lambda}_h$ such that \[\sum_{\ell=1}^r \sum_{\substack{j=1 \\ c_{\ell,j} = 0}}^{a_\ell} |\lambda_{\ell,j}| \leq m - |\mu|.\] For any $\mu \in \cP$, the object $P(\prescript{}{\mu}{\ublambda}(h))$ is the unique direct summand in $F^{(h)}(\mu)P(\ublambda)$ with minimum label. All other labels $\ubnu$ satisfy $|\ubnu_{h,1}| < |\mu|$. 
\end{lem}

We are now ready to show that the weak stratification equips $\cC$ with the structure of a fully stratified category. The idea is to use Lemma \ref{lem:refine_product_order} along with the inductive construction of projective objects in $\cC_{\leq \ublambda^\bullet}$ through upper admissible categorical splitting. 

\begin{lem}\label{lem:standard_exact}
  For each $\ublambda^\bullet \in \Lambda^\bullet$, the standardization and costandardization functors $\cC_{\ublambda^\bullet} \to \cC$ are both exact. Thus, the category $\cC$ is a fully stratified category with respect to the weak stratification.
\end{lem}

\begin{proof}
  The second claim follows from the first and Lemma \ref{lem:refine_product_order}. For brevity, write $\pi$ for the quotient functor $\cC_{\leq \ublambda^{\bullet}} \to \cC_{\ublambda^\bullet}$ and $\pi_!,\pi_*: \cC_{\ublambda^\bullet} \to \cC_{\leq \ublambda^\bullet}$ for the standardization and costandardization functors, respectively. We prove that $\pi_*$ is exact. The exactness of $\pi_!$ is established in the same manner, swapping projectives for injectives and standards for costandards, and thus omitted -- crucially, one needs to use Lemma \ref{lem:restricted_multifock_opposite} for the action of $\eil{i}{\ell}$ and $\fil{i}{\ell}$ on costandard objects.

  It suffices to show that $\pi(P(\ubnu,\ubomega))$ is projective in $\cC_{\ublambda^\bullet}$ for {every} indecomposable projective object $P(\ubnu,\ubomega) \in \cC_{\leq \ublambda^\bullet}$, where $\ubnu \leq \ublambda^{\bullet}$ and $\ubomega \in \Lambda^\circ$. We induct on $p^\circ$ and $q^\circ$. First suppose $\ubomega = \ubempty$. Observe that the subcategory $\cC_{\leq \ublambda^\bullet}$ inherits the structure of a finite highest weight category from the highest weight stratification on $\cC$. In particular, $P(\ubnu,\ubempty)$ admits a filtration by \textit{highest weight} standard objects $\Delta(\ul{\bmu})$ with $\ul{\bmu} \geq \ubnu$. Since $\ubempty$ is a maximal label in $\Lambda^\circ$, each of these labels has the form $\ul{\bmu} = (\ul{\bmu}^\bullet,\ubempty)$ where $\ubnu \leq \ul{\bmu}^\bullet \leq \ublambda^\bullet$. The exact functor $\pi$ sends $\Delta(\ul{\bmu})$ to zero when $\ul{\bmu}^\bullet < \ublambda^\bullet$. On the other hand, note that $\Delta(\ublambda^\bullet,\ubempty) = \pi_!(\Delta_{\ublambda^\bullet}(\ubempty))$, where $\Delta_{\ublambda^\bullet}(\ubempty)$ is the standard object with label $\ubempty$ in the stratum $\cC_{\ublambda^\bullet}$. Hence, the object $\pi(P(\ubnu^\bullet,\ubempty))$ is filtered by \[\pi(\Delta(\ublambda^\bullet,\ubempty)) = \pi \circ \pi_!(\Delta_{\ublambda^\bullet}(\ubempty)) = \Delta_{\ublambda^\bullet}(\ubempty)\] (using $\pi \circ \pi_! \simeq \id$). Since $\ubempty\in \Lambda^\circ$ is maximal, the sections $\Delta_{\ublambda^\bullet}(\ubempty)$ are projective in $\cC_{\ublambda^\bullet}$, as desired.

  Now say $q_h \neq 0$ for $h > n$. Write $\ubomega = (\b{\omega}_{n+1},\b{\omega}_{n+2},\ldots,\b{\omega}_{r})$ and $\ubomega' = (\b{\omega}_{n+1},\ldots,\b{\omega}_{h}',\ldots,\b{\omega}_{r})$, where $\b{\omega}_{h}' = (\emptyset,\omega_{h,2},\ldots,\omega_{h,a_{h}})$. The upper admissible head splitting construction of Subsubsection \ref{subsubsec:upper_admissible_head_splitting} equips the Serre quotient $\ov{\cC}_h$ with the structure of a mixed admissible MFTPC of shape $(p^\bullet,q^\bullet,p^\circ,q^\circ - 1)$. We can also define the weak stratification on $\ov{\cC}_h$ and consider the quotient functor $j: \ov{\cC}_{h,\leq\ublambda^\bullet} \to \ov{\cC}_{h,\ublambda^\bullet}$. By induction, we may assume that $j\ov{P}_h(\ubnu,\ubomega')$ is projective in $\ul{\cC}_{h,\ublambda^\bullet}$, where we write $\ov{P}_h(\ubnu,\ubomega')$ for the indecomposable projective with label $(\ubnu,\ubomega')$ in $\ov{\cC}_{h,\leq\ublambda^\bullet}$.
  
  Let us write $\varphi: \cC \to \ov{\cC}_h$ for the quotient functor with left adjoint $\varphi_!: \ov{\cC}_h \to \cC$. Recall that $\cC_{\ublambda^\bullet}$ has the natural structure of a highest weight category with poset $\Lambda^\circ$. Then, observe that  $\ov{\cC}_{h,\ublambda^\bullet}$ can be identified with the quotient of $\cC_{\ublambda^\bullet}$ corresponding to the poset coideal \[\ov{\Lambda}_h^\circ = \{(\emptyset,\omega_{n+2},\ldots,\omega_{r}) \in \Lambda^\circ\} \subset \Lambda^\circ.\] In particular, if $\varphi': \cC_{\ublambda^\bullet} \to \ov{\cC}_{h,\ublambda^\bullet}$ denotes the quotient functor with left adjoint $\varphi'_!$, the following diagram commutes up to natural isomorphism:
\[\begin{tikzcd}
	{\Delta(\ul{\cC}_{h,\leq \ublambda^\bullet})} && {\Delta(\ul{\cC}_{h,\ublambda^\bullet})} \\
	\\
	{\Delta(\cC_{\leq \ublambda^\bullet})} && {\Delta(\cC_{\ublambda^\bullet})}
	\arrow["j"', from=1-1, to=1-3]
	\arrow["{\varphi_!}", from=1-1, to=3-1]
	\arrow["{\varphi'_!}"', from=1-3, to=3-3]
	\arrow["\pi", from=3-1, to=3-3]
\end{tikzcd}\]
  In particular, since $\varphi_!'$ sends projectives in $\ul{\cC}_{h,\ublambda^\bullet}$ to projectives in $\cC_{\ublambda^\bullet}$, the object 
  \[\pi P(\ubnu,\ubomega') = \pi \varphi_!\ul{P}_h(\ubnu,\ubomega') = \varphi_!'j\ul{P}_h(\ubnu,\ubomega')\]
   is projective in $\cC_{\ublambda^\bullet}$. The same argument as the one we gave for Lemma \ref{lem:produce_proj} shows that $P(\ubnu,\ubomega)$ is a direct summand of $E^{(h)}(\omega_{h,1})P(\ubnu,\ubomega')$. Since $h > n$, we see that
  \[
  \pi(E^{(h)}(\omega_{h,1})P_h(\ubnu,\ubomega')) \simeq E_{\ublambda^\bullet}^{(h)}(\omega_{h,1})\pi(P_h(\ubnu,\ubomega')),
  \]
  which is indeed projective, as the functors $E_{i,\ublambda^\bullet}(\ell)$ send projectives to projectives. We've reduced the claim to the case $q^\circ = 0$. The case $p^\circ \neq 0$ and $q^\circ =0$ is handled similarly, using the \textit{tail} splitting construction and the inductive construction of projectives using the functors $F_i(h)$. 
\end{proof}

We now describe the actions of the categorification functors $\eil{i}{\ell}, \fil{i}{\ell}$ on \[\Delta^\bullet(\ublambda^\bullet, \ublambda^\circ) = \pi^{\ublambda^\bullet}_{!}(P_{\ublambda^\bullet}(\ublambda^\circ)),\] where, as a reminder, $P_{\ublambda^\bullet}(\ublambda^\circ)$ is a projective cover of $L_{\ublambda^\bullet}(\ublambda^\circ) = \pi^{\ublambda^\bullet}(L(\ublambda))$ in the stratum $\cC_{\ublambda^\bullet}$. 

First consider the case $\ell > n$. The object $F_{i,\ublambda^\bullet}(\ell)P_{\ublambda^\bullet}(\ublambda^\circ)$ is a direct sum of indecomposable projectives in $\cC_{\ublambda^\bullet}$. Write $\cL^{F,i,\ell}_{\Xi,\ublambda^\bullet}(\ublambda^\circ) \subset \Lambda_m^\circ$ for the multi-set of labels appearing in this decomposition: \[\fil{i}{\ell}\Delta^\bullet(\ublambda^\bullet,\ublambda^\circ) = \pi_{!}^{\ublambda^\bullet}(F_{i,\ublambda^\bullet}(\ell)P_{\ublambda^\bullet}(\ublambda^\circ)) = \bigoplus_{\ubnu \in \cL^{F,i,\ell}_{\Xi,\ublambda^\bullet}(\ublambda^\circ)} \Delta^\bullet(\ublambda^\bullet,\ubnu),\]
 where we used the isomorphism $\fil{i}{\ell} \circ \pi_!^{\ublambda^\bullet} \simeq \pi_!^{\ublambda^\bullet} \circ F_{i,\ublambda^\bullet}(\ell)$ (a formal consequence of the fact that $\pi^{\ublambda^\bullet} \circ \eil{i}{\ell} \simeq E_{i,\ublambda^\bullet}(\ell) \circ \pi^{\ublambda^\bullet}$). We will similarly define the multi-set $\cL^{E,i,\ell}_{\Xi,\ublambda^\bullet}(\ublambda^\circ)$ so that 
\[\eil{i}{\ell}\Delta^\bullet(\ublambda^\bullet,\ublambda^\circ) = \bigoplus_{\ubnu \in \cL^{E,i,\ell}_{\Xi,\ublambda^\bullet}(\ublambda^\circ)} \Delta^\bullet(\ublambda^\bullet, \ubnu).\] 
Understanding the action of the functors $\eil{i}{\ell}$ and $\fil{i}{\ell}$ for $\ell \leq n$ takes more work.

\begin{lem}\label{lem:bullet_action}
  For any $\ublambda \in \Lambda^\bullet(m-1,m)$, the object $\fil{i}{\ell}\Delta^\bullet(\ublambda)$ admits a $\Delta^\bullet$-filtration with sections $\Delta^\bullet(\ublambda^F(j,\ell)[i])$ for $j = 1,\ldots,a_\ell$, where the labels $\ublambda^F(j,\ell)[i]$ were defined in \ref{cond:mftp4}. Similarly, for any $\ublambda \in \Lambda^\bullet(m,m-1)$, the object $\eil{i}{\ell}\Delta^\bullet(\ublambda)$ has a filtration with sections $\Delta^\bullet(\ublambda^E(j,\ell)[i])$.
\end{lem}

\begin{proof}
  We will prove the claim for $\fil{i}{\ell}$, since the proof for $\eil{i}{\ell}$ is completely analogous. The proof hinges on the fact that \[\Delta^\bullet(\ublambda^\bullet,\ubempty) = \pi^{\ublambda^\bullet}_!(P_{\ublambda^{\bullet}}(\ubempty)) = \pi^{\ublambda^\bullet}_!(\Delta_{\ublambda^\bullet}(\ubempty)) = \Delta(\ublambda^\bullet,\ubempty).\] In particular, Axiom \ref{axiom:RMF3} implies that $\fil{i}{\ell}\Delta^\bullet(\ublambda^\bullet,\ubempty)$ has the desired description for any $\ublambda^\bullet$.

  We proceed by induction on $p^\circ$ and $q^\circ$, having just established the base case $p^\circ = q^\circ = 0$. The idea is to use the inductive construction of $P_{\ublambda^\bullet}(\ublambda^\circ)$ from $P_{\ublambda^\bullet}(\ubempty)$. 

  We first reduce the claim to the case $q^\circ = 0$ by induction. We may assume $q_h \neq 0$ for some $h > n$. The upper admissible head splitting construction from Subsubsection \ref{subsubsec:upper_admissible_head_splitting} equips the corresponding quotient category $\ov{\cC}_h$ with the structure of an MFTPC of shape $(p^\bullet, q^\bullet, p^\circ, q^\circ - 1)$. The left adjoint $\varphi_!$ to the quotient functor $\cC \to \ul{\cC}_h$ commutes with both functors $\eil{i}{\ell}$ and $\fil{i}{\ell}$ for $\ell \leq n$ and sends exact sequences of standardly filtered objects in $\ul{\cC}_h$ to exact sequences of standardly filtered objects in $\cC$. Note that $\varphi_!$ also sends $\Delta^\bullet$-filtered objects in $\ul{\cC}_h$ to $\Delta^\bullet$-filtered objects in $\cC$ and intertwines their labels. Hence, by induction, we may assume that the claim holds in $\cC$ for all $\ublambda^\circ \in \Lambda^\circ$ where $\lambda_{h,1} = \emptyset$. 
  We now induct on $|\lambda_{h,1}|$. Given $\ublambda^\circ \in \Lambda^\circ$, write \[(\ublambda^\circ)' := (\blambda_{n+1},\ldots,\blambda_h',\ldots, \blambda_r), \quad \blambda_h' := (\emptyset, \lambda_{h,2},\ldots,\lambda_{h,a_h}).\] 
  Let us write $\cS$ for the multi-set of labels $\ul{\b{\mu}} \neq \ublambda$ appearing in the decomposition of \[E^{(h)}_{\ublambda^\bullet}(\lambda_{h,1})P_{\ublambda^{\bullet}}((\ublambda^\circ)')\] into indecomposable projective summands. By Lemma \ref{lem:produce_proj}, \begin{align}\label{eqn:direct_sum_thing} \fil{i}{\ell}E^{(h)}(\lambda_{h,1})\Delta^{\bullet}(\ublambda^\bullet,(\ublambda^\circ)') = F_i(\ell)\Delta^\bullet(\ublambda) \oplus \bigoplus_{\ul{\b{\mu}} \in \cS} \fil{i}{\ell}\Delta^{\bullet}(\ublambda^\bullet,\ul{\b{\mu}}).\end{align} On the other hand, we have an isomorphism \[\fil{i}{\ell}E^{(h)}(\lambda_{h,1})\Delta^{\bullet}(\ublambda^\bullet,(\ublambda^\circ)') \simeq E^{(h)}(\lambda_{h,1})\fil{i}{\ell}\Delta^{\bullet}(\ublambda^\bullet,(\ublambda^\circ)').\]
  By induction, we may assume that for any $\ul{\b{\mu}}$ satisfying $|\mu_{h,1}| < |\lambda_{h,1}|$, the object $\fil{i}{\ell}\Delta^{\bullet}(\ublambda^\bullet,\ul{\b{\mu}})$ has a filtration with sections \[\Delta^\bullet((\ublambda^\bullet)^F(j,\ell)[i],\ul{\b{\mu}}).\] 
  In particular, the object $\fil{i}{\ell}\Delta^{\bullet}(\ublambda^\bullet,(\ublambda^\circ)')$ has a $\Delta^\bullet$ filtration with sections \[\Delta^\bullet((\ublambda^\bullet)^F(j,\ell)[i],(\ublambda^\circ)').\] By exactness of $E^{(h)}(\lambda_{h,1})$, we see that $E^{(h)}(\lambda_{h,1})\fil{i}{\ell}\Delta^{\bullet}(\ublambda^\bullet,(\ublambda^\circ)')$ has a $\Delta^\bullet$-filtration with sections $\Delta^\bullet((\ublambda^\bullet)^F(j,\ell)[i],\ul{\b{\mu}})$ for $\ul{\b{\mu}} \in \cS \cup\{\ublambda^\circ\}$. Comparing this filtration with a $\Delta^\bullet$-filtration of (\ref{eqn:direct_sum_thing}) and using the known filtrations of all other summands $\fil{i}{\ell}\Delta^{\bullet}(\ublambda^\bullet,\ul{\b{\mu}})$ (by Lemma \ref{lem:produce_proj}, $|\mu_{h,1}| < |\lambda_{h,1}|$ for all $\ul{\bmu} \in \cS$ so we know these filtrations by induction) yields the desired filtration on $F_i(\ell)\Delta^\bullet(\ublambda)$.
  Here, we use the fact that $\fil{i}{\ell}\Delta^\bullet(\ublambda) \in \Delta(\cC)$, being a summand of the $\Delta^\bullet$-filtered object \[E^{(h)}(\lambda_{h,1})\fil{i}{\ell}\Delta^{\bullet}(\ublambda^\bullet,(\ublambda^\circ)')\] whose other summands are $\Delta^\bullet$-filtered by induction.
  
  Thus, we are left with the case $q_h = 0$. We proceed identically by induction on $p_h$, assuming that some $p_h \neq 0$ for $h > n$ and then inducting on $|\lambda_{h,a_h}|$. We leave the details to the reader.
\end{proof}

\begin{cor}\label{cor:preserve_weak}
  For all $i \in \ZZ$ and $\ell = 1,\ldots,r$, the functors $\eil{i}{\ell}$ and $\fil{i}{\ell}$ preserve $\Delta^\bullet(\cC)$ and $\ov{\nabla}^\bullet(\cC)$.
\end{cor}

\begin{proof}
  Lemma \ref{lem:bullet_action} proves that the functors preserve $\Delta^\bullet(\cC)$. The fact that they preserve $\ov{\nabla}^\bullet(\cC)$ follows from Proposition \ref{prop:homological_criterion}.
\end{proof}

\subsubsection{Weak Tilting Objects}\label{subsubsec:weak_tilting}

In this subsubsection, we introduce and study basic properties of the weak tilting objects, the main characters in the proof of our uniqueness theorem.

\begin{defn}
A \textit{weak tilting object} is simply a tilting object defined with respect to the weak stratification on $\cC$. We write $T^\bullet(\ublambda)$ for the indecomposable weak tilting object with label $\ublambda \in \Lambda$. We will write $\Tilt^\bullet(\cC) := \Delta^\bullet(\cC) \cap \ov{\nabla}^\bullet(\cC)$ for the full subcategory of weak tilting objects.
\end{defn}

Thanks to Corollary \ref{cor:preserve_weak}, all categorification functors $\eil{i}{\ell}$ and $\fil{i}{\ell}$ preserve $\Tilt^\bullet(\cC)$. 
We establish combinatorial lemmas describing the actions of $\eil{i}{\ell}$ and $\fil{i}{\ell}$ on weak tilting objects in $\cC$. Recall the definition of the functors $E^{(h)}(\mu)$ and $F^{(h)}(\mu)$ from Definition \ref{defn:content_string}.

\begin{lem}\label{lem:weak_tilting_action}
Suppose $q_h \neq 0$ for $h \leq n$ and fix $\mu \in \cP$. Moreover, pick $\ublambda \in \ul{\Lambda}_h$ such that \[\sum_{\ell=1}^r \sum_{\substack{j=1 \\ c_{\ell,j} = 1}}^{a_\ell} |\lambda_{\ell,j}| \leq m - |\mu|.\] The object $T^\bullet(\ublambda_{\mu}(h))$ is the unique direct summand in $E^{(h)}(\mu)T^\bullet(\ublambda)$ with maximum label. In fact, all other labels appearing have the form $(\ubnu,\ublambda^\circ)$, where $\ubnu < \ublambda^\bullet_{\mu}(h)$ and $|\nu_{h,1}| < |\mu|$. 

Now suppose $p_h \neq 0$ for $h \leq n$ and fix $\mu \in \cP$. Pick $\ublambda \in \ov{\Lambda}_h$ such that \[\sum_{\ell=1}^r \sum_{\substack{j=1 \\ c_{\ell,j} = 0}}^{a_\ell} |\lambda_{\ell,j}| \leq m - |\mu|.\] The object $T^\bullet(\prescript{}{\mu}{\ublambda}(h))$ is the unique direct summand in $F^{(h)}(\mu)T^\bullet(\ublambda)$ with maximum label. All other labels appearing have the form $(\ubnu,\ublambda^\circ)$, where $\ubnu < \prescript{}{\mu}{\ublambda}$ and $|\ubnu_{h,a_h}| < |\mu|$. 
\end{lem}

\begin{proof}
  In light of Lemma \ref{lem:bullet_action}, this is completely analogous to Lemma \ref{lem:produce_proj}.
\end{proof}

The following lemma describes the action of $\eil{i}{\ell}$ and $\fil{i}{\ell}$ on $T^\bullet(\ublambda)$ when $\ell > n$. 

\begin{lem}\label{lem:circ_action_weak_tilting}
  For $\ell > n$ and any $\ublambda = (\ublambda^\bullet,\ublambda^\circ) \in \Lambda$, we have 
  \begin{align*}
  \eil{i}{\ell}T^\bullet(\ublambda) = \bigoplus_{\ubnu \in \cL_{\Xi,\ublambda^\bullet}^{E,i,\ell}} T^\bullet(\ublambda^\bullet, \ubnu), \quad \fil{i}{\ell}T^\bullet(\ublambda) = \bigoplus_{\ubnu \in \cL_{\Xi,\ublambda^\bullet}^{F,i,\ell}} T^\bullet(\ublambda^\bullet, \ubnu).
  \end{align*}
\end{lem}

\begin{proof}
    The proof is similar to that of Lemma \ref{lem:bullet_action}. We induct on $p^\bullet + q^\bullet$. When $p^\bullet = q^\bullet = 0$, there is a unique stratum in the weak stratification, so $T^\bullet(\ublambda) = \Delta^\bullet(\ublambda)$. Then, the claim is immediate.
    
    We now explain how to reduce the claim to the case $q^\bullet = 0$. Suppose $q_h \neq 0$ for some $h \leq n$. We use lower admissible tail splitting (Subsubsection \ref{subsubsec:lower_admissible_tail_splitting}) to equip subcategory $\ul{\cC}_h$ with the structure of a MFTPC of shape $(p^\bullet,q^\bullet-1,p^\circ,q^\circ)$. We write $\ul{T}^\bullet_h(\ublambda)$ for the indecomposable weak tilting object with label $\ublambda$ in $\ul{\cC}_h$. By induction on $q^\bullet$, we may assume that 
    \[
      \uleil{i}{\ell}\ul{T}_h^\bullet(\ublambda^-(h)) = \bigoplus_{\ubnu \in \cL_{\Xi,\ublambda^\bullet}^{E,i,\ell}} \ul{T}_h^\bullet(\ublambda^-(h)^\bullet, \ubnu).
    \]
    The functor $\uleil{i}{\ell}$ tautologically commutes with the inclusion $\ul{\cC}_h \hookrightarrow \cC$ and moreover, this inclusion sends weak tilting objects to weak tilting objects. In particular, as objects of $\cC$, \[\eil{i}{\ell}{T}^\bullet(\ublambda^-(h)) = \bigoplus_{\ubnu \in \cL_{\Xi,\ublambda^\bullet}^{E,i,\ell}} {T}^\bullet(\ublambda^-(h)^\bullet, \ubnu).\] By Lemma \ref{lem:weak_tilting_action}, the object $T^\bullet(\ublambda)$ is a direct summand of $E_{h}(\lambda_{h,a_h})T^\bullet(\ublambda^-(h))$, and all other indecomposable summands $T^\bullet(\ul{\b{\tau}},\ublambda^\circ)$ have $|\tau_{h,a_h}| < |\lambda_{h,a_h}|$. Let us write $\cS$ for the multi-set of labels $\ul{\b{\tau}} \neq \ublambda^\bullet$ such that $T^\bullet(\ul{\b{\tau}},\ublambda^\circ)$ appears in the decomposition of $E_h(\lambda_{h,a_h})T^\bullet(\ublambda')$. Now observe that \begin{align*} T^\bullet(\ublambda) \oplus \bigoplus_{\ul{\b{\tau}} \in \cS} T^\bullet(\ul{\b{\tau}},\ublambda^\circ) &\simeq \eil{i}{\ell}E^{(h)}(\lambda_{h,a_h}){T}^\bullet(\ublambda^-(h)) \\[5pt] &\simeq E_h(\lambda_{h,a_h})\eil{i}{\ell}{T}^\bullet(\ublambda^-(h)) = \bigoplus_{\ubnu \in \cL_{\Xi,\ublambda^\bullet}^{E,i,\ell}} T^\bullet(\ublambda^\bullet,\ubnu)\oplus \bigoplus_{\ul{\b{\tau}}\in \cS}\bigoplus_{\ubnu \in \cL_{\Xi,\ublambda^\bullet}^{E,i,\ell}} T^\bullet(\ul{\b{\tau}}, \ubnu).\end{align*} Recall that $|\tau_{h,a_h}| < |\lambda_{h,a_h}|$ for all $\ul{\b{\tau}} \in \cS$, so we are done by induction on $|\lambda_{h,a_h}|$. 

    We are left with the case $q^\bullet = 0$ and $p^\bullet \neq 0$. This case is handled completely analogously, by using the Serre subcategory $\ov{\cC}_h$ (for some $h \leq n$ such that $p_h \neq 0$) and swapping the roles of $E$ and $F$. Thus, we omit the proof for brevity.
\end{proof}

\begin{cor}\label{cor:circ_action_weak_tilting}
Suppose $q_h \neq 0$ for some $h > n$ and fix $\mu \in \cP$. Moreover, pick some $\ublambda \in \ov{\Lambda}_h$ such that  \[\sum_{\ell=1}^r \sum_{\substack{j=1 \\ c_{\ell,j} = 1}}^{a_\ell} |\lambda_{\ell,j}| \leq m - |\mu|.\] The object $T^\bullet(\ublambda_{\mu}(h))$ is the unique direct summand in $E^{(h)}(\mu)T^\bullet(\ublambda)$ with minimum label, and all other labels $\ubnu$ have $|\ubnu_{h,1}| < |\lambda_{h,1}|$. 

On the other hand, suppose $p_h \neq 0$ for some $h > n$ and fix $\mu \in \cP$. Pick some $\ublambda \in \ul{\Lambda}_h$ such that \[\sum_{\ell=1}^r \sum_{\substack{j=1 \\ c_{\ell,j} = 0}}^{a_\ell} |\lambda_{\ell,j}| \leq m - |\mu|.\] Then, for any $\mu \in \cP$, the object $T^\bullet(\prescript{}{\mu}{\ublambda}(h))$ is the unique direct summand in $F^{(h)}(\mu)T^\bullet(\ublambda)$ with minimum label, and all other labels $\ubnu$ have $|\ubnu_{h,a_h}| < |\lambda_{h,a_h}|$.
\end{cor}

\begin{proof}
  In light of Lemma \ref{lem:circ_action_weak_tilting}, this is completely analogous to Lemma \ref{lem:produce_proj}.
\end{proof}

Finally, we need to show that these tilting objects are compatible with the upper admissible categorical splitting constructions. To state the result, we set up some notation. Suppose $q_h \neq 0$ for some $h > n$, and consider the Serre quotient $\ov{\cC}_h$ from Subsubsection \ref{subsubsec:upper_admissible_head_splitting}, which we equipped with the structure of a mixed admissible MFTPC of type $\ov{\Xi}_h$. We write $\varphi^h: \cC \to \ov{\cC}_h$ for the corresponding quotient functor and $\varphi^h_!: \ov{\cC}_h \to \cC$ for its left adjoint. For any $\ublambda \in \ov{\Lambda}_h$, we will write $\ov{T}^\bullet_h(\ublambda)$ for the indecomposable weak tilting object in $\ov{\cC}_h$ with label $\ublambda$. 

\begin{prop}\label{prop:weak_tilting_pullback}
  For any $\ublambda \in \ov{\Lambda}_h$, we have $\varphi_!^h\ov{T}^\bullet_h(\ublambda) \simeq T^\bullet(\ublambda)$ as objects of $\cC$. 
\end{prop}

\begin{proof}
  The idea of the proof is similar in spirit to our proof of Lemma \ref{lem:bullet_action}. When $p^\bullet = q^\bullet = 0$, there is a unique stratum, so $T^\bullet(\ublambda)$ and $\ov{T}^\bullet_h(\ublambda)$ coincide with the indecomposable projective objects with label $\ublambda$ in their respective categories. Thus, the claim is immediate.

  We explain how to reduce the claim to the case $q^\bullet = 0$. Suppose $q_j \neq 0$ for some $j \leq n$ and consider the Serre subcategory $\ul{\cC}_j$ from Subsubsection \ref{subsubsec:lower_admissible_tail_splitting}. We equipped $\ul{\cC}_j$ with the structure of an admissible MFTPC of shape $(p^\bullet,q^\bullet-1,p^\circ,q^\circ)$. We construct the analogous subcategory of $\ov{\cC}_h$, which we denote $\cD$ for concision. In particular, $\cD$ has shape $(p^\bullet, q^\bullet-1,p^\circ,q^\circ-1)$. Let us write $\iota_j: \ul{\cC}_j \to \cC$ and $\iota_{\cD}: \cD \to \ul{\cC}_h$ for the respective inclusion functors.
  Similarly, we write $T_j^\bullet(\ublambda)$ and $T_{\cD}^\bullet(\ublambda)$ for the indecomposable weak tilting objects in $\ul{\cC}_j$ and $\cD$, respectively. Finally, write $\psi: \ul{\cC}_j \to \cD$ for the Serre quotient functor with left adjoint $\psi_!$. By the induction hypothesis, we may assume $\psi_!(T_{\cD}^\bullet(\ublambda')) = T_j^\bullet(\ublambda')$ for all $\ublambda' \in \ov{\Lambda}_h \cap \ul{\Lambda}_j$. Thanks to the natural isomorphism $\iota_j \circ \psi_! \simeq \varphi^h_! \circ \iota_{\cD}$ and the fact that $\iota_j$ and $\iota_{\cD}$ send weak tilting objects to weak tilting objects, we deduce that $\varphi^h_!(\ov{T}^\bullet_h(\ublambda')) = T^\bullet(\ublambda')$ as objects of $\cC$.
  
  Recall the notation $\ublambda^-(j) \in \ul{\Lambda}_j$ from Subsection \ref{subsec:categorical_splitting}. Thanks to the natural isomorphisms $\eil{i}{j} \circ \varphi_!^h \simeq \varphi_!^h \circ \uleil{i}{j}$ (since $j \leq n$ and $h > n$), we deduce that \[\varphi_!^h(\ul{E}^{(j)}(\lambda_{j,a_j})\ov{T}_h^\bullet(\ublambda^-(j))) \simeq E^{(j)}(\lambda_{j,a_j})\varphi_!^h(\ov{T}_h^\bullet(\ublambda^-(j))) = E^{(j)}(\lambda_{j,a_j}) T^\bullet(\ublambda^-(j)).\] 
  By Lemma \ref{lem:weak_tilting_action}, we can characterize $\ov{T}^\bullet_h(\ublambda)$ as the unique indecomposable direct summand of $\ul{E}^{(j)}(\lambda_{j,a_j})\ov{T}_h^\bullet(\ublambda^-(j))$ not appearing as a summand in $\ul{E}^{(j)}(\eta)\ov{T}_h^\bullet(\ubnu,\ublambda^\circ)$ for $\eta \in \cP$ and $\ubnu \in \ul{\Lambda}_j^\bullet$ such that $\ubnu_{\eta}(j) < \ublambda^\bullet$. Applying the same characterization to $T^\bullet(\ublambda)$ as an indecomposable direct summand of $E^{(j)}(\lambda_{j,a_j}) T^\bullet(\ublambda^-(j))$ yields the claim. The case where $p^\bullet \neq 0$ but $q^\bullet = 0$ is handled similarly and thus omitted. 
\end{proof}

\subsection{Heisenberg Categorical Actions and Generalized Cyclotomic Quotients} \label{subsec:heisenberg}

We write
\[
\UU := \Delta(\ubempty,\ubempty) = \Delta^\bullet(\ubempty,\ubempty), \quad \TT^\bullet := \bigoplus_{0 \leq \alpha,\beta \leq m} E^\alpha F^\beta\UU \in \Tilt^\bullet(\cC).
\]
In particular, $\UU = T^\bullet(\ubempty,\ubempty)$ and $\TT^\bullet$ are both weak tilting objects in $\cC$. Given two categorifications $\cC$ and $\cD$ of the same type, our goal is to produce an equivalence $\cC \simeq \cD$ by identifying their subcategories of weak tilting objects. This will be done via an extension of \cite[Theorem 5.1]{brundan_losev_webster}. This proof requires that the algebra $\End_{\cC}(\TT^\bullet)$ depends only on the type $\Xi$ and is otherwise independent of the choice of categorification $\cC$. We establish this fact in this subsection.

For our purposes, it is convenient to recast categorical type A actions in terms of module categories over the \textit{(degenerate) Heisenberg category} $\Heis_k$ for some $k \in \ZZ$, following the work of Brundan, Savage, and Webster \cite{BSW_Heisenberg_KM}. 

\begin{defn}\label{defn:heisenberg_k}
    The degenerate Heisenberg category of level $k \in \ZZ$, denoted $\Heis_k$, is the strict monoidal category generated by an object $\uparrow$ and its dual $\downarrow$. The morphisms in this category are generated by
    \[
x = 
\raisebox{-0.25\height}{
  \begin{tikzpicture}[scale=1]
    \filldraw (0,0.2) circle (2pt);  
    \draw[line width=0.5,->] (0,0) -- (0,0.5);  
  \end{tikzpicture}
},
\quad
s = 
\raisebox{-.25\height}{
  \begin{tikzpicture}[scale=1]
    \draw[line width=0.5,->] (0,0) -- (0.5,0.5);
    \draw[line width=0.5,->] (0.5,0) -- (0,0.5);
  \end{tikzpicture}
},
\quad
c = 
\raisebox{-0.25\height}{\begin{tikzpicture}
   \draw[line width = 0.5, ->] (0,0) arc[start angle=180,end angle=360,radius=0.25];
\end{tikzpicture}},
\quad
d = 
\raisebox{-0.25\height}{\begin{tikzpicture}
   \draw[line width =0.5,->] (0,0) arc[start angle=180,end angle=0,radius=0.25];
\end{tikzpicture}},
\]
subject to certain degenerate affine Hecke, adjunction, and inversion relations. We will follow the conventions \cite[Definition 1.1]{brundan_heisenberg} and refer the reader there for an explicit description of these relations. We highlight the following consequence of these relations:
\begin{align*}
  \text{if $k \geq 0$, then $\uparrow \ox \downarrow\ \simeq\ \downarrow \ox \uparrow \oplus\ \unit^{\oplus k}$}, \quad\text{if $k < 0$, then $\uparrow \ox \downarrow \oplus\ \unit^{\oplus (-k)}\ \simeq\ \downarrow \ox \uparrow$},
\end{align*}
where $\unit$ is the monoidal unit.
\end{defn}

For any subset $I \subset \CC$, we say that a weight $\fsl_I$-module $M$ has \textit{central charge} $k$ if each weight $\alpha$ of $M$ satisfies $\sum_{i \in I} \langle h_i,\alpha\rangle = k$. Brundan, Savage, and Webster \cite[Theorem 5.22]{BSW_Heisenberg_KM} showed that any Schurian category $\cC$ equipped with a type A categorification $(\cC,E,F,x,\tau)$ of a $\fsl_I$-module of central charge $k$ can be uniquely equipped with the structure of a $\Heis_k$-module category so that 
    \[
    \uparrow\ \mapsto E, \quad \downarrow\ \mapsto F, \quad \raisebox{-0.25\height}{
  \begin{tikzpicture}[scale=1]
    \filldraw (0,0.3) circle (2pt);  
    \draw[line width=0.5,->] (0,0.5) -- (0,0);  
  \end{tikzpicture}
} \mapsto x \in \End(F), \quad 
\raisebox{-.25\height}{
  \begin{tikzpicture}[scale=1]
    \draw[line width=0.5,->] (0.5,0.5) -- (0,0);
    \draw[line width=0.5,->] (0,0.5) -- (0.5,0);
  \end{tikzpicture}
} \mapsto \tau \in \End(F^2)
    \]
and
\[
\raisebox{-0.25\height}{\begin{tikzpicture}
   \draw[line width = 0.5, ->] (0,0) arc[start angle=180,end angle=360,radius=0.25];
\end{tikzpicture}} \mapsto \varepsilon \in \operatorname{Hom}(\id_{\cC}, FE),\quad
\raisebox{-0.25\height}{\begin{tikzpicture}
   \draw[line width =0.5,->] (0,0) arc[start angle=180,end angle=0,radius=0.25];
\end{tikzpicture}} \mapsto \eta \in \Hom(EF, \id_{\cC}),
\]
where $\varepsilon$ and $\eta$ are the adjunction counit and unit, respectively.

Now consider our restricted MFTPC $\cC$. Recall that $\CC \ox_{\ZZ} K_0(\cC)$ is naturally isomorphic to a direct sum of certain weight spaces in an external tensor product $\cT$ of tensor products of Fock space representations of $\fsl_{\ZZ}^{\oplus r}$. In particular, $\CC \ox_{\ZZ} K_0(\cC)$ has a basis given by the classes of the standard objects $[\Delta(\ublambda)]$. Moreover, if $[\Delta(\ublambda)]$ and $[\Delta(\ubmu)]$ belong to different weight spaces, then  $\Delta(\ublambda)$ and $\Delta(\ubmu)$ are not linked in $\cC$.
 It follows that $\cC = \bigoplus_{i \in I} \cC_i$, where each $\cC_i$ categorifies some weight space in $\CC \ox_{\ZZ} K_0(\cC)$. In particular, we can view $\cC$ as a ``restricted'' $\fsl_{\ZZ}^{\oplus r}$-categorification in the sense of Definition \ref{defn:slI_categorification}. Moreover, an easy computation shows that $\cT$ is an integrable $\fsl_{\ZZ}^{\oplus r}$-module with central charge $p - q$. The key observation is that \cite[Theorem 5.22]{BSW_Heisenberg_KM} is established by checking certain relations among natural transformations in a $\fsl_{I}$-categorification, and these computations are checked ``locally," i.e., involving only small powers of $E$ and $F$. Thus, \cite[Theorem 5.22]{BSW_Heisenberg_KM} easily extends to the setting of \textit{restricted} MFTPCs.

\begin{defn}
We define $k(\Xi) := p - q$. We write $k := k(\Xi)$ unless there is risk of ambiguity. 
\end{defn}

Let us write $\Heis_{k,m}$ for the full subcategory consisting of objects $\uparrow^\alpha \downarrow^\beta$ where $\alpha,\beta \leq m$. It no longer makes sense to ask that $\cC$ is a module category over $\Heis_k$ (higher powers of $E$ and $F$ are not defined), but instead, we can make sense of an ``evaluation functor" $\Heis_{k,m} \to \cC$ at the distinguished standard object $\UU = \Delta(\ubempty,\ubempty)\in \cC$.

\begin{cor}\label{cor:restricted_multifock_heisenberg}
  Suppose $\cC$ is a level $m$ restricted mixed admissible MFTPC. Then, there is a functor \[\Heis_{k,m} \to \cC, \quad \uparrow^{\ox\alpha} \ox \downarrow^{\ox \beta} \mapsto E^\alpha F^\beta \UU\]
corresponding to the following assignments of natural transformations:
    \[
    \raisebox{-0.25\height}{
  \begin{tikzpicture}[scale=1]
    \filldraw (0,0.3) circle (2pt);  
    \draw[line width=0.5,->] (0,0.5) -- (0,0);  
  \end{tikzpicture}
} \mapsto x \in \End(F), \quad 
\raisebox{-.25\height}{
  \begin{tikzpicture}[scale=1]
    \draw[line width=0.5,->] (0.5,0.5) -- (0,0);
    \draw[line width=0.5,->] (0,0.5) -- (0.5,0);
  \end{tikzpicture}
} \mapsto \tau \in \End(F^2),
    \]
and (writing $\varepsilon$ and $\eta$ for the counit and unit, respectively, for the adjunction $(E,F)$),
\[
\raisebox{-0.25\height}{\begin{tikzpicture}
   \draw[line width = 0.5, ->] (0,0) arc[start angle=180,end angle=360,radius=0.25];
\end{tikzpicture}} \mapsto \varepsilon,\quad
\raisebox{-0.25\height}{\begin{tikzpicture}
   \draw[line width =0.5,->] (0,0) arc[start angle=180,end angle=0,radius=0.25];
\end{tikzpicture}} \mapsto \eta.
\]
\end{cor}

In fact, this functor factors through a \textit{cyclotomic quotient} of $\Heis_k$. These cyclotomic quotients $\Heis_k(e|f)$ depend on a choice of polynomials $f(z),e(z) \in \CC[z]$, which encode the minimal polynomials of $x$ and its adjoint. The precise definition is rather lengthy, so we refer the reader to \cite[Section 5.3]{BSW_Heisenberg_KM} for details. In brief, it involves taking the quotient of $\Heis_k$ by a certain tensor ideal $\cI(e|f)$ containing the morphisms $e(\raisebox{-0.25\height}{
  \begin{tikzpicture}[scale=1]
    \filldraw (0,0.2) circle (2pt);  
    \draw[line width=0.5,->] (0,0) -- (0,0.5);  
  \end{tikzpicture}
})$ and $f(\raisebox{-0.25\height}{
  \begin{tikzpicture}[scale=1]
    \filldraw (0,0.3) circle (2pt);  
    \draw[line width=0.5,->] (0,0.5) -- (0,0);  
  \end{tikzpicture}
})$. 

Let us write \[f_\UU(z) = \prod_{\ell=1}^r \prod_{\substack{j=1 \\ c_{\ell,j} = 0}}^{a_\ell} (z-\sigma_{\ell,j}),\quad e_\UU(z) := \prod_{\ell=1}^r \prod_{\substack{j=1 \\ c_{\ell,j} = 1}}^{a_\ell} (z-\sigma_{\ell,j}).\] In particular, $\deg f_\UU = p$ and $\deg e_\UU = q$. Moreover, let $\Heis_{k,m}(e_{\UU}|f_{\UU}) \subset \Heis_k(e_{\UU}|f_{\UU})$ denote the full subcategory consisting of objects $\uparrow^{\alpha}\downarrow^\beta$ for $\alpha,\beta \leq m$. Note that the minimal polynomials of $x_\UU \in \End(F\UU)$ and the adjoint morphism $x^*_\UU \in \End(E\UU)$ divide $f_\UU$ and $e_\UU$, respectively. In fact, we will see that $f_\UU$ and $e_\UU$ actually coincide with the aforementioned minimal polynomials. 

The upshot of the definition in \cite[Section 5.3]{BSW_Heisenberg_KM} and the discussion in \cite[Section 4.2]{BSW_Heisenberg_KM} is that functor $\Heis_{k,m} \to \cC$ from Corollary \ref{cor:restricted_multifock_heisenberg} descends to a functor
\[
  \Phi: \Heis_{k,m}(e_{\UU}|f_{\UU}) \to \cC.
\]

The goal of this subsection is to prove the following proposition, which will allow us to describe certain endomorphism algebras independently of $\cC$. 

\begin{prop}\label{prop:heisenberg_fully_faithful}
Suppose $\cC$ is a level $m$ restricted mixed admissible MFTPC. Then, for any $\alpha,\beta \leq m$, the homomorphism \[\varphi_{\alpha,\beta}: \End_{\Heis_{k}(e_{\UU}|f_{\UU})}(\uparrow^{\alpha}\downarrow^{\beta}) \to \End_{\cC}(E^\alpha F^\beta \UU)\] induced by the functor $\Phi$ is an isomorphism.
\end{prop}

We first prove the following lemma.

\begin{lem}\label{lem:fully_faithful_easy}
  Proposition \ref{prop:heisenberg_fully_faithful} holds if $\alpha = 0$ or $\beta = 0$.
\end{lem}

\begin{proof}
  We handle the case $\alpha = 0$, since the case $\beta = 0$ is completely analogous. Write $\cC'_{0}$ for the Serre subcategory corresponding to the poset ideal \[\Lambda'_{0} := \{\ublambda \in \Lambda \ |\ \lambda_{\ell,j} = \emptyset \text{ whenever $c_{\ell,j} = 1$ for $\ell\leq n$}\}.\] Iteratively applying the lower admissible tail splitting construction, we obtain truncated categorification functors $(E',F')$ on $\cC'_0$ giving it the structure of an admissible MFTPC of shape $(p^\bullet,0,p^\circ,q^\circ)$. If we write $\iota: \cC_0' \to \cC$ for the inclusion functor, observe that $\iota \circ F' \simeq F \circ \iota$ since $F$ preserves the subcategory $\cC'_0$. Since $\iota$ is fully faithful, we see that it induces an isomorphism 
  \[
  \End_{\cC'_0}((F')^\beta \UU) \simeq \End_{\cC}(F^\beta \UU).
  \]   
  In turn, we define the Serre quotient $\cC_0$ of $\cC'_0$ corresponding to the poset coideal 
  \[
  \Lambda_{0} := \{\ublambda \in \Lambda_{0}' \ |\ \lambda_{\ell,j} = \emptyset \text{ whenever $c_{\ell,j} = 1$ for $\ell> n$}\} \subset \Lambda_{0}'
  \]
  Iteratively applying the upper admissible head splitting construction, we obtain truncated categorification functors $(E^\sharp,F^\sharp)$ on $\cC_0$ giving it the structure of an admissible MFTPC of shape $(p^\bullet,0,p^\circ,0)$. Let $\varphi_!: \cC_0 \to \cC_0'$ denote the left adjoint to the quotient functor. Due to Lemma \ref{lem:f_commute_quotient} and the full faithfulness of $\varphi_!$, the functor $\varphi_!$ induces an isomorphism \[\End_{\cC_0}((F^\sharp)^\beta \UU) \simeq \End_{\cC'_0}((F')^\beta \UU) \simeq \End_{\cC}(F^\beta \UU).\] 
  
  By the definition of an MFTPC, there exist homomorphisms $\dAH_\beta \to \End_{\cC}(F^\beta \UU)$ and $\dAH_\beta \to \End_{\cC_0}((F^\sharp)^\beta \UU)$ induced by the natural transformations $x \in \End(F)$ and $\tau \in \End(F^2)$. Thanks to the axioms \ref{axiom:RMF2_bullet} and \ref{axiom:RMF2_circ}, these morphisms factor through the degenerate cyclotomic Hecke algebra $\dAH_\beta(f_{\UU})$, i.e., the quotient of $\dAH_\beta$ by the two-sided ideal generated by $f_{\UU}(x_1)$. 

  Since the action on $F^\sharp$ is inherited from the action on $F$, the following diagram commutes:
\[\begin{tikzcd}
	& {\dAH_{\beta}(f_{\UU})} \\
	{\End_{\cC_0}((F^\sharp)^\beta \UU)} && {\End_{\cC}(F^\beta \UU)}
	\arrow[from=1-2, to=2-1]
	\arrow[from=1-2, to=2-3]
	\arrow["\sim", from=2-1, to=2-3]
\end{tikzcd}\]
  In fact, by \cite[Lemma 5.13]{BSW_Heisenberg_KM}, there is an isomorphism \[\End_{\Heis_k(e_{\UU}|f_{\UU})}(\downarrow^\beta) \simeq \dAH_\beta(f_{\UU})\] and under this identification, $\varphi_{0,\beta}$ corresponds to the homomorphism $\dAH_\beta(f_{\UU}) \to \End_{\cC}(F^\beta \UU)$ above. Thus, it suffices to show that the map $\dAH_\beta(f_{\UU}) \to \End_{\cC_0}((F^\sharp)^\beta \UU)$ is an isomorphism.

  Applying Corollary \ref{cor:restricted_multifock_heisenberg} to the MFTPC $\cC_0$ and factoring through the appropriate cyclotomic quotient, we produce a functor \[\Phi_0: \Heis_{p,m}(1|f_{\UU}) \to \cC_0, \quad \uparrow^{\ox \alpha} \ox \downarrow^{\ox \beta}  \ \mapsto\ (E^\sharp)^\alpha (F^\sharp)^\beta \UU.\] Moreover, arguing as before, the induced homomorphism \[\varphi^\sharp_{0,\beta}: \End(\downarrow^\beta) \to \End_{\cC_0}((F^\sharp)^\beta \UU)\] coincides with $\dAH_\beta(f_{\UU}) \to \End_{\cC_0}((F^\sharp)^\beta \UU)$ under the identification $\End(\downarrow^\beta) \simeq \dAH_\beta(f_{\UU})$. Hence, it remains to show that $\varphi_{0,\beta}^\sharp$ is an isomorphism. Observe that the following diagram commutes, where all vertical maps are the maps induced by the functor $\Phi_0$:
\[\begin{tikzcd}
	{\End(\downarrow^{\ox \beta})} && {\Hom(\unit, \uparrow^{\ox \beta} \ox \downarrow^{\ox \beta})} && {\Hom(\unit, \unit^{\oplus (\beta! p^\beta)})} \\
	\\
	{\End((F^\sharp)^\beta \UU)} && {\Hom(\UU, (E^\sharp)^\beta (F^\sharp)^\beta \UU)} && {\Hom(\UU,\UU^{\oplus (\beta!p^\beta)})}
	\arrow["\sim", from=1-1, to=1-3]
	\arrow[from=1-1, to=3-1, "\varphi_{0,\beta}^\sharp", swap]
	\arrow["\sim", from=1-3, to=1-5]
	\arrow[from=1-3, to=3-3]
	\arrow[from=1-5, to=3-5]
	\arrow["\sim"', from=3-1, to=3-3]
	\arrow["\sim"', from=3-3, to=3-5]
\end{tikzcd}\]
In the diagram above, the leftmost pair of horizontal arrows are induced by the adjunction $(\downarrow,\uparrow)$ and $(F^\sharp,E^\sharp)$ and the rightmost pair of horizontal arrows are obtained by repeatedly applying the relations $\uparrow \ox \downarrow \ \simeq \ \downarrow \ox \uparrow \oplus \unit^{\oplus p}$ and $EF \simeq FE \oplus \id^{\oplus p}$ and as well as $\uparrow \ox \unit = 0$ and $E\UU = 0$. The rightmost vertical map is clearly an isomorphism, so we are done. 
\end{proof}

\begin{rem}
  By Lemma \ref{lem:fully_faithful_easy}, the map \[\End_{\Heis_{k}(e_{\UU}|f_{\UU})}(\downarrow) \simeq \CC[z]/(f_{\UU}(z)) \to \End_{\cC}(F \UU)\] given by $z \mapsto x_{\UU}$ is an isomorphism. Thus, the minimal polynomial of $x_\UU \in \End_{\cC}(F\UU)$ is $f_\UU(z)$.
\end{rem}

\begin{proof}[Proof of Proposition \ref{prop:heisenberg_fully_faithful}]
  Thanks to adjunction, it suffices to show that the induced map \[\varphi_{\alpha,\beta}': \Hom_{\Heis_{k}(e_\UU|f_\UU)}(\unit, \uparrow^{\beta}\downarrow^\alpha \uparrow^\alpha \downarrow^\beta) \to \Hom_{\cC}(\UU, E^\beta F^\alpha E^\alpha F^\beta \UU)\] is an isomorphism for any $\alpha,\beta \leq m$. First suppose $k \geq 0$. In this case, we can use the isomorphisms $EF \simeq FE \oplus \id^{\oplus k}$ (resp. $\uparrow \ox \downarrow \ \simeq \ \downarrow \ox \uparrow \oplus \unit^{\oplus k}$) to move all copies of $E$ (resp. $\uparrow$) to the right, yielding isomorphisms
  \[
    \uparrow^{\beta}\downarrow^\alpha \uparrow^\alpha \downarrow^\beta\ \simeq\ \bigoplus_{\gamma \leq m}(\downarrow^{\gamma} \ox \uparrow^{\gamma})^{\oplus n_{\gamma}}, \quad E^{\beta}F^\alpha E^\alpha F^\beta\UU \ \simeq\ \bigoplus_{\gamma \leq m}(E^{\gamma} F^{\gamma} \UU)^{\oplus n_{\gamma}},
  \]  
  for some nonnegative integers $n_{\gamma}$ (note that the aforementioned relations preserve the \textit{difference} between the number of occurrences of $E$ and $F$ or $\uparrow$ and $\downarrow$). Then, the following diagram commutes, where the horizontal arrows are the aforementioned isomorphisms and the vertical arrows are the maps induced by $\Phi_0$:
\[\begin{tikzcd}
	{\Hom(\unit,\uparrow^{\beta}\downarrow^{\alpha}\uparrow^{\alpha}\downarrow^{\beta})} && {\displaystyle\bigoplus_{\gamma \leq m} \Hom(\unit, \uparrow^{\gamma}\downarrow^\gamma)^{\oplus n_{\gamma}}} \\
	\\
	{\Hom_{\cC}(\UU, E^\beta F^\alpha E^\alpha F^\beta \UU)} && {\displaystyle\bigoplus_{\gamma \leq m} \Hom_{\cC}(\UU, E^\gamma F^\gamma\UU)^{\oplus n_{\gamma}}}
	\arrow["\sim", from=1-1, to=1-3]
	\arrow[from=1-1, to=3-1]
	\arrow[from=1-3, to=3-3]
	\arrow["\sim"', from=3-1, to=3-3]
\end{tikzcd}\]
The right vertical arrow is an isomorphism by Lemma \ref{lem:fully_faithful_easy}, so we are done.

The case where $k \leq 0$ is handled identically, except we use the relations $FE \simeq EF \oplus \id^{\oplus |k|}$ (resp. $\downarrow \ox \uparrow\ \simeq \ \uparrow \ox \downarrow \oplus \unit^{\oplus |k|}$) to bring all copies of $F$ in the word $E^\beta F^\alpha E^\alpha F^\beta$ (resp. all copies of $\downarrow$ in the word $\uparrow^\beta \downarrow^\alpha \uparrow^\alpha \downarrow^\beta$) to the right.
\end{proof}

\begin{cor}\label{cor:endomorphism_alg_independent}
    If $\cC$ is a level $m$ restricted mixed admissible MFTPC, then the algebra
    \[
    A(\Xi) := \End_{\cC}(\TT^\bullet) \simeq \bigoplus_{0 \leq \alpha,\beta,\alpha',\beta' \leq m} \Hom_{\Heis_{k}(e_{\UU}|f_{\UU})}(\uparrow^{\alpha}\downarrow^{\beta}, \uparrow^{\alpha'}\downarrow^{\beta'})
    \] is \textit{independent} of the categorification $\cC$. 
    In particular, any MFTPC $\cC$ of type $\Xi$ admits a functor 
    \[
    \rho_{\cC}: \cC \to \rmodcat{A(\Xi)}, \quad M \mapsto \Hom_{\cC}(\TT^\bullet,M).
    \]
\end{cor}

\begin{ex}
If $\Xi$ consists of only the pair $((0,\sigma),(0,1))$, where $\sigma \in \ZZ$, then $\End_{\cC}(E^\alpha F^\beta \UU)$ is isomorphic to the walled Brauer algebra $B_{\alpha,\beta}(\sigma)$ (i.e., of dimensions $\alpha,\beta$ and parameter $\sigma$). In fact, $\rmodcat{A(\Xi)}$ is Ringel dual to the abelian envelope of the Deligne category $\cD_\sigma$.
\end{ex}

\begin{rem}\label{rem:same_argument}
A completely analogous argument shows that
\[
    \End_{\cC}(F^\alpha E^\beta \UU) \simeq \End_{\Heis_k(e_{\UU}|f_{\UU})}(\downarrow^\alpha\uparrow^\beta)
\]
is also independent of the categorification $\cC$.
\end{rem}

\begin{defn}
  For $\alpha,\beta \leq m$, we define $H_{\alpha,\beta}(\Xi) := \End_{\Heis_k(e_{\UU}|f_{\UU})(\uparrow^\alpha \downarrow^\beta)}$
\end{defn}

The category $\rmodcat{A(\Xi)}$ admits a tautological categorical type A action given as follows. For each $\alpha,\beta \geq 0$, we have $\fk$-algebra homomorphisms \[H_{\alpha,\beta}(\Xi) \to H_{\alpha+1,\beta}(\Xi), \quad H_{\alpha,\beta}(\Xi) \to H_{\alpha,\beta+1}(\Xi)\] given by \[\varphi \mapsto \id_{\uparrow} \ox \varphi \in \End(\uparrow \ox \uparrow^{\ox \alpha} \ox \downarrow^{\ox \beta}), \quad \varphi \mapsto \varphi \ox \id_{\downarrow} \in \End(\uparrow^{\ox \alpha} \ox \downarrow^{\ox \beta} \ox \downarrow),\] respectively. Then, we simply define the functor $E: \rmodcat{A(\Xi)} \to \rmodcat{A(\Xi)}$ as the direct sum of functors $E: \rmodcat{H_{\alpha,\beta}(\Xi)} \to \rmodcat{H_{\alpha+1,\beta}(\Xi)}$ given by induction: \[M \mapsto M\ox_{H_{\alpha,\beta}(\Xi)}H_{\alpha+1,\beta}(\Xi).\] We similarly define $F: \rmodcat{A(\Xi)} \to \rmodcat{A(\Xi)}$ as the direct sum of the functors \[M \mapsto  M\ox_{H_{\alpha,\beta}(\Xi)}H_{\alpha,\beta+1}(\Xi) .\] The natural transformations $x \in \End(F)$ and $\tau \in \End(F^2)$ act via right multiplication by \[\id_{\uparrow^{\ox \alpha} \ox \downarrow^{\ox \beta}} \ox \raisebox{-0.25\height}{
   \begin{tikzpicture}[scale=1]
    \filldraw (0,0.3) circle (2pt);  
    \draw[line width=0.5,->] (0,0.5) -- (0,0);  
  \end{tikzpicture}
} \in H_{\alpha,\beta+1}(\Xi), 
\quad 
\id_{\uparrow^{\ox \alpha} \ox \downarrow^{\ox \beta}} \ox \raisebox{-.25\height}{
  \begin{tikzpicture}[scale=1]
    \draw[line width=0.5,->] (0.5,0.5) -- (0,0);
    \draw[line width=0.5,->] (0,0.5) -- (0.5,0);
  \end{tikzpicture}
} \in H_{\alpha,\beta+2}(\Xi),\] respectively.
Immediately from the defining relations of $\Heis_k$, we deduce that the data $(E,F,x,\tau)$ defines a categorical type A action on $\rmodcat{A(\Xi)}$. Moreover, we see that the functor 
\[
\rho_{\cC}: \cC \to \rmodcat{A(\Xi)}, \quad M \mapsto \Hom_{\cC}(\TT^\bullet,M)
\]
is tautologically strongly equivariant with respect to the (restricted) categorical actions defined on both $\cC$ and $\rmodcat{A(\Xi)}$. 

\subsection{Double Centralizer Theorem for Weak Tilting Objects} We now aim to prove the following weak tilting ``double centralizer'' theorem à la Losev--Webster \cite[Theorem 5.1]{losev_webster}. This theorem is the core technical step in the proof of our uniqueness theorem, as it reduces the proof to identifying the full subcategories of weak tilting objects in two MFTPCs of the same type. 

\begin{defn}
Let us write $\Tilt^\bullet(\cC,m)$ for the full subcategory of weak tilting objects in $\cC$ with labels belonging to the subset $\Lambda(m,m) \subset \Lambda$. 
\end{defn}

\begin{thm}\label{thm:mixed_fully_faithful}
  The functor \[\rho_{\cC}: \cC \to \rmodcat{A(\Xi)}, \quad M \mapsto \Hom_{\cC}(\TT^\bullet, M)\] is fully faithful when restricted to $\Tilt^\bullet(\cC,m)$.
\end{thm}

Our proof strategy is to use the structural results on weak tilting objects that we established in Subsubsection \ref{subsubsec:weak_tilting} along with the categorical splitting construction to inductively ``peel away'' the lowest weight Fock space factors in our MFTPC, reducing the claim to the case $q^\bullet = q^\circ = 0$. In this latter case, we will see that the proof of \cite[Theorem 5.1]{losev_webster} generalizes with little difficulty.

\subsubsection{Inductive Base Case}\label{subsubsec:induction_base_case}

In this subsubsection, we prove the following proposition, which establishes the base case $q^\bullet = q^\circ= 0$ for our proof of Theorem \ref{thm:mixed_fully_faithful}.

\begin{prop}\label{prop:base_case_fully_faithful}
  Assume $q^\bullet = q^\circ = 0$. The functor \[\rho_{\cC}: \cC \to \rmodcat{A(\Xi)}, \quad M \mapsto \Hom_{\cC}(\TT^\bullet, M)\] is fully faithful on $\Tilt^\bullet(\cC,m)$.
\end{prop}

In fact, Proposition \ref{prop:base_case_fully_faithful} is essentially a direct consequence of \cite[Section 5]{losev_webster}. Since the label $\ubempty$ is maximal in $\Lambda^\circ$, we have $\Delta_{\ubempty}(\ubempty) = P_{\ubempty}(\ubempty)$ and hence \[\UU = \Delta^\bullet(\ubempty, \ubempty) = \Delta(\ubempty,\ubempty).\] Since $q^\bullet = q^\circ = 0$, the label $(\ubempty,\ubempty)$ is incomparable to any other label in $\Lambda$, so the object $\UU$ (and hence any direct summand of $\TT^\bullet$) is both projective and injective. 

\begin{lem}\label{lem:no_kill_cokernel}
  Suppose $q^\bullet = q^\circ = 0$, and take any $\ublambda^\bullet \in \Lambda^\bullet$ with $k$ boxes and any $\ublambda^\circ \in \Lambda^\circ$ with $k'$ boxes, where $k + k' \leq m$. For any simple $L$ appearing in the socle of $\Delta^\bullet(\ublambda^\bullet,\ublambda^\circ)$, we have $({F}^\bullet)^k ({F}^\circ)^{k'} L \neq 0$. 
\end{lem}

\begin{proof}
  Completely analogous to the proof of \cite[Proposition 5.2]{losev_webster}.
\end{proof}

\begin{cor}\label{cor:projective_cover_I}
  Suppose $q^\bullet = q^\circ = 0$. For $M \in \Delta^\bullet(\cC)$ such that the labels belonging to any $\Delta^\bullet$-filtration of $M$ lie in $\Lambda(m,m)$, the injective hull of $M$ is a summand of $\TT^\bullet$. 
\end{cor}

\begin{proof}
  Lemma \ref{lem:no_kill_cokernel} shows that there is a nonzero map from the socle of $M$ to $\TT^\bullet$. Under the assumption $q^\bullet = q^\circ = 0$, the object $\TT \in \cC$ is both projective and injective. Hence, the aforementioned map factors through the injective hull of the socle.
\end{proof}

\begin{lem}\label{lem:surjection_costandard_kernel}
  Suppose $q^\bullet = q^\circ = 0$. For any indecomposable projective $P \in \cC$ whose cosocle has label belonging to $\Lambda(m,m)$, there exists $\beta \leq m$ and an injection $P \hookrightarrow {F}^\beta \UU$ with standardly filtered cokernel.
\end{lem}

\begin{proof}
  Completely analogous to \cite[Lemma 5.3]{losev_webster}.
\end{proof}

\begin{proof}[Proof of Proposition \ref{prop:base_case_fully_faithful}]
  In view of Corollary \ref{cor:projective_cover_I} and Lemma \ref{lem:surjection_costandard_kernel}, the proof follows by the Five Lemma argument used to prove \cite[Theorem 5.1]{losev_webster}.
\end{proof}

\subsubsection{Full Faithfulness Theorem for Mixed Admissible Categorifications}\label{subsubsec:full_faithfulness}

In this subsubsection, we use the categorical splitting construction along with properties of weak tilting objects to inductively reduce the proof of Theorem \ref{thm:mixed_fully_faithful} to Proposition \ref{prop:base_case_fully_faithful}. We introduce an auxiliary definition which relates $\rmodcat{A(\ul{\Xi}_h)}$ with $\rmodcat{A(\Xi)}$, a required step for our inductive proof strategy.

\begin{defn} 
Define the universal Heisenberg category $\ul{\Heis}$ as the diagrammatic strict monoidal category generated by objects $\uparrow$ and $\downarrow$ and morphisms $x,s,c,d$ as in Definition \ref{defn:heisenberg_k} subject to the degenerate affine Hecke and adjunction relations from \cite[Definition 1.1]{brundan_heisenberg}, but \textbf{not} the inversion relation.
\end{defn}

In particular, we have a natural full functor $\ul{\Heis} \to \Heis_k$ for any $k \in \ZZ$. We define
\[R := \bigoplus_{\alpha,\beta,\alpha',\beta' \in \ZZ_{\geq 0}} \Hom_{\ul{\Heis}}(\uparrow^\alpha \ox \downarrow^\beta, \uparrow^{\alpha'} \otimes \downarrow^{\beta'}).\] The surjective $\fk$-algebra homomorphism $\zeta: R \twoheadrightarrow A(\Xi)$ yields a fully faithful pullback functor \[\zeta^*: \rmodcat{A(\Xi)} \to \rmodcat{R}.\]

First suppose $q_h \neq 0$ for some $h \leq n$. In Subsubsection \ref{subsubsec:lower_admissible_tail_splitting}, we showed that the Serre subcategory $\ul{\cC}_h$ inherits the structure of a level $m$ restricted MFTPC of type $\ul{\Xi}_h$. 

\begin{lem}\label{lem:universal_heisenberg_pullback}
The following diagram commutes:
  \begin{equation}\label{eqn:diagram0}
    \begin{tikzcd}[scale=0.5]
    \ul{\cC}_h && \cC \\ \\
    \rmodcat{A(\ul{\Xi}_h)} && \rmodcat{A(\Xi)} \\ \\
    & \rmodcat{R}&
    \arrow["\iota", from=1-1, to = 1-3]
    \arrow["\rho_{\ul{\cC}_{h}}", from = 1-1, to = 3-1, swap]
    \arrow["\rho_\cC", from = 1-3, to = 3-3]
    \arrow["\zeta^*", from = 3-1, to = 5-2, swap]
    \arrow["\zeta^*", from = 3-3, to = 5-2]
    \end{tikzcd}
  \end{equation}
  In particular, if $\rho_{\ul{\cC}_h}$ is fully faithful on $\Tilt^\bullet(\ul{\cC}_h(m))$, then $\rho_{\cC}$ induces an isomorphism 
  \[
    \Hom_{\cC}(T^\bullet(\ul{\bmu}'), T^\bullet(\ublambda')) \xto{\sim} \Hom_{A(\Xi)}(\rho_{\cC}(T^\bullet(\ul{\bmu}')), \rho_{\cC}(T^\bullet(\ublambda')))
  \]
  for any $\ul{\bmu}',\ublambda' \in \ul{\Lambda}_h(m,m)$. 
\end{lem}

\begin{proof}
  For any $M \in \ul{\cC}_h$, we have a sequence of $\fk$-vector space isomorphisms
  \begin{align*}
    \Hom_{\ul{\cC}_h}(\ul{F}^\alpha \ul{E}^\beta \UU, M) &\simeq \Hom_{\ul{\cC}_h}(\ul{E}^\beta \UU, \ul{E}^\alpha M) \\[5pt] &\simeq \Hom_{\cC}(\iota \ul{E}^\beta \UU, \iota \ul{E}^\alpha M) \simeq \Hom_{\cC}(E^\beta \UU, E^\alpha \iota M) \simeq \Hom_{\cC}(F^\alpha E^\beta \UU,\iota M),
  \end{align*}
  where we used the natural isomorphism $E\iota \simeq \iota \ul{E}$. In particular, we have a $\fk$-vector space isomorphism $\pi_{\ul{\cC}_h}(M) \simeq \pi_{\cC}(\iota M)$. This isomorphism tautologically intertwines the $R$-actions on both sides, so that (\ref{eqn:diagram0}) indeed commutes. The second claim follows from the full faithfulness of $\zeta^*$ and the fact that $\iota$ sends indecomposable weak tilting objects in $\ul{\cC}_h$ to indecomposable weak tilting objects in $\cC$ (with the same label).
\end{proof}

Now suppose $q_h \neq 0$ for some $h > n$ and write $\ov{\cC}_h$ for the corresponding Serre quotient of $\cC$ from Subsubsection \ref{subsubsec:upper_admissible_head_splitting}. We equipped this quotient with the structure of an admissible MFTPC of type $\ov{\Xi}_h$. Let $\varphi: \cC \to \ov{\cC}_h$ denote the quotient functor with left adjoint $\varphi_!: \ov{\cC}_h \to \cC$.

\begin{lem}\label{lem:mixed_universal_heisenberg_pullback}
  The following diagram commutes:
  \begin{equation}\label{eqn:diagram3}
    \begin{tikzcd}[scale=0.5]
    \Delta(\ov{\cC}_h) && \Delta(\cC) \\ \\
    \rmodcat{A(\ov{\Xi}_h)} && \rmodcat{A(\Xi)} \\ \\
    & \rmodcat{R}&
    \arrow["\varphi_!", from=1-1, to = 1-3]
    \arrow["\rho_{\ov{\cC}_{h}}", from = 1-1, to = 3-1, swap]
    \arrow["\rho_\cC", from = 1-3, to = 3-3]
    \arrow["\zeta^*", from = 3-1, to = 5-2, swap]
    \arrow["\zeta^*", from = 3-3, to = 5-2]
    \end{tikzcd}
  \end{equation}
  In particular, if $\rho_{\ov{\cC}_h}$ is fully faithful on $\Tilt^\bullet(\ov{\cC}_h,m)$, then $\rho_{\cC}$ induces an isomorphism 
  \[
    \Hom_{\cC}(T^\bullet(\ublambda'), T^\bullet(\ubnu')) \xto{\sim} \Hom_{A(\Xi)}(\rho_{\cC}(T^\bullet(\ublambda')), \rho_{\cC}(T^\bullet(\ubnu')))
  \]
  for any $\ublambda',\ubnu' \in \ov{\Lambda}_h(m)$. 
\end{lem}

\begin{proof}
Fix some $M \in \Delta(\ov{\cC}_h)$. Thanks to Lemma \ref{lem:f_commute_quotient} and the fact that $\varphi_!$ is fully faithful, we have natural $\fk$-vector space isomorphisms
\begin{align*}
\Hom_{\ov{\cC}_h}(\ov{E}^\alpha \ov{F}^\beta \UU, M) \simeq \Hom_{\ov{\cC}_h}(\ov{F}^\beta \UU, \ov{F}^\alpha M) \simeq \Hom_{\cC}(\varphi_!\ov{F}^\beta \UU, \varphi_!\ov{F}^\alpha M) \simeq \Hom(E^\alpha F^\beta \UU, \varphi_!M).
\end{align*}
Clearly, these isomorphisms intertwine the induced $R$-actions, so we are done. The second claim follows from the full faithfulness of $\zeta^*$ and Proposition \ref{prop:weak_tilting_pullback}.
\end{proof}

We are finally ready to complete our proof of Theorem \ref{thm:mixed_fully_faithful}. 

\begin{proof}[Proof of Theorem \ref{thm:mixed_fully_faithful}]
  We proceed by induction on $q^\bullet$ and $q^\circ$. The base case $q^\bullet = q^\circ = 0$ was handled in Proposition \ref{prop:base_case_fully_faithful}. First, we reduce the claim to the case $q^\bullet = 0$. Suppose $q_h \neq 0$ for some $h \leq n$. Consider the highest weight subcategory $\ul{\cC}_h$ from Subsubsection \ref{subsubsec:lower_admissible_tail_splitting}. This subcategory is equipped with the structure of a restricted MFTPC of shape $(p^\bullet,q^\bullet-1,p^\circ,q^\circ)$. Recall the notation $\ublambda^-(h)$ from Notation \ref{notation:add_partition}. For any $\ublambda \in \Lambda(m,m)$, Lemma \ref{lem:weak_tilting_action} shows that $T^\bullet(\ublambda)$ is a direct summand of $E^{(h)}(\lambda_{h,a_h})T^\bullet(\ublambda^-(h))$. Hence, we want to show that the induced map 
  \[
    \Hom_{\cC}(T^\bullet(\ul{\bmu}), E^{(h)}(\lambda_{h,a_h})T^\bullet(\ublambda^-(h))) \to \Hom_{A(\Xi)}(\rho_\cC(T^\bullet(\ul{\bmu})), E^{(h)}(\lambda_{h,a_h})\rho_{\cC}(T^\bullet(\ublambda^-(h))))
  \]
  is an isomorphism for any $\ul{\bmu} \in \Lambda(m,m)$ . By adjunction, it suffices to prove that the map 
  \[
    \Hom_{\cC}(T^\bullet(\ul{\bmu}), T^\bullet(\ublambda')) \to \Hom_{A(\Xi)}(\rho_\cC(T^\bullet(\ul{\bmu})), \rho_{\cC}(T^\bullet(\ublambda')))
  \]
  is an isomorphism for any $\ublambda' \in \ul{\Lambda}_h(m,m)$ and any $\ul{\bmu} \in \Lambda(m,m)$. Using Lemma \ref{lem:weak_tilting_action} again, we can prove that the induced map \[
    \Hom_{\cC}(E^{(h)}(\lambda_{\mu,a_h})T^\bullet(\ul{\bmu}'), T^\bullet(\ublambda')) \to \Hom_{A(\Xi)}(E^{(h)}(\lambda_{h,a_h})\rho_\cC(T^\bullet(\ul{\bmu}')), \rho_{\cC}(T^\bullet(\ublambda')))
  \] is an isomorphism. Applying adjunction once more along with the fact that the functor $E$ preserves the subcategory $\ul{\cC}_h$, we reduce the induction step to the claim that the induced map \[\Hom_{\cC}(T^\bullet(\ul{\bmu}'), T^\bullet(\ublambda')) \to \Hom_{A(\Xi)}(\rho_\cC(T^\bullet(\ul{\bmu}')), \rho_{\cC}(T^\bullet(\ublambda')))\] is an isomorphism for $\ublambda',\ul{\bmu}' \in \ul{\Lambda}_h(m,m)$. We are done by induction and Lemma \ref{lem:universal_heisenberg_pullback}.
  
  It remains to consider $q^\bullet = 0$ and $q^\circ \neq 0$. We proceed similarly by induction on $q^\circ$, using upper admissible head splitting instead of lower admissible tail splitting, Corollary \ref{cor:circ_action_weak_tilting} instead of Lemma \ref{lem:weak_tilting_action} and Lemma \ref{lem:mixed_universal_heisenberg_pullback} instead of Lemma \ref{lem:universal_heisenberg_pullback}. We leave the details to the reader.
\end{proof}

\subsection{Uniqueness of Mixed Admissible Categorifications}\label{subsec:mixed_admissible}

We finally prove the uniqueness theorem that lies at the heart of the technical machinery in this paper.

\begin{thm}\label{thm:mixed_admissible_uniqueness}
  Let $\cC$ and $\cD$ be level $m$ restricted mixed admissible MFTPCs of type $\Xi$. Then, there exists a strongly equivariant equivalence $\cC(m,m) \simeq \cD(m,m)$ that intertwines the labels of simple objects.
\end{thm}

Take $\cC$ and $\cD$ from the statement of Theorem \ref{thm:mixed_admissible_uniqueness}. We will use subscripts $\cC$ and $\cD$ to indicate which category an object, morphism, or functor is associated to. For instance, we write $\UU_{\cC}$ (resp. $\UU_{\cD}$) for the standard object in $\cC$ (resp. $\cD$) with minimal label -- this object is both simple and weak tilting. By Theorem \ref{thm:mixed_fully_faithful}, we have fully faithful functors 
\[
\rho_{\cC}: \Tilt^\bullet(\cC,m) \to \rmodcat{A(\Xi)}, \quad \rho_{\cD}: \Tilt^\bullet(\cD,m) \to \rmodcat{A(\Xi)}.
\]

\begin{proof}[Proof of Theorem \ref{thm:mixed_admissible_uniqueness}]

  We first explain how an equivalence $\Tilt^\bullet(\cC,m) \simeq \Tilt^\bullet(\cD,m)$ would give rise to an equivalence of categories $\cC(m,m) \simeq \cD(m,m)$. Thanks to Proposition \ref{prop:weak_tilting_pullback}, we can identify $\Tilt^\bullet(\cC(m,m))$ (resp. $\Tilt^\bullet(\cD(m,m))$) with $\Tilt^\bullet(\cC,m)$ (resp. $\Tilt^\bullet(\cD,m)$) through an appropriate sequence of inclusion functors and left adjoints to quotient functors. Thus, we would have an equivalence $\Tilt^\bullet(\cC(m,m)) \simeq \Tilt^\bullet(\cD(m,m))$. Such an equivalence would imply an equivalence of Ringel duals $\cC(m,m)' \simeq \cD(m,m)'$ taken with respect to the weak tilting objects. Then, \cite[Corollary 4.11]{brundan_stroppel_semiinfinite} would allow us to recover an equivalence $\cC(m,m) \simeq \cD(m,m)$. In fact, if the equivalence $\Tilt^\bullet(\cC,m) \simeq \Tilt^\bullet(\cD,m)$ intertwined labels of weak tilting objects, then the resulting equivalence $\cC(m,m) \simeq \cD(m,m)$ would also intertwine labels of simple objects. Finally, the functors $\rho_{\cC}$ and $\rho_{\cD}$ are tautologically strongly equivariant. 
  
  Thus, in light of Theorem \ref{thm:mixed_fully_faithful}, it suffices to show that
  \begin{align}\tag{$\heartsuit$}\label{claim:heart}
    \rho_{\cC}(T^\bullet_{\cC}(\ublambda)) = \rho_{\cD}(T^\bullet_{\cD}(\ublambda))
  \end{align}
  for each $\ublambda \in \Lambda(m,m)$. 
  To prove (\ref{claim:heart}), we first reduce to the case $p^\bullet = q^\bullet = 0$ by induction on $p^\bullet$ and $q^\bullet$. Suppose $q_h \neq 0$ for some $h \leq n$. As before, consider the full subcategory $\ul{\cC}_h$ along with its truncated categorification functors $(\ul{E},\ul{F})$ that realize it as a (restricted) MFTPC of shape $(p,q-1)$. We similarly consider $\ul{\cD}_h \subset \cD$. The commutativity of diagram (\ref{eqn:diagram0}) implies that the following diagram commutes:
\begin{equation}\label{eqn:diagram2}\begin{tikzcd}
	& {\ul{\cC}_h} && \cC \\
	{\rmodcat{A(\ul{\Xi}_h)}} && {\rmodcat{R}} && {\rmodcat{A(\Xi)}} \\
	& {\ul{\cD}_h} && \cD
	\arrow["\iota"{description}, from=1-2, to=1-4]
	\arrow["{\rho_{\ul{\cC}_h}}"', from=1-2, to=2-1]
	\arrow["{\rho_{\cC}}", from=1-4, to=2-5]
	\arrow["{\zeta_{k-1}^*}"{description}, from=2-1, to=2-3]
	\arrow["{\zeta_k^*}"{description}, from=2-5, to=2-3]
	\arrow["{\rho_{\ul{\cD}_h}}", from=3-2, to=2-1]
	\arrow["\iota"{description}, from=3-2, to=3-4]
	\arrow["{\rho_{\cD}}"', from=3-4, to=2-5]
\end{tikzcd}\end{equation}
  
  By induction, we assume that $\rho_{\ul{\cC}_h}(T^\bullet_{\ul{\cC}_h}(\ublambda')) = \rho_{\ul{\cD}_h}(T^\bullet_{\ul{\cD}_h}(\ublambda'))$ for all $\ublambda' \in \ul{\Lambda}_h(m,m)$. But then, the commutativity of the diagram (\ref{eqn:diagram2}) along with the full faithfulness of $\zeta_k^*$ implies that \[\rho_{\cC}(T^\bullet_{\cC}(\ublambda')) = \rho_{\cC}(\iota T^\bullet_{\ul{\cC}_h}(\ublambda')) = \rho_{\cD}(\iota T^\bullet_{\ul{\cD}_h}(\ublambda')) = \rho_{\cD}(T^\bullet_{\cD}(\ublambda'))\] for all $\ublambda' \in \ul{\Lambda}_h(m,m)$. Now consider an arbitrary $\ublambda \in \Lambda(m,m)$ and write $\ublambda' := \ublambda^-(h)$. Thanks to Lemma \ref{lem:weak_tilting_action} and the fact that the functors $\rho_{\cC}$ and $\rho_{\cD}$ are strongly equivariant, we can characterize both $\rho_{\cC}(T^\bullet_{\cC}(\ublambda))$ and $\rho_{\cD}(T^\bullet_{\cD}(\ublambda))$ as the unique indecomposable direct summand of 
  {\small
  \begin{align*}
    \rho_{\cC}(E^{(h)}(\lambda_{h,a_h})T^\bullet_{\cC}(\ublambda')) = E^{(h)}(\lambda_{h,a_h})\rho_{\cC}(T^\bullet_{\cC}(\ublambda')) = E^{(h)}(\lambda_{h,a_h})\rho_{\cD}(T^\bullet_{\cD}(\ublambda')) = \rho_{\cD}(E^{(h)}(\lambda_{h,a_h})T^\bullet_{\cD}(\ublambda'))
  \end{align*}
  }
  not appearing as an indecomposable direct summand in 
  \[E^{(h)}(\mu)\rho_{\cC}(T^\bullet_{\cC}(\ubnu)) = E^{(h)}(\mu)\rho_{\cD}(T^\bullet_{\cD}(\ubnu))\]
  for any partition $\mu \in \cP$ and $\ubnu \in \ul{\Lambda}_h(m,m)$ such that $\ubnu_\ell = \ublambda_\ell$ for $\ell \neq h$ and $\ubnu_\mu(h) < \ublambda$ (note there are only finitely many such pairs $(\mu,\ubnu)$). Thus, we have reduced the claim to the case $q^\bullet = 0$. 
  
  We now proceed by induction on $p^\bullet$, assuming $p_h \neq 0$ for some $h \leq n$. As before, we consider the subcategories $\ov{\cC}_h$ and $\ov{\cC}_h$ along with their truncated categorification functors $\ov{E}$ and $\ov{F}$. For each $m' \leq m$, let us write \[\rho_{\cC,m'}': \cC \to \rmodcat{A(\Xi)}, \quad M \mapsto \bigoplus_{0 \leq \beta \leq m'} \Hom(F^{\beta}\UU,M).\] We similarly define the functor $\rho_{\cD,m'}': \cD \to \rmodcat{A(\Xi)}$. A similar inductive argument to the one we made in the case $q \geq 1$ shows that $\pi_{\cC,m'}'(T^\bullet_{\cC}(\ublambda)) = \pi_{\cD,m'}'(T^\bullet_{\cD}(\ublambda))$ for each $\ublambda \in \Lambda(m',m')$.
  
  On the other hand, if $\ublambda \not \in \Lambda(m',m')$, we claim that $\rho_{\cC,m'}'(T^\bullet_{\cC}(\blambda)) = \rho_{\cD,m'}'(T^\bullet_{\cD}(\blambda)) = 0$. Indeed, note that the label $\ublambda$ of any costandard object appearing in a costandard filtration of $E^\beta T^\bullet_{\cC}(\blambda)$ has $|\blambda_\ell| > 0$ for at least one index $\ell = 1,\ldots,r$. Hence, \[\rho_{\cC,m
  }'(T^\bullet_{\cC}(\ublambda)) = \bigoplus_{0 \leq \beta \leq m'}\Hom_{\cC}(\UU, E^\beta T^\bullet_{\cC}(\ublambda)) = 0.\]
  Thanks to the relation $EF \simeq FE \oplus \id^{\oplus p}$ and the observation $E\UU = 0$, note that we can write $\rho_{\cC} = \bigoplus_{m' \leq m} a_{m'}\rho'_{\cC,m'}$ for some integers $a_{m'} \geq 0$. We can also decompose $\rho_{\cD} = \bigoplus_{m' \leq m} a_{m'}\rho'_{\cD,m'}$ for the \textit{same} integers $a_{m'}$, so we are done.

  It remains to consider $p^\bullet = q^\bullet = 0$. This case proceeds identically by induction on $p^\circ$ and $q^\circ$. One should use upper admissible splitting instead of lower admissible splitting, Corollary \ref{cor:circ_action_weak_tilting} instead of Lemma \ref{lem:weak_tilting_action} and diagram (\ref{eqn:diagram3}) instead of diagram (\ref{eqn:diagram0}). We leave details to the reader.
\end{proof}



\subsection{Relating Restricted Categorifications and Full Categorifications}\label{subsec:restricted_vs_full}

We finally clarify the relationship between restricted and full admissible categorifications in the following three ways.

\begin{enumerate}
\item We construct restricted mixed admissible categorifications from full categorifications.
\item We ``glue'' restricted mixed admissible categorifications to form \textit{pre-MFTPCs}, which are $\fsl^{\oplus r}_{\ZZ}$-categorifications satisfying some weaker conditions than those of honest MFTPCs.
\item We explain how to recover full categorifications from the restricted categorifications constructed in (1) using the gluing construction from (2). As a consequence, we prove a uniqueness result for \textit{full} MFTPCs. 
\end{enumerate}

\subsubsection{Constructing Restricted Categorifications}\label{subsubsec:restricted_from_full}

Let $\cC$ be a mixed admissible multi-Fock tensor product categorification of type $\Xi$ and fix some $m \geq 0$. Write $\Lambda = \cP^{a_1} \times \cdots \times \cP^{a_r}$ for the poset underlying $\cC$. Moreover, define \[\Lambda_m := \cP^{a_1}(m) \times \cdots \times \cP^{a_r}(m),\] which is a poset coideal in the upper finite poset ideal \[\Lambda_m' := \cP^{a_1}(m) \times \cdots \times \cP^{a_n}(m) \times \cP^{a_{n+1}} \times \cdots \times \cP^{a_r} \subset \Lambda.\] Define $\cC_m'$ as the Serre subcategory of $\cC$ corresponding to $\Lambda_m'$. Note that $\cC_m'$ is a Schurian category from the axioms of a mixed semi-infinite highest weight category. In fact, the induced stratification on $\cC_m'$ realizes it as an upper finite highest weight category. Finally, we write $\cC_m$ for the Serre quotient of $\cC_m'$ corresponding to the poset coideal $\Lambda_m$. This category naturally inherits the structure of a finite highest weight category with poset $\Lambda_m$. We demonstrate how to equip $\cC_m$ with the structure of a level $m$ restricted mixed admissible MFTPC of type $\Xi$.

For $j,k \geq 0$, we define the Serre subcategories $\hat{\cC}^\bullet(j,k) \subset \cC$ corresponding to the poset ideal \[\Lambda^\bullet(j,k) := \left\{\ublambda \in \Lambda\ \bigg|\ \sum_{\ell=1}^n \sum_{i=1}^{p_\ell} |\lambda_{\ell,i}| < j, \quad \sum_{\ell=1}^n \sum_{i=1}^{q_\ell} |\lambda_{\ell,i+p_\ell}| \leq k\right\} \subset \Lambda.\] When $j,k \leq m$, we also define the Serre subcategories $\cC_m^\bullet(j,k) \subset \cC_m$ given by the poset ideal \[\Lambda_m^\bullet(j,k) := \Lambda^\bullet(j,k) \cap \Lambda_m \subset \Lambda_m.\]
For each $\ell = 1,\ldots,n$, observe that the functors $F_i(\ell)$ and $E_i(\ell)$ restrict to biadjoint functors \[\hat{F}_i(\ell): \hat{\cC}^{\bullet}(j,k) \to \hat{\cC}^{\bullet}(j+1,k-1), \quad \hat{E}_i(\ell): \hat{\cC}^\bullet(j+1,k-1) \to \hat{\cC}^\bullet(j,k).\] From the axioms for $F_i(\ell)$ and $E_i(\ell)$, we see that these functors preserve the Serre subcategory of $\cC$ corresponding to the poset ideal $\Lambda_m' \setminus \Lambda_m \subset \Lambda_m'$. Thus, for $j,k \leq m$, the functors $\hat{F}_i(\ell)$ and $\hat{E}_i(\ell)$ descend to biadjoint functors \[\bar{F}_i(\ell): \cC_m^\bullet(j,k) \to \cC_m^\bullet(j+1,k-1), \quad \bar{E}_i(\ell): \cC_m^\bullet(j+1,k-1) \to \cC_m^\bullet(j,k).\] Moreover, from this construction, it is clear that the functors \[\bar{F}_i(\ell): \cC_m^\bullet(j,k) \to \cC_m^\bullet(j+1,k-1), \quad \bar{E}_i(\ell): \cC_m^\bullet(j+1,k-1) \to \cC_m^\bullet(j,k)\] are given by restricting the functors \[\bar{F}_i(\ell): \cC_m^\bullet(m-1,m) \to \cC_m^\bullet(m,m-1), \quad \bar{E}_i(\ell): \cC_m^\bullet(m,m-1) \to \cC_m^\bullet(m-1,m)\] to the Serre subcategories $\cC_m^\bullet(j,k) \subset \cC_m^\bullet(m-1,m)$ and $\cC_m^\bullet(j+1,k-1) \subset \cC_m^\bullet(m,m-1)$, respectively. These latter functors satisfy (RMF3$^\bullet$) by definition. Moreover, the functor \[\bar{F}^\bullet := \bigoplus_{\ell = 1}^n \bigoplus_{\ell \in \ZZ} \bar{F}_i(\ell): \cC_m^\bullet(m-1,m) \to \cC_m^\bullet(m,m-1)\] inherits endomorphisms $\bar{x}^\bullet \in \End(\bar{F}^\bullet)$ and $\tau^\bullet \in \End((\bar{F}^\bullet)^2)$ satisfying (RMF2$^\bullet$) from the endomorphisms $x \in \End(F)$ and $\tau \in \End(F^2)$, respectively.

For $j,k \geq 0$, define the Serre quotients $\hat{\cC}^\circ(j,k)$ of $\cC$ corresponding to the poset coideals \[\Lambda^\circ(j,k) := \left\{\ublambda \in \Lambda \ \bigg|\ \sum_{\ell=n+1}^r \sum_{i=1}^{q_\ell} |\lambda_{\ell,i}| < j, \quad \sum_{\ell=n+1}^r \sum_{i=1}^{p_\ell} |\lambda_{\ell,i+q_\ell}| < k\right\}.\] When $j,k \leq m$, we define the Serre quotients $\cC_m^\circ(j,k)$ of $\cC_m$ corresponding to the poset coideals \[\Lambda_m^\circ(j,k) := \Lambda^\circ(j,k) \cap \Lambda_m \subset \Lambda_m.\]
For each $\ell = n+1,\ldots,r$, we see that the functors $F_i(\ell)$ and $E_i(\ell)$ restrict to biadjoint functors 
\[
  \hat{F}_i(\ell): \pi_!\Delta(\hat{\cC}^\circ(j,k)) \to \pi_!\Delta(\hat{\cC}^\circ(j-1,k+1)), \quad \hat{E}_i(\ell): \pi_!\Delta(\hat{\cC}^\circ(j-1,k+1)) \to \pi_!\Delta(\hat{\cC}^\circ(j,k)).
\]
Since $F_i(\ell)$ and $E_i(\ell)$ preserve the Serre subcategory $\cC_m'$, they descend to biadjoint functors 
\[
  \bar{F}_i(\ell): \pi_!\Delta({\cC}_m^\circ(j,k)) \to \pi_!\Delta({\cC}_m^\circ(j-1,k+1)), \quad \bar{E}_i(\ell): \pi_!\Delta({\cC}_m^\circ(j-1,k+1)) \to \pi_!\Delta({\cC}_m^\circ(j,k)),
\]
where $\pi_!$ here denotes the left adjoint to the quotient functor $\cC_m \to \cC_m(j,k)$. Every object in $\cC_m^\circ(j,k)$ is the cokernel of some map of projective, hence standardly filtered, objects. In particular, since $\pi_!$ is right exact, we conclude $\bar{F}_i(\ell)$ and $\bar{E}_i(\ell)$ uniquely extend to biadjoint functors \[  \bar{F}_i(\ell): \pi_!{\cC}_m^\circ(j,k) \to \pi_!{\cC}_m^\circ(j-1,k+1), \quad \bar{E}_i(\ell): \pi_!{\cC}_m^\circ(j-1,k+1) \to \pi_!{\cC}_m^\circ(j,k).\] Moreover, it is clear that these functors are restrictions of the functors \[\bar{F}_i(\ell): \pi_!{\cC}_m^\circ(m,m-1) \to \pi_!{\cC}_m^\circ(m-1,m), \quad \bar{E}_i(\ell): \pi_!{\cC}_m^\circ(m-1,m) \to \pi_!{\cC}_m^\circ(m,m-1),\] which satisfy (RMF3$^\circ$) by definition. The endomorphisms $x \in \End(F^\circ)$ and $\tau \in \End((F^\circ)^2)$ descend to endomorphisms $\bar{x}^\circ \in \End(\bar{F}^\circ)$ and $\bar{\tau}^\circ \in \End((\bar{F}^\circ)^2)$ satisfying (RMF2$^\circ$), where \[\bar{F}^\circ := \bigoplus_{\ell=n+1}^r \bigoplus_{i \in \ZZ} \bar{F}_i(\ell): \pi_!\cC_m^\circ(m-1,m) \to \pi_! \cC_m^\circ(m,m-1).\] 
Altogether, the functors $\bar{F}_i(\ell)$, $\bar{E}_i(\ell)$ along with the endomorphisms $\bar{x}^\bullet,\bar{x}^\circ,\bar{\tau}^\bullet, \bar{\tau}^\circ$ equip the subquotient $\cC_m$ with the desired structure of a level $m$ restricted MFTPC of type $\Xi$.


\subsubsection{Gluing Restricted Categorifications}\label{subsubsec:gluing}

Suppose we have a family $\mathfrak{F}$ of admissible restricted multi-Fock tensor product categorifications $\cC_m$ for $m \geq 0$, each of type $\Xi$. That is, each $\cC_m$ is a level $m$ restricted mixed admissible MFTPC of type $\Xi$. We will write \[\Lambda_m = \cP^{a_1}(m) \times \cdots \times \cP^{a_r}(m), \quad \Lambda = \bigcup_{m \geq 0} \Lambda_m\] for the poset of $\cC_m$. We explain how to glue these categorifications together to form a full MFTPC.

For $j,k \geq 0$ with $j + k \leq m$, we write $\cC_{m}^{j,k}$ for the Serre subquotient of $\cC_m$ corresponding to {\small\begin{align*}\Lambda^{j,k} := \bigg\{\ublambda \in \Lambda \ \bigg| \ \sum_{\ell=1}^n \sum_{i=1}^{p_\ell} |\lambda_{\ell,i}| \leq j, \quad \sum_{\ell=1}^n\sum_{i=1}^{q_\ell} |\lambda_{\ell, i+p_\ell}| \leq j \quad \sum_{s=n+1}^r \sum_{i=1}^{q_s} |\lambda_{s,i}| \leq k, \quad \sum_{s=n+1}^r \sum_{i=1}^{p_s} |\lambda_{s, i+q_s}| \leq k\bigg\},\end{align*}} which is a poset coideal in the poset ideal  \[\Lambda^j_m := \left\{\ublambda \in \Lambda_m \ \bigg| \ \sum_{\ell=1}^n\sum_{i=1}^{p_\ell} |\lambda_{\ell,i}| \leq j, \quad \sum_{\ell=1}^n \sum_{i=1}^{q_\ell} |\lambda_{\ell, i+p_\ell}| \leq j\right\} \subset \Lambda_m.\]

Observe we have equivalences $\cC_m^{j,k} \simeq \cC_{m'}^{j,k}$ whenever $m' \geq m \geq j + k$. Indeed, Theorem \ref{thm:mixed_admissible_uniqueness} yields highest weight equivalences $\cC_m(m,m) \simeq \cC_{m'}(m,m)$, which restrict to equivalences of subquotients $\cC_m^{j,k} \simeq \cC_{m'}^{j,k}$. Thus, we define the stable category $\cC(\mathfrak{F})^{j,k} := \lim_{m \to \infty} \cC_m^{j,k}$. In other words, $\cC(\mathfrak{F})^{j,k}$ is simply equivalent to $\cC_m^{j,k}$ for any $m \gg 0$. 

For fixed $j \geq 0$ and $k' \geq k$, we have a Serre quotient functor \[\pi_{j}(k',k): \cC(\mathfrak{F})^{j,k'} \to \cC(\mathfrak{F})^{j,k}\] (given by the quotient $\cC_m^{j,k'} \to \cC_m^{j,k}$ for any $m \gg 0$). This quotient functor is fully faithful on standardly filtered objects and in particular induces an isomorphism \[\Hom_{\cC(\mathfrak{F})^{j,k}}(P^{j,k}(\ublambda), P^{j,k}(\ubnu)) \simeq \Hom_{\cC(\mathfrak{F})^{j,k-1}}(P^{j,k-1}(\ublambda), P^{j,k-1}(\ubnu)),\] where $P^{j,k}(\ublambda)$ is the projective cover of the simple $L(\ublambda)$ in $\cC(\mathfrak{F})^{j,k}$. Thus, we can define \[A^{j,k} := \bigoplus_{\ublambda,\ubnu \in \Lambda^{j,k}} H(\ublambda,\ubnu), \quad H(\ublambda,\ubnu) := \lim_{k \to \infty} \Hom_{\cC(\mathfrak{F})^{j,k}}(P^{j,k}(\ublambda),P^{j,k}(\ubnu)),\] where multiplication is the opposite of the composition of morphisms. 

In particular, there is an equivalence of finite abelian categories $\cC(\mathfrak{F})^{j,k} \simeq A^{j,k}\modcat$. Now, $A^{j,k-1} \subset A^{j,k}$ is a direct summand of $A^{j,k}$ corresponding to some idempotent $e_{k-1} \in A^{j,k}$. In this notation, the functor $M \mapsto e_{k-1}M$ corresponds to the quotient functor \[\cC(\mathfrak{F})^{j,k} = A^{j,k}\modcat \to A^{j,k-1}\modcat  = \cC(\mathfrak{F})^{j,k-1}.\]

In turn, we define the $\fk$-algebra 
\[
  A^{j,\infty} := \bigoplus_{\ublambda,\ubnu \in \Lambda^j} H(\ublambda,\ubnu), \quad \Lambda^j := \bigcup_{k \geq 0} \Lambda^{j,k}.
\]
The algebra $A^{j,\infty}$ can be equipped with the structure of a locally finite-dimensional locally unital $\fk$-algebra via the distinguished orthogonal idempotents $e_{\ublambda}$ (for $\ublambda \in \Lambda^j$) defined by \[e_{\ublambda}A^{j,\infty}e_{\ubnu} = H(\ublambda,\ubnu).\] In this notation, observe that \[A^{j,k} = \bigoplus_{\ublambda,\ubnu \in \Lambda^{j,k}} e_{\ublambda} A^{j,\infty} e_{\ubnu}.\] 

\begin{defn}
  Define $\cC(\mathfrak{F})^{j,\infty}$ as the category of locally finite-dimensional $A^{j,\infty}$-modules.
\end{defn}

For each $k \geq 0$, we have a natural functor $\pi^j_k: \cC(\mathfrak{F})^{j,\infty} \to \cC(\mathfrak{F})^{j,k}$ given by \[M \mapsto \bigoplus_{\ublambda \in \Lambda^{j,k}} e_{\ublambda}M.\] We observe that this functor can actually be identified with a Serre quotient functor whose kernel is the subcategory of objects killed by the idempotent $e_k := \sum_{\ublambda \in \Lambda^{j,k}} e_{\ublambda}$. Given any sequence of modules $M_k \in A^{j,k}\modcat$ with $M_k = e_kM_{k+1}$, we can define the $\fk$-vector space \[M := \bigcup_{k \geq 0} M_k\] and equip it with the structure of a locally finite-dimensional $A^{j,\infty}$-module by prescribing the action of $A^{j,k} = e_kA^{j,\infty}e_k \subset A^{j,\infty}$ to be the same as its natural action on $M_{k'}$ for $k' \geq k$. 

\begin{lem}
    The Schurian category $\cC(\mathfrak{F})^{j,\infty}$ has the natural structure of an upper finite highest weight category with poset $\Lambda^j = \bigcup_{m \geq 0} \Lambda^j_m.$ In fact, $\cC(\mathfrak{F})^{j,k}$ is identified with the Serre quotient of $\cC(\mathfrak{F})^{j,\infty}$ corresponding to the poset coideal $\Lambda^{j,k} \subset \Lambda^j$.
\end{lem}

\begin{proof}
    \textbf{Step 1.} First, we claim that the upper finite poset $\Lambda^j$ indexes the simple objects in $\cC(\mathfrak{F})^{j,\infty}$. Given $\ublambda \in \Lambda^{j}$, write $L^{j,k}(\ublambda)$ for the simple object in $\cC(\mathfrak{F})^{j,k}$ with label $\ublambda$. If $\ublambda \not \in \Lambda^{j,k}$, we agree that $L^{j,k}(\ublambda) = 0$. In particular, we have $e_kL^{j,k+1}(\ublambda) = L^{j,k}(\ublambda)$. Thus, consider the $A^{j,\infty}$-module \[L^{j,\infty}(\ublambda) := \bigcup_{k \geq 0} L^{j,k}(\ublambda).\] Any proper submodule $N \subset L^{j,\infty}(\ublambda)$ gives rise to a nontrivial proper submodule $e_kN \subset L^{j,k}(\ublambda)$, and hence $L^{j,\infty}(\ublambda)$ is simple. 
    Conversely, if $M \in \cC(\mathfrak{F})^{j,\infty}$ is simple, then since $e_k$ is a Serre quotient functor, we must have $e_kM = e_kL^{j,\infty}(\ublambda)$ for $k \gg 0$ and some $\ublambda \in \Lambda_m^j$. We deduce $M = L^{j,\infty}(\ublambda)$. Thus, the category $\cC(\mathfrak{F})^{j,\infty}$ admits a stratification by the upper finite poset $\Lambda^j$. 
    
    Note that the induction functor \[M \mapsto A^{j,\infty}e_k \ox_{A_{j,k}} M, \quad \cC(\mathfrak{F})^{j,k} \to \cC(\mathfrak{F})^{j,\infty}\] is left adjoint to the quotient functor $N \mapsto e_kN$. For each $\ublambda \in \Lambda^{j,k}$, we now define the objects \[\Delta^{j,\infty}(\ublambda) := A^{j,\infty}e_k \ox_{A^{j,k}} \Delta^{j,k}(\ublambda), \quad P^{j,\infty}(\ublambda) := A^{j,\infty}e_k \ox_{A^{j,k}} P^{j,k}(\ublambda),\] where $\Delta^{j,k}(\ublambda)$ and $P^{j,k}(\ublambda)$ are the standard and indecomposable projective objects, respectively, with label $\ublambda$ in $\cC(\mathfrak{F})^{j,k}$. 
    These objects clearly belong to $\cC(\mathfrak{F})^{j,\infty}$. In fact, it is not hard to see that these objects do not depend on the choice of $k \geq 0$ provided $\ublambda \in \Lambda^{j,k}$. Indeed, for any $k' \geq k$, the induction functor $A^{j,k'}e_k \ox_{A^{j,k}} -$ is precisely the left adjoint to the quotient functor $\cC(\mathfrak{F})^{j,k'} \to \cC(\mathfrak{F})^{j,k}$. In particular, induction sends standard objects to standard objects. Hence, \[A^{j,\infty}e_{k'} \ox_{A^{j,k'}} \Delta^{j,k'}(\ublambda) = A^{j,\infty}e_{k'} \ox_{A^{j,k'}} A^{j,k'}e_k \ox_{A^{j,k}} \Delta^{j,k}(\ublambda) = A^{j,\infty}e_k \ox_{A^{j,k}} \Delta^{j,k}(\ublambda).\] The same argument holds for $P^{j,\infty}(\ublambda)$ as well.
  
    \textbf{Step 2.} We claim that $\Delta^{j,\infty}(\ublambda)$ is standard with respect to the stratification. A completely formal argument shows that the manifestly indecomposable object $\Delta^{j,\infty}(\ublambda)$ is projective in $\cC(\mathfrak{F})^{j,\infty}_{\leq\ublambda}$. This uses the following three facts. First, the Serre quotient functor $\cC(\mathfrak{F})^{j,\infty} \to \cC(\mathfrak{F})^{j,k}$ restricts to an exact functor $\cC(\mathfrak{F})^{j,\infty}_{\leq \ublambda} \to \cC(\mathfrak{F})^{j,k}_{\leq \ublambda}$. Second, $\Delta^{j,k}(\ublambda) = e_k\Delta^{j,\infty}(\ublambda)$ is projective in $\cC(\mathfrak{F})^{j,k}_{\leq \ublambda}$. Finally, the quotient functor is the exact right adjoint of induction. Now, note that \[A^{j,\infty} = \bigcup_{k' \geq 0} e_{k'}A^{j,\infty} \implies \Delta^{j,\infty}(\ublambda) = \bigcup_{k' \geq 0} e_{k'}A^{j,\infty}e_k \ox_{A^{j,k}} \Delta^{j,k}(\ublambda) = \bigcup_{k' \geq 0} e_{k'}\Delta^{j,\infty}(\ublambda).\]   
    Thus, we obtain a surjective morphism \[\Delta^{j,\infty}(\ublambda) = \bigcup_{k' \geq 0} e_{k'}\Delta^{j,\infty}(\ublambda)= \bigcup_{k' \geq 0} \Delta^{j,k'}(\ublambda) \twoheadrightarrow \bigcup_{k' \geq 0} L^{j,k'}(\ublambda) = \bigcup_{k' \geq 0} e_{k'}L^{j,\infty}(\ublambda) = L^{j,\infty}(\ublambda).\] By \cite[Lemma 3.1]{brundan_stroppel_semiinfinite}, we see that $\Delta^{j,\infty}(\ublambda)$ is standard with respect to the $\Lambda^j$-stratification.
    
    \textbf{Step 3.} It remains to show that this stratification gives a highest weight structure. A formal argument also shows that $P^{j,\infty}(\ublambda)$ is projective in $\cC(\mathfrak{F})^{j,\infty}$. Moreover, the right-exactness of induction yields a surjection $P^{j,\infty}(\ublambda) \twoheadrightarrow \Delta^{j,\infty}(\ublambda)$ of $A^{j,\infty}$-modules. We need to show that the kernel $N$ of this surjection is filtered by modules $\Delta^{j,\infty}(\ubnu)$ with $\ubnu \geq \ublambda$. For each $k \geq 0$, we have a short exact sequence of $A^{j,k}$-modules \[0 \to e_kN \to P^{j,k}(\ublambda) \to \Delta^{j,k}(\ublambda) \to 0.\] Since each $\cC(\mathfrak{F})^{j,k}$ is a highest weight category, each $e_kN$ is filtered by objects $\Delta^{j,k}(\ubnu)$ for $\ubnu \in \Lambda^{j,k}$ with $\ubnu \geq \ublambda$. In fact, for $k' \geq k$, the induction functor $\cC(\mathfrak{F})^{j,k} \to \cC(\mathfrak{F})^{j,k'}$ (i.e., the right adjoint to the quotient functor) sends short exact sequences of standardly filtered objects to short exact sequences, so the multiset of labels $\ubnu$ appearing in the standard filtration of $e_kN$ does not depend on the choice of $k \geq 0$ as long as $\ublambda \in \Lambda_{j,k}$. Let us write $\ubnu_1,\ldots,\ubnu_s$ for these labels, so that $e_kN$ admits a filtration with sections $\Delta^{j,k}(\ubnu_1),\ldots,\Delta^{j,k}(\ubnu_s)$ for $k \gg 0$. Thus, $N = \bigcup_{k \geq 0} e_kN$ admits a filtration by the modules $\Delta^{j,\infty}(\ubnu_t) = \bigcup_{k \geq 0} e_k\Delta^{j,\infty}(\ubnu_t)$ for $t = 1,\ldots,s$.  
\end{proof}


For each $j' > j \geq 0$, there is a natural surjection $A^{j'} \to A^j$, which induces a fully faithful pullback functor $A^j\modcat \to A^{j'}\modcat$. This functor sends locally finite-dimensional $A^j$-modules to locally finite-dimensional $A^{j'}$-modules and thus restricts to a fully faithful functor \[\iota(j,j'): \cC(\mathfrak{F})^j \to \cC(\mathfrak{F})^{j'}.\]
Thus, we define the union \[\cC'(\mathfrak{F}) := \bigcup_{j \geq 0} \cC(\mathfrak{F})^{j,\infty},\] where the functor $\cC(\mathfrak{F})^{j,\infty} \to \cC(\mathfrak{F})^{j',\infty}$ is the embedding $\iota(j,j')$. The following result is manifest.

\begin{lem}
For $j' \geq j$, the functor $\iota(j,j')$ identifies $\cC(\mathfrak{F})^{j,\infty}$ as the Serre subcategory of the upper-finite highest weight category $\cC(\mathfrak{F})^{j',\infty}$ corresponding to the poset ideal $\Lambda^j \subset \Lambda^{j'}$. In particular, the simple objects in $\cC'(\mathfrak{F})$ are labelled by $\Lambda = \bigcup_{j \geq 0} \Lambda^j$.
\end{lem}

It remains to equip $\cC'(\mathfrak{F})$ with the structure of a type $\Xi$ MFTPC. For $m_1,m_2 \leq j$, define \[\Lambda^j(m_1,m_2) := \left\{\ublambda \in \Lambda \ \bigg|\ \sum_{\ell=1}^n\sum_{i=1}^{p_\ell} \lambda_{\ell,i}\leq m_1, \quad \sum_{\ell=1}^n \sum_{i=1}^{q_\ell} \lambda_{\ell,i+p_\ell} \leq m_2\right\} \subset \Lambda^j\] along with the corresponding Serre subcategories $\cC(\mathfrak{F})^{j,\infty}(m_1,m_2) \subset \cC(\mathfrak{F})^{j,\infty}$. 

We first construct biadjoint functors \[F^{(j)}: \cC(\mathfrak{F})^{j,\infty}(j-1,j) \to \cC(\mathfrak{F})^{j,\infty}(j,j-1), \quad E^{(j)}: \cC(\mathfrak{F})^{j,\infty}(j,j-1) \to \cC(\mathfrak{F})^{j,\infty}(j-1,j)\] as follows. First, for any $M \in \Delta(\cC(\mathfrak{F})^{j,\infty}(j-1,j))$ or $N \in \Delta(\cC(\mathfrak{F})^{j,\infty}(j,j-1))$, there exists $k_0$ such that $e_kM \in \cC(\mathfrak{F})^{j,k}(j-1,j)$ and $e_kN \in \cC(\mathfrak{F})^{j,k}(j,j-1)$ for $k \geq k_0$. 

First, one checks that $e_{k-1}F = Fe_{k-1}$ as functors $\cC^{j,k}(j-1,j) \to \cC^{j,k-1}(j,j-1)$. Similarly, we have $e_{k-1}E = Ee_{k-1}$. In turn, we define \[F^{(j)}M := \bigcup_{k \geq k_0} Fe_kM, \quad E^{(j)}N:= \bigcup_{k \geq k_0} Ee_kN.\]
Since every projective object in $\cC(\mathfrak{F}^{j,\infty})(j-1,j)$ (resp. $\cC(\mathfrak{F}^{j,\infty})(j,j-1)$) admits a finite standard filtration, we can uniquely extend the exact functors $F^{(j)}$ and $E^{(j)}$ to the desired biadjoint functors 
 \[F^{(j)}: \cC(\mathfrak{F})^{j,\infty}(j-1,j) \to \cC(\mathfrak{F})^{j,\infty}(j,j-1), \quad E^{(j)}: \cC(\mathfrak{F})^{j,\infty}(j,j-1) \to \cC(\mathfrak{F})^{j,\infty}(j-1,j).\]

Now, note that \[\cC'(\mathfrak{F}) = \bigcup_{j \geq 1} \cC(\mathfrak{F})^{j,\infty}(j-1,j) = \bigcup_{j \geq 1} \cC(\mathfrak{F})^{j,\infty}(j,j-1),\] so taking the direct limit of the functors $F^{(j)}$ and $E^{(j)}$, we obtain induced biadjoint functors \[F: \cC'(\mathfrak{F}) \to \cC'(\mathfrak{F}), \quad E: \cC'(\mathfrak{F}) \to \cC'(\mathfrak{F}).\]
From the definitions, these functors satisfy (MF4). Moreover,the natural transformations $x_m := x_m^\bullet + x_m^\circ,\tau_m := \tau_m^\bullet + \tau^\circ_m$ from the restricted categorification structure extend to natural transformations $x \in \End(F)$ and $\tau \in \End(F^2)$ satisfying (MF2) and (MF3). Thus, we have equipped $\cC'(\mathfrak{F})$ with the structure of a MFTPC of type $\Xi$. 

\subsubsection{Recovering Full Categorifications from Restricted Categorifications}

Finally, suppose we have an admissible (full) multi-Fock tensor product categorification $\cC$. As before, in the upper admissible case, we make the additional assumption that $\cC$ is a highest weight category. From this full categorification, we apply the construction of Subsubsection \ref{subsubsec:restricted_from_full} to obtain a family $\mathfrak{F}$ of restricted admissible categorifications $\cC_m$ for $m \in \ZZ$. The gluing construction in Subsubsection \ref{subsubsec:gluing} yields a MFTPC $\cC'(\mathfrak{F})$ from this family of restricted MFTPCs. First, we relate $\cC$ and $\cC'(\mathfrak{F})$.

\begin{thm}\label{thm:recover_from_gluing}
There is a strongly equivariant equivalence $\cC \simeq \cC'(\mathfrak{F})$ intertwining labels of simples.
\end{thm}

\begin{proof} 
  Let us write $\cC^j$ for the Schurian Serre subcategory of $\cC$ corresponding to the upper finite poset ideal $\Lambda^j \subset \Lambda$. From basic Morita theory for locally unital algebras (e.g., see \cite[Theorem 2.4]{Brundan_Davidson}), we have an equivalence of categories $\cC^j \simeq A_{\cC}^{j}\modcat_{\lfd}$, where $A_{\cC}^j$ is the locally finite-dimensional locally unital algebra \[A_{\cC}^{j} := \bigoplus_{\ublambda,\ubnu \in \Lambda^j} \Hom_{\cC}(P^j(\ublambda),P^j(\ubnu)),\] where we write $P^j(\ublambda)$ for the projective cover of $L(\ublambda)$ in the Schurian subcategory $\cC^j$. On the other hand, if $\ublambda,\ubnu \in \Lambda^{j,k}$ for some $k \geq 0$, then we actually have \[H(\ublambda,\ubnu) = \Hom_{\cC^{j,k}}(P^{j,k}(\ublambda),P^{j,k}(\ubnu)) \simeq \Hom_{\cC^j}(P^j(\ublambda),P^j(\ubnu)),\] so we have an isomorphism of locally unital algebras \[A^{j,\infty} = \bigoplus_{\ublambda,\ubnu \in \Lambda^j} H(\ublambda,\ubnu) \simeq A^j_{\cC}\] intertwining labels of distinguished idempotents on both sides. It follows that we have an equivalence $\cC^j \simeq \cC(\mathfrak{F})^{j,\infty}$ intertwining labels of simples. The isomorphisms $A^{j,\infty} \simeq A^j_{\cC}$ naturally intertwine the quotient maps $A^{j+1,\infty} \to A^{j,\infty}$ and $A^{j+1}_\cC \to A^{j}_\cC$ and hence the aforementioned equivalences intertwine the embeddings $\cC^{j,\infty}(\mathfrak{F}) \hookrightarrow \cC^{j+1,\infty}(\mathfrak{F})$ and $\cC^j \hookrightarrow \cC^{j+1}$, respectively. Thus, we take the direct limit of these equivalences to obtain the desired equivalence $\cC \simeq \cC'(\mathfrak{F})$. It is not hard to see that this functor is also strongly equivariant -- we omit the details for brevity.
\end{proof}

\begin{cor}\label{cor:full_uniqueness}
    Suppose $\cC$ and $\cD$ are (full) MFTPCs of the same admissible type $\Xi$. Then, there exists a strongly equivariant equivalence of categories $\cC \simeq \cD$ intertwining labels of simple objects.
\end{cor}

\begin{proof}
  In view of Theorem \ref{thm:mixed_admissible_uniqueness}, the glued categorifications $\cC' := \cC'(\mathfrak{F})$ and $\cD' := \cD'(\mathfrak{F})$ depend only on the type $\Xi$ and not on the underlying categories $\cC$ and $\cD$. Thus, we obtain a manifest strongly equivariant equivalence $\cC' \simeq \cD'$ intertwining labels of simple objects which yields the desired equivalence $\cC \simeq \cD$ in light of Theorem \ref{thm:recover_from_gluing}.
\end{proof}



\section{Constructing Admissible Fock Space Tensor Product Categorifications}
\label{sec:construction}

In this section, we construct restricted admissible multi-Fock space tensor product categorifications by considering truncations of classical parabolic categories $\cO$. We will first construct these categorifications in case where $r = 1$, i.e., with a single external tensor factor. Taking Deligne tensor products of these Fock tensor product categorifications gives us constructions of restricted \textit{multi-Fock} tensor product categorifications. We will glue these categorifications together to obtain full MFTPCs. As a corollary, we obtain Kazhdan--Lusztig type formulas for multiplicities of simple objects in standard objects in full admissible categorifications.

\subsection{Basic Notation}\label{subsec:construction_notation}
 
We first establish some notational conventions for this section and set up a common framework for constructing our admissible categorifications. Throughout Subsections \ref{subsec:construction_notation} to \ref{subsec:constructing_restricted_fock}, we will fix a sequence of integers $\b{\sigma} := (\sigma_1,\ldots,\sigma_a)$ and a $01$-sequence $\b{c} := (c_1,\ldots,c_a)$. We will assume that the sequence $\b{c}$ is \textit{admissible} with swapping index $\kappa$ (see Definition \ref{defn:admissible_sequence}). In turn, we set $\beta := a - \kappa$. We will write \[\sigma_i^\sharp := \begin{cases} \sigma_i & \text{$\b{c}$ is lower admissible} \\ -\sigma_i & \text{$\b{c}$ is upper admissible}. \end{cases}\] for each $i = 1,\ldots,\kappa$ and \[\sigma_i^\flat := \begin{cases} -\sigma_{i+\kappa} & \text{$\b{c}$ is lower admissible} \\ \sigma_{i+\kappa} & \text{$\b{c}$ is upper admissible}\end{cases}\] for each $i = 1,\ldots,\beta$. \textit{Note the signs in these definitions -- they will be crucial later}. 

Moreover, throughout this section, we fix some very large positive integer $N$ -- we will specify how large $N$ needs to be later. For concision, we will write $\sigma_i^\sharp(N) := \sigma_i^\sharp + N$ and $\sigma_j^\flat(N) := \sigma_j^\flat + N$ for each $i = 1,\ldots,\kappa$ and $j = 1,\ldots,\beta$. Finally, set
\[
d^\sharp(N) := \sum_{i=1}^\kappa \sigma_i^\sharp(N), \quad d^\flat(N) := \sum_{i=1}^\beta \sigma_i^\flat(N), \quad d(N) := d^\sharp(N) + d^\flat(N).
\]
Consider the (ordered) decomposition
\begin{align}\label{eqn:super_decomp}
    V(N) := \fk^{d(N)} = \fk^{\sigma_1^\sharp(N)} \oplus \cdots \oplus  \fk^{\sigma_\kappa^\sharp(N)} \oplus \fk^{\sigma_1^\flat(N)} \oplus \cdots \oplus \fk^{\sigma_\beta^\flat(N)}
\end{align}
and fix an ordered basis for $V(N)$ by concatenating bases for the summands in (\ref{eqn:super_decomp}) in order. We will identify elements of $\fg(N)$ with $d(N) \times d(N)$ matrices according to this ordered basis.

Let $\fb(N) \subset \fg(N)$ denote the Borel subalgebra of upper triangular matrices. Also define the Cartan subalgebra $\fh(N) \subset \fg(N)$ spanned by the unit diagonal matrices. Pick a basis 
\[\{\varepsilon^\sharp_i(1),\ldots,\varepsilon^\sharp_i(\sigma^\sharp(N))\}_{i=1}^\kappa \cup \{\varepsilon^\flat_i(1),\ldots,\varepsilon^\flat_i(\sigma^\flat(N))\}_{i=1}^\beta\] for $\fh(N)^*$ dual to the basis of unit diagonal matrices.
Identify elements $\nu \in \fh(N)^*$ with tuples \[(\b{\nu}^\sharp_1,\ldots,\b{\nu}^\sharp_\kappa,\b{\nu}^\flat_1,\ldots,\b{\nu}^\flat_\beta),\] where each $\b{\nu}_i^\sharp = (\nu_{i,1}^\sharp,\ldots,\nu^\sharp_{i,\sigma_i^\sharp(N)}) \in \fk^{\sigma_i^\sharp(N)}$ and $\b{\nu}_i^\flat = (\nu_{i,1}^\flat,\ldots,\nu^\flat_{i,\sigma_i^\flat(N)}) \in \fk^{\sigma_i^\flat(N)}$.

We define the \textit{weight lattice} $\fh(N)^*_{\ZZ}$ to be the set of weights $\nu$ for which all numbers $\nu_{ij}^\sharp$ and $\nu_{ij}^\flat$ above are integers. The space $\fh(N)^*$ has a symmetric pairing $(-,-)$ defined so that the basis  $\{\varepsilon_i^\sharp(j),\varepsilon_k^\flat(\ell)\}$ is orthonormal with respect to $(-,-)$.

The root system $R \subset \fh(N)_{\ZZ}^*$ of $\fg(N)$ is the union of the sets of weights \begin{align*}
R^+_0 &:= \{\varepsilon^\sharp_i(j) - \varepsilon^\sharp_k(\ell) \ |\ (i < k) \text{ or } (i = k, j < \ell)\} \cup \{\varepsilon^\flat_i(j) - \varepsilon^\flat_k(\ell) \ |\ (i < k) \text{ or } (i = k, j < \ell)\} \\[5pt]
R^-_0 &:= \{\varepsilon^\sharp_i(j) - \varepsilon^\sharp_k(\ell) \ |\ (i > k) \text{ or } (i = k, j > \ell)\} \cup \{\varepsilon^\flat_i(j) - \varepsilon^\flat_k(\ell) \ |\ (i > k) \text{ or } (i = k, j > \ell)\}
\\[5pt] R^+_1 &:= \{\varepsilon^\sharp_i(j) - \varepsilon^\flat_k(\ell)\} \\[5pt]
R^-_1 &:= \{\varepsilon^\flat_k(\ell) - \varepsilon^\sharp_i(j)\}.
\end{align*}
The roots belonging to $R^+_0 \cup R^+_1$ are called \textit{positive}, and all other roots are called \textit{negative}. The \textit{Weyl vector} $\rho \in \frac{1}{2}\fh(N)^*_{\ZZ}$ is defined to be the half-sum of the positive roots.

Let $\fl(N) \subset \fg(N)$ denote the Levi subalgebra consisting of endomorphisms stabilizing the decomposition (\ref{eqn:super_decomp}) and define the parabolic subalgebra $\fp(N) = \fl(N) + \fb(N)$. 

\begin{defn}
The category $\cO^{\b{c},\b{\sigma}}(N)$ is the full subcategory of $\fg(N)$-modules $M$ that are locally finite-dimensional over $\fp(N)$ and where the Cartan subalgebra $\fh(N)$ acts semisimply, such that all weights lie in $\fh(N)_{\ZZ}^*$ and all weight spaces are finite-dimensional.
\end{defn}

\begin{defn}
A \textit{parabolic highest weight} is a weight $\nu \in \fh(N)^*_{\ZZ}$ such that each $\b{\nu}_i^\sharp$ and $\b{\nu}_i^\flat$ is a \textit{non-increasing} sequence of integers. Let us write $\sDecorate{\Upsilon} \subset \fh(N)^*_{\ZZ}$ for the set of parabolic highest weights. We will be interested in the set $\cTheta := \rho + \sDecorate{\Upsilon}$ of \textit{$\rho$-shifted parabolic highest weights}.
\end{defn}

\begin{rem}\label{rem:strict_decrease}
For any $\vartheta \in \fh(N)^*$, observe that $\vartheta \in \cTheta$ if and only if $\vartheta \in \frac{1}{2}\fh(N)^*_{\ZZ}$ and
\begin{itemize}
\item for each $i = 1,\ldots,\kappa$ and $j = 1,\ldots, \sigma^\sharp_i(N)-1$, we have $\vartheta_{ij}^\sharp - \vartheta_{i,{j+1}}^\sharp > 0$, and
\item for each $i = 1,\ldots,\beta$ and $j = 1,\ldots, \sigma^\flat_i(N)-1$, we have  $\vartheta_{ij}^\flat - \vartheta_{i,{j+1}}^\flat > 0$. 
\end{itemize}
\end{rem}

The irreducible representations of $\fl(N)$ are in bijection with $\sDecorate{\Upsilon}$. For any $\vartheta \in \cTheta$, we write $\sDecorate{X}(\vartheta-\rho)$ for the corresponding irreducible $\fl(N)$-module. Define the \textit{parabolic Verma module}
\[
\cDelta(\vartheta) := U(\fg(N)) \ox_{\fp(N)} X_{N}(\vartheta-\rho)
\] by inflating $\sDecorate{X}(\vartheta-\rho)$ to a $\fp(N)$-module and then applying the parabolic induction functor. Observe the $\rho$-shift in the definition. The following proposition is standard.

\begin{prop}
The Verma module $\sDecorate{\Delta}(\vartheta)$ has a unique simple quotient $\cIrr(\vartheta)$. Moreover, the assignment $\vartheta \mapsto \cIrr(\vartheta)$ gives a bijection between $\cTheta$ and the set of simple objects in $\cParO$. 
\end{prop}

We recall the linkage principle for $\rho$-shifted parabolic highest weights. Identify $\nu \in \fh(N)^*$ with $d(N)$-tuples of complex numbers by concatenating the sequences $\b{\nu}_1^\sharp,\ldots,\b{\nu}_\kappa^\sharp,\b{\nu}_1^\flat,\ldots,\b{\nu}_\beta^\flat$ in order. Under this identification, we allow the symmetric group $\mathfrak{S}_{d(N)}$ to act on $\fh(N)^*$ naturally.

For $\alpha \in R^+$, the reflection $s_\alpha$ acts on $\nu \in \fh(N)^*$ via $s_\alpha\nu := \nu - (\nu,\alpha^\vee)\alpha.$ Identifying $\fh(N)^*$ with $\fk^{d(N)}$, we identify $s_\alpha$ with a certain transposition in $\mathfrak{S}_{d(N)}$, which we also denote $s_\alpha$.

For $\vartheta \in \cTheta$ and $\alpha \in R^+$ satisfying $(\vartheta,\alpha^\vee) > 0$, Remark \ref{rem:strict_decrease} implies that there exists a unique \[w \in \mathfrak{S}_{\sigma_1^\sharp(N)} \times \cdots \times \mathfrak{S}_{\sigma_\kappa^{\sharp}(N)}\times \mathfrak{S}_{\sigma_1^\flat(N)} \times \cdots \mathfrak{S}_{\sigma_\beta^{\flat}(N)} \subset \mathfrak{S}_{d(N)}\] such that $w(s_\alpha \vartheta) \in \cTheta$. We will write $(s_\alpha \vartheta)_+ := w(s_\alpha \vartheta)$.

\begin{rem}
When $\alpha = \varepsilon^\sharp_i(j) - \varepsilon^\sharp_k(\ell) \in R_0^+$, we have $(\vartheta,\alpha^\vee) > 0$ if and only if $\vartheta^\sharp_{ij} > \vartheta^\sharp_{k\ell}$. On the other hand, when $\alpha = \varepsilon^\sharp_i(j) - \varepsilon^\flat_k(\ell) \in R_1^+$, we have $(\vartheta,\alpha^\vee) > 0$ if and only if $\vartheta^\sharp_{ij} > \vartheta^\flat_{k\ell}$ . 
\end{rem}

\begin{defn}
Let $A_{\vartheta} := \{\alpha \in R^+ \ |\ (\vartheta,\alpha^\vee) > 0\}$. We say that $\vartheta'$ is \textit{linked} to $\vartheta$, denoted $\vartheta' \uparrow \vartheta$ if there exists $\alpha \in A_{\vartheta}$ such that $(s_\alpha\vartheta)_+ = \vartheta'$ and $\vartheta'$ precedes $\vartheta$ in the dominance order on weights. We define the \textit{linkage order} on $\cTheta$ to be the transitive closure $\preceq$ of the relation $\uparrow$. The linkage order $\preceq$ is a partial order on $\cTheta$. 
\end{defn}

The following result is a standard consequence of the Jantzen sum formula.

\begin{lem}
    For any $\vartheta,\zeta\in \cTheta$, we have $[\cDelta(\vartheta):\cIrr(\zeta)] \neq 0$ only if $\zeta \preceq \vartheta$. Moreover, the category $\cParO$ is a direct sum of finite highest weight categories, each of whose posets is given by a poset component of $(\Theta_N^{\b{c},\b{\sigma}},\preceq)$. 
\end{lem}

\subsection{Combinatorics of Virtual Multipartitions}\label{subsec:virtual}

Fix a positive integer $m$ and choose the integer $N$ so that $\sigma^\sharp_i(N) > m$ for all $i = 1,\ldots,\kappa$ and $\sigma^\flat_i(N) > m$ for all $i = 1,\ldots,\beta$. Our ultimate goal will be to embed the poset $\cP^a(m)$ from Corollary \ref{cor:restricted_partitions} as a poset ideal in $(\cTheta, \preceq)$ when $N$ is sufficiently large. To do so, we will need to reinterpret $\rho$-shifted super-parabolic highest weights as \textit{virtual multipartitions}, in a construction inspired by \cite[Section 2.3]{losev_vv_conjecture}.

\begin{defn}\label{defn:virtual_partition}
A \textit{virtual partition} $\mu^v$ is a diagram consisting of a finite set of horizontal rows of unit squares in the upper half-plane $\RR \times \RR_{\geq 0}$ satisfying the following properties.
\begin{enumerate}
\item The rows are stacked on top of each other, with the bottom-most row lying on the $x$-axis.
\item Each row extends infinitely far to the left but terminates to the right.
\item The \textit{virtual length} of a row is the $x$-coordinate of the right-most box in the row. The virtual lengths of the rows must be integers that form a non-increasing sequence when read from the bottom row to the top row.
\end{enumerate}
We will write $\mu^v = (\mu^v_1,\ldots,\mu^v_r)$, where $r$ is the number of rows in $\mu^v$ and $\mu^v_i$ is the virtual length of the $i$th row in $\mu^v$ when read from the bottom to the top. In particular, if $\mu^v_r > 0$, then $\mu^v$ is a partition in the conventional sense.

A \textit{$(\b{c},\b{\sigma})$-virtual multipartition} is a collection $\vartheta^v = (\vartheta^{\sharp,v}_{1},\ldots,\vartheta^{\sharp,v}_{\kappa}, \vartheta^{\flat,v}_{1},\ldots,\vartheta^{\flat,v}_{\beta})$ of virtual partitions such that each $\vartheta^{\sharp,v}_i$ has exactly $\sigma_i^\sharp(N)$ rows and each $\vartheta^{\flat,v}_k$ has exactly $\sigma_k^\flat(N)$ rows. Let us write $\cP^{\mathrm{vir}}$ for the collection of all $(\b{c},\b{\sigma})$-virtual multipartitions.
\end{defn}

We will assign each $\vartheta \in \cTheta$ to a $(\b{c},\b{\sigma})$-virtual multipartition \[\vartheta^v = (\vartheta^{\sharp,v}_{1},\ldots,\vartheta^{\sharp,v}_{\kappa}, \vartheta^{\flat,v}_{1},\ldots,\vartheta^{\flat,v}_{\beta}) \in \cP^{\mathrm{vir}},\]  where \[\vartheta_i^{\sharp,v} := \left(\vartheta_{ij}^\sharp - \frac{d(N)-1}{2} + j - 1\right)_{j=1}^{\sigma^\sharp_i(N)}, \quad \vartheta_k^{\flat,v} := \left(\vartheta_{k\ell}^\flat - \frac{d(N)-1}{2} + \ell - 1\right)_{\ell=1}^{\sigma^\flat_k(N)}\] for each $i =1,\ldots,\kappa$ and $k = 1,\ldots,\beta$. Each coordinate defined above is an integer, so $\vartheta^v$ is indeed a virtual multipartition. The assignment $\vartheta \mapsto \vartheta^v$ gives a bijection $\cTheta \simeq \cP^{\mathrm{vir}}$. 

On the other hand, we have a \textit{virtual realization map} $V_{\b{c}}: \cP^a(m) \to \cP^{\mathrm{vir}}$. For any \[\blambda = (\lambda_1^\sharp,\ldots,\lambda_\kappa^\sharp,\lambda_1^\flat,\ldots,\lambda_\beta^\flat) \in \cP^a(m),\] with the notation $\ell_i^\sharp := \ell(\lambda_i^\sharp)$ and $\ell_i^\flat := \ell(\lambda_i^\flat)$, we define \[V_{\b{c}}(\blambda) = (V(\lambda_1^\sharp),\ldots, V(\lambda_\kappa^\sharp),V(\lambda_1^\flat),\ldots, V(\lambda_\beta^\flat)),\] where
\begin{align*}
V(\lambda_i^\sharp) := \begin{cases} (\lambda_{i,1}^\sharp + \sigma_i^\sharp,\lambda_{i,2}^\sharp + \sigma_i^\sharp,\ldots,\lambda_{i,\ell^\sharp_i}^\sharp + \sigma_i^\sharp,\sigma_i^\sharp ,\ldots,\sigma_i^\sharp) & \text{if $\b{c}$ is lower admissible} \\
(N,\ldots,N,-\lambda_{i,\ell^\sharp_i}^\sharp + N,\ldots,-\lambda_{i,2}^\sharp + N,-\lambda_{i,1}^\sharp + N) & \text{if $\b{c}$ is upper admissible}
\end{cases} \\[5pt]
V(\lambda_i^\flat) := \begin{cases} (N,\ldots,N,-\lambda_{i,\ell^\flat_i}^\flat + N,\ldots,-\lambda_{i,2}^\flat + N,-\lambda_{i,1}^\flat + N) & \text{if $\b{c}$ is lower admissible} \\
(\lambda_{i,1}^\flat + \sigma_i^\flat,\lambda_{i,2}^\flat + \sigma_i^\flat,\ldots,\lambda_{i,\ell^\flat_i}^\flat + \sigma_i^\flat,\sigma_i^\flat ,\ldots,\sigma_i^\flat) & \text{if $\b{c}$ is upper admissible} 
\end{cases} 
\end{align*}

Let us write $\weight(\blambda) \in \cTheta$ for the weight corresponding to $V(\blambda) \in \cP^{\mathrm{vir}}$ under the identification $\cP^{\mathrm{vir}} \simeq \cTheta$. Thus, we obtain an embedding $\weight: \cP^a(m) \hookrightarrow \cTheta$. We will prove that, for $N \gg 0$, the restriction of $\preceq$ to the image of $\weight$ realizes $\cP^a(m)$ as a poset ideal in $\cTheta$. 

\begin{lem}\label{lem:multipartition_good_swaps}
For $\blambda \in \cP^a(m)$, $\alpha \in R_0^+$, we have $(s_\alpha\weight(\blambda))_+ \in \weight(\cP^a(m'))$ for some $m' \geq 0$.  
\end{lem}

\begin{proof}
This lemma is essentially \cite[Proposition A.6.1]{VV} -- we adapt their proof to our language for the reader's convenience. Let us write \[\vartheta := (s_\alpha\weight(\blambda))_+ = (\vartheta_1^\sharp,\ldots,\vartheta_\kappa^\sharp,\vartheta_1^\flat,\ldots,\vartheta_\beta^\flat)\] as well as the corresponding virtual multipartition\[\vartheta^v := (s_\alpha\weight(\blambda))^v_+ = (\vartheta_1^{v,\sharp},\ldots,\vartheta_\kappa^{v,\sharp},\vartheta_1^{v,\flat},\ldots,\vartheta_\beta^{v,\flat}).\] 
We consider the case $\alpha = \varepsilon_i^\sharp(j) - \varepsilon_k^\sharp(\ell)$ and $\b{c}$ is lower admissible -- all other cases are identical. For each $t = 1,\ldots,\kappa$, we show that the virtual partition $\vartheta^{v,\sharp}_t$ corresponding to $\vartheta^{\sharp}_t$ has the form
\[
\vartheta^{v,\sharp}_t = (\eta_{t,1}^\sharp + \sigma_t^\sharp, \eta_{t,2}^\sharp + \sigma_t^\sharp,\ldots, \eta_{t,r}^\sharp + \sigma_t^\sharp,\sigma_t^\sharp, \ldots,\sigma_t^\sharp),
\]
where $(\eta_{t,1}^\sharp,\ldots,\eta_{t,r}^\sharp)$ is a partition. Whenever $t \neq i,k$, this claim is immediate since $\vartheta^{v,\sharp}_t = V(\blambda_t^\sharp)$ in that case. Thus, suppose $t = i$ or $j$. In this case, we claim it suffices to show \begin{align}\label{eqn:inequality_good_swap}\vartheta^\sharp_{t,\sigma_t^\sharp(N)} \geq \sigma_t^\sharp + \frac{d(N)-1}{2} - \sigma_t^\sharp(N) + 1 = \frac{d(N)+1}{2} + N.\end{align} Indeed, the inequality (\ref{eqn:inequality_good_swap}) is equivalent to the condition that $\vartheta^{\sharp,v}_{t,\sigma_t^\sharp} \geq \sigma_t^\sharp$ and hence $\vartheta^{\sharp,v}_{t,s} \geq \vartheta^{\sharp,v}_{t,\sigma_t^\sharp} \geq \sigma_t^\sharp$ for all $s = 1,\ldots,\sigma_t^\sharp$. By definition, observe that \begin{align*}
    \{\vartheta^{\sharp}_{i,s}\}_{s=1}^{\sigma_i^\sharp(N)} = \{\weight(\blambda)^\sharp_{i,s}\}_{s=1}^{\sigma_i^\sharp(N)} \setminus \{\weight(\blambda)^\sharp_{i,j}\} \cup \{\weight(\blambda)^\sharp_{k,\ell}\}.
\end{align*}
Then, we see that 
\[
\vartheta^{\sharp}_{i,\sigma_i^\sharp} = \min_s(\vartheta^{\sharp}_{i,s}) \geq \min(\{\weight(\blambda)^{\sharp}_{i,s}\}_s \cup \{\weight(\blambda)^\sharp_{k,\ell}\}).
\]
We know that $\weight(\blambda)^{\sharp}_{i,s}$ and $\weight(\blambda)^\sharp_{k,\ell}$ are both bounded below by $\frac{d(N) + 1}{2} + N$ since the weight $\weight(\blambda)$ comes from a multipartition. Hence, we see that $\vartheta^\sharp_{i,\sigma_i^\sharp}$ satisfies (\ref{eqn:inequality_good_swap}). An identical argument shows that $\vartheta^\sharp_{j,\sigma_j^\sharp}$ satisfies this inequality as well. 
\end{proof}

\begin{lem}\label{lem:multipartition_bad_swaps}
    Suppose $N >  m + 1 + \max_{i=1}^\kappa(|\sigma_i^\sharp|) + \max_{i=1}^\beta(|\sigma_i^\flat|)$. For any $\blambda \in \cP^a(m)$ and $\alpha \in R_1^+$, we have $(s_\alpha \weight(\blambda))_+ \in \weight(\cP^a(m'))$ for some $m' \geq 0$. 
\end{lem}

\begin{proof}
    This proof is similar to the proof of Lemma \ref{lem:multipartition_good_swaps}. Suppose $\alpha = \varepsilon_i(j)^\sharp - \varepsilon_k(\ell)^\flat$. We consider the case where $\b{c}$ is lower admissible -- the proof is the same in the upper admissible case. Let us write \[\vartheta := (s_\alpha\weight(\blambda))_+ = (\vartheta_1^\sharp,\ldots,\vartheta_\kappa^\sharp,\vartheta_1^\flat,\ldots,\vartheta_\beta^\flat)\] as well as the corresponding virtual multipartition\[\vartheta^v := (s_\alpha\weight(\blambda))^v_+ = (\vartheta_1^{v,\sharp},\ldots,\vartheta_\kappa^{v,\sharp},\vartheta_1^{v,\flat},\ldots,\vartheta_\beta^{v,\flat}).\] 
    By the same argument that we made for Lemma \ref{lem:multipartition_good_swaps}, we must establish
    \begin{align*}
        \vartheta_{i,\sigma_i^\sharp(N)}^\sharp \geq N + \frac{d(N) + 1}{2}, \quad \vartheta_{k,1}^\flat \leq N + \frac{d(N)-1}{2}.
    \end{align*}
    For the latter inequality, we use the fact that
    \[
    \{\vartheta_{k,s}^\flat\}_{s=1}^{\sigma_k^\flat(N)} \subset \{\weight(\blambda)_{k,s}^\flat\}_{s=1}^{\sigma_k^\flat(N)} \cup \{\weight(\blambda)_{ij}^\sharp\}
    \]
    so that
    \[
    \vartheta_{k,1}^\flat = \max(\{\vartheta_{k,s}^\flat\}_{s=1}^{\sigma_k^\flat(N)}) \leq \max(\{\weight(\blambda)_{k,s}^\flat\}_{s=1}^{\sigma_k^\flat(N)} \cup \{\weight(\blambda)_{ij}^\sharp\}).
    \]
    On the one hand, we already know that $\weight(\blambda)_{k,s}^\flat \leq N + \frac{d(N) -1}{2}$ since this weight comes from a multipartition. On the other hand, \[\weight(\blambda)_{ij}^\sharp = \lambda_{ij}^\sharp + \sigma_i^\sharp + \frac{d(N) - 1}{2} - j + 1 \leq m + \sigma_i^\sharp + \frac{d(N) - 1}{2} + 1 \leq N + \frac{d(N)- 1}{2}\] by the assumed bound on $N$. The inequality $\vartheta_{i,\sigma_i^\sharp(N)}^\sharp \geq N + \frac{d(N) + 1}{2}$ is proved similarly.
\end{proof}

Recall the inverse dominance order $\leq$ on $\cP^a$ associated with $(\b{\sigma},\b{c})$ from Definition \ref{defn:inv_dominance}. First, we show that the inverse dominance order refines the linkage order on $\weight(\cP^a(m))$. In turn, we construct a Serre subquotient of $\cO^{\b{c},\b{\sigma}}(N)$ corresponding to $\weight(\cP^a(m))$ that is highest weight with respect to the inverse dominance order.

\begin{prop}\label{prop:construct_refine}
Pick some $\blambda \in \cP^a(m)$ and $\alpha \in R^+$. Write $\b{\mu} \in \cP^a(m')$ (for some $m' \geq 0$) for the genuine multipartition corresponding to $(s_\alpha\weight(\blambda))_+$. If $\b{c}$ is lower admissible and $(\weight(\blambda),\alpha) > 0$, then $\b{\mu} \leq \b{\lambda}$. On the other hand, if $\b{c}$ is upper admissible and $(\weight(\blambda),\alpha) < 0$, then $\b{\mu} \geq \b{\lambda}$.
\end{prop}

\begin{proof}
Given a virtual multipartition $\vartheta^v = (\vartheta^{v,\sharp}_1,\ldots,\vartheta^{v,\sharp}_\kappa,\vartheta^{v,\flat}_1,\ldots,\vartheta^{v,\flat}_\beta)$, we say that a box $b \in \vartheta^v$ has \textit{virtual content} $\mathrm{ct}(b) := x - y$, where $(x,y)$ are the coordinates of the top right corner of $b$.

We prove the claim when $\b{c}$ is lower admissible -- the argument is exactly the same in the upper admissible case. Consider arbitrary $\vartheta \in \cTheta$ such that $(\vartheta,\alpha) > 0$. By definition, $(s_\alpha\vartheta)_+$ is obtained from $\vartheta$ by permuting its coordinates. We shall understand the action of this permutation on the corresponding virtual multipartition. 

Write $s_\alpha = \varepsilon_i^?(j) - \varepsilon_k^\invques(\ell)$, where $(?,\invques) = (\sharp,\sharp), (\sharp,\flat)$, or $(\flat,\flat)$. The condition that $(\vartheta,\alpha) > 0$ means that $\vartheta_{ij}^? > \vartheta_{k\ell}^\invques$. Moreover, the weight $s_\alpha\vartheta$ is obtained from $\vartheta$ by swapping $\vartheta_{ij}^?$ and $\vartheta^\invques_{k\ell}$. On the level of ``improper" virtual multipartitions (i.e., a set of ``improper virtual partitions'' where the lengths of the rows may not be non-increasing), we see that \begin{align*}(s_\alpha \vartheta)^{v,?}_i &= (\vartheta^{v,?}_{i,1},\ldots,\vartheta^{v,?}_{i,j-1}, \vartheta^{v,?}_{ij} - (\vartheta_{ij}^? - \vartheta_{k\ell}^\invques), \vartheta_{i,j+1}^{v,?},\ldots,\vartheta_{i,\sigma_i^\sharp(N)}^{v,?}) \\[5pt] 
(s_\alpha \vartheta)^{v,\invques}_k &= (\vartheta^{v,\invques}_{k,1},\ldots,\vartheta^{v,\invques}_{k,\ell-1}, \vartheta^{v,\invques}_{k\ell} + (\vartheta_{ij}^? - \vartheta_{k\ell}^\invques), \vartheta_{k+1}^{v,\invques},\ldots,\vartheta_{k,\sigma_k^\invques(N)}^{v,\invques}).\end{align*}
Thus, we removed $\vartheta_{ij}^? - \vartheta_{k\ell}^\invques$ boxes from $\vartheta^{v,?}_i$ and added them to $\vartheta^{v,\invques}_k$. Observe that the removed boxes and added boxes have the same virtual contents: the left-most removed box has content \[\vartheta^{v,?}_{ij} - (\vartheta_{ij}^? - \vartheta_{k\ell}^\invques) + 1 - j = \vartheta_{k\ell}^\invques - \frac{d(N)-1}{2} = \vartheta_{k\ell}^{v,\invques} + 1 - k,\] which is the content of the left-most added box.

Let us also comment on how to obtain $(s_\alpha\vartheta)_+^v$ from the improper virtual multipartition $(s_\alpha\vartheta)^v$. For each row $(s_\alpha\vartheta)_{i,s}^{v,?}$ where $(s_\alpha\vartheta)_{i,s}^{v,?} > (s_\alpha\vartheta)_{i,s-1}^{v,?}$, we move the last $(s_\alpha\vartheta)_{i,s-1}^{v,?} - (s_\alpha\vartheta)_{i,s}^{v,?}$ from the row $s$ and add them to the row $s-1$. Of particular importance, the virtual contents of the boxes are preserved under this procedure, since each box stays in the same diagonal that it originally belonged to. Repeating this procedure eventually results in a genuine virtual partition. We perform an identical procedure to $(s_\alpha\vartheta)_{k}^{v,\invques}$ to obtain a genuine virtual partition. Then, the resulting virtual multipartition is precisely $(s_\alpha\vartheta)_+^v$. 


Suppose now that $\vartheta^v = V(\blambda)$ for some $\blambda \in \cP^a(m')$. We relate the contents of boxes in $\blambda$ to contents of boxes in $\vartheta^v$. Removing or adding a box with content $\gamma$ in some $\blambda_i^\sharp$ is the same thing as removing or adding a box in $\vartheta^{v,\sharp}_i$ with virtual content $\gamma + \sigma_i^\sharp$. On the other hand, removing (resp. adding) a box with content $\gamma$ in some $\blambda_k^\flat$ is the same thing as \textit{adding} (resp. \textit{removing}) a box in $\vartheta^{v,\flat}_k$ with virtual content $-\gamma -\sigma_k^\flat$. Thus, for each $b \geq 0$ and each $N = 1,\ldots,a$, we have
\[
\sum_{k=1}^N (-1)^{c_k}|\mu_k|_{(-1)^{c_k}(b-\sigma_k)} \leq \sum_{k=1}^N (-1)^{c_k}|\lambda_k|_{(-1)^{c_k}(b-\sigma_k)},
\]
with equality for $N = a$. In other words, $\b{\mu} \leq_S \b{\lambda}$ in the inverse dominance order. 
\end{proof}

\begin{cor}
Let us make the same assumptions as in Lemma \ref{lem:multipartition_bad_swaps}. In the lower (resp. upper) admissible case, the subset $\weight(\cP^a(m))$ is a poset ideal (resp. coideal) in $\cTheta$ with respect to the linkage order $\preceq$.
\end{cor}

\begin{defn}\label{defn:partition_subcategory}
Suppose $\b{c}$ is lower (resp. upper) admissible. Let $\cParOm{m}$ denote the Serre subcategory (resp. quotient) of $\cParO$ corresponding to the ideal (resp. coideal) $\weight(\cP^a(m))$.
\end{defn} 

For each $\blambda \in \cP^a(m)$, write $\Delta(\blambda)$ for the corresponding standard object in $\cParOm{m}$. 

\begin{cor}
    Under the same assumptions as in Lemma \ref{lem:multipartition_bad_swaps}, the subcategory $\cParOm{m}$ is a highest weight category with poset $(\cP^a_m, \leq_S)$, where $\leq_S$ is the inverse dominance order.
\end{cor}

\subsection{Constructing Restricted Admissible Fock Space Tensor Product Categorifications}\label{subsec:constructing_restricted_fock}

Now, we construct a restricted categorical action on the highest weight category $\cParOm{m}$ from Definition \ref{defn:partition_subcategory}. Continue to make the same assumptions on $N$ as in the statement of Lemma \ref{lem:multipartition_bad_swaps}.

Note that $\cParO$ is preserved under the functors $F := V(N) \ox -$ and $E := V(N)^* \ox -$, so Example \ref{ex:categorical_type_a_glX} equips $\cParO$ with a type A categorical action $(E,F,x,\tau)$, where $x \in \End(F)$ is given by the action of the so-called tensor Casimir $\Omega \in U(\fg(N))$ and $\tau \in \End(F^2)$ is given by the symmetric braiding.
For any weight \[\vartheta = (\vartheta_1^\sharp,\ldots,\vartheta_\kappa^\sharp,\vartheta_1^\flat,\ldots,\vartheta_\beta^\flat) \in \cTheta,\] observe that $F\cDelta(\vartheta)$ has a Verma filtration with sections $\Delta(\vartheta^{\sharp}(F,i,j))$ and $\Delta(\vartheta^{\flat}(F,r,\ell))$ for $i = 1,\ldots,\kappa$, $j = 1,\ldots,\sigma_i^\sharp(N)$, $r = 1,\ldots,\beta$, and $\ell = 1,\ldots,\sigma_r^\flat(N)$, where
\begin{align*}
    \vartheta(\sharp,F,i,j) &:= (\vartheta^\sharp_{1},\ldots,\vartheta^\sharp_{i-1},\vartheta^\sharp_i[F,j], \vartheta^\sharp_i, \ldots,\vartheta^\sharp_\kappa,\vartheta^\flat_1,\ldots,\vartheta^\flat_\beta), \\[5pt]
    \vartheta^\sharp_i[F,j] &:= (\vartheta_{i,1}^\sharp,\ldots,\vartheta_{i,j-1}^\sharp, \vartheta_{i,j}^\sharp+1,\vartheta^\sharp_{i,j+1},\ldots,\vartheta^\sharp_{i,\sigma_i^\sharp(N)})
\end{align*}
if such a weight $\vartheta^\sharp(F,i,j)$ is indeed in $\cTheta$; otherwise, we agree that the corresponding section does not appear in the filtration. Similarly, 
\begin{align*}
    \vartheta(\flat,F,r,\ell) &:= (\vartheta^\sharp_{1},\ldots,\vartheta^\sharp_{\kappa},\vartheta_1^\flat,\ldots,\vartheta_{r-1}^\flat, \vartheta^\flat_{r}[F,\ell], \vartheta^\flat_{r+1},\ldots,\vartheta^\flat_\beta), \\[5pt]
    \vartheta^\flat_{r}[F,\ell] &:= (\vartheta_{r,1}^\flat,\ldots,\vartheta_{r,\ell-1}^\flat, \vartheta_{r,\ell}^\flat+1,\vartheta^\flat_{r,\ell+1},\ldots,\vartheta^\flat_{r,\sigma_r^\flat(N)})
\end{align*}
if such a weight $\vartheta(\flat, F,r,\ell)$ is indeed in $\cTheta$; otherwise, we agree that the corresponding section does not appear in the filtration. These sections appear in the following order: all sections $\Delta(\vartheta(\sharp,F,i,j))$ occur higher than the sections $\Delta(\vartheta(\flat,F,r,\ell))$, a section $\Delta(\vartheta(\sharp,F,i,j))$ occurs higher than $\Delta(\vartheta(\sharp,F,i',j'))$ if and only if $i < i'$ or $i = i'$ and $j < j'$, and a section $\Delta(\vartheta(\flat,F,r,\ell))$ occurs higher than $\Delta(\vartheta(\flat,F,r',\ell'))$ if and only if $r < r'$ or $r = r'$ and $\ell < \ell'$. Any weight that appears in this filtration will be called \textit{$F$-valid}.

Similarly, the module $E\cDelta(\vartheta)$ has a filtration with sections $\Delta(\vartheta(\sharp,E,i,j))$ and $\Delta(\vartheta(\flat,E,r,\ell))$ for $i = 1,\ldots,\kappa$, $j = 1,\ldots,\sigma_r^\sharp(N)$, $r = 1,\ldots,\beta$, and $\ell = 1,\ldots,\sigma_r^\flat(N)$, where the label
$\vartheta(\sharp,E,i,j)$ (resp. $\vartheta(\flat, E, r,\ell)$) is defined similarly to $\vartheta(\sharp, F, i,j)$ (resp. $\vartheta(\flat,F,r,\ell)$), except that the ``adjusted" coordinate $\theta_{i,j}^\sharp + 1$  (resp. $\theta_{r,\ell}^\flat + 1$) in $\vartheta_i^\sharp[F,j]$ (resp. $\vartheta_r^\flat[F,\ell]$) should be replaced by $\theta_{i,j}^\sharp - 1$ (resp. $\theta_{r,\ell}^\flat$). The ordering of these sections is the same as the ordering of the analogous sections in the filtration of $F\cDelta(\vartheta)$. Any weight appearing in this filtration will be called \textit{$E$-valid}.


\begin{notation}
    For any $\blambda \in \cP^a(m)$, we write $\blambda(\sharp,F,i,j) := \weight(\blambda)(\sharp,F,i,j)$. We similarly define $\blambda(\flat,F,r,\ell)$, $\blambda(\sharp,E,i,j)$, and $\blambda(\flat,E,r,\ell)$.
\end{notation}

\subsubsection{Lower Admissible Case} First assume that $\b{c}$ is lower admissible and fix some $\b{\lambda} \in \cP^a(m)$. We consider the restriction of the functors $E$ and $F$ to the highest weight subcategory $\cParOm{m}$.

For each $i = 1,\ldots,\kappa$, observe that any $F$-valid section $\blambda(\sharp,F,i,j)$ corresponds to the multipartition obtained from $\blambda$ by adding a single box to $\lambda_i^\sharp$. Similarly, for each $k = 1,\ldots,\beta$, any $F$-valid section $\blambda(\flat,F,j,\ell)$ corresponds to the multipartition obtained from $\blambda$ by removing a single box from $\lambda_j^\flat$, \textit{with the exception} of the weight $\blambda(\flat,F,j,1)$, which does \textit{not} come from a multipartition.

Let us compute the eigenvalues of $x$ on $\cDelta(\blambda)$. By definition, we must compute the action of \[\frac{1}{2}(\mathrm{coprod}(C_2) - C_2 \ox 1 - 1 \ox C_2)\] on $V(N) \ox \cDelta(\blambda)$, where $C_2$ is the quadratic Casimir for $\fg(N)$. Note $C_2$ acts on $V(N)$ by $d(N)$. 

Recall that the quadratic Casimir $C_2$ acts on any highest weight $\fg(N)$-module with \textit{$\rho$-shifted highest weight} by the scalar $(\vartheta,\vartheta) - (\rho,\rho)$. In particular, if $\blambda^+ = \blambda(\sharp,F,i,j)$ is an $F$-valid weight, then $x$ acts on the section $\cDelta(\blambda^+)$ by \begin{align*}\frac{1}{2}((\blambda^+,\blambda^+) - (\blambda,\blambda)-d(N)) &=  \frac{1}{2}((\varepsilon_i^\sharp(j), \varepsilon_i^\sharp(j))(2\lambda_{ij} + 2\sigma_i^\sharp + d(N) - 1 - 2j + 2 + 1) - d(N)) \\[5pt] &= \mathrm{ct}(\lambda_i^+ - \lambda_i) + \sigma_i^\sharp = \mathrm{ct}(\lambda_i^+ - \lambda_i) + \sigma_i,\end{align*} where $\mathrm{ct}(\lambda_i^+ - \lambda_i)$ is the content of the box added to $\lambda_i$. Similarly, if $\blambda^- = \blambda(\flat,F,j,\ell)$ is an $F$-valid weight \textit{corresponding to a multipartition}, then $x$ acts on the section $\cDelta(\blambda^-)$ by
\begin{align*}
\frac{1}{2}((\blambda^-,\blambda^-) - (\blambda,\blambda)-d(N)) &=  \frac{1}{2}((-2\lambda_{j,\sigma_j^\flat(N) - \ell + 1} + 2N + d(N) - 1 - 2\ell + 2 + 1) - d(N)) \\[5pt] &= -\mathrm{ct}(\lambda_j - \lambda_j^-) - \sigma_j^\flat = -\mathrm{ct}(\lambda_j - \lambda_j^-) + \sigma_{j + \kappa}
\end{align*}

Finally, $x$ acts on the sections $\cDelta(\blambda(\flat,F,j,1))$ by 
\begin{align*}
 \frac{1}{2}((2N + d(N) - 1 - 2 + 2 + 1) - d(N)) &= N.
\end{align*}

We claim the sections $\cDelta(\blambda(\flat,F,k,1))$ are the only sections in the aforementioned filtration where $x$ acts by $N$. Indeed, the eigenvalues of all other sections are bounded above by \[m + \max_{i=1,\ldots,\kappa} |\sigma_i^\sharp| + \max_{j=1,\ldots,\beta} |\sigma_j^\flat|,\] which is strictly less than $N$ (by Lemma \ref{lem:multipartition_bad_swaps}). Consider the functor $\ul{F} := \bigoplus_{\alpha \neq N} F_\alpha$. Recall the poset ideals $\cP^a(m_1,m_2) \subset \cP^a(m)$ and write $\cParOsub{m}{N,m_1,m_2} \subset \cParOm{m}$ for the corresponding Serre subcategory. We see that $\ul{F}$ restricts to functors \[\ul{F}: \cParOsub{m}{N,m_1,m_2} \to \cParOsub{m}{N,m_1+1,m_2-1}.\] The endomorphisms $x \in \End(F)$ and $\tau \in \End(F^2)$ restrict to endomorphisms $\ul{x} \in \End(\ul{F})$ and $\ul{\tau} \in \End(\ul{F}^2)$ still satisfying the degenerate affine Hecke relations (since they restrict to endomorphisms of each eigenfunctor).

Similarly, any $E$-valid section $\blambda(\flat,E,k,\ell)$ corresponds to the multipartition obtained from $\blambda$ by adding a single box to $\lambda_k^\flat$ for each $k = 1,\ldots,\beta$. On the other hand, for each $i = 1,\ldots,\kappa$, any $E$-valid section $\blambda(\sharp,E,i,j)$ corresponds to the multipartition obtained from $\blambda$ by removing a single box from $\lambda_i^\sharp$, \textit{with the exception} of the weight $\blambda(\sharp,E,i,\sigma_i^\sharp(N))$, which does \textit{not} come from a multipartition. The eigenvalue of $x$ on each $\Delta(\blambda(\flat,E,k,\ell))$ is $-\operatorname{ct}(\lambda_k^+ - \lambda_k) - \sigma_k^\flat$, where $\lambda_k^+$ is the partition obtained from $\lambda_k$ by adding a box in row $\sigma_k^\flat(N) - \ell + 1$. Similarly, the  eigenvalue of $x$ on each $\Delta(\blambda(\sharp,E,i,j))$ for $j \neq \sigma_i^\sharp(N)$ is $\operatorname{ct}(\lambda_i - \lambda_i^-) + \sigma_i^\sharp$. Finally, the eigenvalue of $x$ on the sections $\cDelta(\blambda(\sharp,E,i,\sigma_i(N)))$ is given by
\[
-\frac{1}{2}(2\sigma_i^\sharp - 1 - 2(\sigma_i^\sharp + N) + 2 - 1) = N
\]
as well. As before, the only sections in the filtration of $E\cDelta(\blambda)$ on which $x$ acts by $N$ are these ``bad'' sections $\cDelta(\blambda(\sharp,E,i,\sigma_i(N)))$. Thus, $\ul{E} := \bigoplus_{\alpha \neq N} E_\alpha$ restricts to functors \[\ul{E}: \cParOsub{m}{N,m_1,m_2} \to \cParOsub{m}{N,m_1-1,m_2+1}.\] 

\begin{thm}\label{thm:construction_lower_admissible}
With the data $(\ul{E},\ul{F},\ul{x},\ul{\tau})$, the category $\cParOm{m}$ is a level $m$ restricted lower admissible (multi-)Fock space tensor product categorification of type $\Xi = ((\b{\sigma},\b{c}))$.
\end{thm}

\begin{proof}
    Clear from the explicit descriptions of $\ul{F}$ and $\ul{E}$ on standard objects.
\end{proof}

\begin{cor}
    The categorification $\cParOm{m}$ is independent of $N \gg 0$ in the following sense: for $N,N' \gg 0$, there exists a strongly equivariant equivalence $\cParOsub{m}{N} \simeq \cParOsub{m}{N'}$ intertwining the labels of simple objects.
\end{cor}

\begin{proof}
An immediate consequence of Theorem \ref{thm:mixed_admissible_uniqueness} -- note $\cC(m,m) = \cC$ for a level $m$ restricted admissible MFTPC where $|\Xi| = 1$.
\end{proof}

\subsubsection{Upper Admissible Case}

It remains to handle the case where the sequence $\b{c}$ is upper admissible. We approach this case using Ringel duality. Let us define
\[
\b{c}^\vee := (c_a,c_{a-1},\ldots,c_1),\quad \b{\sigma}^\vee := (\sigma_a,\sigma_{a-1},\ldots,\sigma_1).
\]
Finally, we write $\ul{0}_a$ for a sequence $(0,\ldots,0)$ of $a$ zeros.

\begin{lem}\label{lem:slz_categorifications}
In the language of \cite{brundan_losev_webster}, the category $\cParO$ with the categorification data $(E,F,x,\tau)$ is a $\fsl_{\ZZ}$-tensor product categorification of type $(\b{\sigma} + N,\ul{0}_a)$ and level $a$. Similarly, $\cO^{\b{c}^\vee,\b{\sigma}^\vee}(N)$ has the natural structure of an $\fsl_{\ZZ}$-tensor product categorification of type $(\b{\sigma}^\vee + N,\ul{0}_a)$ and level $a$.
\end{lem}

\begin{proof}
This lemma is just a special case of \cite[Theorem 3.10]{brundan_losev_webster}.
\end{proof}

If $\cC$ is a direct sum of finite highest weight categories, then its Ringel dual $\cC^\vee$ is the direct sum of the Ringel duals of its highest weight summands. Finally, for each $m_1,m_2 \leq m$, we write $\cParOsub{m}{m_1,m_2}$ for the Serre quotient of $\cParOm{m}$ corresponding to the poset coideal $\cP^a(m_1,m_2)$. 

\begin{prop}
    Let us write $\cParO^\vee$ for the Ringel dual of $\cParO$. There exists an equivalence of highest weight categories $\cParO^\vee \simeq \cO^{\b{c}^\vee,\b{\sigma}^\vee}(N)$. Moreover, each highest weight subcategory $\cO^{\b{c}^\vee,\b{\sigma}^\vee}_m(m_1,m_2)$ can be identified with the Ringel dual of the quotient category $\cParOsub{m}{m_1,m_2}$.
\end{prop}

\begin{proof}
     Thanks to \cite[Theorem B]{brundan_losev_webster} and Lemma \ref{lem:slz_categorifications}, it suffices to equip $\cParO^\vee$ with the structure of an $\fsl_{\ZZ}$-tensor product categorification of type $(\b{\sigma}^\vee + N,\b{c}^\vee)$ and level $a$. This procedure can be done similarly to Lemma \ref{lem:restricted_multifock_ringel}. The second part follows since Ringel duality interchanges Serre subcategories corresponding to poset ideals with Serre quotient categories corresponding to poset coideals in the opposite order.
\end{proof}

In turn, Lemma \ref{lem:restricted_multifock_ringel} transports the lower admissible restricted categorification data $(\ul{E},\ul{F},x,\tau)$ on $\cO^{\b{c}^\vee,\b{\sigma}^\vee}_m(N)$ to upper admissible categorification data $(\ov{E},\ov{F},x,\tau)$ on the quotient $\cParOsub{m}{N}$. 

\begin{thm}
    With restricted categorical action given by the data $(\ov{E},\ov{F},\ov{x},\ov{\tau})$, the category $\cParOm{m}$ is a level $m$ restricted upper admissible (multi-)Fock space tensor product categorification of type $\Xi = ((\b{\sigma},\b{c}))$.
\end{thm}

\begin{proof}
    An immediate consequence of Lemma \ref{lem:restricted_multifock_ringel}.
\end{proof}

\begin{cor}
    The category $\cParOm{m}$ is independent of $N$ whenever $N$ is sufficiently large.
\end{cor}



\begin{defn}\label{defn:lower_admissible_stable}
Suppose $\b{c}$ is lower (resp. upper) admissible. We define $\cO^{\b{c},\b{\sigma}}_m := \cParOm{m}$ for sufficiently large $N$. For $m_1,m_2 \leq m$, we write $\cParOsub{m}{m_1,m_2}$ for the highest weight subcategory (resp. quotient) corresponding to the poset ideal (resp. coideal) $\cP^a(m_1,m_2)$. The standard, costandard, simple, and tilting objects in $\cO^{\b{c},\b{\sigma}}_m$ are denoted \[\Delta^{\b{c},\b{\sigma}}_m(\blambda), \quad \nabla^{\b{c},\b{\sigma}}_m(\blambda), \quad L^{\b{c},\b{\sigma}}_m(\blambda), \quad T^{\b{c},\b{\sigma}}_m(\blambda),\] respectively, for $\blambda \in \cP^a(m)$.
\end{defn}

\subsection{External Tensor Products of Restricted Categorifications}\label{subsec:external_tensor_product_restricted}

Fix an admissible sequence 
\[
\Xi := ((\b{\sigma}_1,\b{c}_1),(\b{\sigma}_2,\b{c}_2),\ldots,(\b{\sigma}_r,\b{c}_r)).
\]
We will explicitly construct a restricted multi-Fock tensor product categorification of type $\Xi$ using the construction of the previous subsection.

For each $\ell = 1,\ldots,r$, let us write $(E(\ell),F(\ell),x(\ell),\tau(\ell))$ for the (restricted) categorification data on the (stable) category $\cO^{\b{c}_\ell,\b{\sigma}_\ell}_m$ that we introduced in Definition \ref{defn:lower_admissible_stable}.

In turn, define \[\cO_m^{\Xi} := \cO^{\b{c}_1,\b{\sigma}_1}_m \boxtimes \cO^{\b{c}_2,\b{\sigma}_2}_m \boxtimes \cdots \boxtimes \cO^{\b{c}_r,\b{\sigma}_r}_m.\] Observe that $\cO_m^{\Xi}$ is a highest weight category with poset $\cP^{a_1}_{m} \times \cdots \times \cP^{a_r}_m$ and standard objects \[\Delta^{\Xi}_m(\ublambda) := \Delta^{\b{c}_1,\b{\sigma}_1}_m(\blambda_1) \boxtimes \cdots \boxtimes \Delta^{\b{c}_r,\b{\sigma}_r}_m(\blambda_r)\] for $\ublambda = (\blambda_1,\ldots,\blambda_r) \in \cP^{a_1}_m \times \cdots \cP^{a_r}_m$. This follows by the bi-exactness of the canonical functor \[\cC \times \cD \to \cC \boxtimes \cD\] for finite abelian categories, see e.g., \cite[Proposition 1.11.2]{EGNO}.

We extend the functors $F(\ell)$ and $E(\ell)$ to functors on $\cO^\Xi_m$ (given by acting on the $\ell$th external tensor factor). Writing $F := \bigoplus_{\ell = 1}^r F(\ell)$ and $E := \bigoplus_{\ell = 1}^r E(\ell)$, we define $x \in \End(F)$ (resp. $\tau \in \End(F^2)$) as the direct sum of the endomorphisms $x(\ell)$ (resp. $\tau(\ell)$). Thanks to this setup, the following theorem is straightforward (and thus its proof is omitted for brevity).

\begin{thm}\label{thm:construct_mftp}
    The data $(E,F,x,\tau)$ endows $\cO^\Xi_m$ with the structure of a level $m$ restricted mixed admissible MFTPC of type $\Xi$.
\end{thm}

As a consequence of the structural properties of $\cO^{\Xi}_m$ and Theorem \ref{thm:mixed_admissible_uniqueness}, we obtain several structural results on mixed admissible (restricted) MFTPCs. First and foremost, we obtain explicit formulas for multiplicities of simples in standards.

\begin{prop}\label{prop:mftpc_multiplicities}
    Let $\cC$ be a level $m$ restricted mixed admissible MFTPC of type $\Xi$. For any labels $\ublambda,\ubnu \in \cP^{a_1}(m) \times \cdots \times \cP^{a_r}(m)$, 
    \[[\Delta(\ublambda):L(\ubnu)] = \prod_{\ell=1}^r [\Delta^{\b{c}_\ell,\b{\sigma}_\ell}_{m,N}(\blambda_\ell):L^{\b{c}_\ell,\b{\sigma}_\ell}_{m,N}(\b{\nu}_\ell)],\] 
    where $\Delta^{\b{c}_\ell,\b{\sigma}_\ell}_{m,N}(\blambda_\ell)$ and $L^{\b{c}_\ell,\b{\sigma}_\ell}_{m,N}(\b{\nu}_\ell)$ are the Verma and simple modules in $\cO_m^{\b{c}_\ell,\b{\sigma}_\ell}(N)$ for $m,N \gg 0$. 
    In particular, the multiplicity $[\Delta(\ublambda):L(\ubnu)]$ is the product of parabolic Kazhdan--Lusztig polynomials evaluated at 1.
\end{prop}

\begin{proof}
    Without loss of generality, assume $\b{c}_1,\ldots,\b{c}_n$ are lower admissible, while $\b{c}_{n+1},\ldots,\b{c}_r$ are upper admissible. Pick $m \gg 0$ so that \[\ublambda \in \cP^{a_1}(m) \times \cdots \times \cP^{a_n}(m) \times \cP^{a_{n+1}} \times \cdots \times \cP^{a_r}.\] Let $\cC_m'$ denote the Serre subcategory of $\cC$ corresponding to the poset ideal on the right-hand side. Then, the multiplicity $[\Delta(\ublambda):L(\ubnu)]$ can be computed in $\cC_m'$, which is in fact a Schurian category. Using the construction of Subsubsection \ref{subsubsec:restricted_from_full}, we obtain a Serre quotient $\cC_m$ of $\cC_m'$ with the structure of level $m$ restricted MFTPCs of type $\Xi$. In turn, \cite[Lemma 2.27]{brundan_stroppel_semiinfinite} tells us that the multiplicity $[\Delta(\ublambda):L(\ubnu)]$ can be computed in $\cC_m$ for $m \gg 0$. Theorem \ref{thm:mixed_admissible_uniqueness} gives us an equivalence $\cC_m \simeq \cO_{m}^\Xi$. Then, the proposition follows from the bi-exactness of the Deligne tensor product functor.
\end{proof}

\begin{defn}\label{defn:stable_mftpc}
    Applying the construction of Subsubsection \ref{subsubsec:gluing} to the level $m$ restricted MFTPCs $\cO^{\Xi}_m$, we obtain the MFTPC $\cO^{\Xi}_\infty$, which we call the \textit{stable parabolic category $\cO$ of type $\Xi$}.
\end{defn}

Thanks to Theorem \ref{thm:recover_from_gluing}, any genuine MFTPC of type $\Xi$ is strongly equivariantly equivalent to $\cO^{\Xi}_\infty$ (respecting labels of simple objects). Indeed, the category $\cC'(\mathfrak{F})$ from the statement of that theorem is equivalent to $\cO^{\Xi}$ because of Theorem \ref{thm:mixed_admissible_uniqueness}.

\section{Structure of Complex Rank Parabolic Category \strO}\label{sec:complex_rank_category_o}

We now study structural properties of $\ots$, in parallel with classical properties of the BGG category $\cO$. For instance, we interpolate the Jantzen determinant formula (Theorem \ref{thm:jantzen_determinant}) and use it to prove a complex rank analog of the Jantzen sum formula (Theorem \ref{thm:jantzen_sum}). In turn, we show that certain subcategories of $\ots$ are upper finite highest weight categories (Corollary \ref{cor:ots_highest_weight}) and define the categories $\extots$ (Definition \ref{defn:extots}) from the statement of Conjecture \ref{conj:categorical_conjecture}. 

\subsection{Dual Verma Modules}\label{subsec:duality}

We first introduce a contravariant \textit{duality functor} that interpolates the usual BGG duality to the complex rank setting. We will use the following notation throughout: let $w_0 \in \mathfrak{S}_n$ denote the longest element, and for each permutation $w \in \mathfrak{S}_n$, we write $w\bt := (t_{w(1)},\ldots,t_{w(n)})$ and $w\bs = (s_{w(1)},\ldots,s_{w(n)})$. 

\begin{notation}
    For any $\beta \in \fz_{\bt}^*$ and $M \in \widehat{\cO}^{\bt}$, write $M^\beta \in \Ind\cL_{\bt}$ for the $\beta$-weight subobject in $M$. 
\end{notation}

This duality functor will be defined on the Serre subcategory $\tildeots \subset \ots$ consisting of modules $M$ whose weight subobjects $M^{\beta}$ are \textit{finite length} objects in $\cL^{\bt}$ (for all $\beta \in \fz_{\bt}^*$). In particular, the Verma modules $\mts(\blambda)$ and their simple quotients $\lts(\blambda)$ belong to $\tildeots$.

Given $W \in \tildeots$, we consider the decomposition $W = \bigoplus_{\beta \in \fz_{\bt}^*} W^\beta$ into $\fz_{\bt}$-weight subobjects. Then, we define the $U_{\bt}$-module \[W^\circledast := \bigoplus_{\beta \in \fz_{\bt}^*} (W^\beta)^*,\] where the action of $\fg_{\bt}$ is given as the direct sum of morphisms $\rho^*: \fg^{\alpha}_{\bt} \ox (W^{\beta})^* \to (W^{\beta-\alpha})^*$ (for a root $\alpha$, we write $\fg^{\alpha}_{\bt}$ for the $\alpha$-root space in $\fg_{\bt}$) characterized by the following diagram:
\[\begin{tikzcd}
	{\fg_{\bt}^{\alpha} \ox (W^\beta)^* \ox W^{\beta - \alpha}} && {(W^\beta)^* \ox \fg^{\alpha}_{\bt} \ox W^{\beta - \alpha}} \\
	\\
	{(W^{\beta-\alpha})^* \ox W^{\beta-\alpha}} & \unit & {(W^\beta)^* \ox W^\beta}
	\arrow["\sim", from=1-1, to=1-3]
	\arrow["{\rho^* \ox \id_{W^{\beta-\alpha}}}"', from=1-1, to=3-1]
	\arrow["{\id_{(W^\beta)^*} \ox \rho}", from=1-3, to=3-3]
	\arrow["\ev"', from=3-1, to=3-2]
	\arrow["\ev", from=3-3, to=3-2]
\end{tikzcd}\]
where $\rho: \fg_{\bt}^\alpha \ox W^{\beta - \alpha} \to W^{\beta}$ is the original action map. In general, the resulting module $W^\circledast$ does not lie in $\tildeots$. We introduce the following categorical analog of the Chevalley involution.

\begin{defn}\label{defn:different_choice_parabolic}
For each $w \in \mathfrak{S}_n$, we define a functor $\phi_w: U_{\bt}\modcat \to U_{w\bt}\modcat$ by using the natural equivalence $P_w: \Ind\cL_{\bt} \simeq \Ind\cL_{w\bt}$. That is, given some $U_{\bt}$-module $M$ and $w \in W$, the $U_{w\bt}$-module structure on $\phi_w M = P_wM$ is given by identifying $V^{w\bt}_{ij}$ with $P_wV^{\bt}_{w^{-1}i, w^{-1}j}$.
\end{defn}

Then, for any $W \in \tildeots$, the module $\phi_{w_0}(W^{\circledast}) \in U_{w_0\bt}\modcat$ lies in $\widetilde{\cO}^{w_0\bt}_{w_0\bs}$.

\begin{defn}\label{defn:duality}
We define the duality functor \[\DD: \tildeots \to \widetilde{\cO}^{w_0\bt}_{w_0\bs}, \quad W \mapsto \phi_{w_0}(W^\circledast).\] 
\end{defn}

\begin{notation}
For any $\blambda \in \cP^{2n}$, we adopt the notation
\[
\blambda^* := (\lambda_{2n-1},\lambda_{2n},\ldots,\lambda_3, \lambda_4, \lambda_1,\lambda_2).
\]
\end{notation}

The following lemma is standard.

\begin{lem}\label{lem:duality_properties}
    The following properties about the duality functors hold.
\begin{enumerate}[(a)]
\item The composition of $\DD: \tildeots \to \cO^{w_0\bt}_{w_0\bs}$ with $\DD: \widetilde{\cO}^{w_0\bt}_{w_0\bs} \to \tildeots$ is isomorphic to the identity endofunctor of $\tildeots$. In particular, the duality functors furnish a contravariant equivalence \[\tildeots \simeq \widetilde{\cO}^{w_0\bt}_{w_0\bs}.\]
\item For any $\blambda \in \cP^{2n}$, we have \[\DD\lts(\blambda) = L^{w_0\bt}_{w_0\bs}(\blambda^*).\] In particular, the module $\DD\mts(\blambda)$ has simple socle $L^{w_0\bt}_{w_0\bs}(\blambda^*)$. 
\end{enumerate}
\end{lem}

\begin{defn}
    For each $\blambda \in \cP^{2n}$, we define $\checkmts(\blambda):=\DD M^{w_0\bt}_{w_0\bs}(\blambda^*) \in \ots$. We call this module a \textit{dual Verma module}. It has simple socle $\lts(\blambda)$.
\end{defn}

The following fact from the classical BGG category $\cO$ easily generalizes to our setting.

\begin{lem}\label{lem:homological_verma_filtration}
    For any $\blambda,\b{\nu} \in \cP^{2n}$, we have \[\dim \Ext^n_{\ots}(\mts(\blambda), \checkmts(\b{\nu})) = \delta_{n,0}\delta_{\blambda,\b{\nu}}.\] In particular, if an $M \in \ots$ has a finite filtration by Verma modules, then $\Ext^1(M, \checkmts(\b{\nu})) = 0$ for every $\b{\nu} \in \cP^{2n}$, and the multiplicity of $\mts(\b{\nu})$ in that filtration is $\dim \Hom(M, \checkmts(\b{\nu}))$. 
\end{lem}

\subsection{Interpolated Equivariant Jantzen Determinant Formula}\label{subsec:jantzen_determinant}

We now interpolate Jantzen's determinant formula for the Shapovalov form of parabolic Verma modules. As suggested in \cite[Remark 4.6]{Etingof_complex_rank_2}, this formula will provide us with a criterion for the simplicity of $\mts(\blambda)$. 

\subsubsection{Classical Story}\label{subsubsec:classical} Let us first recall Jantzen's formula for the parabolic Shapovalov morphism in classical Lie theory. This will require some notational setup. 

Fix positive integers $m_1,\ldots,m_n \in \ZZ$ and set $m := m_1 + \cdots + m_n$. Consider the parabolic category $\ul{\cO}^{\b{m}}$ associated with the usual Lie algebra $\ul{\fgl}_m$ and the composition $m = m_1 + \cdots + m_n$. Write $\ul{\fl}_{\b{m}}$ and $\ul{\fp}_{\b{m}}$ for the corresponding Levi and parabolic subalgebras, respectively. Let $\ul{\fz}_{\b{m}} \simeq \fk^n$ denote the center of $\ul{\fl}_{\b{m}}$.

\begin{notation}
For $i = 1,\ldots,n$ and $a = 1,\ldots,m_i$, we write $\varepsilon_{i}^a \in \fk^{m_1} \times \cdots \times \fk^{m_n}$ for the corresponding coordinate vector. In turn, for $i,j \in\{1,\ldots,n\}$, $p \in \{1,\ldots,m_i\}$ and $q \in \{1,\ldots,m_j\}$, we write $\alpha_{i,j}^{p,q} := \varepsilon_i^p - \varepsilon_j^q$, which we regard as a root for $\ul{\fgl}_{m}$. We write $R^+$ for the set of positive roots, i.e., roots $\alpha_{i,j}^{p,q}$ where $i < j$ or $i = j$ and $p < q$. Then, let $R_{\b{m}}^+$ denote the set of positive roots for the Levi subalgebra $\ul{\fl}_{\b{m}}$, i.e., the positive roots $\alpha_{i,j}^{p,q}$ where $i = j$. 

Let $W_{\b{m}}$ denote the Weyl group of $\ul{\fl}_{\b{m}}$, and let ${\varrho}$ denote the half sum of all positive roots $\beta \in R^+$. More explicitly, define \[M_i := m_1 + \cdots + m_{i-1}\] so that for each $i = 1,\ldots,n$ and $p = 1,\ldots,m_i$, the $(m_1 + \cdots + m_{i-1} + p)$th coordinate of $\varrho$ is \[\varrho_{i,p} := \frac{1}{2}((m_i + \cdots + m_n) - (m_1 + \cdots + m_{i-1}) + 1) - p = \frac{m + 1}{2} - M_i - p.\] We write $\tau \cdot \eta$ for the $\varrho$-shifted action of $\tau \in W_{\b{m}}$ on $\eta \in \fk^{M}$. 
\end{notation}

\begin{defn}
For any $\eta \in \fk^{m}$ with trivial stabilizer under the dot action of $W_{\b{m}}$, there exists unique $\tau \in W_{\b{m}}$ such that $\tau \cdot \eta$ is an $\ul{\fl}_{\b{m}}$-dominant weight. We define $\mathrm{sgn}(\eta) := \mathrm{sgn}(\tau)$. Moreover, we write $\eta_+ := \tau \cdot \eta$. For any $\eta \in \fk^M$ with nontrivial stabilizer, we declare $\mathrm{sgn}(\eta) = 0$. 
\end{defn}  

For each $\blambda \in \cP^{2n}$ where $m_i \geq \ell(\lambda_{2i-1}) + \ell(\lambda_{2i})$ for each $i =1,\ldots,n$, write $\ul{X}^{\b{m}}(\blambda,\bs)$ for the finite-dimensional irreducible $\ul{\fl}_{\b{m}}$-module with highest weight \[\Theta_{\b{m}}^{\b{s}}(\blambda) := (\phi_{m_1,s_1}(\blambda)_1,\phi_{m_2,s_2}(\blambda)_2,\ldots,\phi_{m_n,s_n}(\blambda)_n),\] where $\phi_{m_i,s_i}(\blambda)_i \in \fk^{m_i}$ is defined by 
    \begin{align*}
        \phi_{m_i,s_i}(\blambda)_i &= (\lambda_{2i-1,1} + s_i, \lambda_{2i-1,2} + s_i,\ldots,\lambda_{2i-1,\ell(\lambda_{2i-1})} + s_i, s_i, \ldots,
        \\[5pt] &\quad \ \ \ldots, s_i, s_i, -\lambda_{2i,\ell(\lambda_{2i})} + s_i, \ldots, -\lambda_{2i,2} + s_i, -\lambda_{2i,1} + s_i).
\end{align*}  
In turn, we define $\ul{M}^{\b{m}}(\blambda,\bs)$ as the parabolic Verma module with highest weight $\Theta_{\bs}^{\b{m}}(\blambda)$, i.e., \[\ul{M}^{\b{m}}(\blambda,\bs) := U(\ul{\fgl}_M) \ox_{U(\ul{\fp}_{\b{m}})} X^{\b{m}}(\blambda,\bs).\]
More generally, for any $\ul{\fl}_{\b{m}}$-dominant weight $\vartheta \in \fk^M$, we write $\ul{M}^{\b{m}}(\vartheta)$ for the parabolic Verma module with highest weight $\vartheta$.

Finally, recall that the parabolic Verma module $\ul{M}^{\b{m}}(\blambda,\bs)$ carries the contravariant Shapovalov form, which is $\ul{\fl}_{\b{m}}$-equivariant. Thus, it restricts to pairings on each $\ul{\fl}_{\b{m}}$-isotypic component of $\ul{M}^{\b{m}}(\blambda,\bs)$. We write $\ul{D}^{\b{m}}_{\bs}(\blambda)_{\bmu} \in \fk(s_1,\ldots,s_n)$ for the determinant of the Shapovalov form in the $\ul{X}^{\b{m}}(\bmu,\bs)$-isotypic component of $\ul{M}^{\b{m}}(\blambda,\bs)$. Refer to Appendix \ref{appendix:equivariant_jantzen} for the proof of the following proposition.

\begin{prop}\label{prop:equivariant_jantzen}
    Up to a nonzero scalar, the following formula holds for $m_i \gg 0$:
    {\small
    \[
        \ul{D}^{\b{m}}_{\bs}(\blambda)_{\bmu} = \prod_{\alpha = \alpha^{p,q}_{i,j} \in R^+ \setminus R^+_{\b{m}}} \ \prod_{k \geq 1} (\phi_{m_i,s_i}(\blambda)_{i,p} - \phi_{m_j,s_j}(\blambda)_{j,q} - (M_i+p) + (M_j+q) - k)^{\tilde{r}_{\b{m}}(\blambda, k \alpha,\bmu)},
    \]
    }
    where we write \[\tilde{r}_{\b{m}}(\blambda,k\alpha,\bmu) = \mathrm{sgn}(\Theta_{\b{m}}^{\b{s}}(\blambda) - k\alpha)r_{\b{m}}((\Theta_{\b{m}}^{\b{s}}(\blambda) - k\alpha)_+,\bmu),\] and $r_{\b{m}}((\Theta_{\b{m}}^{\b{s}}(\blambda) - k\alpha)_+,\bmu)$ denotes the multiplicity of the $X^{\b{m}}(\bmu,\bs)$-isotypic component in the parabolic Verma module $\ul{M}^{\b{m}}((\Theta_{\b{m}}^{\b{s}}(\blambda) - k\alpha)_+)$. 
\end{prop}

\subsubsection{Interpolated Dominance Order}

To prove the complex rank analog of the Shapovalov morphism, we first introduce an interpolated variant of the dominance order.

For sufficiently large positive integers $N$, we define 
\[
    \Phi_{N}(\blambda) := (\phi_{N}(\blambda)_1,\ldots,\phi_{N}(\blambda)_n) \in (\fk^N)^n,
\]
where, for each $i = 1,\ldots,n$, we have
\begin{align*}
    \phi_{N}(\blambda)_i &:= (\lambda_{2i-1,1},\lambda_{2i-1,2},\ldots,\lambda_{2i-1,\ell(\lambda_{2i-1})}, 0,\ldots,0,-\lambda_{2i,\ell(\lambda_{2i})}, \ldots,-\lambda_{2i,2},-\lambda_{2i,1}) \in \fk^N.
\end{align*}

\begin{defn}
    We define the \textit{interpolated dominance order} $\trianglelefteq$ on $\cP^{2n}$ by declaring $\bmu \trianglelefteq \blambda$ if $\Phi_N(\bmu) \leq \Phi_N(\blambda)$ in the dominance order on weights of the usual Lie algebra $\fgl_{Nn}$ for $N \gg 0$. Note that this definition does not depend on the choice of $N \gg 0$. 
\end{defn}

It is easy to see that $\trianglelefteq$ defines an interval-finite partial order on $\cP^{2n}$. 

\begin{lem}\label{lem:isotypic_weights}
    Assume $t_i \not \in \ZZ$ for all $i \in \{1,\ldots,n\}$. As an object of $\Ind \cL_{\bt}$, we have
    \[
    M_{\b{s}}^t(\b{\lambda}) \simeq \Sym(\mathfrak{u}^-_{\bt}) \ox X^{\bt}_{\bs}(\b{\lambda}).
    \]
    In particular, the $\fl_{\b{t}}$-isotypic components of $M_{\b{s}}^t(\b{\lambda})$ are given by $X_{\b{s}}(\b{\lambda}')$, where $\b{\lambda}'$ is obtained from $\b{\lambda}$ through a sequence of the following manipulations (applied at each step to some $\b{\mu} \in \cP^{2n}$):
    \begin{itemize}
        \item[(i)] add a box to $\mu_{2j}$ and add a box to $\mu_{2i-1}$ for some $i > j$;
        \item[(ii)] add a box to $\mu_{2j}$ and remove a box from $\mu_{2i}$ for some $i > j$;
        \item[(iii)] remove a box from $\mu_{2j-1}$ and add a box to $\mu_{2i-1}$ for some $i > j$;
        \item[(iv)] remove a box from $\mu_{2j-1}$ and remove a box from $\mu_{2i}$ for some $i > j$.
    \end{itemize}
\end{lem}

\begin{proof}
    The first claim follows from the explicit definition of $\mts(\blambda)$ as an induced module and the PBW theorem, which gives us an isomorphism \[U_{\bt} \simeq \Sym(\mathfrak{u}^-_{\bt}) \ox U(\mathfrak{p}_{\bt})\] of right $U(\fp_{\bt})$-modules (this is not an isomorphism of algebras). Writing $V_{ij} := V_i \boxtimes V_j^* \in \cL_{\bt}$, it follows that the isotypic components of $\mts(\blambda)$ are the direct summands of \[V_{i_1,j_1} \ox \cdots \ox V_{i_r,j_r} \ox X_{\bs}^{\bt}(\blambda)\] for arbitrary $r \geq 0$ and integers $i_\ell,j_\ell \in \{1,\ldots,n\}$ satisfying $i_\ell > j_\ell$.  The claim follows from the explicit decomposition, given by \cite[Corollary 7.1.2]{comes_wilson}, of the tensor products \[V_i \ox X^{t_i}(\blambda_i), \quad V_i^* \ox X^{t_i}(\blambda_i). \qedhere\] 
\end{proof}

\begin{cor}
    If $\lts(\bmu)$ is a composition factor of $\mts(\blambda)$, then $\bmu \trianglelefteq \blambda$. 
\end{cor}

\begin{defn}\label{defn:dominance_ordering_map}
    We define a map $\omega_{\bs}: \cP^{2n} \to \fz_{\bt}^*$ via \[\blambda \mapsto (|\lambda_1| - |\lambda_2| + s_1, |\lambda_3| - |\lambda_4| + s_2, \ldots, |\lambda_{2n-1}| - |\lambda_{2n}| + s_n).\] In particular, $\fz_{\bt}$ acts on $X^{\bt}_{\bs}(\blambda)$ by $\omega_{\bs}(\blambda)$. 
\end{defn}

To avoid confusion, we write $\trianglelefteq_{\mathrm{dom}}$ for the dominance order on $\fz_{\bt}^*$. In particular, $\beta_1 \trianglelefteq_{\mathrm{dom}} \beta_2$ if and only if $\beta_2 - \beta_1$ is a $\ZZ_{\geq 0}$-linear combination of simple roots. The following lemma is obvious.

\begin{lem}
    For any $\bmu,\blambda \in \cP^{2n}$, we have that $\bmu \trianglelefteq \blambda$ implies $\omega_{\bs}(\bmu) \trianglelefteq_{\mathrm{dom}} \omega_{\bs}(\blambda)$.
\end{lem}

\subsubsection{Shapovalov Morphism}

Fix some $\blambda,\bmu \in \cP^{2n}$ with $\bmu \trianglelefteq_{\bt} \blambda$ and set $\beta := \omega_{\bs}(\bmu) - \omega_{\bs}(\blambda) \in \fz_{\bt}^*$. Write $U_{\bt}[\beta]$ for the $\beta$-weight subobject in $U_{\bt}$. Thus, the $X_{\bs}^{\bt}(\bmu)$-isotypic component $Y_{\bs}^{\bt}(\blambda)_{\bmu}$ of $\mts(\blambda)$ is a summand of $U_{\bt}[\beta] \ox X_{\bs}^{\bt}(\blambda) \subset \mts(\blambda)$. The action map $U_{\bt}[\beta] \ox \mts(\blambda) \to \mts(\blambda)$ restricts to an $\fl_{\bt}$-equivariant morphism \[U_{\bt}[-\beta] \ox U_{\bt}[\beta] \ox X_{\bs}^{\bt}(\blambda) \to X_{\bs}^{\bt}(\blambda),\] which induces (by adjunction) an $\fl_{\bt}$-equivariant morphism \[U_{\bt}[\beta] \ox X_{\bs}^{\bt}(\blambda) \to U_{\bt}[\beta] \ox X_{\bs}^{\bt}(\blambda).\] In particular, we obtain the Shapovalov morphism $\tilde{B}_{\bs}^{\bt}(\blambda)_{\bmu} \in \End_{\fl_{\bt}}(Y_{\bs}^{\bt}(\blambda)_{\bmu})$. We may write \[Y_{\bs}^{\bt}(\blambda)_{\bmu} := \fk^{\oplus r(\blambda,\bmu)} \ox X_{\bs}^{\bt}(\bmu)\] for some $r(\blambda,\bmu) \geq 0$ independent of $\bt$. Altogether, we have a linear map $B_{\bs}^{\bt}(\blambda)_{\bmu} \in \End(\fk^{\oplus r(\blambda,\bmu)})$. 

\begin{defn}
    We define the \textit{Shapovalov determinant}
    \[
        D_{\bs}^{\bt}(\blambda)_{\bmu} := \det B_{\bs}^{\bt}(\blambda)_{\bmu} \in \fk(s_1,\ldots,s_n),
    \]
    which we regard as a rational function in variables $s_1,\ldots,s_n$, well-defined up to a nonzero scalar.
\end{defn}

From the construction, we see that $\mts(\blambda)$ is simple if and only if $D_{\bs}^{\bt}(\blambda)_{\bmu} \neq 0$ for all $\bmu \trianglelefteq_{\bt} \blambda$. We will compute $D_{\bs}^{\bt}(\blambda)_{\bmu}$ (for fixed $\blambda$ and $\bmu$) by interpolation. 

Fix indeterminates $x_1,\ldots,x_n$ and repeat the constructions of Subsection \ref{subsec:categoryO} in the category
\[
    \cL_{\b{x}} := \rep{\GL_{x_1}} \boxtimes \cdots \boxtimes \rep{\GL_{x_n}}
\]
over the field of rational functions $\fk(x_1,\ldots,x_n)$. In particular, we have the Verma modules $\mxs(\blambda)$ for $\blambda \in \cP^{2n}$, and for each $\bmu \in \cP^{2n}$, we have a Shapovalov morphism $B_{\bs}^{\bx}(\blambda)_{\bmu} \in \End(\fk^{\oplus r(\blambda,\bmu)})$ and thus a Shapovalov determinant $D_{\bs}^{\bx}(\blambda)_{\bmu} \in \fk(x_1,\ldots,x_n,s_1,\ldots,s_n)$.

The Shapovalov morphism $B_{\bs}^{\bx}(\blambda)_{\bmu}$ is defined explicitly as a $\fk(x_1,\ldots,x_n)$-linear combination of braidings, evaluations and coevaluations. In fact, each coefficient is actually \textit{polynomial} in the variables $(x_1,\ldots,x_n)$. Hence, for $(t_1,\ldots,t_n) \in \fk^n$ with $t_i \not \in \ZZ$ for each $i$, the morphism $B_{\bs}^{\bt}(\blambda)_{\bmu}$ is given by the same linear combination of braidings, evaluations, and coevaluations, with the aforementioned polynomial coefficients evaluated at $(t_1,\ldots,t_n)$. Moreover, the determinant of $B_{\bs}^{\bx}(\blambda)_{\bmu}$ is also a polynomial function in $(x_1,\ldots,x_n)$, and we deduce the following lemma.

\begin{lem}\label{lem:evaluate_non_integer}
Up to a nonzero scalar, the determinant $D_{\bs}^{\bt}(\blambda)_{\bmu} \in \fk(s_1,\ldots,s_n)$ is given by evaluating the polynomial $D_{\bs}^{\bx}(\blambda)_{\bmu} \in \fk(s_1,\ldots,s_n)(x_1,\ldots,x_n)$ at $(t_1,\ldots,t_n) \in \fk^n$ with $t_i \not \in \ZZ$ for each $i$. 
\end{lem}

\subsubsection{Interpolated Shapovalov Determinant Formula}

Continue to fix $\blambda,\bmu \in \cP^{2n}$ with $\bmu \trianglelefteq_{\bt} \blambda$. Moreover, keep the notation for the classical Lie theoretic constructions from Subsubsection \ref{subsubsec:classical}.

Write $U(\ul{\fgl}_{M})[\beta]$ for the $\beta$-isotypic component of $U(\ul{\fgl}_M)$ with respect to the $\fz_{\b{m}}$-action, where $\beta$ has the same meaning as in the previous subsubsection. The categorical definition of the Shapovalov morphism \[B^{\bt}_{\bs}(\blambda)[\beta]: U_{\bt}[\beta] \ox X^{\bt}_{\bs}(\blambda) \to U_{\bt}[\beta] \ox X^{\bt}_{\bs}(\blambda)\] that we provided in complex rank can also be used to define the usual Shapovalov morphism \[U(\ul{\fgl}_{M})[\beta] \ox X^{\b{m}}(\blambda,\bs) \to U(\ul{\fgl}_{M})[\beta] \ox X^{\b{m}}(\blambda,\bs)\] in $\fz_{\b{m}}$-weight $\beta$. Write $\cS$ for the multiset of multipartitions so that \begin{align}\label{eqn:isotypic1}U_{\bt}[\beta] \ox X_{\bs}^{\bt}(\blambda) = \bigoplus_{\bnu \in \cS} X_{\bs}^{\bt}(\bnu)\end{align} is the $\fl_{\bt}$-isotypic decomposition of $U_{\bt}[\beta] \ox X_{\bs}^{\bt}(\blambda)$. 

When the integers $m_i$ are sufficiently large, we have the parallel $\ul{\fl}_{\b{m}}$-isotypic decomposition 
\begin{align}\label{eqn:isotypic2} U(\ul{\fgl}_{M})[\beta] \ox X^{\b{m}}(\blambda,\bs) \simeq \bigoplus_{\bnu \in \cS} X^{\b{m}}(\bnu,\bs).\end{align} 
Write $\cL_{\bx}^{\mathrm{poly}} \subset \cL_{\bx}$ for the subcategory consisting of the same objects as $\cL_{\bx}$ but whose morphisms are the morphisms that can be expressed as $\fk[x_1,\ldots,x_n]$-linear combinations of braidings, evaluations, and coevaluations.  
In particular, the morphism $B^{\bx}_{\bs}(\blambda)[\beta]$ belongs to $\cL_{\bx}^{\mathrm{poly}}$. For any $m_1,\ldots,m_n \in \ZZ$, we have a functor $F_{\b{m}}: \cL_{\bx} \to \mathrm{Rep}(\ul{\fl}_{\b{m}})$ given by evaluating at $(x_1,\ldots,x_n) = (m_1,\ldots,m_n)$ and then applying the semisimplification functor. When the integers $m_1,\ldots,m_n$ are sufficiently large, we have $F_{\b{m}}(X^{\bx}_{\bs}(\bnu)) = X^{\b{m}}(\bnu,\bs)$ for each $\bnu \in \cS$. We conclude that $F_{\b{m}}$ sends $B_{\bs}^{\bx}(\blambda)_{\bmu}$ to the \textit{restriction} of the Shapovalov morphism \[U(\ul{\fgl}_{M})[\beta] \ox X^{\b{m}}(\blambda,\bs) \to U(\ul{\fgl}_{M})[\beta] \ox X^{\b{m}}(\blambda,\bs)\] to the $X^{\b{m}}(\bmu,\b{s})$-isotypic component. Thus, we deduce the following elementary lemma.

\begin{lem}\label{lem:evaluate_integer}
    For $m_i \gg 0$, the evaluation of the rational function $D_{\bs}^{\bx}(\blambda)_{\bmu} \in \fk(s_1,\ldots,s_n)(x_1,\ldots,x_n)$ at $x_i = m_i \in \ZZ$ equals, up to $\fk^\times$, the classical Shapovalov determinant $\ul{D}_{\bs}^{\b{m}}(\blambda)_{\bmu} \in \fk(s_1,\ldots,s_n)$.
\end{lem}

Thus, we seek a formula for $D_{\bs}^{\bx}(\blambda)_{\bmu}$ interpolating Jantzen's formula from Proposition \ref{prop:equivariant_jantzen}.

We need a different perspective on the labelling set $\cP^{2n}$ for Verma modules in $\ots$. We identify this set with a subset of the abelian group $\cW := (\fk^\NN)^{2n}$. We always use the convention that $\NN = \{1,2,\ldots\}$.
Assign any multipartition $\blambda \in \cP^{2n}$ to the interpolated weight $\Psi(\blambda) := (\Psi(\blambda)_1,\Psi(\blambda)_2,\ldots,\Psi(\blambda)_{2n})$, where for each $i = 1,\ldots,n$, we have infinite sequences
\begin{align*}
\Psi(\blambda)_{2i-1} &:= (\lambda_{2i-1,1}, \lambda_{2i-1,2},\ldots), \quad \Psi(\blambda)_{2i} := (-\lambda_{2i,1}, -\lambda_{2i,2},\ldots),
\end{align*}
and we agree that $\eta_j = 0$ if $j$ exceeds the length of a partition $\eta$. 

\begin{notation}\label{notation:interpolated_weights}
We introduce the \textit{interpolated Weyl vector}. Recalling $t = t_1 + \cdots + t_n$, we write
\[
T_i := t_1 + \cdots + t_{i-1}, \quad T_i' := t_i + \cdots + t_n, \quad R_i := \frac{T_i' - T_i + 1}{2} = \frac{t + 1}{2} - T_i
\]
and define $\rho_{\bt} := (\bar\rho_1,\ldots,\bar\rho_{2n}) \in \cW$, where 
\[
\bar\rho_{2i-1} = (R_i-1, R_i - 2, R_i -3 ,\ldots), \quad \bar\rho_{2i} = (R_{i+1}, R_{i+1} + 1, R_{i+1} + 2,\ldots)
\]
for each $i = 1,\ldots,n$. Given $\blambda \in \cP^{2n}$, we define its \textit{$\bar{\rho}$-shifted weight} to be
$\Psi_{\bt}(\blambda)_{\bar{\rho}} := \Psi_{\bt}(\blambda) +\bar{\rho}$. For brevity, let's use the notation \[\Psi_{\bt}(\blambda)_{\bar{\rho}} := (\psi_{\bt}(\blambda)_1,\ldots,\psi_{\bt}(\blambda)_{2n}).\] 
\end{notation}

We now introduce the following notational shorthand. Set 
\[\SS := \{(i,j,p,q,k) \ |\ \lceil i/2\rceil < \lceil j/2\rceil\} \subset \{1,\ldots,2n\} \times \{1,\ldots,2n\} \times \NN \times \NN \times \NN.\]
For any $\zeta = (i,j,p,q,k)\in \SS$, define a multipartition $\blambda(\zeta) \in \cP^{2n}$ as follows. First, pick integers $m_r \geq \ell(\lambda_{2r-1}) + \ell(\lambda_{2r}) + p + q$ for each $r = 1,\ldots,n$. Recall the notation $\Theta_{\b{m}}^{\b{s}}(\blambda) \in \fk^{m_1 + \cdots + m_n}$ from before and set $\alpha := \alpha_{\lceil i/2\rceil,\lceil j/2 \rceil}^{p',q'}$, where $p' = p$ if $i$ is odd and $p' = m_{i/2} - p + 1$ if $i$ is even; the number $q'$ is defined similarly. We let $\blambda(\zeta)$ denote the unique multipartition so that \[\Theta_{\b{m}}^{\b{s}}(\blambda(\zeta)) = (\Theta_{\b{m}}^{\bs}(\blambda) - k\alpha)_+.\] Note that $\blambda(\zeta)$ does not depend on the choice of $m_1,\ldots,m_{n}$ provided that these integers are sufficiently large. Moreover, $\blambda_r = \blambda(\zeta)_r$ for all $r \neq i,j$. We also set \[\mathrm{sgn}(\blambda(\zeta)) := \mathrm{sgn}(\Theta_{\b{m}}^{\bs}(\blambda) - k\alpha),\] which is independent of the choice of $m_i \gg 0$ as well. 

Without further ado, we arrive at the following interpolated Jantzen determinant formula.

\begin{thm}\label{thm:jantzen_determinant}
    Assume $t_i \not \in \ZZ$ for any $i = 1,\ldots,n$. Up to a nonzero scalar, we have 
    \[
        D_{\bs}^{\bt}(\blambda)_{\bmu} = \prod_{\zeta = (i,j,p,q,k) \in \SS} (\psi_{\bt}(\blambda)_{i,p} - \psi_{\bt}(\blambda)_{j,q} + s_i - s_j - k)^{\tilde{r}(\blambda(\zeta),\bmu)},
    \]
    where \[\tilde{r}(\blambda(\zeta),\bmu) = \mathrm{sgn}(\blambda(\zeta))r(\blambda(\zeta),\bmu)\] and $r(\blambda(\zeta),\bmu)$ is the multiplicity of the $X^{\bt}_{\bs}(\bmu)$-isotypic component of $\mts(\blambda(\zeta))$.
\end{thm}

Necessarily, $\blambda(\zeta) \triangleleft_{\bt} \blambda$ because $(\Theta_{\b{m}}^{\b{s}}(\blambda) - k\alpha)_+ < \Theta_{\b{m}}^{\b{s}}(\blambda)$ in the dominance order. Hence, only finitely many of the exponents $\tilde{r}(\blambda(\zeta),\bmu)$ can be nonzero because there are only finitely many $\bnu \in \cP^{2n}$ with $\bmu \trianglelefteq_{\bt} \bnu \trianglelefteq_{\bt} \blambda$. 

\begin{proof}[Proof of Theorem \ref{thm:jantzen_determinant}]
    Let's write $R(x_1,\ldots,x_n) \in \fk(s_1,\ldots,s_n)(x_1,\ldots,x_n)$ for the rational function 
    such that $R(t_1,\ldots,t_n)$ is the expression given in statement of Theorem \ref{thm:jantzen_determinant}. For sufficiently large integers $m_1,\ldots,m_n$, the multiplicity $r_{\b{m}}((\Theta_{\b{m}}^{\bs}(\blambda)-k\alpha_{i,j}^{p,q})_+,\bmu)$ from the statement of Proposition \ref{prop:equivariant_jantzen} stabilize, and its stable limit coincides with the multiplicity $r(\blambda(\zeta),\bmu)$ for a uniquely determined $\zeta \in \SS$. In light of the expression for $\ul{D}^{\b{m}}_{\bs}(\blambda)_{\bmu}$ derived in Proposition \ref{prop:equivariant_jantzen}, we see that $R(m_1,\ldots,m_n) = \ul{D}^{\b{m}}_{\bs}(\blambda)_{\bmu}$ when $m_i \gg 0$. On the other hand, Lemma \ref{lem:evaluate_integer} tells us that $\ul{D}^{\b{m}}_{\bs}(\blambda)_{\bmu}$ is the evaluation of $D^{\b{x}}_{\bs}(\blambda)_{\bmu} \in \fk(s_1,\ldots,s_n)(x_1,\ldots,x_n)$ at $(x_1,\ldots,x_n) = (m_1,\ldots,m_n)$ when $m_i \gg 0$. It follows that $R(x_1,\ldots,x_n) = D^{\bx}_{\bs}(\blambda)_{\bmu}$. Lemma \ref{lem:evaluate_non_integer} thus concludes the proof.
\end{proof}

Not only does Theorem \ref{thm:jantzen_determinant} give us a concrete computable criterion for the simplicity of $\mts(\blambda)$, but it also helps characterize when $\lts(\bmu)$ might be a composition factor of $\mts(\blambda)$ and hence gives rise to a highest weight structure on $\ots$. The construction of this highest weight structure parallels the construction for the classical BGG category $\cO$.

\subsection{Equivariant Jantzen Sum Formula}\label{subsec:jantzen_sum}

In this subsection, we use Theorem \ref{thm:jantzen_determinant} to derive a Jantzen sum formula for Verma modules $\mts(\blambda)$. Just like in the classical setting, this sum formula will help us determine an appropriate highest weight order for $\ots$. 

Let $R := \fk[[y]]$ denote the ring of power series over $\fk$ in the indeterminate $y$. We deform $\cO^{\bt}$ to an $R$-linear abelian category $\cO^{\bt}_R$ defined over the semisimple $R$-linear abelian category 
\[
\cL_{\bt,R} := \rep{\GL_{t_1(1+y)}} \boxtimes \cdots \boxtimes \rep{\GL_{t_n(1+y)}}
\]
in the same way that $\cO^{\bt}$ was defined over $\cL_{\bt}$. This category contains Verma modules $M^{\bt,R}_{\bs}(\blambda)$ for $\blambda \in \cP^{2n}$, which are flat deformations of the Verma modules $\mts(\blambda)$. 

Moreover, we can make sense of the Shapovalov morphisms $B_{\bs}^{\bt,R}(\blambda)_{\bmu}$ just as before. In particular, their determinants give rise to Laurent series \[D_{\bs}^{\bt,R}(\blambda)_{\bmu}(y) = \det B_{\bs}^{\bt,R}(\blambda)_{\bmu} \in R[y^{-1}](s_1,\ldots,s_n).\]
Let $B_{\bs}^{\bt,R}(\blambda) \in \End(M^{\bt,R}_{\bs}(\blambda))$ denote the direct sum of the morphisms $B_{\bs}^{\bt,R}(\blambda)_{\bmu}$. In turn, we may define the \textit{Jantzen filtration}  
\[
M^{\bt}_{\bs}(\blambda) = M^{\bt}_{\bs}(\blambda)_0 \supset M^{\bt}_{\bs}(\blambda)_1 \supset M^{\bt}_{\bs}(\blambda)_2 \supset \cdots 
\]
where $M^{\bt}_{\bs}(\blambda)_1$ is the maximal proper submodule of $M^{\bt}_{\bs}(\blambda)$, and for each $j \geq 0$, the submodule $\mts(\blambda)_j$ is the specialization of \[M^{\bt,R}_{\bs}(\blambda)_j = B_{\bs}^{\bt,R}(\blambda)^{-1}(y^jM^{\bt,R}_{\bs}(\blambda))\] to the closed point. Taking $\bmu$-isotypic components (using the semisimplicity of the $\fl_{\bt}$-action), we have a filtration 
    \[
        \mts(\blambda)_{0,\bmu} \supset \mts(\blambda)_{1,\bmu} \supset \mts(\blambda)_{2,\bmu} \supset \cdots,
    \]
    where $\mts(\blambda)_{j,\bmu}$ denotes the $X_{\bs}^{\bt}(\bmu)$-isotypic component of $\mts(\blambda)_j$. Taking multiplicity spaces, we can view this as a classical Jantzen filtration of a finite-dimensional $\fk$-vector space. The following lemma is simply an explicit interpretation of a classical result of Jantzen on such filtrations. Let $\mathrm{ord}_y(f)$ denote the $y$-adic valuation of any $f \in \mathrm{Frac}(R)$.

\begin{lem}\label{lem:classical_jantzen}
    Assume $t_i \not \in \ZZ$ for any $i = 1,\ldots,n$. For any $\blambda,\bmu \in \cP^{2n}$, we have
    \[
        \mathrm{ord}_y(D_{\bs}^{\bt,R}(\blambda)_{\bmu}(y)) = \sum_{j \geq 1} r_j(\mts(\blambda),\bmu),
    \]
    where $r_j(\mts(\blambda),\bmu)$ is the multiplicity of the $X^{\bt}_{\bs}(\bmu)$-isotypic component of $\mts(\blambda)_j$.  
\end{lem}

Our complex rank Jantzen sum formula will be an identity of formal \textit{characters} of certain objects in $\ots$. These characters should be viewed as \textit{Levi-equivariant} version of the usual notion of characters from the classical BGG category $\cO$.

\begin{defn}
    For any $M \in \ots$, we define its \textit{equivariant character} $\mathrm{ch}(M)$ as the formal series 
    \[
        \mathrm{ch}(M) := \sum_{\bmu \in \cP^{2n}} \mathrm{mult}(M,\bmu)e^{\bmu},
    \]  
    where $\mathrm{mult}(M,\bmu)$ is the multiplicity of the $X^{\bt}_{\bs}(\bmu)$-isotypic component of $M$.
\end{defn}

The statement of the Jantzen sum formula will require some additional notational setup. First, we need to partition the set of indices $\{1,\ldots,2n\}$ based on the nonvanishing of certain factors
\[
\psi_{\bt}(\blambda)_{i,p} - \psi_{\bt}(\blambda)_{j,q} + s_i - s_j - k
\]
in the Jantzen determinant formula (from Theorem \ref{thm:jantzen_determinant}). In particular, these factors will not contribute to the order of vanishing of the determinant $D_{\bs}^{t,R}(\blambda)_{\bmu}(y)$. 

\begin{notation}\label{notation:set_partition}
Recall the notation $T_i := t_1 + \cdots + t_{i-1}$ for $i=1,\ldots,n+1$. Define the constants \[\gamma_{2i-1} := s_i - T_i, \quad \gamma_{2i} := s_i - T_{i+1}\] for $i = 1,\ldots,n$. In turn, form the set partition
\[
\{1,\ldots,2n\} = \Gamma_1 \sqcup \Gamma_2 \sqcup \cdots \sqcup \Gamma_r
\]
so that $\gamma_i - \gamma_j \in \ZZ$ if and only if $i$ and $j$ both belong to $\Gamma_\ell$ for some $\ell = 1,\ldots,r$. Then, we enumerate the elements of $\Gamma_\ell$ in increasing order: \[\Gamma_\ell = \{i_\ell(1) < i_\ell(2) < \cdots < i_\ell(a_\ell)\},\]
where $a_\ell = \# \Gamma_\ell$. In particular, if $i \in \Gamma_\ell$ and $j \in \Gamma_{\ell'}$, we observe that \[\psi_{\bt}(\blambda)_{i,p} - \psi_{\bt}(\blambda)_{j,q} + s_i - s_j - k \neq 0\] for \textit{any} $\blambda \in \cP^{2n}$ and $k \geq 1$ unless $\ell = \ell'$. 
\end{notation}

In turn, we want to introduce some notation to describe the tuples $\zeta = (i,j,p,q,k) \in \SS$ where the corresponding factor \[\psi_{\bt}(\blambda)_{i,p} - \psi_{\bt}(\blambda)_{j,q} + s_i - s_j - k\] \textit{does} contribute to the order of vanishing of $D_{\bs}^{t,R}(\blambda)_{\bmu}(y)$. In particular, there must exist $\ell = 1,\ldots,r$ such that $i = i_\ell(j_1)$ and $j = i_\ell(j_2)$ for some $1 \leq j_1 < j_2 \leq a_\ell$.

\begin{notation}
For each $\ell = 1,\ldots,r$, let us write \[\XX_\ell := \{(j_1,j_2,p,q) \in \NN^{4} \ |\ 1 \leq j_1 < j_2 \leq a_\ell\}.\]
Fix some $\chi = (j_1,j_2,p,q) \in \XX_{\ell}$ and set \[\delta_{\ell,\chi}(\blambda) := \psi_{\bt}(\blambda)_{i_\ell(j_1),p} - \psi_{\bt}(\blambda)_{i_\ell(j_2),q} + s_{i_\ell(j_1)} - s_{i_\ell(j_2)}.\] Note that $\delta_{\ell,\chi}(\blambda) \in \ZZ$ by the definition of the set partition $\Gamma_1 \sqcup \cdots \sqcup \Gamma_r$. Then, we set 
\[
\widetilde{\XX}_\ell^+(\blambda) := \{\chi \in \XX_\ell \ |\ \delta_{\ell,\chi}(\blambda) > 0\}.
\]
For any $\chi = (j_1,j_2,p,q) \in \widetilde{\XX}_\ell^+(\blambda)$, we define $\zeta^{\ell,\chi} := (i_\ell(j_1),i_\ell(j_2),p,q,\delta_{\ell,\chi}(\blambda))$ and set \[\blambda^{\ell,\chi} := \blambda(\zeta^{\ell,\chi}).\]
Finally, we define
\[
\XX_\ell^+(\blambda) := \{\chi = (j_1,j_2,p,q) \in \widetilde{\XX}_\ell^+(\blambda) \ |\ T_{\lceil i_\ell(j_2)/2 \rceil} - T_{\lceil i_\ell(j_1)/2 \rceil} \neq 0\}.
\]
\end{notation}

\begin{thm}\label{thm:jantzen_sum}
    Assume $t_i \not \in \ZZ$ for any $i = 1,\ldots,n$. The following complex rank Jantzen sum formula holds for any $\blambda \in \cP^{2n}$:
    \[
        \sum_{j \geq 1} \mathrm{ch}(\mts(\blambda)_j) = \sum_{\ell=1}^r \ \sum_{\chi \in \XX_\ell^+(\blambda)} \mathrm{sgn}(\blambda^{\ell,\chi})\mathrm{ch}(\mts(\blambda^{\ell,\chi}))
    \]
\end{thm}

\begin{proof}
    On the one hand, Lemma \ref{lem:classical_jantzen} tells us that 
    \begin{align}\label{eqn:order_vanishing}
        \sum_{\bmu \in \cP^{2n}}\mathrm{ord}_y(D_{\bs}^{\bt,R}(\blambda)_{\bmu}(y))e^{\bmu} = \sum_{j \geq 1}\sum_{\bmu \in \cP^{2n}} r_j(\mts(\blambda),\bmu)e^{\bmu} = \sum_{j \geq 1} \mathrm{ch}(\mts(\blambda)_j).
    \end{align}
    On the other hand, Theorem \ref{thm:jantzen_determinant} gives us the following explicit formula:
    \[
        D_{\bs}^{\bt,R}(\blambda)_{\bmu}(y) = \prod_{\zeta = (i,j,p,q,k) \in \SS} ((T_{\lceil j/2 \rceil}-T_{\lceil i/2 \rceil})y + \psi_{\bt}(\blambda)_{i,p} - \psi_{\bt}(\blambda)_{j,q} + s_i - s_j - k)^{\tilde{r}(\blambda(\zeta),\bmu)}.
    \]
    Each factor \begin{align}\label{eqn:factor}((T_{\lceil j/2 \rceil}-T_{\lceil i/2 \rceil})y + \psi_{\bt}(\blambda)_{i,p} - \psi_{\bt}(\blambda)_{j,q} + s_i - s_j - k)\end{align} contributes to the order of vanishing of $D_{\bs}^{\bt,R}(\blambda)_{\bmu}(y)$ only if 
    \[
        0 = \psi_{\bt}(\blambda)_{i,p} - \psi_{\bt}(\blambda)_{j,q} + s_i - s_j - k,
    \]
    whence $\zeta = (i_\ell(j_1),i_\ell(j_2),p,q,k)$ for some $\ell = 1,\ldots,r$. Letting $\chi = (j_1,j_2,p,q)$, we must have that $\delta_{\ell,\chi}(\blambda) = k > 0$, i.e., $\chi \in \widetilde{\XX}^+_\ell(\blambda)$. Under these assumptions, the factor (\ref{eqn:factor}) contributes to the order of vanishing if and only if $(T_{\lceil j/2 \rceil}-T_{\lceil i/2 \rceil}) \neq 0$, i.e., $\chi \in \XX^+_\ell(\blambda)$ and $\blambda(\zeta) = \blambda^{\ell,\chi}$. Thus,
    \[
        \mathrm{ord}_y(D_{\bs}^{\bt,R}(\blambda)_{\bmu}(y)) = \sum_{\ell=1}^r \sum_{\chi \in \XX^+_\ell(\blambda)} \tilde{r}(\blambda^{\ell,\chi},\bmu) = \sum_{\ell=1}^r \sum_{\chi \in \XX^+_\ell(\blambda)} \mathrm{sgn}(\blambda^{\ell,\chi}){r}(\blambda^{\ell,\chi},\bmu),
    \]
    where ${r}(\blambda^{\ell,\chi},\bmu)$ is the multiplicity of the $X_{\bs}^{\bt}(\bmu)$-isotypic component in $\mts(\blambda^{\ell,\chi})$. It follows that 
    \begin{align*}
        \sum_{\bmu \in \cP^{2n}} \mathrm{ord}_y(D_{\bs}^{\bt,R}(\blambda)_{\bmu}(y))e^{\bmu} &= \sum_{\bmu \in \cP^{2n}}\sum_{\ell=1}^r \sum_{\chi \in \XX^+_\ell(\blambda)} \mathrm{sgn}(\blambda^{\ell,\chi}){r}(\blambda^{\ell,\chi},\bmu)e^{\bmu} \\[5pt]
        &= \sum_{\ell=1}^r \sum_{\chi \in \XX^+_\ell(\blambda)} \sum_{\bmu \in \cP^{2n}} \mathrm{sgn}(\blambda^{\ell,\chi}){r}(\blambda^{\ell,\chi},\bmu)e^{\bmu} \\[5pt] 
        &= \sum_{\ell=1}^r \ \sum_{\chi \in \XX_\ell^+(\blambda)} \mathrm{sgn}(\blambda^{\ell,\chi})\mathrm{ch}(\mts(\blambda^{\ell,\chi})).
    \end{align*}
    Comparing to Equation (\ref{eqn:order_vanishing}) completes the proof.
\end{proof}

\subsection{Linkage Order}\label{subsec:interpolated_linkage_order}

Using the Jantzen sum formula, we introduce a partial order $\preceq$ on $\cP^{2n}$ describing a necessary condition for $\lts(\bmu)$ to be a composition factor of $\mts(\blambda)$. Eventually, we will show that the poset $(\cP^{2n},\preceq)$ gives $\ots$ the structure of an upper-filtered highest weight category. We will define the partial order $\preceq$ as a product of the inverse dominance orders from Definition \ref{defn:inv_dominance}. Moreover, the way that $\preceq$ is factored as a product will depend on the set partition $\{1,\ldots,2n\} = \Gamma_1 \sqcup \cdots \sqcup \Gamma_r$ from Notation \ref{notation:set_partition}.

\begin{notation}\label{notation:bijection}
Consider the bijection \begin{equation}\label{eqn:bijection} \cP^{2n} \xto{\sim} \cP^{a_1} \times \cdots \times \cP^{a_r}, \quad \blambda \mapsto (\ublambda_1, \ublambda_2,\ldots,\ublambda_r),\end{equation} where $\ublambda_\ell := (\lambda_{i_\ell(1)},\ldots,\lambda_{i_\ell(a_\ell)})$. We now define the sequences
\[
\b{\sigma}_\ell = (\sigma_{\ell,1},\ldots,\sigma_{\ell,a_\ell}), \quad \sigma_{\ell,k} := \gamma_{i_\ell(k)} - \gamma_{i_\ell(1)}
\]
and
\[
\b{c}_\ell = (c_{\ell,1},\ldots,c_{\ell,a_\ell}), \quad c_{\ell,k} = \begin{cases} 0 & \text{if $i_\ell(k)$ is odd} \\ 1 & \text{if $i_\ell(k)$ is even}.\end{cases}
\]
\end{notation}

For each $\ell = 1,\ldots,r$, recall the definition of the inverse dominance order of type $(\b{c}_\ell,\b{\sigma}_\ell)$ on $\cP^{a_\ell}$ from Definition \ref{defn:inv_dominance}.

\begin{defn}
    We define the \textit{linkage order} $\preceq$ on $\cP^{2n}$ by declaring $\bnu \preceq \blambda$ whenever $\ubnu_\ell \leq \ublambda_\ell$ in the inverse dominance order of type $(\b{c}_\ell,\b{\sigma}_\ell)$ on $\cP^{a_\ell}$ for each $\ell = 1,\ldots,r$. 
\end{defn}

The following lemma is the key ingredient in showing that the linkage order $\preceq$ gives a necessary condition for $[\mts(\blambda):\lts(\bnu)]$ to be nonzero.

\begin{lem}\label{lem:one_step_down}
    For any $\ell = 1,\ldots,r$ and any $\chi \in \XX_\ell^+(\blambda)$, we have $\blambda^{\ell,\chi} \preceq \blambda$ in the linkage order.
\end{lem}

\begin{proof}
    For concision, set $\bnu := \blambda^{\ell,\chi}$. We will prove the following stronger claim: we have $\ubnu_\ell \leq \ublambda_\ell$ in the inverse dominance order of type $(\b{c}_\ell,\b{\sigma}_\ell)$ on $\cP^{a_\ell}$, and for each $\ell' \in \{1,\ldots,r\}$ with $\ell' \neq \ell$, we have $\ubnu_{\ell'} = \ublambda_{\ell'}$. Let us write $\chi = (j_1,j_2,p,q)$. Pick large positive integers $m_1,\ldots,m_n$ so that \[\Theta_{\bs}^{\b{m}}(\bnu) = (\Theta_{\bs}^{\b{m}}(\blambda) - \delta^{\ell,\chi}(\blambda)\alpha_{\ell,\chi})_+,\] where, using the notation \[p' := \begin{cases}p & \text{if $i_\ell(j_1)$ is odd} \\ m_{i_\ell(j_1)/2} - p + 1 & \text{else}, \end{cases} \quad q' := \begin{cases}q & \text{if $i_\ell(j_2)$ is odd} \\ m_{i_\ell(j_2)/2} - q + 1 & \text{else}. \end{cases},\] we have \[\alpha_{\ell,\chi} := \alpha^{p',q'}_{\lceil i_\ell(j_1)/2\rceil,\lceil i_\ell(j_2)/2\rceil} \in R^+ \setminus R^+_{\b{m}}.\] When $\ell' \neq \ell$ and the integers $m_i$ are all sufficiently large, subtracting $\alpha$ from $\Theta_{\bs}^{\b{m}}(\blambda)$ does not affect the entries outside of the ``blocks'' corresponding to the indices in $\Gamma_\ell$. In particular, we have $\ubnu_{\ell'} = \ublambda_{\ell'}$ when $\ell' \neq \ell$. 

    Now fix some very large integer $N \gg 0$ and define the sequence
    \[
        \Theta^{N}(\blambda)_\ell := (\theta^{N}(\blambda)_{\ell,1},\ldots,\theta^N(\blambda)_{\ell,a_\ell}) \in (\fk^N)^{a_\ell}, 
    \]
    where each $\theta^N(\blambda)_{\ell,k} \in \fk^N$ is defined by \[\theta^N(\blambda)_{\ell,k} := \begin{cases} (\lambda_{i_\ell(k),1} + \sigma_{\ell,k},\lambda_{i_\ell(k),2} + \sigma_{i_\ell(k)},\ldots,\lambda_{i_\ell(k),N} + \sigma_{\ell,k}) & \text{if $i_\ell(k)$ is odd} \\ 
    (-\lambda_{i_\ell(k),N} + \sigma_{\ell,k},\ldots,-\lambda_{i_\ell(k),2} + \sigma_{\ell,k},-\lambda_{i_\ell(k),1} + \sigma_{\ell,k}) & \text{if $i_\ell(k)$ is even}.
    \end{cases} \] In the definition above, we understand that $\lambda_{i_\ell(k),j} = 0$ if $j > \ell(\lambda_{i_\ell(k)})$. We similarly define the sequence $\Theta^N(\bnu)_\ell \in (\fk^N)^{a_\ell}$. 
    
    We now regard these sequences as \textit{virtual $a_\ell$-multipartitions} in the sense of Definition \ref{defn:virtual_partition}. We can now adapt the argument in the proof of Proposition \ref{prop:construct_refine} \textit{mutatis mutandis} to this situation, noting that the virtual multipartition $\Theta^N(\bnu)_\ell$ is obtained from $\Theta^N(\blambda)_\ell$ by some sequence of virtual box moves, each of which does not affect the total number of boxes of a given virtual content. Translating the statement about virtual contents into the language of contents of boxes in $\ubnu_\ell$ and $\ublambda_\ell$, just as in Proposition \ref{prop:construct_refine}, then completes the proof. We leave details to the reader.
\end{proof}

We now apply the Jantzen sum formula in tandem with the previous lemma to relate $\preceq$ with composition multiplicities. We make the following observation.
For any $\bnu \in \cP^{2n}$, any $V \in \ots$ whose $\fz_{\bt}^*$-weights are bounded above by a finite set of weights has a finite filtration $\mathscr{F}_{\bnu}$ given by \[0 = V_0 \subset V_1 \subset V_2 \subset \cdots \subset V_M = V\] such that for each $i = 1,\ldots,M$, either (1) $V_{i}/V_{i-1} \simeq \lts(\bnu_i)$ with $\bnu \trianglelefteq \bnu_i$ or (2) the $X^{\bt}(\b{\eta})$-isotypic component of $V_{i}/V_{i-1}$ is zero for any $\b{\eta} \trianglerighteq \bnu$ -- for instance, the proof of \cite[Lemma 2.1.9]{Kumar} readily generalizes to our setting. For any $\bnu' \trianglerighteq \bnu$, we define $[V:\lts(\bnu')]_{\mathscr{F}_{\bnu}}$ as the number of times that $\bnu'$ shows up in the filtration $\mathscr{F}_{\bnu}$. An easy argument using characters shows that this number is independent of the choice of $\bnu \trianglelefteq_{\bt} \bnu'$ as well as the filtration $\mathscr{F}_{\bnu}$. Moreover, any finite filtration of $V$ can be refined to a filtration of the form $\mathscr{F}_{\bnu}$ for some $\bnu \in \cP^{2n}$. Thus, we see that $[V:\lts(\bnu')]_{\mathscr{F}_{\bnu}}$ coincides with the usual composition multiplicity $[V:\lts(\bnu')]$. We deduce (e.g., following the proof of \cite[Lemma 2.1.12]{Kumar}) that for any $\bnu \in \cP^{2n}$, \begin{align}\label{eqn:composition_character}\mathrm{ch}(\mts(\bnu)) = \sum_{\bnu' \in \cP^{2n}} [\mts(\bnu):\lts(\bnu')] \mathrm{ch}(\lts(\bnu')).\end{align} 

\begin{prop}\label{prop:linkage}
    If $[\mts(\blambda):\lts(\bmu)] \neq 0$ for $\blambda,\bmu \in \cP^{2n}$, then $\bmu \preceq \blambda$ in the linkage order.
\end{prop}

\begin{proof}
    We induct on the number $N$ of multipartitions $\bnu \in \cP^{2n}$ with $\bmu \trianglelefteq \bnu \trianglelefteq \blambda$, with trivial base case $N = 1$. Otherwise, we have $\bmu \triangleleft \blambda$, so $\lts(\bmu)$ must be a composition factor of the maximal submodule $\mts(\blambda)_1$. Then, the linear independence of the characters $\mathrm{ch}(\lts(\bnu'))$, Equation (\ref{eqn:composition_character}), and the Jantzen sum formula from Theorem \ref{thm:jantzen_sum} jointly imply that $\lts(\bmu)$ must be a composition factor in $\mts(\blambda^{\ell,\chi})$ for some $\chi \in \XX_\ell^+(\blambda)$. Since $\blambda^{\ell,\chi} \triangleleft\blambda$, the number of $\bnu \in \cP^{2n}$ with $\bmu \trianglelefteq \bnu \trianglelefteq \blambda^{\ell,\chi}$ is strictly less than $N$, so we are done by Lemma \ref{lem:one_step_down} and induction.
\end{proof}

\begin{cor}\label{cor:linkage_duality}
    For $\blambda,\bmu \in \cP^{2n}$, we have $[\checkmts(\blambda): \lts(\bmu)] \neq 0$ only if $\bmu \preceq \blambda$.
\end{cor}

\begin{proof}
    Since $\DD$ is an equivalence, we have \[[\checkmts(\blambda): \lts(\bmu)]  = [\DD M^{w_0\bt}_{w_0\bs}(\blambda^*): \DD L_{w_0 \bs}^{w_0 \bt}(\bmu^*)] = [M^{w_0\bt}_{w_0\bs}(\blambda^*): L_{w_0 \bs}^{w_0 \bt}(\bmu^*)].\] Let $\preceq_{w_0}$ denote the linkage order on $\cO^{w_0\bt}_{w_0\bs}$. Then, it is easy to check that $\bmu^* \preceq_{w_0} \blambda^*$ if and only if $\bmu \preceq \blambda$. The result follows from Proposition \ref{prop:linkage}.
\end{proof}

\subsection{Projectives, Highest Weight Structure, and Restricted Category \strO}\label{subsec:extended_category_o}

We are now ready to put the category $\ots$ in the framework of multi-Fock tensor product categorifications.




\subsubsection{Projectives}\label{subsubsec:projectives_truncations}

We first describe a ``filtration" of $\ots$ by subcategories with a projective generating family. We record some facts regarding objects with finite filtrations by Verma modules. We omit their proofs, since their classical proofs directly generalize to our complex rank setting.

\begin{lem}\label{lem:verma_filtration_summand}
    The direct summand of any object with a finite filtration by Verma modules also has a finite filtration by Verma modules.
\end{lem}

\begin{cor}\label{cor:kernel_verma_filtration}
    Suppose $M \to \mts(\blambda)$ is a surjective morphism and $M$ has a finite filtration by Verma modules. Then, the kernel of this morphism also has a finite filtration by Verma modules. 
\end{cor}

\begin{lem}\label{lem:verma_projective}
    For each $\blambda \in \cP^{2n}$, the Verma module $\mts(\blambda)$ is the projective cover of $\lts(\blambda)$ in the Serre subcategory $\ots(\preceq \blambda)$ spanned by the simple objects $\lts(\bmu)$ for $\bmu \preceq \blambda$.
\end{lem}

\begin{proof}
    The functor $\Hom_{\tildeots}(\mts(\blambda),-)$ sends any $M \in \ots$ to the $\fl_{\bt}$-submodule of its $X_{\bs}^{\bt}(\blambda)$-isotypic component that is annihilated by $\fp_{\bt}$. On objects in $\tildeots(\preceq\blambda)$, this submodule coincides with the entire $X_{\bs}^{\bt}(\blambda)$-isotypic component. Thanks to the semisimplicity of $\Ind \cL_{\bt}$, the functor sending $M$ to its $X_{\bs}^{\bt}(\blambda)$-isotypic component is exact, so we are done. 
\end{proof}

Recall that we identify $\fz_{\bt}^* \simeq \fk^n$ with weights for the classical Lie algebra $\ul{\fgl}_n$. We write $\alpha_1,\ldots,\alpha_{n-1}$ for the simple roots of $\ul{\fgl}_n$ under this identification. Moreover, for any $\beta \in \fz_{\bt}^*$ and $X \in \widehat{\cO}^{\bt}_{\bs}$, recall that $X^\beta \in \Ind \cL_{\bt}$ denotes the $\beta$-weight subobject in $X$. 

\begin{defn}
For any poset ideal $I \subset \fz_{\bt}^*$, we write $\ots[I]$ (resp. $\hatots[I]$) for the full subcategory of objects $M \in \ots$ (resp. $M \in \hatots$) such that $M^\gamma \neq 0$ for any $\gamma \in \fz_{\bt}^*$ implies $\gamma \in I$.
\end{defn}

\begin{prop}\label{prop:verma_quotient_of_projective}
    For any \textup{finitely generated} ideal $I \subset \fz_{\bt}^*$, the subcategory $\ots[I]$ has enough projectives. Moreover, for any $\blambda \in \omega_{\bs}^{-1}(I)$, the kernel of the projective cover $P^I_{\bs}(\blambda) \twoheadrightarrow \mts(\blambda)$ has a finite filtration by Verma modules $\mts(\bmu)$ with $\blambda \preceq \blambda$. 
\end{prop}

\begin{proof}
    Our proof closely follows the approach of \cite[Theorem 3.2(1)]{soergel98}.
    
    Write $U(\fp^{\bt})\modcat_I$ for the full subcategory of $U(\fp^{\bt})\modcat$ consisting of objects whose $\fz_{\bt}^*$-weights lie in the poset ideal $I$, and let $\iota_I: U(\fp^{\bt})\modcat_I \to U(\fp^{\bt})\modcat$ denote the inclusion functor. We write $\iota_I^!: U(\fp^{\bt})\modcat \to U(\fp^{\bt})\modcat_I$ for its left adjoint, which sends any object $M \in U(\fp^{\bt})\modcat$ to its largest quotient contained in $U(\fp^{\bt})\modcat_I$. Now, define
    \[
        Q_{\bs}^I(\blambda) := U_{\bt} \ox_{U(\fp_{\bt})} \iota_I^!(U(\fp_{\bt}) \ox_{U(\fl_{\bt})} X_{\bs}^{\bt}(\blambda)).
    \]
    By definition, there is a surjection $Q_{\bs}^I(\blambda) \to \mts(\blambda)$.
    Since the ideal $I$ is finitely generated, the parabolic subalgebra $\fp_{\bt}$ acts locally finitely on $Q_{\bs}^I(\blambda)$, and as an object of $\Ind \cL^{\bt}$, every irreducible constituent of $Q_{\bs}^I(\blambda)$ has the form $X_{\bs}^{\bt}(\bmu)$ for some $\bmu \in \cP^{2n}$. Hence, we have $Q_{\bs}^I(\blambda) \in \widehat{\cO}^{\bt}_{\bs}$.  In fact, since $I$ is a poset ideal, all $\fz_{\bt}$-weights of $Q_{\bs}^I(\blambda)$ belong to $I$ and hence $P_{\bs}^I(\blambda) \in  \widehat{\cO}^{\bt}_{\bs}[I]$ as well. By the usual Frobenius reciprocity argument, we see that $Q_{\bs}^I(\blambda)$ represents the exact functor $\widehat{\cO}^{\bt}_{\bs} \to \Ind \cL_{\bt}$ sending any $M \in \widehat{\cO}^{\bt}_{\bs}$ to its $X_{\bs}^{\bt}(\blambda)$-isotypic component. It follows that $Q_{\bs}^I(\blambda)$ is projective as an object in $\widehat{\cO}^{\bt}_{\bs}$. 

    We will now show that $Q_{\bs}^I(\blambda)$ has a finite filtration by Verma modules. Thanks to the exactness of the functor $U_{\bt} \ox_{U(\fp_{\bt})} -$, it suffices to show that $\iota_I^!(U(\fp_{\bt}) \ox_{U(\fl_{\bt})} X_{\bs}^{\bt}(\blambda))$ has a filtration by objects of the form $X_{\bs}^{\bt}(\b{\nu})$, viewed as a $U(\fp_{\bt})$-module via inflation. As a module over $U(\fl_{\bt})$, we certainly have a decomposition
    \[\iota_I^!(U(\fp_{\bt}) \ox_{U(\fl_{\bt})} X_{\bs}^{\bt}(\blambda)) = X_{\bs}^{\bt}(\blambda) \oplus \bigoplus_{\b{\nu} \in \cP_{\blambda}} X_{\bs}^{\bt}(\b{\nu})\]
    for some subset $\cP_{\blambda} \subset \cP^{2n}$. This subset is in fact finite, because $I$ is finitely generated. Moreover, the weights $\omega_{\bs}(\b{\nu})$ are pairwise distinct and $\omega_{\bs}(\blambda) \triangleleft_{\mathrm{dom}} \omega_{\bs}(\b{\nu})$ for all $\b{\nu} \in \cP_{\blambda}$ by definition. Now, we will pick an ordering $\b{\nu}_1,\ldots,\b{\nu}_v$ such that $\omega_{\bs}(\b{\nu}_i) \triangleleft_{\mathrm{dom}} \omega_{\bs}(\b{\nu}_j)$ only if $j > i$. For $j = 1,\ldots,v$, let us write $N_j \subset \iota_I^!(U(\fp_{\bt}) \ox_{U(\fl_{\bt})} X_{\bs}^{\bt}(\blambda))$ for the $U(\fp_{\bt})$-submodule generated by the subobjects $X_{\bs}(\b{\nu}_i)$ for $i = 1,\ldots,j$. Then, our choice of ordering implies that 
     \[0 \subset N_1 \subset N_2 \subset \cdots \subset N_v \subset \iota_I^!(U(\fp_{\bt}) \ox_{U(\fl_{\bt})} X_{\bs}^{\bt}(\blambda))\]
    is the desired filtration.
    
    By Corollary \ref{cor:kernel_verma_filtration}, the kernel of the surjection $Q_{\bs}^I(\blambda) \to \mts(\blambda)$ has a {finite} filtration by Verma modules $\mts(\b{\nu})$ where $\omega_{\bs}(\b{\nu}) \triangleright_{\mathrm{dom}} \omega_{\bs}(\blambda)$. Thus, this projective object $Q_{\bs}^I(\blambda)$ must belong to $\ots[I]$ itself. Let $P_{\bs}^I(\blambda)$ be any indecomposable summand with a nonzero surjection to $\mts(\blambda)$. By Lemma \ref{lem:verma_filtration_summand}, we see that $P_{\bs}^I(\blambda)$ has a finite Verma filtration as well. Any $\mts(\bnu)$ in the kernel of the resulting surjection to $\mts(\blambda)$ must have $\bnu \succeq \blambda$ in the linkage order. Indeed, by Lemma \ref{lem:homological_verma_filtration}, this follows from the easy observation that $[\check{M}_{t}^s(\bnu): L(\blambda)] \neq 0$ only if $\bnu \succeq \blambda$. 
    We remark that, by construction, the summand $P_{\bs}^I(\blambda)$ is in fact projective in $\widehat{\cO}_{\bs}^{\bt}[I]$.
\end{proof}

For any poset ideal $I \subset \fz_{\bt}^*$ with respect to the dominance order, observe that $\omega_{\bs}^{-1}(I) \subset \cP^{2n}$ is a poset ideal with respect to the linkage order $\preceq$, since the map $\omega_{\bs}$ is order preserving.

\begin{cor}\label{cor:ots_I_schurian}
    For any finitely generated poset ideal $I \subset \fz_{\bt}^*$, the subcategory $\ots[I]$ is Schurian.
\end{cor}

\begin{proof}
 Consider the locally unital algebra \[A := \bigoplus_{\blambda,\b{\nu} \in \omega_{\bs}^{-1}(I)} \Hom_{\ots}(P^I_{\bs}(\blambda),P^I_{\bs}(\b{\nu})).\] Note that $A$ is locally finite-dimensional because each $P^I_{\bs}(\blambda)$ has finite composition multiplicities. Moreover, each $P^I_{\bs}(\blambda)$ is a compact object in $\hatots[I]$, as any $P^I_{\bs}(\blambda)$ is finitely generated as a $U_{\bt}$-module. By its definition, the category $\hatots[I]$ is closed under small coproducts. Thus, by \cite[Lemma 2.4]{brundan_stroppel_semiinfinite}, we deduce that $\hatots[I]$ is equivalent to $A^{\op}\modcat$. In turn, \cite[Lemma 2.13]{brundan_stroppel_semiinfinite} implies that $\ots[I]$ is Schurian.
\end{proof}

\begin{defn}
For an upper finite poset ideal $J \subset \cP^{2n}$ (with respect to the linkage order $\preceq$), we write $\ots(J)$ for the Serre subcategory spanned by the simples $\lts(\blambda)$ for $\blambda \in J$.
\end{defn}

Fix an upper finite poset ideal $J \subset \cP^{2n}$. For each \textit{finite} subset $S \subset J$, we write $I_S \subset \fz_{\bt}^*$ for the poset ideal generated by $\omega_{\bs}(\blambda)$ for $\blambda \in S$. In turn, we define $J_S := \omega_{\bs}^{-1}(I_S) \cap J$. 

Define $P^S_{\bs}(\blambda) := \iota^!_S P^{I_S}_{\bs}(\blambda)$, where $\iota^!$ is left adjoint to the inclusion functor $\iota: \ots(J_S) \to \ots[I_S]$. That is, the object $P^S_{\bs}(\blambda)$ is the largest quotient of $P^{I_S}(\blambda)$ belonging to $\hatots(J_S)$. Then, $P^S_{\bs}(\blambda)$ admits a surjection to $\mts(\blambda)$ for all $\blambda \in J_S$, and in fact, is the projective cover of $\lts(\blambda)$ in $\ots(J_S)$. Note that $P^S_{\bs}(\blambda)$ is also projective in $\hatots(J_S)$, since it is also the maximal quotient of the projective object $P^{I_S}_{\bs}(\blambda) \in \hatots[I_S]$ belonging to $\hatots(J_S)$. 

\begin{lem}\label{lem:projectives_in_finite}
    For each finite subset $S \subset J$, the Serre subcategory $\ots(J_S)$ is an upper finite highest weight category with poset $(J_S, \preceq)$ and standard objects $\mts(\blambda)$ for $\blambda \in J_S$.
\end{lem}

\begin{proof}
    By \cite[Lemma 2.24]{brundan_stroppel_semiinfinite}, the Serre subcategory $\ots(J_S) \subset \ots[I_S]$ is also Schurian. Thanks to Lemma \ref{lem:verma_projective}, it remains to check is that the kernel of the surjection $P^S_{\bs}(\blambda) \twoheadrightarrow \mts(\blambda)$ is filtered by Verma modules $\mts(\bmu)$ where $\bmu \succ \blambda$. 

    A completely analogous argument to the proof of Lemma \cite[Lemma 3.44]{brundan_stroppel_semiinfinite} tells us that $\Ext_{\ots[I_S]}^1(\mts(\bnu),\mts(\bmu)) = 0$ unless $\bmu \preceq \bnu$.
    Thus, we may rearrange the Verma module filtration of $P^{I_S}_{\bs}(\blambda)$ so that all sections $\mts(\bmu)$ with $\bmu \not \in J$ appear below the sections with $\bmu \in J$, i.e., there is a short exact sequence 
    \[
        0 \to N_1 \to P^{I_S}_{\bs}(\blambda) \to N_2 \to 0
    \] 
    such that $N_1$ has a filtration by Verma modules $\mts(\bmu)$ with $\bmu \not\in J$ and $N_2$ has a filtration by Verma modules $\mts(\bmu)$ such that $\bmu \in J$ and $\mts(\blambda)$ appears at the top of $N_2$. Then, it is not hard to see that $N_2 = P^S_{\bs}(\blambda)$. Moreover, by Proposition \ref{prop:verma_quotient_of_projective}, all Verma modules $\mts(\bmu)$ appearing below $\mts(\blambda)$ in $P^{I_S}_{\bs}(\blambda)$ already satisfy $\bmu \succ \blambda$, so we are done.
\end{proof}

Now suppose we have an inclusion $S \subset S'$ for some finite subset $S' \subset J$. Then, there is a surjective map $P^{S'}_{\bs}(\blambda) \twoheadrightarrow P^S_{\bs}(\blambda)$. Define the inverse limit
\[
    P^J_{\bs}(\blambda) := \varprojlim_{S \subset J} P^S_{\bs}(\blambda). 
\]

\begin{lem}\label{lem:enough_projectives}
    We have an equality $P^J_{\bs}(\blambda) = P^{S_0}_{\bs}(\blambda)$ for some finite subset $S_0 \subset J$. In particular, there is a surjection $P^J_{\bs}(\blambda)\twoheadrightarrow \mts(\blambda)$ with kernel filtered by Verma modules $\mts(\bmu)$ with $\bmu \succ \blambda$. Moreover, $P^J_{\bs}(\blambda)$ is projective in $\hatots(J)$.
\end{lem}

\begin{proof}
    Let $S_0 \subset J$ denote the set of labels $\bmu \in J$ for which $\bmu \succeq \blambda$. Note that this set is finite because we assumed that $J$ is upper finite. For any finite subset $S \subset J$ containing $S_0$, we claim that $P^{}_{\bs}(\blambda) = P^{S_0}_{\bs}(\blambda)$. We consider the canonical surjective morphism $P^{S}_{\bs}(\blambda) \twoheadrightarrow P^{S_0}_{\bs}(\blambda)$. By arguments from the proof of Lemma \ref{lem:projectives_in_finite}, the kernel of this epimorphism is filtered by Verma modules $\mts(\bmu)$ with $\bmu \in J_{S} \setminus J_{S_0}$ such that $\bmu \succ \blambda$. By the very construction of $S_0$, however, no such labels $\bmu$ exist, so the epimorphism $P^{S}_{\bs}(\blambda) \twoheadrightarrow P^{S_0}_{\bs}(\blambda)$ is actually an isomorphism. 

    It remains to show that $P^J_{\bs}(\blambda) = P^{S_0}_{\bs}(\blambda)$ is projective in $\hatots(J)$. Suppose we have an epimorphism $g: X \twoheadrightarrow Y$ in $\hatots(J)$ along with a morphism $f: P^J_{\bs}(\blambda) \to Y$. Write $Y' \subset Y$ for the image of $f$. Note that $Y'$ is finitely generated as a $U_{\bt}$-module since $P^J_{\bs}(\blambda)$ is. In particular, we can pick a finitely generated subobject $X'\subset X$ such that $g$ restricts to a surjection $X' \twoheadrightarrow Y'$. This means that all $\fz_{\bt}^*$-weights of $X'$ belong to some finitely generated poset ideal $I \subset \fz_{\bt}^*$. Then, we consider the poset ideal $K := (\omega_s^{-1}(I) \cap J) \cup J_{S_0}$ and the corresponding Serre subcategory $\ots(K)$. In particular, the objects $X'$ and $Y'$ both belong to $\ots(K)$. 
    Write $P^K_{\bs}(\blambda)$ for the largest quotient of $P^{I_S}_{\bs}(\blambda)$ belonging to $\ots(K)$. Since $S_0 \subset K$, the argument given in the first paragraph of this proof generalizes to show that $P^K_{\bs}(\blambda)$ actually equals $P^{J}_{\bs}(\blambda)$. Thus, $P^J_{\bs}(\blambda)$ is also projective in $\ots(K)$. Then, there exists $h: P^J_{\bs}(\blambda) \to X' \subset X$ such that $g \circ h = f$, as desired.
\end{proof}

\begin{cor}\label{cor:ots_highest_weight}
    For any upper finite poset ideal $J \subset \cP^{2n}$, the subcategory $\ots(J)$ is Schurian. In particular, it is an upper finite highest weight category with poset $(J,\preceq)$. 
\end{cor}

\begin{proof}
   In light of Lemma \ref{lem:enough_projectives}, the proof is completely analogous to that of Corollary \ref{cor:ots_I_schurian}.
\end{proof}

\subsubsection{Highest Weight Structure and the Extended Category} \label{subsubsec:highest_weight}

We finally define the restricted category $\extots$ advertised in the introduction and Conjecture \ref{conj:categorical_conjecture}. In brief, the category $\ots$ is too ``large" to have the structure of an MFTPC (it fails to satisfy (MF2) from Definition \ref{defn:multi_fock}), so we will instead work with a natural subcategory.

\begin{defn}\label{defn:extots}
    Define the subcategory $\extots \subset \ots$ as the full subcategory consisting of objects $M \in \ots$ for which there exists a \textit{finite} subset $E_M \subset \cP^{2n}$ such that for every composition factor $\lts(\blambda)$ of $M$, the label $\blambda$ is bounded above (in the linkage order $\preceq$) by some label in $E_M$.
\end{defn}
In particular, any simple object and any object with a finite filtration by Verma modules belongs to $\extots$. Moreover, this restricted subcategory $\extots$ is clearly an abelian subcategory of $\ots$. Moreover, when the poset $\cP^{2n}$ factors as a product of \textit{admissible} posets, the subcategory $\extots$ can be described rather explicitly.

Recall the bijection \[\cP^{2n} \simeq \cP^{a_1} \times \cdots \times \cP^{a_r}, \quad \blambda \mapsto (\ublambda_1,\ldots,\ublambda_r)\] from Notation \ref{notation:bijection}. Assume each $\cP^{a_\ell}$ is either upper or lower admissible with respect to the inverse dominance order. For $j \geq 0$ and $\ell = 1,\ldots,r$, we write $\Lambda_{\ell}(j) \subset \cP^{a_\ell}$ for the subset given by \[\Lambda_{\ell}(j) := \begin{cases} \cP^{a_\ell} &\text{if $c_\ell$ is upper admissible} \\ \cP^{a_\ell}(j) & \text{if $c_\ell$ is lower admissible}. \end{cases}\] In turn, we define $\cP^{2n}(j) \subset \cP^{2n}$ to be the \textit{upper finite poset ideal} \[\Lambda_1(j) \times \Lambda_2(j) \times \cdots \times \Lambda_r(j) \subset \cP_{a_1} \times \cdots \times \cP^{a_r} \simeq \cP^{2n}.\]  Write $(\ots)^j \subset \ots$ for the Serre subcategory spanned by the simple objects labelled by $\blambda \in \cP^{2n}(j)$. By Corollary \ref{cor:ots_highest_weight}, this Serre subcategory is an upper finite highest weight category with poset $(\cP^{2n}(j), \preceq)$ and standard objects $\mts(\blambda)$ for $\blambda \in \cP^{2n}(j)$. Then, the following lemma, which demonstrates that $\extots$ satisfies (MF2) from Definition \ref{defn:multi_fock}, is manifest.

\begin{lem}\label{lem:extots_admissible}
    When each $\cP^{a_\ell}$ is upper or lower admissible, every object in $\extots$ is isomorphic to an object in $(\ots)^j$ for some $j \geq 0$. That is, 
    \[
        \extots = \bigcup_{j\geq 0} (\ots)^j.
    \]
\end{lem}


\section{Categorical Actions and Multiplicities in Complex Rank Category \strO}\label{sec:categorical_actions}

In this section, we study a natural categorical type A action on the category $\extots$. In the \textit{admissible} case, we will show that this categorical action gives $\extots$ the structure of a multi-Fock tensor product categorification (MFTPC) in the sense of Definition \ref{defn:multi_fock}. In turn, we will leverage the results of Sections \ref{sec:uniqueness} and \ref{sec:construction} to explicitly describe $\extots$ in terms of (stable) parabolic categories $\cO$ in integer rank. When $\fk =\CC$, these structural results also give the desired proof of Conjecture \ref{conj:categorical_conjecture} in the admissible case.

\subsection{Categorical Type A Action and MFTPCs}\label{subsec:categorical_action}

Observe that the endofunctors
\[F: M \mapsto V \ox M, \quad E: M \mapsto V^* \ox M\] of $\fg_{\bt}\modcat$ restrict to endofunctors of $\extots$. Thus, the general construction given in Example \ref{ex:categorical_type_a_glX} yields a categorical type A action $(E,F,x,\tau)$ on $\extots$, where the natural transformation $x \in \End(F)$ is given by the action of the tensor Casimir (see Example \ref{ex:categorical_type_a_glX} for the definition), and $\tau \in \End(F^2)$ comes from the symmetric braiding. 

The usual argument, along with \cite[Corollary 7.1.2]{comes_wilson}, shows that $F\mts(\blambda) = V \ox \mts(\blambda)$ (resp. $E\mts(\blambda) = V^* \ox \mts(\blambda)$) admits a filtration with sections $\mts(\blambda')$, where $\blambda'$ is obtained from $\blambda$ by either adding (resp. removing) a single box to some partition $\lambda_{2i-1}$ or removing (resp. adding) a single box from a partition $\lambda_{2i}$ for $i = 1,\ldots,n$. 

For $i = 1,\ldots,n$, recall the constants \[T_i := t_1 + \cdots + t_{i-1}, \quad \gamma_{2i-1} := s_i - T_i, \quad \gamma_{2i} := s_i - T_{i+1}.\] 

\begin{lem}\label{lem:eigenvalues}
    The eigenvalues of $x$ on $FM_{\bds{s}}^t(\blambda)$ are given by 
    \[
    p+\gamma_{2i-1}, \quad -q + \gamma_{2i},
    \]
    where $i$ ranges from $1$ to $n$, $p$ ranges through the contents of all addable boxes for $\lambda_{2i-1}$, and $q$ ranges through contents of all removable boxes for $\lambda_{2i}$. 
\end{lem}

\begin{proof}
    This immediately follows from our calculation of the action of $C_2$ on Verma modules in Proposition \ref{prop:central_character} and the explicit descriptions of $F\mts(\blambda)$ and $E\mts(\blambda)$. 
\end{proof}

Recall the set partition
\[
\{1,\ldots,2n\} = \Gamma_1 \sqcup \Gamma_2 \sqcup \cdots \sqcup \Gamma_r
\]
from Notation \ref{notation:set_partition}, along with the sequences $\b{\sigma}_\ell$ and $\b{c}_\ell$ from Notation \ref{notation:bijection}. Define \[\Xi_{\bs}^{\bt} := ((\b{\sigma}_1,\b{c}_1),(\b{\sigma}_2,\b{c}_2),\ldots,(\b{\sigma}_r,\b{c}_r)).\] 
Observe that $\Xi_{\bs}^{\bt}$ is an admissible type if and only if each poset factor $\cP^{a_\ell}$ from Notation \ref{notation:bijection} is either upper or lower admissible. 

By Lemma \ref{lem:eigenvalues}, we have an eigenfunctor decomposition
\[
F = \bigoplus_{\ell=1}^r \bigoplus_{m \in \ZZ} F_{\gamma_{i_\ell(1)} + m}, \quad E = \bigoplus_{\ell=1}^r \bigoplus_{m \in \ZZ} E_{\gamma_{i_\ell(1)} + m},
\]
where the subscripts denote the corresponding eigenvalue. In this notation, $F_{\gamma_{i_\ell(1)} + m}\mts(\blambda)$  (resp. $E_{\gamma_{i_\ell(1)+m}}\mts(\blambda)$) admits a filtration with sections $\mts(\blambda')$, where $\blambda'$ is obtained from $\blambda$ by either adding (resp. removing) a box of content $m - \sigma_{\ell,k}$ to $\lambda_{i_\ell(k)}$ if $c_{\ell,k} = 0$ (if such a partition exists) or removing (resp. adding) a box of content $-m+\sigma_{\ell,k}$ from $\lambda_{i_\ell(k)}$ if $c_{\ell,k} = 1$. 

We define 
\[
F(\ell) := \bigoplus_{m \in \ZZ} F_{\gamma_{i_\ell(1)} + m}, \quad E(\ell) = \bigoplus_{m \in \ZZ} E_{\gamma_{i_\ell(1)} + m}.
\]
Decompose $x = \tilde{x}(1) + \cdots + \tilde{x}(r)$, so that the endomorphisms $x(\ell) := \tilde{x}(\ell) - \gamma_{i_\ell(1)}\id$ are pairwise commuting endomorphisms of $F$ with eigenvalues belonging to $\ZZ$. In particular, each $x(\ell)$ restricts to an endofunctor of $F(\ell)$ acting on $F_{\gamma_{i_\ell(1)} + m}$ with eigenvalue $m \in \ZZ$. Finally, the endofunctor $\tau$ restricts to endofunctors $\tau(\ell)$ of each $F(\ell)^2$. 

\begin{thm}\label{thm:ots_mftpc}
    Assume $t_i \not \in \ZZ$ for $i = 1,\ldots,n$ and that $\Xi_{\bs}^{\bt}$ is an admissible type. Then, the datum $(E,F,x,\tau)$ equips the category $\extots$ with the structure of a (full) MFTPC of type $\Xi_{\bs}^{\bt}$.
\end{thm}

\begin{proof}
    Combine Corollary \ref{cor:ots_highest_weight} and Lemma \ref{lem:extots_admissible} with the explicit descriptions of the functors $F_{\gamma_{i_\ell}(1)+m}, E_{\gamma_{i_\ell}(1)+m}$ on standard objects.
\end{proof}


\subsection{Structural Results for Admissible Categories \texorpdfstring{$\extots$}{O}} \label{subsec:structural_results}
Under combinatorial conditions on the parameters $\bt$ and $\bs$ (namely, so that the sequences $\b{c}_\ell$ are admissible), we use our results from Sections \ref{sec:uniqueness} and \ref{sec:construction} to completely understand the structure of $\extots$.

\begin{ex}\label{ex:admissibility_two_blocks}
In the case $n = 2$, the admissibility condition is extraneous. Let us write $s := s_2 - s_1$. As explained in Subsection \ref{subsec:categorical_action}, the triple $\{t_1,t_2,s\}$ yields a set partition \[\{1,2,3,4\} = \Gamma_1 \sqcup \Gamma_2 \sqcup \cdots \sqcup \Gamma_r.\]
Admissibility means the subsets $\Gamma_\ell$ cannot equal $\{1,2,3\}, \{2,3,4\}$, or $\{1,2,3,4\}$, which were already made impossible by the non-integrality of $t_1$ and $t_2$.
\end{ex}

In particular, if the parameters $(\bt,\bs)$ are admissible, then $\extots$ is a \textit{mixed admissible} MFTPC of type $\Xi_{\bs}^{\bt}$. Thus, we can use the results of Sections \ref{sec:uniqueness} and \ref{sec:construction} to completely understand the structure of $\extots$. For example, recall the stable category $\cO^{\Xi_{\bs}^{\bt}}_\infty$ from Definition \ref{defn:stable_mftpc}. We can now make rigorous the idea that parabolic category $\cO$ in complex rank is a ``stable limit'' of ordinary parabolic categories $\cO$. The following corollary is a direct consequence of Theorem \ref{thm:recover_from_gluing}.

\begin{cor}
    Assume $t_i \not \in \ZZ$ for $i = 1,\ldots,n$ and that the type $\Xi_{\bs}^{\bt}$ is admissible. Then, there is a strongly equivariant equivalence \[\extots \simeq \cO^{\Xi_{\bs}^{\bt}}_\infty\] intertwining the labels of simples.
\end{cor}

Moreover, Proposition \ref{prop:mftpc_multiplicities} allows us to compute multiplicities of simple objects in Verma modules, giving the following answer to Problem \ref{question:etingof} in the admissible case. Recall the category $\cO_m^{\b{c}_\ell,\b{\sigma}_\ell}(N)$ from Definition \ref{defn:partition_subcategory}.

\begin{cor}\label{cor:admissible_multiplicities}
     Assume $t_i \not \in \ZZ$ for $i = 1,\ldots,n$ and that the parameters $\{\bt,\bs\}$ are admissible. For any multipartitions $\bnu,\blambda \in \cP^{2n}$, we have 
    \begin{align}\label{eqn:multiplicity} [\mts(\blambda):\lts(\bnu)] = \prod_{\ell=1}^r [\Delta^{\b{c}_\ell,\b{\sigma}_\ell}_{m,N}(\ublambda_\ell):L^{\b{c}_\ell,\b{\sigma}_\ell}_{m,N}(\ubnu_\ell)],\end{align} where $\Delta^{\b{c}_\ell,\b{\sigma}_\ell}_{m,N}(\blambda_\ell)$ and $L^{\b{c}_\ell,\b{\sigma}_\ell}_{m,N}(\b{\nu}_\ell)$ are the Verma and simple modules in $\cO_m^{\b{c}_\ell,\b{\sigma}_\ell}(N)$ for $m,N \gg 0$. 
\end{cor}

Finally, take $\fk = \CC$ and recall the stratification $W_I$ of parameter space $\scr{X}_n$ defined in Subsection \ref{subsec:stratification}. The key observation is that $\Xi_{\bs}^{\bt}$ is an admissible type (i.e., each sequence $\b{c}_\ell$ is admissible) if and only if the point $(\bt,\bs) \in \scr{X}_n$ belongs to an admissible facet $W_I \subset \scr{X}_n$ for some $I \subset \cH$.

Thus, the results in this section allow us to prove Conjecture \ref{conj:categorical_conjecture} in the admissible case. Moreover, Corollary \ref{cor:admissible_multiplicities} yields explicit formulas for the stable multiplicities on each admissible facet.

\begin{cor}
    Conjecture \ref{conj:categorical_conjecture} holds for admissible facets $W_I$. In particular, multiplicities can be directly computed by the formula (\ref{eqn:multiplicity}) for points $(\bt,\bs) \in \scr{X}$ on admissible facets.
\end{cor}

\begin{rem}
    For \textit{interpolatable} admissible facets $W_I$, our proof of Conjecture \ref{conj:categorical_conjecture} gives an explicit interpretation of the stable multiplicities from Proposition \ref{prop:ultraproduct_multiplicity}. On the other hand, for non-interpolatable admissible facets $W_I$ (where our ultraproduct approach was ineffective), our new proof yields an explicit formula for multiplicities on $W_I$ as well as an interpretation in terms of stable representation theory.
\end{rem}

\begin{ex}
Assume $n = 2$. By Example \ref{ex:admissibility_two_blocks}, the results in this subsection apply to any set of parameters $\{t_1,t_2,s_1,s_2\}$ provided $t_1,t_2 \not \in \ZZ$. In particular, we can completely describe the structure of $\extots$ in terms of usual parabolic categories $\cO$ when $n  = 2$ and $t_1,t_2 \not \in \ZZ$. Moreover, we have a complete proof of Conjecture \ref{conj:categorical_conjecture} in the case $n = 2$.
\end{ex}

\section{Applications to Stable Modular Representation Theory}\label{sec:applications}

We use the validity of Conjecture \ref{conj:categorical_conjecture} to interpret stable multiplicities (Corollary \ref{cor:admissible_modular}) from the modular representation theory of the hyperalgebra of $\GL_m$ in terms of parabolic Kazhdan--Lusztig polynomials. This section uses a construction of Deligne categories through ultraproducts of representation categories $\GL_m$ in large positive characteristic. For the sake of brevity, we refer the reader to \cite{harman16} for a detailed treatment of this construction.

\subsection{Basic Notation}

Let $\FF := \ov{\FF}_p$ denote the algebraic closure of $\FF_p$ for some prime $p > 0$. Consider the algebraic group $\mathrm{GL}_m$ along with its Lie algebra $\fgl_m$.

We write $\mathrm{Dist}(\mathrm{GL}_m)$ for the \textit{hyperalgebra} of $\mathrm{GL}_m$. More precisely, let $\FF[\GL_m]$ denote the $\FF$-algebra of functions on $\GL_m$, and let $I \subset \FF[\GL_m]$ denote the ideal of functions vanishing at the identity. Then, we define the filtered Hopf algebra \[\mathrm{Dist}(\mathrm{GL}_m) := \bigcup_{k \geq 0} F^k\mathrm{Dist}(\mathrm{GL}_m), \quad F^k\mathrm{Dist}(\mathrm{GL}_m) := \bigcup_{0 \leq n \leq k} (\FF[\GL_m]/I^{n+1})^*.\]
For a composition $m_1 + \cdots + m_n = m$, consider the associated parabolic subgroup $P_{p,\b{m}} \subset \GL_m(\FF)$ and Levi subgroup $L_{p,\b{m}} \simeq G_{m_1}(\FF) \times \cdots \times G_{m_n}(\FF)$. Let $U_+$ denote the unipotent radical of $P$ and $U_-$ denote the unipotent radical of the opposite parabolic. Then,
\begin{align*}
\mathrm{Dist}(\GL_m) &\simeq \mathrm{Dist}(U_-) \ox \mathrm{Dist}(L_{p,\b{m}}) \ox \mathrm{Dist}(U_+).
\end{align*}

We recall the finite-dimensional representation theory of these algebras. Consider the $P$-dominant integral weights $\Lambda^{\fp} \subset \ZZ^m$ consisting of sequences $(\lambda_1,\ldots,\lambda_n)$, where each entry $\lambda_i = (\lambda_{i,1},\ldots,\lambda_{i,m_i}) \in \ZZ^{m_i}$ is a dominant integral weight for $\fgl_{m_i}(\FF)$. These weights parameterize finite-dimensional irreducible representations $X^{p,\b{m}}(\lambda)$ of the reductive algebraic group \[L_{p,\b{m}} := \GL_{m_1} \times \cdots \times \GL_{m_n}.\]
Let $T \subset \GL_m$ denote the torus of diagonal matrices and let $\Rep{\mathrm{Dist}(L_{p,\b{m}}),T}$ denote the category of finite-dimensional representations of $\mathrm{Dist}(L_{p,\b{m}})$, equipped with a rational action of $T$ that is compatible with the existing $\mathrm{Dist}(T)$-action. Then, we have equivalences of categories
\[
\Rep{\mathrm{Dist}(L_{p,\b{m}}),T} \simeq \Rep{L_{p,\b{m}}} \simeq \Rep{\GL_{m_1}(\FF)} \boxtimes \cdots \boxtimes \Rep{\GL_{m_n}(\FF)}.
\]
For any $\lambda \in \Lambda^{\fp}$, we define the \textit{parabolic hyperalgebra Verma module}
\[
\Delta^{p,\b{m}}(\lambda) := \mathrm{Dist}(\GL_m(\FF)) \ox_{\mathrm{Dist}(P)} X^{p,\b{m}}(\lambda).
\]
By the usual theory of highest weight modules, these Verma modules have unique simple quotients $L^{p,\b{m}}(\lambda)$ respectively. These Verma modules are parabolic analogs of hyperalgebra Verma modules introduced by Haboush \cite{haboush}.



 Note that $\mathrm{Dist}(G)$ is an object in $\Ind \Rep{\mathrm{Dist}(L_{p,\b{m}}),T}$ since all eigenvalues of $\mathrm{Dist}(T)$ are integers. It follows that $\Delta^{p,\b{m}}(\lambda)$ and $L^{p,\b{m}}(\lambda)$ are also objects in $\Ind \Rep{\mathrm{Dist}(L_{p,\b{m}}),T}$.

Now suppose $m_i \gg 0$. Given a multipartition $\blambda \in \cP^{2n}$, we consider the dominant weight \begin{align*}\theta^{p,\b{m}}(\blambda) &:= (\theta^{p,\b{m}}(\blambda)_1,\ldots,\theta^{p,\b{m}}(\blambda)_n), \\[5pt] \theta^{p,\b{m}}(\blambda)_i &:= (\lambda_{2i-1,1},\ldots,\lambda_{2i-1,\ell(\lambda_{2i-1})},0,\ldots,0,-\lambda_{2i,\ell(\lambda_{2i})},\ldots,-\lambda_{2i,1}).\end{align*} 

\subsection{Hyperalgebra Verma Modules and Complex Rank Category \texorpdfstring{$\cO$}{O}}

Fix a set partition \[\{1,\ldots,2n\} = \Gamma_1 \sqcup \cdots \sqcup \Gamma_r\] with $\Gamma_\ell = \{i_\ell(1),\ldots,i_\ell(a_\ell)\}$. Moreover, for each $\ell = 1,\ldots,r$, we define a sequence \[\b{c}_{\ell} := (c_{\ell,1},\ldots,c_{\ell,a_\ell}) \in \{0,1\}^{a_\ell}\] where $c_{\ell,k} = 0$ if $i_\ell(k)$ is odd and $c_{\ell,k} = 1$ otherwise. Finally, we pick $\b{\sigma}_\ell = (\sigma_{\ell,1},\ldots,\sigma_{\ell,a_\ell}) \in \ZZ^{a_\ell}$. 

Now consider an increasing sequence of prime numbers $(p_1,p_2,\ldots)$ indexed by $\NN$, and let $\FF_j := \ov{\FF}_{p_j}$ for each $j \in \NN$. For each $j \in \NN$, we pick positive integers $m_{j,1},\ldots,m_{j,n}$ and set $m_j := m_{j,1} + \cdots + m_{j,n}$. Let's assume that $p_j \gg m_j$ for each $j \gg 0$. Then, define
\[
\gamma_{j,2i-1} := -m_{j,1} - \cdots - m_{j,i-1} \in \ZZ/p\ZZ, \quad \gamma_{j,2i} := -m_{j,1} - \cdots - m_{j,i} \in \ZZ/p\ZZ
\]
for each $i \in \{1,\ldots,n\}$. Our choices of $m_{j,i}$ must satisfy the following constraints.
\begin{enumerate}
\item For each fixed $i = 1,\ldots,n$, the sequence $(m_{j,i})_{j \in \NN}$ is strictly increasing with $j$. 
\item For each $\ell = 1,\ldots,r$ and each $k = 1,\ldots,a_\ell$, we require that \[\gamma_{j,i_\ell(k)} - \gamma_{j,i_\ell(1)} = \sigma_{\ell,k} \in \ZZ/p_j\ZZ. \]
\item On the other hand, for any pair $\ell < \ell'$, we require that \[\gamma_{j,i_\ell(1)} - \gamma_{j,i_{\ell'}(1)} = (\ell'-\ell)j \in \ZZ/p_j\ZZ.\] 
\end{enumerate}

Let $\mathscr{F}$ be a non-principal ultrafilter on $\NN$. Any algebraically closed field with cardinality continuum is isomorphic to $\CC$ and thus furnishes us with a (non-canonical) isomorphism \[\Psi: \prod_{\substack{j \in \NN \\ \mathscr{F}}} \ov{\FF}_{p_j} \xto{\sim} \CC.\] Thus, writing $\b{m}_j := (m_{j,1},\ldots,m_{j,n})$, we define \[\widetilde{\cC} := \prod_{\substack{j \in \NN \\ \mathscr{F}}} \Ind \Rep{L_{p_j,\b{m}_j}}, \quad V := \prod_{\substack{j \in \NN \\ \mathscr{F}}} V_j \in \widetilde{\cC},\] where $V_j$ is the vector representation of $L_{p_j,\b{m}_j}$. In turn, write $\cC \subset \widetilde{\cC}$ for the subcategory tensor generated by $V$. Finally, for each $i = 1,\ldots,n$, let $t_i := \Psi\left(\prod_{\substack{j \in \NN \\ \mathscr{F}}} m_{j,i}\right) \in \CC$ 

Then, by the ultraproduct construction of Deligne categories, there exists a symmetric monoidal equivalence $\cC \simeq \cL_{\b{t}}$, where $\b{t} := (t_1,\ldots,t_n)$.


Note that $U_{\bt}$ has a natural PBW filtration with filtered pieces $F^kU_{\bt}$, defined as the image of $\bigoplus_{0 \leq u \leq k} \fgl_{\bt}^{\ox u}$ under the natural multiplication map. 

\begin{prop}\label{prop:ultraproduct_universal_enveloping}
For each $k \geq 0$, we have isomorphisms \[F^kU_{\bt}\simeq \ultraprod F^kU(\fgl_{m_j}(\FF_j)) \simeq \ultraprod F^k\mathrm{Dist}(\GL_{m_j})\] of objects in $\cC$. In particular, we have algebra isomorphisms \[U_{\bt} \simeq \bigcup_{k \geq 0}\ultraprod F^k\mathrm{Dist}(\GL_{m_j}) \simeq \bigcup_{k \geq 0}\ultraprod F^kU(\fgl_{m_j}(\FF_j)).\]
\end{prop}

\begin{proof}
    The first isomorphism is straightforward, using $\fgl_t = \ultraprod \fgl_{m_j}(\FF_j).$ The isomorphism \[\ultraprod F^kU(\fgl_{m_j}(\FF_j)) \simeq \ultraprod F^k\mathrm{Dist}(\GL_{m_j})\] follows since natural maps $F^kU(\fgl_{m_j}(\FF_j)) \hookrightarrow F^k\mathrm{Dist}(\GL_{m_j})$ are isomorphisms for $k \leq p_j$. 
\end{proof}

In turn, for any $\blambda \in \cP^{2n}$, we define the object \begin{align*}\Delta^{\bt}(\blambda) &:= \bigcup_{k \geq 0}\prod_{\substack{j \in \NN \\ \mathscr{F}}} F^k\Delta^{p_j,\b{m}_j}(\theta^{p_j,\b{m}_j}(\blambda)) \in \widetilde{\cC}.\end{align*}
Strictly speaking, in order for $\theta^{p_j,\b{m}_j}(\blambda)$ to be well-defined, we require $\ell(\blambda_{2i-1}) + \ell(\blambda_{2i}) \leq m_{j,i}$ for each $j \in \NN$ and $i \in \{1,\ldots,n\}$. Our assumption that $\{m_{i,j}\}_{j \in \NN}$ is strictly increasing means that these bounds hold for $\mathscr{F}$-many indices $j \in \NN$.

By \L os's theorem and Proposition \ref{prop:ultraproduct_universal_enveloping}, each $\Delta^{\bt}(\blambda)$ has the structure of a $U_{\bt}$-module. 

\begin{prop}
    For $\blambda \in \cP^{2n}$, we have a $U_{\bt}$-module isomorphism 
    $M^{\bt}_{\b{0}}(\blambda) \simeq \Delta^{\bt}(\blambda)$.
\end{prop}

\begin{proof}
    Completely analogous to Proposition \ref{prop:verma_ultraproduct} in Appendix \ref{sec:ultraproducts}.
\end{proof}

\begin{lem}\label{lem:technical_lemma}
    Let $\cC$ be an $\fk$-abelian category where $[M:L] < \infty$ for any $M \in \cC$ and any simple object $L$. Let $P \in \cC$ be an object satisfying $\dim \Hom(P,L) = 1$ for some simple $L \in \cC$ and such that the functor $\Hom(P,-)$ preserves short exact sequences of subquotients of some object $M \in \cC$. Then, we have the equality $[M:L] = \dim \Hom(P,M)$. 
\end{lem}

\begin{proof}
A straightforward induction on $[M:L]$. 
\end{proof}

\begin{cor}\label{cor:modular_multiplicities}
    For $j \gg 0$ and $\blambda,\bmu \in \cP^{2n}$, we have equalities 
    \begin{align*}
        [M^{\bt}_{\b{0}}(\blambda):L^{\bt}_{\b{0}}(\bmu)] &= [\Delta^{p_j,\b{m}_j}(\theta^{p_j,\b{m}_j}(\blambda)):L^{p_j,\b{m}_j}(\theta^{p_j,\b{m}_j}(\bmu))].
    \end{align*}
\end{cor}

\begin{proof}
    Let us write $Z \simeq \GG_m^n$ for the center of $L_{p_j,\b{m}_j}$. Let $\iota^!$ denote the functor taking any $L_{p_j,\b{m}_j}$-module $M$ to its largest quotient whose $Z$-weights are less than \[\omega(\lambda) = (|\lambda_1| - |\lambda_2|,\ldots,|\lambda_{2n-1}| - |\lambda_{2n}|)\] in the dominance order for weights of $Z$.  We may construct objects 
    \[
        Q^j_{\blambda}(\bmu) := \mathrm{Dist}(\GL_{m_j}(\FF_j)) \ox_{\mathrm{Dist}(P_{p_j,\b{m}_j})} \iota^!(\mathrm{Dist}(P_{p_j,\b{m}_j}) \ox_{\mathrm{Dist}(L_{p_j,\b{m}_j})} X^{p_j,\b{m}_j}(\theta^{p_j,\b{m}_j}(\bmu))).
    \]
    Thus, we see that $Q^j_{\blambda}(\bmu)$ represents the functor \[\cC \to \mathrm{Vect}_{\FF_j}, \quad M \mapsto \Hom_{L_{p_j,\b{m}_j}}(X^{p_j,\b{m}_j}(\theta^{p_j,\b{m}_j}(\bmu)),M),\] where $\cC$ is the full subcategory of objects in $\Ind \mathrm{Rep}(L_{p_j,\b{m}_j})$ where $Z$ acts semisimply with finite-dimensional $Z$-weight spaces of weight bounded above by $\omega(\blambda)$ in the dominance order. Let $P_{\blambda}^j(\bmu)$ be an indecomposable summand of $Q_{\blambda}^j(\bmu)$ with a nonzero surjection to $\Delta^{p_j,\b{m}_j}(\theta^{p_j,\b{m}_j}(\bmu))$. 
    
    Note that $P_{\blambda}^j(\bmu)$ admits a natural $\ZZ$-grading $P_{\blambda}^j(\bmu)_d$ by the heights of $Z$-weights. Then, $\bigoplus_{d \in \ZZ} \prod_{\scr{F}} P_{\blambda}^j(\bmu)_d$ is a projective cover of $L^{\bt}(\bmu)$ in the full subcategory of $U_{\bt}$-modules semisimple over $\fz_{\bt}$, where the $\fz_{\bt}$-weights are bounded above by $\omega_{\bt}(\blambda)$. In particular, \[\dim \Hom_{U(\fg_{\bt})}\left(\bigoplus_{d \in \ZZ}\ultraprod P_{\blambda}^j(\bmu)_d, L^{\bt}(\bmu)\right) = 1.\] On the other hand, we see that $L^{\bt}(\bmu)$ can be realized as the $\ZZ$-graded ultraproduct of the modules $L^{p_j,\b{m}_j}(\theta^{p_j,\b{m}_j}(\bmu))$ (where the grading comes from weights of $Z$). Thus, an easy generalization of Lemma \ref{lem:ultraproduct_dimension} from Appendix \ref{sec:ultraproducts} to our setting shows that for $j \gg 0$, we have
    \begin{align*}
        \dim \Hom_{\mathrm{Dist}(\GL_{m_j})}(P_{\blambda}^j(\bmu), L^{p_j,\b{m}_j}(\theta^{p_j,\b{m}_j}(\bmu))) = \dim \Hom_{U(\fg_{\bt})}\left(\bigoplus_{d \in \ZZ}\ultraprod P_{\blambda}^j(\bmu)_d, L^{\bt}(\bmu)\right) = 1.
    \end{align*}
    Now, we wish to apply Lemma \ref{lem:technical_lemma} to show that \begin{align}\label{eqn:dim_equality} \dim \Hom(P_{\blambda}^j(\bmu),\Delta^{p_j,\b{m}_j}(\theta^{p_j,\b{m}_j}(\blambda))) = [\Delta^{p_j,\b{m}_j}(\theta^{p_j,\b{m}_j}(\blambda)):L^{p_j,\b{m}_j}(\theta^{p_j,\b{m}_j}(\bmu))]\end{align} for $j \gg 0$. The only nontrivial hypothesis that needs to be checked is that the functor \[\Hom_{L_{p_j,\b{m}_j}}(X^{p_j,\b{m}_j}(\theta^{p_j,\b{m}_j}(\bmu)),-)\] preserves short exact sequences of $L_{p_j,\b{m}_j}$-subquotients of $\Delta^{p_j,\b{m}_j}(\theta^{p_j,\b{m}_j}(\blambda))$. Letting $L_j^0$ denote the (semisimple) derived subgroup of $L_{p_j,\b{m}_j}$ so that $L_{p_j,\b{m}_j} = L_j^0Z$, it suffices to show that the $\omega(\mu)$-weight space in $\Delta^{p_j,\b{m}_j}(\theta^{p_j,\b{m}_j}(\blambda))$ is semisimple as a $L_j^0$-module when $j \gg 0$. This fact simply follows from the observation that the (finitely many) composition factors of this weight space have $L_j^0$-highest weights belonging to the fundamental alcove for $L_j^0$ when $p_j \gg 0$.
    
    Having now established (\ref{eqn:dim_equality}), the rest of the argument proceeds similarly to the proof of Proposition \ref{prop:ultraproduct_multiplicity} presented in Appendix \ref{sec:ultraproducts}, i.e., through another application of Lemma \ref{lem:ultraproduct_dimension}.
\end{proof}

Now define \[\gamma_{2i-1} := -t_{1} - \cdots - t_{i-1} \in \CC, \quad \gamma_{2i} := -t_1 - \cdots - t_i \in \CC.\] By \L os's theorem, our assumptions on the choices of $m_{j,i}$ imply that 
\begin{enumerate}
\item For each $\ell = 1,\ldots,r$ and each $k = 1,\ldots,a_\ell$, we have \[\gamma_{i_\ell(k)} - \gamma_{i_\ell(1)} = \sigma_{\ell,k} \in \ZZ. \]
\item On the other hand, for any pair $\ell \neq \ell'$, the difference $\gamma_{i_\ell(1)} - \gamma_{i_{\ell'}(1)}$ is non-integral.
\end{enumerate}

\begin{cor}\label{cor:admissible_modular}
    Assume that each sequence $\b{c}_\ell$ is admissible. For $j \gg 0$, the numbers 
    \begin{align*}
        [\Delta^{p_j,\b{m}_j}(\theta^{p_j,\b{m}_j}(\blambda)):L^{p_j,\b{m}_j}(\theta^{p_j,\b{m}_j}(\bmu))]
    \end{align*}
    stabilize and can be computed as the product of parabolic Kazhdan--Lusztig polynomials evaluated at $1$, i.e., via the formula given in Corollary \ref{cor:admissible_multiplicities}.
\end{cor}

Recall the facets $W_I$ in parameter space $\scr{X}_n$ from Subsection \ref{subsec:stratification}. Equivalently, the set partition $\Gamma_1 \sqcup \cdots \sqcup \Gamma_r$ and the collection of sequences $\b{\sigma}_\ell$ correspond to some subset $I \subset \cH$ (in the notation of Section \ref{sec:main_conjecture}) such that the facet $W_I$ is admissible if and only if the sequences $\b{c}_\ell$ are all admissible. Proposition \ref{prop:ultraproduct_multiplicity} along with the validity of Conjecture \ref{conj:categorical_conjecture} in the admissible case give us an explicit interpretation of the stable multiplicities \[[\Delta^{p_j,\b{m}_j}(\theta^{p_j,\b{m}_j}(\blambda)):L^{p_j,\b{m}_j}(\theta^{p_j,\b{m}_j}(\bmu))]\] when $W_I$ is interpolatable. However, when $W_I$ is not interpolatable, our proof of Conjecture \ref{conj:categorical_conjecture} gives an interpretation of the aforementioned multiplicities that is inaccessible via model theory.

\begin{rem}\label{rem:algebraic}
    Fix some subset $\{i_1,\ldots,i_w\} \subset \{1,\ldots,n\}$. For each $u = 1,\ldots,w$, let $q_u(x) \in \ZZ[x]$ be polynomials satisfying the following condition: there exists a sequence of primes $p_{1} < p_{2} < \cdots$ such that each $q_u(x)$ has a root modulo $p_j$ for each $j \in \NN$. If we pick each $m_{j,i_u}$ to be a root of $q_u(x)$ modulo $p_{j}$, then each $t_{i_u}$ is an algebraic number, namely a complex root of $q_u(x)$. Thus, in the notation of Section \ref{sec:main_conjecture}, Corollary \ref{cor:admissible_modular} also gives us a formula for multiplicities at points $(\bt,\bs) \in \scr{X}_n$ where $\bs = (0,\ldots,0)$ and where the coordinates $t_i$ may be algebraic numbers.
\end{rem}

\newpage

\section*{Index of Notation}

\begin{itemize} \setlength{\itemsep}{0.5em}
    \item $\cD_t$,  $\cL_t$, $\fgl_{\bt}, \fl_{\bt}, \f{u}^+_{\bt}, \f{u}^-_{\bt}, U_{\bt}$, $X_{\bs}^{\bt}(\blambda), \mts(\blambda), \lts(\blambda),\widehat{\cO}^{\bt}_{\bs}, \ots$ \dotfill Subsection \ref{subsec:categoryO}
    \item $\scr{X}_n$, interpolatable points \dotfill Subsection \ref{subsec:interpolatable}
    \item $H_k(i,j), W_I$, interpolatable/admissible strata \dotfill Subsection \ref{subsec:stratification}
    \item $\cC_{\leq \lambda}, \cC_\lambda$ \dotfill Subsection \ref{subsec:standardly_stratified}
    \item $\leq$ (inverse dominance order on multipartitions), $\cP^a(m), \cP^a(m_1,m_2)$ \dotfill Subsection \ref{subsec:multipartition_poset}
    \item $\cF_{s}, \cF_s^{\vee}$ \dotfill Subsection \ref{subsec:fock}
    \item $\Xi$, MFTPC \dotfill Subsection \ref{subsec:ftp}
    \item $\Lambda^\bullet, \Lambda^\circ, \Lambda^\bullet(j,k),\Lambda^\circ(j,k), \Lambda(j,k)$ \dotfill Subsection \ref{subsec:restricted_multi_fock}
    \item $\cC^\bullet(j,k), \cC^\circ(j,k),\cC(j,k)$ \dotfill Subsection \ref{subsec:restricted_multi_fock}
    \item restricted MFTPC \dotfill Subsection \ref{subsec:restricted_multi_fock}
    \item $p^\bullet, p^\circ, q^\bullet, q^\circ, p, q$ \dotfill Start of Section \ref{sec:uniqueness}
    \item $\ul{\Lambda}_h, \ov{\Lambda}_h, \ul{\Xi}_h, \ov{\Xi}_h$, $\ublambda_\mu(h), \ublambda^-(h)$, $\prescript{}{\mu}{\ublambda}(h)$, $\prescript{-}{}{\ublambda}(h)$ \dotfill Subsection \ref{subsec:categorical_splitting}
    \item weak stratification, $\Delta^\bullet(\ublambda), T^\bullet(\ublambda)$, $E^{(\ell)}(\mu), F^{(\ell)}(\mu)$ \dotfill Subsection \ref{subsec:structure}
    \item $\UU, \TT^\bullet, \Heis_k, k(\Xi), \Heis_{k,m}$, $\Heis_k(e|f), A(\Xi), H_{\alpha,\beta}(\Xi)$ \dotfill Subsection \ref{subsec:heisenberg}
    \item $\ul{\Heis}$ \dotfill Subsubsection \ref{subsubsec:full_faithfulness}
    \item $\cC(\mathfrak{F})^{j,k}, \cC(\mathfrak{F})^{j,\infty}, \cC(\mathfrak{F})$ \dotfill Subsubsection \ref{subsubsec:restricted_from_full}
    \item $\cO^{\b{c},\b{\sigma}}(N)$, $\Delta_{N}^{\b{c},\b{\sigma}}(\vartheta), L_{N}^{\b{c},\b{\sigma}}(\vartheta)$ \dotfill Subsection \ref{subsec:construction_notation}
    \item virtual multipartition, $\weight(\blambda)$ \dotfill Subsection \ref{subsec:virtual}
    \item $\cO_m^{\b{c},\b{\sigma}}$ \dotfill Subsection \ref{subsec:constructing_restricted_fock}
    \item $\cO_m^\Xi$, $\cO_\infty^\Xi$ \dotfill Subsection \ref{subsec:external_tensor_product_restricted}
    \item $\DD, \check{M}^{\bt}_{\bs}(\blambda)$ \dotfill Subsection \ref{subsec:duality}
    \item $\trianglelefteq_{\bt}, \omega_{\bs}(\blambda), \trianglelefteq_{\mathrm{dom}}, B_{\bs}^{\bt}(\blambda)_{\bmu}, D_{\bs}^{\bt}(\blambda)_{\bmu}, T_i, \psi_{\bt}(\blambda)_i$ \dotfill Subsection \ref{subsec:jantzen_determinant}
    \item $\SS, \blambda(\zeta), \mathrm{sgn}(\blambda(\zeta)), \tilde{r}(\blambda(\zeta),\bmu)$ \dotfill Subsection \ref{subsec:jantzen_determinant}
    \item $\gamma_{2i-1},\gamma_{2i}, \Gamma_\ell, \b{\sigma}_\ell,\b{c}_\ell, D_{\bs}^{t,R}(\blambda)_{\bmu}(y), \mts(\blambda)_j$, equivariant character, $\XX_\ell^+, \blambda^{\ell,\chi}$ \dotfill Subsection \ref{subsec:jantzen_sum}
    \item linkage order $\preceq$ for $\ots$ \dotfill Subsection \ref{subsec:interpolated_linkage_order}
    \item $\extots$ \dotfill Subsubsection \ref{subsubsec:highest_weight}
    \item $\Xi_{\bs}^{\bt}$ \dotfill Subsection \ref{subsec:categorical_action}
    \item $\mathrm{Dist}(\GL_m), \Delta^{p,\b{m}}(\lambda), L^{p,\b{m}}(\lambda),\theta^{p,\b{m}}(\blambda)$ \dotfill Section \ref{sec:applications}
\end{itemize}

\newpage

\appendix

\section{Ultraproduct Constructions}\label{sec:ultraproducts}

The goal of this appendix is to provide a proof of Proposition \ref{prop:ultraproduct_multiplicity} using ultraproducts. We keep the notation from Subsection \ref{subsec:interpolatable}. We refer the reader to \cite{Schoutens2010Ultraproducts} for background on ultraproducts and \cite[Section 2]{hu2024introductiondelignecategories} for the details related to Deligne categories.

Fix a non-principal ultrafilter $\scr{F}$ on the natural numbers $\NN$. Consider the ultraproduct category
\[
\cC' := \ultraprod {\mathrm{Rep}}_{\ov{\QQ}}(\GL_j),
\]
where we use the convention that, when $j < 0$, the symmetric tensor category ${\mathrm{Rep}}_{\ov{\QQ}}(\GL_j)$ is the category of super-representations of the supergroup $\GL(0|-j)$ such that the negative identity matrix acts by parity. 
By \L os's Theorem, the category $\cC'$ has the natural structure of an abelian symmetric monoidal category over the field $\prod_{\scr{F}} \ov{\QQ}$. By a classical result of Steinitz, we can fix an isomorphism 
\[
\phi: \ultraprod \ov{\QQ} \simeq \CC,
\] 
which is determined up to the action of $\mathrm{Aut}(\CC)$. Write 
$t := \phi\left(\ultraprod j\right) \in \CC,$
which is transcendental over $\ov{\QQ}$.
Let $V$ denote the ultraproduct of the vector representations $V_j \in {\mathrm{Rep}}_{\ov{\QQ}}(\GL_j)$ and finally, let $\cC$ denote the tensor subcategory of $\cC'$ generated by $V$. 

\begin{thm}[Deligne]
    There exists a symmetric monoidal equivalence $\cD_t \simeq \cC$ given by $V_t \mapsto V$. In particular, the category $\cC$ does not depend on the choice of ultrafilter $\scr{F}$.
\end{thm}

Let $\mathscr{F}$ be a non-principal ultrafilter on $\NN$. Define \[\widetilde{\cC}_{\scr{F}} := \prod_{\substack{j \in \NN \\ \mathscr{F}}} \Ind \mathrm{Rep}_{\ov{\QQ}}({\fl(j)}), \quad V := \prod_{\substack{j \in \NN \\ \mathscr{F}}} V(j) \in \tilde{\cC}_{\scr{F}}\] and let $\cC \subset \widetilde{\cC}_{\scr{F}}$ denote the Serre subcategory of $\widetilde{\cC}_{\scr{F}}$ tensor generated by $V$. The ultraproduct construction of Deligne categories yields exists a symmetric monoidal equivalence \[\cC \simeq \cL_{\b{t}} = \cD_{t_1} \boxtimes \cdots \boxtimes \cD_{t_n}.\] In particular, $\cC$ does not depend on the choice of $\scr{F}$.

The universal enveloping algebra $U_{\bt}$ admits a PBW filtration with pieces $F^kU_{\bt}$ given by the image of $\bigoplus_{0 \leq u \leq k} \fgl_{\bt}^{\ox u}$ under the multiplication map. The following lemma is straightforward.

\begin{lem}\label{lem:ultraproduct_universal_enveloping}
For each $k \geq 0$, we have an isomorphism \[F^kU_{\bt} \simeq \ultraprod F^kU(\fgl(j))\] of objects in $\cC$. In particular, we have an algebra isomorphism \[U_{\bt} \simeq \bigcup_{k \geq 0}\ultraprod F^kU(\fgl(j)).\]
\end{lem}

For any $\blambda \in \cP^{2n}$ and $j \gg 0$, the Verma module $\Delta^j_{\bt,\bs}(\blambda)$ admits a $\ZZ$-grading by the heights of $\fz(j)$-weights, regarded as roots for $\fgl_n$. For each $d \in \ZZ$, we let $\Delta^j_{\bt,\bs}(\blambda)_d$ denote the corresponding homogeneous component, which is a finite-dimensional $\ov{\QQ}$-supervector space. In turn, we define
\[
    \Delta^{\bt}_{\bs}(\blambda) := \bigoplus_{d \in \ZZ} \ultraprod \Delta^j(\blambda)_d,
\]
which admits the natural structure of a $U_{\bt}$-module -- to see this, simply note that the filtration on $\Delta^j_{\bt,\bs}(\blambda)$ induced by the grading is the same as the PBW filtration (up to some degree shift), so $\Delta^{\bt}_{\bs}(\blambda)$ has the natural structure of a filtered $U_{\bt}$-module with respect to the PBW filtrations. In other words, there is a canonical identification 
\[
    \Delta^{\bt}_{\bs}(\blambda) = \bigcup_{k \geq 0} \ultraprod F^k\Delta^j_{\bt,\bs}(\blambda),
\]
where $F^k\Delta^j(\blambda)$ is the PBW filtration. 

\begin{prop}\label{prop:verma_ultraproduct}
    For $\blambda \in \cP^{2n}$, we have $U_{\bt}$-module isomorphisms $\mts(\blambda) \simeq \Delta_{\bs}^{\bt}(\blambda).$
\end{prop}

\begin{proof}
    For the sake of brevity, we sketch the proof, leaving the details to the reader. Let $X^j(\blambda)$ denote the simple supermodule in $\mathrm{Rep}_{\ov{\QQ}}(\fl(j))$ with highest weight $\theta^j(\blambda)$. Observe \[\ultraprod X^j(\blambda) \simeq X^{\bt}_{\bs}(\blambda) \in \cL_{\bt}.\] Indeed, the simple object $X^{\bt}_{\bs}(\blambda)$ can be inductively constructed a direct summand of some mixed tensor powers of $V = \prod_{\mathscr{F}} V(j)$ corresponding to a suitable sequence of eigenvalues of the Casimir operator. For sufficiently large, hence $\mathscr{F}$-many, indices $j$, there is an identical characterization of $X^j(\blambda)$ as a direct summand of some mixed tensor power of the vector representation $V(j)$, so the observation follows from \L os's theorem. 

    On the other hand, thanks to \L os's theorem, it is not hard to see that $\Delta_{\bs}^{\bt}(\blambda)$ is generated as a $U_{\bt}$-module in filtration degree zero, i.e., by \[F^0\Delta_{\bs}^{\bt}(\blambda) = \ultraprod F^0\Delta^j_{\bt,\bs}(\blambda) = X_{\bs}^{\bt}(\blambda).\] \L os's theorem also gives us that the nilpotent radical $\mathfrak{u}^+_{\bt}$ of $\mathfrak{p}_{\bt}$ annihilates $F^0\Delta_{\bs}^{\bt}(\blambda)$. Hence, we obtain a surjective $U_{\bt}$-module homomorphism $\mts(\blambda) \to \Delta_{\bs}^{\bt}(\blambda).$ On the other hand, by the same proof as Lemma \ref{lem:ultraproduct_universal_enveloping}, note that \[U(\mathfrak{u}_{\bt}^-) \simeq \bigcup_{k \geq 0}\ultraprod F^kU(\mathfrak{u}^-_{p_j,\b{m}_j}).\]
    By \L os's theorem, $Y_{\bs}^{\bt}(\blambda)$ is isomorphic to the free $U(\mathfrak{u}_{\bt}^-)$-module $U(\mathfrak{u}_{\bt}^-) \ox X_{\bs}^{\bt}(\blambda)$ with \[F^kY_{\bs}^{\bt}(\blambda) = F^kU(\mathfrak{u}_{\bt}^-) \ox X_{\bs}^{\bt}(\blambda).\] The module $\mts(\blambda)$ also satisfies $F^k\mts(\blambda) = F^kU(\mathfrak{u}_{\bt}^-) \ox X_{\bs}^{\bt}(\blambda)$. The surjection $\mts(\blambda) \to Y^{\bt}_{\bs}(\blambda)$ respects filtrations and thus induces a surjective endomorphism of $F^kU(\mathfrak{u}_{\bt}^-) \ox X_{\bs}^{\bt}(\blambda)$ for each $k \geq 0$, which must be an isomorphism. We conclude that $\mts(\blambda) \simeq Y_{\bs}^{\bt}(\blambda)$.  
\end{proof}

\begin{lem}\label{lem:ultraproduct_dimension}
    For each $j \in \NN$, let $M_j,N_j$ be $U(\fgl(j))$-modules where $\fz(\fl(j))$ acts semisimply. In particular, $M_j$ and $N_j$ inherits $\ZZ$-gradings by the heights of weights of $\fz(\fl(j))$. Let \[M := \bigoplus_{d \in \ZZ} \ultraprod M_{j,d}, \quad N := \bigoplus_{d \in \ZZ} \ultraprod N_{j,d},\] and assume that the isomorphism types of these $\fg_{\bt}$-modules do not depend on the choice of $\scr{F}$. If $M$ is finitely generated as a $\fg_{\bt}$-module, then for $j \gg 0$, we have
    \begin{align}\label{eqn:ultraproduct_dimension}
        \dim \Hom_{U(\fgl(j))}\left(M_j,N_j\right)  = \dim \Hom_{U_{\bt}}(M,N).
    \end{align}
\end{lem}

\begin{proof}
    Suppose $M$ is generated in degrees $d_1,\ldots,d_r$. In particular, we can identify $\Hom_{U_{\bt}}(M,N)$ as the subspace $W \subset \bigoplus_{i = 1}^r \Hom_{\fl_{\bt}}(M_{d_i},N_{d_i})$ of functions commuting with the action of $\fg_{\bt}$. 
    
    On the other hand, thanks to \L os's theorem, we see that $M_j$ is also generated in these degrees for $\scr{F}$-many indices $j$. In particular, for $\scr{F}$-many $j \in \NN$, we identify $\Hom_{U(\fgl(j))}(M_j,N_j)$ with the subspace $W_j$ of \[\bigoplus_{i = 1}^r \Hom_{\fl_{\bt}}(M_{j,d_i},N_{j,d_i})\] consisting of functions commuting with the action of $\fgl(j)$. By the ultraproduct construction of Deligne's category, we have \[\Hom_{\fl_{\bt}}(M_d, N_d) = \ultraprod \Hom_{\fl(j)}(M_{j,d}, N_{j,d})\] for any $d \in \ZZ$.  \L os's theorem thus implies that \[\ultraprod \Hom_{U(\fgl(j))}(M_j,N_j) = \Hom_{U_{\bt}}(M,N).\] Since we can specify the dimension of a vector space using first-order sentences in the two-sorted language of vector spaces over a field, \L os's theorem implies that the equality (\ref{eqn:ultraproduct_dimension}) holds for $\scr{F}$-many indices $j$. The independence of $M$ from the choice of $\scr{F}$ means (\ref{eqn:ultraproduct_dimension}) holds for $j \gg 0$.
\end{proof}

We can now prove Proposition \ref{prop:ultraproduct_multiplicity}. We recall the statement for the reader's convenience. 

\begin{prop}\label{prop:ultraproduct_multiplicity_appendix}
    Continue to use the notation from the previous paragraph. For any $\blambda,\bmu \in \cP^{2n}$, the multiplicities $[\Delta^j(\blambda):L^j(\bmu)]$ stabilize as $j \to \infty$. In fact, we have 
    \[
        [M^{\bt}_{\bs}(\blambda): L^{\bt}_{\bs}(\bmu)] = \lim_{j \to \infty} [\Delta_{\bt,\bs}^j(\blambda):L_{\bt,\bs}^j(\bmu)].
    \]
    In particular, the right-hand side is independent of the choice of $(\mathbf{m}_j,\b{\varsigma}_j)$.
\end{prop}

\begin{proof}
    As before, fix a non-principal ultrafilter $\scr{F}$ on $\NN$. For each $j \in \NN$, let $\cI_j(\blambda)$ be the poset ideal in $\ov{\QQ}^n$ (equipped with the dominance order for weights of $\fgl_n$) generated by \[\omega(\blambda, \b{\varsigma}_j) := (|\lambda_1| - |\lambda_2| + \varsigma_{j,1}, |\lambda_3| - |\lambda_4| + \varsigma_{j,2}, \ldots, |\lambda_{2n-1}| - |\lambda_{2n}| + \varsigma_{j,2n}) \in \ov{\QQ}^n.\] Write $\fz(j)$ for the center of $\fl(j)$. Let $U(\fp(j))\modcat_{\blambda}$ denote the full subcategory of $U(\fp(j))\modcat$ consisting of objects whose $\fz(j)$-weights belong to $\cI_j(\blambda)$, and let $\iota_{j,\blambda}: U(\fp(j))\modcat_{\blambda} \to U(\fp(j))\modcat$ denote the inclusion functor with left adjoint $\iota_{j,\blambda}^!$. Then, define the $\fgl(j)$-supermodule \[Q_{j,\blambda}(\bmu) := U(\fgl(j)) \ox_{U(\fp(j))} \iota_{j,\blambda}^!(U(\fp(j)) \ox_{U(\fl(j))} X^j(\bmu)).\] 

    In the proof of Proposition \ref{prop:verma_quotient_of_projective}, we constructed a projective object $Q_{\bs}^{\bt}(\bmu,\blambda)$ in the Serre subcategory of consisting of objects in $\widehat{\cO}^{\bt}$ whose $\fz_{\bt}$-weights are bounded above by $\omega_{\bs}(\blambda)$, along with a surjection $Q_{\bs}^{\bt}(\bmu,\blambda) \twoheadrightarrow \mts(\bmu)$. The object $Q_{\bs}^{\bt}(\bmu,\blambda)$ admits a grading $Q_{\bs}^{\bt}(\bmu,\blambda)_d$ by heights of $\fz_{\bt}$-weights. Similarly, the module $Q_{j,\blambda}(\bmu)$ also admits a grading $Q_{j,\blambda}(\bmu)_d$ from the heights of $\fz(j)$-weights. From their explicit descriptions, we have that 
    \[
        Q_{\bs}^{\bt}(\bmu,\blambda)_d = \ultraprod F^kQ_{j,\blambda}(\mu)_d \implies Q_{\bs}^{\bt}(\bmu,\blambda) = \bigoplus_{d \in \ZZ} \ultraprod Q_{j,\blambda}(\mu)_d.
    \]

    In turn, let $P_{j,\blambda}(\bmu)$ denote the unique indecomposable summand of $Q_{j,\blambda}(\bmu)$ that admits a nonzero map to $\Delta^j(\bmu)$. Thus, $P_{\blambda,j}(\bmu)$ is a projective cover of $L^j(\bmu)$ and \[\dim_{\ov{\QQ}} \Hom_{\fgl(j)}(P_{\blambda,j}(\bmu), \Delta^j(\blambda)) = [\Delta^j(\blambda):L^j(\bmu)].\] 
    Similarly, we take $P_{\bs}^{\bt}(\bmu,\blambda)$ to be an indecomposable summand of $Q_{\bs}^{\bt}(\bmu,\lambda)$ admitting a surjective map to $\Delta^j(\bmu)$, so that \[\dim_{\CC}\Hom_{\fg_t}(P_{\bs}^{\bt}(\bmu,\blambda),\mts(\blambda)) = [\mts(\blambda):\lts(\bmu)].\]  
    By \L os's theorem, we must have 
    \[
        P_{\bs}^{\bt}(\bmu,\blambda)_d = \ultraprod P_{j,\blambda}(\bmu)_d.
    \]
    Applying Lemma \ref{lem:ultraproduct_dimension} with $M_j = P_{j,\blambda}(\bmu)$ and $N_j = \Delta^j(\blambda)$, we deduce that 
    \[
        \dim_{\CC} \Hom_{U_{\bt}}(P_{\bs}^{\bt}(\bmu,\blambda), \mts(\blambda)) = \dim_{\ov{\QQ}} \Hom_{\fgl(j)}(P_{j,\blambda}(\bmu), \Delta^j(\blambda))
    \] for $j \gg 0$, whence the proposition follows.
\end{proof}

\section{Interpolated Central Characters}\label{sec:central_characters}

We now compute \textit{central characters} of Verma modules in the complex rank setting. Our notion for the center of $U_{\bt}$ is the invariant subalgebra $Z(U_{\bt}) := \Hom(\unit, U_{\bt})$. For each $k \geq 1$, we define the central element $C_k \in Z(U_{\bt})$ as follows. For each positive integer $i$, we write
\[c_i := \id_{\fgl_{\bt}^{\ox (i-2)}} \ox \id_V \ox \coev_{V^*} \ox \id_{V^*}: \fgl_{\bt}^{\ox(i-1)} \to \fgl_{\bt}^{\ox i},\] 
where we use the identification $\fgl_{\bt} = V \ox V^*$. In turn, we define $C_k$ as the composition
\[
C_k := \unit \xto{c_{k-1} \circ c_{k-2} \circ \cdots \circ c_1} \fgl_{\bt}^{\ox k} \xto{m_k} U_{\bt},
\]
where $m_k$ is the natural morphism $\fgl_{\bt}^{\ox k}  \to U_{\bt}$. 
For any $\blambda \in \cP^{2n}$, the central element $C_k$ acts on $\mts(\blambda)$ by a scalar, which we denote $\chi_{\bs,k}(\blambda,\bt) \in \fk$.



The following description of central characters in the case $n= 1$ follows by a standard interpolation argument. Refer to \cite[pg. 20]{Etingof_complex_rank_2} for an explicit description of these characters. 

\begin{lem}\label{lem:interpolate_casimir}
    Assume $n = 1$ and fix a bipartition $\blambda := (\lambda^\sharp,\lambda^\flat)$, where $\ell(\lambda^\sharp) = \ell^\sharp$ and $\ell(\lambda^\flat) = \ell^\flat$. Moreover, fix a large positive integer $m$ and define \[\theta_m(\blambda) := (\lambda^\sharp_1,\lambda^\sharp_2,\ldots,\lambda^\sharp_{\ell^\sharp}, 0,\ldots,0, -\lambda^\flat_{\ell^\flat},-\lambda^\flat_{\ell^\flat - 1},\ldots, -\lambda^\flat_1) \in \ZZ^{m}.\] 
    Then, there exists a polynomial $p_{k,\blambda}(x)$ for which
    \begin{itemize}
    \item[(a)] for any $t \in \fk$, the scalar $p_{k,\blambda}(t) \in \fk$ describes the action of $C_k \in U_t$ on $X(\blambda) \in \cD_t$, and
    \item[(b)] for sufficiently large integers $m$, the scalar $p_{k,\blambda}(m)$ describes the action of $\ul{C}_k^{(m)} \in U(\ul{\fgl}_m)$ on the irreducible $\ul{GL}_m$-module $\ul{X}_m(\theta_m(\blambda))$ with highest weight $\theta_m(\blambda)$. 
    \end{itemize}
\end{lem}

Recall (from Notation \ref{notation:interpolated_weights}) the interpolated Weyl vector $\rho_{\bt}$ and the $\rho$-shifted highest weight \[\Psi(\blambda)_{\rho_{\bt}} = (\psi_{\bt}(\blambda)_1,\ldots,\psi_{\bt}(\blambda)_{2n}).\]

\begin{prop}\label{prop:central_character}
    For each $k \geq 0$, the eigenvalue $\chi_{\bs,k}(\blambda,\bt)$ of $C_k$ on $\mts(\blambda)$ is given by
    \small
    \begin{align*}
     \sum_{i=1}^{2n} \Bigg(\sum_{j=1}^{\infty} \left(\psi_{\bt}(\blambda)_{i,j}^k -\psi_{\bt}(\b{\emptyset})_{i,j}^k\right)  + P_{i,k}(t_1,\ldots,t_n)\Bigg),
    \end{align*}
    \normalsize
    where $P_{i,k,\bs}(z_1,\ldots,z_n)$ is the polynomial interpolation of the following multivariate expression defined for positive integer tuples $(m_1,\ldots,m_n)$ by
    \[
    P_{i,k,\bs}(m_1,\ldots,m_n) :=\sum_{j=1}^{m_i} \left( s_i + \frac{(m_i + m_{i+1} + \cdots + m_n) - (m_1 + \ldots + m_{i-1}) + 1}{2} - j\right)^k.
    \]
\end{prop}

\begin{proof}
    The idea is to use a categorical algorithm (i.e., applicable in any symmetric tensor category) to reduce to the computation of central characters of simples in the Deligne categories $\cD_{t_i}$. Then, we apply Lemma \ref{lem:interpolate_casimir} to interpolate the classical central character formula to complex rank. 

    For concision, write $E_{ij} := V_i \ox V_j^* \subset \fgl_{\bt}$. Observe that the morphism $C_k$ equals the composition
    \[
    \unit \to \bigoplus_{i_1,i_2,\ldots,i_k} E_{i_1,i_2} \ox E_{i_2,i_3} \ox \cdots\ox E_{i_k,i_1} \xto{m} U_{\bt},
    \]
    where the first morphism is the direct sum of morphisms $\varphi_{i_1,i_2,\ldots,i_k}$ given by
    \begin{align*}
    \unit &\to V_{i_1} \ox V_{i_1}^* \to V_{i_1} \ox (V_{i_2}^* \ox V_{i_2}) \ox V_{i_1}^* \to V_{i_1} \ox V_{i_2}^* \ox V_{i_2} \ox (V_{i_3}^* \ox V_{i_3}) \ox V_{i_1}^* \to \cdots \\[5pt] &\to V_{i_1} \ox V_{i_2}^* \ox V_{i_2} \ox V_{i_3}^* \ox \cdots \ox V_{i_{k-1}}^* \ox V_{i_{k-1}} \ox (V_{i_k}^* \ox V_{i_k}) \ox V_{i_1}^* = E_{i_1,i_2} \ox E_{i_2,i_3} \ox \cdots \ox E_{i_{k},i_1}.
    \end{align*}
    The $j$th arrow above is given by coevaluation for $V_{i_j}^*$.

    For each $i = 1,\ldots,n$, we define a morphism $C_k(i): \unit \to U(\fgl_{t_i})$ in the same manner as $C_k$. Thus, if $i_1 = i_2 = \cdots = i_k = i$, then $m \circ \varphi_{i,i,\ldots,i}$ is the same as the composition \[\psi_{k,i} := \unit \xto{C_k(i)} U(\fgl_{t_{i}}) \hookrightarrow U_{\bt}.\]
    For any sequence $\{i_1,\ldots,i_k\}$, we claim that the composition $m \circ \varphi_{i_1,\ldots,i_k} \in Z(U_{\bt})$ acts on the highest weight subobject $\unit \ox X_{\bs}^{\bt}(\blambda)$ according to \[\sum_{a=1}^k \sum_{i=1}^n Q_{a,i}(t_1,\ldots,t_n)\psi_{a,i} \in Z(U_{\bt}),\] where $Q_{a,i}(x_1,\ldots,x_n) \in \fk[x_1,\ldots,x_n]$. This claim follows by iteratively applying the commutator to ``move'' factors $E_{i_\ell, i_{\ell+1}}$ to the right whenever $i_\ell > i_{\ell+1}$, using the following facts:
    \begin{enumerate}[(1)]
        \item we have an equality of morphisms $m \circ ([,] \oplus \tau) = m$, where $[,]: \fg_{\bt} \ox \fg_{\bt} \to \fg_{\bt}$ is the commutator and $\tau \in \End(\fg_{\bt} \ox \fg_{\bt})$ is the symmetric braiding,
        \item the sequence of coevaluations used to produce $\varphi_{i_1,\ldots,i_r}$ can be composed in any order,
        \item we have an equality $\ev_{V_{i_j}} \circ \coev_{V_{i_j}^*} = \dim V_{i_j} \id_{\unit} = t_{i_j}\id_{\unit}$,
        \item and $E_{i_{\ell},i_{\ell+1}}$ acts by zero on $\unit \ox X_{\bs}^{\bt}(\blambda)$ whenever $i_{\ell+1} > i_\ell$.
    \end{enumerate}
    In particular, we see that the action of $C_k$ on $\unit \ox X_{\bs}^{\bt}(\blambda)$ coincides with the action of the morphism 
    \[
        \sum_{a=1}^k \sum_{i=1}^n Q_{a,i}(t_1,\ldots,t_n)\psi_{a,i} \in Z(U_{\bt})
    \]
    for some \textit{polynomials} $Q_{a,i}(x_1,\ldots,x_n) \in \fk[x_1,\ldots,x_n]$. 

    Recall that underlines denote classical (integer rank) Lie theoretic objects. For any integer $m \geq n$ and composition $m_1 + m_2 + \cdots + m_n = m$, we pick the corresponding principal Levi subalgebra of $\ul{\fgl}_m$ and use this subalgebra to fix embeddings $Z(U(\ul{\fgl}_{m_i})) \hookrightarrow Z(U(\ul{\fgl}_{m}))$.
    We write $\ul{\psi}_{a,i}$ for the image of $\ul{C}_a^{(m_i)} \in Z(U(\ul{\fgl}_{m_i}))$ in $Z(U(\ul{\fgl}_m))$. 
    Observe that the commutator algorithm from complex rank is purely \textit{categorical} and also expresses the action of $\ul{C}_k^{(m)}$ in terms of the $\ul{\psi}_{a,i}$. That is, $\ul{C}_k^{(m)}$ acts on the highest weight space in a parabolic Verma module for $\ul{\fgl}_m$ by 
    \[
        \sum_{a=1}^k \sum_{i=1}^n Q_{a,i}(m_1,\ldots,m_n)\ul{\psi}_{a,i}\in \fk,
    \]
    where the polynomials $Q_{a,i}$ are the \textit{same} as the polynomials we found in complex rank.

    For each $i = 1,\ldots,n$, let us write $\blambda_i = (\lambda_i^\sharp,\lambda_i^\flat)$. Assuming that $m_i \gg 0$, define 
    {\small
    \[
    \theta(\blambda_i,s_i) = (\lambda^\sharp_{i,1} + s_i, \lambda^{\sharp}_{i,2} + s_i,\ldots, \lambda^{\sharp}_{i,\ell_i^\sharp} + s_i, s_i, \ldots, s_i, -\lambda^\flat_{i,\ell_i^\flat} + s_i, -\lambda^\flat_{i,\ell_i^\flat-1} + s_i, \ldots, -\lambda^\flat_{i,1} + s_i) \in \ZZ^{m_i} + s_i,
    \]
    }
    where $\ell_i^\sharp = \ell(\blambda_i^\sharp)$ and $\ell_i^\flat = \ell(\blambda_i^\flat)$. Similarly define
    \[
    \theta(\blambda_i) = (\lambda^\sharp_{i,1}, \lambda^{\sharp}_{i,2},\ldots, \lambda^{\sharp}_{i,\ell_i^\sharp}, 0, \ldots, 0, -\lambda^\flat_{i,\ell_i^\flat}, -\lambda^\flat_{i,\ell_i^\flat-1}, \ldots, -\lambda^\flat_{i,1}) \in \ZZ^{m_i},
    \]
    For a dominant $\ul{\fgl}_{m_i}$-weight $\theta$, we write $\ul{\psi}_{a,i}(\theta)$ for the scalar by which $\ul{\psi}_{a,i}$ acts on the irreducible $\ul{\fgl}_{m_i}$-module of highest weight $\theta$. A straightforward computation shows that
    \[
    \ul{\psi}_{a,i}(\theta(\blambda_i,s_i)) = \sum_{\ell=0}^{a} \binom{a}{\ell} \ul{\psi}_{a-\ell,i}(\theta(\blambda_i)) s_i^\ell.
    \] 
    In Lemma \ref{lem:interpolate_casimir}, we produced a polynomial $p_{a-\ell,\blambda_i}(x) \in \fk[x]$ such that $\ul{\psi}_{a-\ell,i}(\theta(\blambda_i)) = p_{a-\ell,\blambda_i}(m_i)$ when $m_i \gg 0$. Thus, whenever $m_i \gg 0$, the central element $C^{(m)}_k$ acts on the parabolic Verma module $M(\b{\theta})$ of highest weight $\b{\theta} := (\theta(\blambda_1,s_1), \theta(\blambda_2,s_2), \ldots, \theta(\blambda_n, s_n))$ by the scalar 
    \[
        \xi_{\bs,k}(\blambda)(m_1,\ldots,m_n) = \sum_{a=1}^k \sum_{i=1}^n \sum_{\ell=0}^a \binom{a}{\ell} Q_{a,i}(m_1,\ldots,m_n)p_{a-\ell,\blambda_i}(m_i)s_i^\ell.
    \]

    On the other hand, we know that $\psi_{a,i} \in Z(U_{\bt})$ acts on $\unit \ox X_{\bs}^{\bt}(\blambda)$ by the same scalar by which $\psi_{a,i} \in Z(U_{t_i})$ acts on the object $X(\blambda_i) \ox \unit_{s_i} \in \cD_{t_i}$. Write $\psi_{a,i}(\blambda_i,s_i)$ for this scalar. Moreover, write $\psi_{a,i}(\blambda_i)$ for the scalar by which $\psi_{a,i}$ acts on $X(\blambda_i)$. In this case, we can also compute
    \[
        \psi_{a,i}(\blambda_i,s_i) = \sum_{\ell=0}^{a} \binom{a}{\ell}  \psi_{a-\ell,i}(\blambda_i)s_i^\ell.
    \] 
    By Lemma \ref{lem:interpolate_casimir}, we see that $\psi_{a-\ell,i}(\blambda_i) = p_{a-\ell,\blambda_i}(t_i)$ for the \textit{same} polynomials $p_{a-\ell,\blambda_i}(x)$ from earlier. In particular, we have shown that 
    \[
        \chi_{\bs,k}(\blambda,\bt) = \sum_{a=1}^k \sum_{i=1}^n \sum_{\ell=0}^a \binom{a}{\ell}Q_{a,i}(t_1,\ldots,t_n)p_{a-\ell,\blambda_i}(t_i)s_i^\ell.
    \]
    In summary, we see that there is a polynomial $R(x_1,\ldots,x_n) \in \fk[x_1,\ldots,x_n]$ for which \[R(t_1,\ldots,t_n) =  \chi_{\bs,k}(\blambda,\bt)\] for all $t_1,\ldots,t_n \in \fk$ and \[R(m_1,\ldots,m_n) = \xi_{\bs,k}(\blambda)(m_1,\ldots,m_n)\] for $m_1,\ldots,m_n \gg 0$. We also have another expression for $\xi_{\bs,k}(\blambda)(m_1,\ldots,m_n)$: writing $\ul{\rho}(m)$ for the Weyl vector for $\ul{\fgl}_m$, then we also know that (for $m_1,\ldots,m_n \gg 0$), \[\xi_{\bs,k}(\blambda)(m_1,\ldots,m_n) = \sum_{r} ((\b{\theta}_r + \ul{\rho}(m)_r)^k - \ul{\rho}(m)_r^k).\] In coordinates, the entries of $\ul{\rho}(m)$ can be expressed as some polynomials in $m_1,\ldots,m_n$, i.e., there exists some polynomial $P(x_1,\ldots,x_n) \in \fk[x_1,\ldots,x_n]$ such that 
    \[P(m_1,\ldots,m_n) = \sum_r ((\b{\theta}_r + \ul{\rho}(m)_r)^k - \ul{\rho}(m)_r^k)\]
    for all integers $m_1,\ldots,m_n$. In other words, we have an equality of polynomials \[R(m_1,\ldots,m_n) = \xi_{\bs,k}(\blambda)(m_1,\ldots,m_n) = P(m_1,\ldots,m_n)\] for sufficiently large integers $m_1,\ldots,m_n$, so
    \[
    \chi_{\bs,k}(\blambda,\bt) = R(t_1,\ldots,t_n) = P(t_1,\ldots,t_n)
    \]
    for $t_1,\ldots,t_n \in \fk$. The expression in the proposition statement is the polynomial interpolation of \[P(m_1,\ldots,m_n) = \sum_r ((\b{\theta}_r + \ul{\rho}(m)_r)^k - \ul{\rho}(m)_r^k). \qedhere\]
\end{proof}

\section{Equivariant Jantzen Determinant Formula}\label{appendix:equivariant_jantzen}

We give a proof of an ``equivariant" Jantzen determinant formula (Proposition \ref{prop:equivariant_jantzen}).

\begin{prop}
    Up to a nonzero scalar, the following formula holds for $m_i \gg 0$:
    {\small
    \[
        \ul{D}^{\b{m}}_{\bs}(\blambda)_{\bmu} = \prod_{\alpha = \alpha^{p,q}_{i,j} \in R^+ \setminus R^+_{\b{m}}} \ \prod_{k \geq 1} (\phi_{m_i}(\blambda)_{i,p} - \phi_{m_j}(\blambda)_{j,q} + (s_i - M_i-p) - (s_j - M_j-q) - k)^{\tilde{r}_{\b{m}}(\blambda, k \alpha,\bmu)},
    \]
    }
    where we write \[\tilde{r}_{\b{m}}(\blambda,k\alpha,\bmu) = \mathrm{sgn}(\Theta_{\b{m}}^{\b{s}}(\blambda) - k\alpha)r_{\b{m}}((\Theta_{\b{m}}^{\b{s}}(\blambda) - k\alpha)_+,\bmu),\] and $r_{\b{m}}((\Theta_{\b{m}}^{\b{s}}(\blambda) - k\alpha)_+,\bmu)$ denotes the multiplicity of the $X^{\b{m}}(\bmu,\bs)$-isotypic component in the parabolic Verma module $\ul{M}^{\b{m}}((\Theta_{\b{m}}^{\b{s}}(\blambda) - k\alpha)_+)$. 
\end{prop}

\begin{proof}
    For any $\eta \in \fk^m$, Jantzen \cite{jantzen_determinant}\footnote{Strictly speaking, Jantzen's formula is stated for the case of a semisimple Lie algebra, but his formula readily generalizes to the case of $\fgl_{M}$. The reader can also refer to \cite{kac_kazhdan} for another derivation of Jantzen's formula.} gave the following formula for the determinant of the Shapovalov form on $M^{\b{m}}(\blambda,\bs)$ restricted to the $\eta$-weight space:
    {\small
    \[
        \ul{D}_{\bs}^{\b{m}}(\blambda)[\eta] = \prod_{\alpha = \alpha^{p,q}_{i,j} \in R^+ \setminus R^+_{\b{m}}} \ \prod_{k \geq 1} (\phi_{m_i}(\blambda)_{i,p} - \phi_{m_j}(\blambda)_{j,q} + (s_i - M_i-p) - (s_j - M_j-q) - k)^{\tilde{E}_{\b{m}}(\blambda, k \alpha,\eta)}
    \]
    }
    where \[\tilde{E}_{\b{m}}(\blambda, k \alpha,\eta) = \mathrm{sgn}(\Theta_{\b{m}}^{\b{s}}(\blambda) - k\alpha)E_{\b{m}}((\Theta_{\b{m}}^{\b{s}}(\blambda) - k\alpha)_+,\eta)\] and $E_{\b{m}}((\Theta_{\b{m}}^{\b{s}}(\blambda) - k\alpha)_+,\eta)$ is the dimension of the $\eta$-weight space in $M^{\b{m}}((\Theta_{\b{m}}^{\bs}(\blambda) - k\alpha)_+)$. 

    We seek an \textit{equivariant} version of Jantzen's formula, giving the determinant in the $X^{\b{m}}(\bmu,\bs)$-isotypic component. More generally, we hope to find the determinant in the $X^{\b{m}}(\vartheta)$-isotypic component for any $\fl_{\b{m}}$-dominant weight $\vartheta$. This equivariant formula takes the form
    {\small
    \[
        \ul{D}_{\bs}^{\b{m}}(\blambda)_{\vartheta} = \prod_{\alpha = \alpha^{p,q}_{i,j} \in R^+ \setminus R^+_{\b{m}}} \ \prod_{k \geq 1} (\phi_{m_i}(\blambda)_{i,p} - \phi_{m_j}(\blambda)_{j,q} + (s_i - M_i-p) - (s_j - M_j-q) - k)^{\tilde{r}_{\b{m}}(\blambda, k \alpha,\vartheta)}
    \]
    }
    for some exponents $\tilde{r}_{\b{m}}(\blambda, k \alpha,\vartheta)$. We must have 
    \[
       \tilde{E}_{\b{m}}(\blambda, k \alpha,\eta) = \sum_{\vartheta} \tilde{r}_{\b{m}}(\blambda, k \alpha,\vartheta) \dim X^{\b{m}}(\vartheta)_{\eta},
    \]
    where $\dim X^{\b{m}}(\vartheta)_{\eta}$ is the dimension of the $\eta$-weight space in $X^{\b{m}}(\vartheta)$. Thus, setting \[r_{\b{m}}((\Theta_{\b{m}}^{\b{s}}(\blambda) - k\alpha)_+,\vartheta) := \mathrm{sgn}(\Theta_{\b{m}}^{\b{s}}(\blambda) - k\alpha)\tilde{r}_{\b{m}}(\blambda,k\alpha,\vartheta),\] it follows that 
    \begin{align*}
        \sum_{\eta} E_{\b{m}}((\Theta_{\b{m}}^{\b{s}}(\blambda) - k\alpha)_+,\eta)e^{\eta} &= \sum_{\eta} \sum_{\vartheta} r_{\b{m}}((\Theta_{\b{m}}^{\b{s}}(\blambda) - k\alpha)_+,\vartheta) \dim X^{\b{m}}(\vartheta)_{\eta}e^{\eta} \\[5pt]
        &= \sum_{\vartheta} r_{\b{m}}((\Theta_{\b{m}}^{\b{s}}(\blambda) - k\alpha)_+,\vartheta)\mathrm{ch}(X^{\b{m}}(\vartheta))
    \end{align*}
    The left-hand side is $\mathrm{ch}(M^{\b{m}}((\Theta_{\b{m}}^{\bs} - k\alpha)_+))$, so it follows that $r_{\b{m}}((\Theta_{\b{m}}^{\b{s}}(\blambda) - k\alpha)_+,\vartheta)$ must be the multiplicity of $X^{\b{m}}(\vartheta)$ as a summand of $M^{\b{m}}((\Theta_{\b{m}}^{\bs} - k\alpha)_+)$. This completes the proof.
\end{proof}

\bibliography{parabolic_multiplicities}
\bibliographystyle{alpha}

\end{document}